\documentclass{article}
\usepackage{amsmath}
\usepackage{amssymb}
\usepackage{amsthm}
\usepackage{bm}
\usepackage{enumitem}
\usepackage{ifthen}
\usepackage{mathtools}
\usepackage{scalerel}
\usepackage{thmtools}
\usepackage{tikz-cd}
\usepackage{xparse} 
\usepackage{xstring}

\usepackage[cmtip,all,barr]{xy}
\newdir_{ (}{{ }*!/-.5em/@_{(}}

\newcommand \Ar {\mathrm{Ar}}

\DeclareMathOperator \arin {\stackrel{\Ar}{\in}}

\DeclareMathOperator \Atl {\mathcal{A}\mathrm{tl}}

\DeclareMathOperator \Atlas {\mathrm{Atlas}}

\newcommand \Bun {\mathrm{Bun}}

\newcommand \BunProd {\mathrm{Bun}-\mathrm{prod}}

\newcommand \Cat [1] {\mathbf{#1}}

\newcommand \catname [1] {\mathcal{#1}}
\newcommand \catseqname [1] {\bm{\mathcal{#1}}}

\newcommand \Ck {\mathrm{C^k}}

\newcommand {\compose} [1] [def]
   { \ifthenelse {\equal{#1}{def}}
       {\mathbin \circ}
       {\mathbin{\overset{#1} {\circ}}}
   }

\newcommand {\composeh} [1] [def]
   { \ifthenelse {\equal{#1}{def}}
       {\mathbin \odot}
       {\mathbin{\overset{#1} {\odot}}}
   }

\newcommand {\composet} [1] [def]
   { \ifthenelse {\equal{#1}{def}}
       {\mathbin \cdot}
       {\mathbin{\overset{#1} {\cdot}}}
   }

\DeclareMathOperator \defeq {\stackrel{\mathrm{def}}{=}}

\DeclareMathOperator \domain {\mathrm{domain}}

\DeclareMathOperator \Domain {\seqname{domain}}

\newcommand \false {\mathrm{False}}

\newcommand \Fib   {\mathrm{Fib}}

\newcommand \full [2] [] {\underset{\mathrm {{#1}full}}{#2}}

\newcommand \fullcref [1]
   { \ifthenelse {\equal{\nameref{#1}}{}}
       {\cref{#1}}
       {\cref{#1} (\nameref{#1})}
   }

\newcommand \funcname [1] {\mathit{#1}}
\newcommand \funcseqname [1] {\bm{#1}}

\DeclareMathOperator \Functor {\mathop{\mathcal{F}}}

\DeclareMathOperator \head {\mathrm{head}}

\DeclareMathOperator \Hom {\mathrm{Hom}}

\DeclareMathOperator \Id {\mathrm{Id}}

\DeclareMathOperator \ID {\mathbf{Id}}

\DeclareMathOperator \isCk {\mathrm{isCk}}

\DeclareMathOperator \isAtl {\mathrm{isAtl}}

\DeclareMathOperator \isLCS {\mathrm{isLCS}}

\DeclareMathOperator \iso {\stackrel{\sim}{=}}

\DeclareMathOperator \join {\mathrm{join}}

\newcommand \LCS {\mathrm{LCS}}

\DeclareMathOperator \length {\mathrm{length}}

\DeclareMathOperator \lengtho {\mathrm length0}

\newcommand \M {\mathrm{M-atlas}}

\newcommand \Man {\mathrm{Man}}

\newcommand \maps {\!\colon}

\newcommand{\mathvarname}[1]{%
  \begingroup\noexpandarg
  \StrLen{#1}[\temp]%
  \ifnum\temp>1
    \mathrm{#1}%
  \else
    #1%
  \fi
  \endgroup
}

\newcommand \maxfull [2] [] {\underset{\mathrm {{#1}max-full}}{#2}}
\newcommand \maximal [2] [] {\underset{\mathrm {{#1}max}}{#2}}

\newcommand \minimal [2] [] {\underset{{#1}\mathrm {min}}{#2}}

\newcommand \Mod {\mathrm{Mod}}

\newcommand \Ob {\mathrm{Ob}}
\DeclareMathOperator \objin {\stackrel{\Ob}{\in}}

\newcommand \onto {\epi}

\newcommand \op [2] [] {\underset{{#1}\mathrm {op}}{#2}}

\newcommand \optriv [2] [] {\underset{{#1}\mathrm {op-triv}}{#2}}

\newcommand \pagecref [1]
   { \ifthenelse {\equal{\nameref{#1}}{}}
       {\cref{#1} on \cpageref{#1}}
       {\cref{#1} (\nameref{#1}) on \cpageref{#1}}
   }

\newcommand \Pagecref [1]
   { \ifthenelse {\equal{\nameref{#1}}{}}
       {\Cref{#1} on \cpageref{#1}}
       {\Cref{#1} (\nameref{#1}) on \cpageref{#1}}
   }

\DeclareMathOperator \range {\mathrm{range}}

\DeclareMathOperator \Range {\seqname{range}}

\DeclareMathOperator \seqin {\stackrel{()}{\in}}

\newcommand \Set {\mathbf{Set}}

\newcommand \singcat [2] [] {\underset{{#1}\mathcal{S}ing}{#2}}
\newcommand \Singcat [2] [] {\underset{\bm{{#1}\mathcal{S}ing}}{#2}}

\newcommand \strict [2] [] {\underset{\mathrm {{#1}strict}}{#2}}

\newcommand \subcat [1] [] {\overset{\mathrm{{#1}cat}}{\subseteq}}
\newcommand \SUBCAT [1] [] {\overset{\mathbf{{#1}cat}}{\subseteq}}

\newcommand \submod [1] [] {\overset{\mathrm{{#1}mod}}{\subseteq}}

\DeclareMathOperator \SUBSETEQ {\stackrel{()}{\subseteq}}

\DeclareMathOperator \tail {\mathrm{tail}}

\newcommand \toiso {\,\to/{>}->>/^{\iso}}

\newcommand \Top {\mathrm{Top}}
\newcommand \Topcat {\mathcal{T}\mathrm{op}}

\newcommand \topname [1] {\mathfrak{#1}}

\newcommand \Topology {\mathfrak{Top}}

\newcommand \Triv [2] [] {\underset{{#1}\mathbf {triv}}{#2}}

\newcommand \trivcat [2] [] {\underset{{#1}\mathcal{T}\mathrm{riv}}{#2}}
\newcommand \Trivcat [2] [] {\underset{{#1}\bm{\mathcal{T}\mathrm{riv}}}{#2}}

\newcommand \true {\mathrm{True}}
\newcommand \truthcat {\mathcal{T}}
\newcommand \truthset {\mathbb{T}}
\newcommand \truthspace {\mathrm{Truthspace}}
\newcommand \truthtop {\mathfrak{Truthtop}}

\newcommand \unioncat [1] [] {\overset{\mathrm{{#1}cat}}{\cup}}
\newcommand \UNIONCAT [1] [] {\overset{\mathbf{{#1}cat}}{\cup}}

\ExplSyntaxOn

\int_gzero_new:N \g_style_quant_parens_int
\int_gzero_new:N \g_style_quant_subscr_int
\int_gzero_new:N \g_style_set_subscr_int

\NewDocumentCommand{\equant}{mm}
  {
    \quant:nnn {\exists} {#1} {#2}
  }

\NewDocumentCommand{\uquant}{mm}
  {
    \quant:nnn {\forall} {#1} {#2}
  }

\NewDocumentCommand \setupquant {m}
  {
    \keys_set:nn {shmuel / quant} {#1}
  }

\keys_define:nn {shmuel / quant}
 {
  subscript            .choices:nn =
    {
      {
        none,
        stacked,
        multiple
      }
      {
        \int_gset:Nn \g_style_quant_subscr_int {\l_keys_choice_int-1}
      }
    },
  subscript            .default:n = multiple,
  subscript            .initial:n = none,
  parentheses          .choices:nn =
    {
      {
        none,
        single,
        multiple
      }
      {
        \int_gset:Nn \g_style_quant_parens_int {\l_keys_choice_int - 1}
      }
    },
  parentheses          .default:n = multiple,
  parentheses          .initial:n = none,
  separater            .tl_set:N = \g_style_quant_sep_tl,
  separater            .default:n = {.},
  separater            .initial:n = {}
 }

\cs_new:Npn \quant:nnn #1 #2 #3
  {
    \clist_set:Nn \l_tmpa_clist {#2}
    \int_case:nn
      {\g_style_quant_subscr_int}
      {
        {0}
        {
          \int_compare:nTF {\g_style_quant_parens_int = 2}
            {\tl_set:Nn \l_tmpa_tl {\right ) \left ( #1}}
            {\tl_set:Nn \l_tmpa_tl {#1}}
          \int_compare:nT {\g_style_quant_parens_int > 0} {\left (}
          #1
          \clist_use:Nn \l_tmpa_clist {\l_tmpa_tl}
          \int_compare:nT {\g_style_quant_parens_int > 0} {\right )}
          \g_style_quant_sep_tl #3
        }
        {1}
        {
          \fp_set:Nn
            \l_tmpa_fp
            {ceil{\clist_count:N{\l_tmpa_clist} - 1} * .1 + 1}
          \int_compare:nT {\g_style_quant_parens_int > 0} {\left (}
          \scaleobj{\fp_to_decimal:N \l_tmpa_fp}{#1} \sb
             { \substack { \clist_use:Nn \l_tmpa_clist { \\ } } }
          \int_compare:nT {\g_style_quant_parens_int > 0} {\right )}
          #3
        }
        {2}
        {
          \clist_map_inline:Nn
            \l_tmpa_clist
            {
              {
                \int_compare:nT
                  {\g_style_quant_parens_int > 0}
                  {\left (}
                \scaleobj{1.2}{#1} \sb
                {##1}
                \int_compare:nT
                  {\g_style_quant_parens_int > 0}
                  {\right )}
              }
            }
          #3
        }
      }
  }

\NewDocumentCommand{\set}{mos}
  {
    \IfBooleanTF {#3}
    {
      \bool_set_true:N \l_tmpa_bool
    }
    {
      \bool_set_false:N \l_tmpa_bool
    }
    \set_of:nnn {#1} {#2} {\l_tmpa_bool}
  }

\NewDocumentCommand \setupset {m}
  {
    \keys_set:nn {shmuel / set} {#1}
  }

\keys_define:nn {shmuel / set}
  {
    separater .tl_set:N = \g_style_set_sep_tl,
    subscript            .choices:nn =
      {
        {
          stacked,
          multiple
        }
        {
          \int_gset:Nn \g_style_set_subscr_int {\l_keys_choice_int-1}
        }
      },
    separater .initial:n = {\mid},
    subscript .initial:n = stacked
  }

\cs_new:Npn \set_of:nnn #1 #2 #3
  {
    \IfValueTF {#2}
    {
    }
    {
    }
    \clist_set:Nn \l_tmpa_clist {#2}
    \tl_gset:Nn \g_tmpa_tl {\clist_use:Nn \l_tmpa_clist {\land}}
    \IfValueTF {#2}
    {
      \bool_if:nTF {#3}
      {
        \scalerel*{\{}{#1\g_tmpa_tl}
        #1
        \clist_set:Nn \l_tmpa_clist {#2}
        \scalerel*{\mathbin{\g_style_set_sep_tl}}{#1\g_tmpa_tl}
        \\
        \clist_set:Nn \l_tmpa_clist {#2}
        \g_tmpa_tl
        \clist_set:Nn \l_tmpa_clist {#2}
        \scalerel*{\}}{#1\g_tmpa_tl}
      }
      {
        \left \{
        #1
        \clist_set:Nn \l_tmpa_clist {#2}
        \scalerel*{\mathbin{\g_style_set_sep_tl}}{#1\g_tmpa_tl}
        \g_tmpa_tl
        \right \}
      }
    }
    {
      \left \{ #1 \right \}
    }
  }

\NewDocumentCommand{\seqname}{m}
  {
    \seqname:n {#1}
  }

\cs_new:Npn \seqname:n #1
  {
    \int_compare:nTF {\tl_count:n{#1} > 1}
      {
        {\bm{\mathrm{#1}}}
      }
      {
        {\bm {#1}}
      }
  }

\NewDocumentCommand{\intersection}{om}
  {
    \unint_of:nnn \bigcap {#1} {#2}
  }

\NewDocumentCommand{\union}{om}
  {
    \unint_of:nnn \bigcup {#1} {#2}
  }

\cs_new:Npn \unint_of:nnn #1 #2 #3
  {
    \IfValueTF {#2}
    {
      \clist_set:Nn \l_tmpa_clist {#2}
      \int_case:nn
        {\g_style_set_subscr_int}
        {
          {0}
          {
            \fp_set:Nn
              \l_tmpa_fp
              {ceil{\clist_count:N{\l_tmpa_clist} - 1} * .1 + 1}
            \scaleobj{\fp_to_decimal:N \l_tmpa_fp}{#1} \sb
            {\substack { \clist_use:Nn \l_tmpa_clist { \\ } }}
            #3
          }
          {1}
          {
            #1 \sb
            {\clist_use:Nn \l_tmpa_clist {,}}
            #3
          }
        }
    }
    {
      #1 #3
    }
  }

\ExplSyntaxOff

\newtheorem{theorem}{Theorem}[section]
\newtheorem{lemma}[theorem]{Lemma}
\newtheorem{corollary}[theorem]{Corollary}

\theoremstyle{definition}
\newtheorem{definition}[theorem]{Definition}
\newtheorem{example}[theorem]{Example}

\theoremstyle{remark}
\newtheorem{remark}[theorem]{Remark}

\numberwithin{equation}{section}

\newcommand \Alpha   A
\newcommand \Beta    B
\newcommand \Epsilon E
\newcommand \Zeta    Z
\newcommand \Eta     H
\newcommand \Iota    I
\newcommand \Kappa   K
\newcommand \Mu      M
\newcommand \Nu      N
\newcommand \Omicron O
\newcommand \Rho     P
\newcommand \Tau     T
\newcommand \Chi     S

\usepackage[colorlinks,hidelinks,draft=false]{hyperref}

\usepackage{cleveref}
\hypersetup {
   colorlinks,
   pdfinfo={
      Author={Shmuel (Seymour J.) Metz},
      Keywords={fiber bundles,manifolds},
      Subject={Topology},
      Title={Local Coordinate Spaces: a proposed unification of manifolds with fiber bundles, and associated machinery}
   }
}
\usepackage[final]{showlabels}
\showlabels{cite}
\showlabels{cref}
\showlabels{crefrange}

\begin{document}

\title%
{%
   Local Coordinate Spaces:
  a proposed unification of manifolds and fiber bundles,
  and associated machinery\thanks%
  {
  I wish to gratefully thank Walter Hoffman (z"l),
  Milton Parnes, Dr. Stanley H. Levy, the Mathematics department of
  Wayne State University, the Mathematics department of the State
  University of New York at Buffalo and others who guided my education.
  }
}
\author{Shmuel (Seymour J.) Metz
}
\maketitle




\begin{abstract}
This paper presents a unified view of manifolds and fiber bundles,
which, while superficially different, have strong parallels. It
introduces the notions of an m-atlas and of a local coordinate space,
and shows that special cases are equivalent to fiber bundles and
manifolds. Along the way it defines some convenient notation, defines
categories of atlases, and constructs potentially useful functors.
\end{abstract}

\setupquant {parentheses=multiple,subscript=stacked}
\setupset   {}

\tikzset{
    rot90/.style={anchor=south, rotate=90, inner sep=.2mm}
}

\section {Introduction}
\label{sec:intro}
Historically, the concept of pseudo-groups allowed unifying manifolds
and manifolds with boundary. The definitions of fiber bundles and
manifolds have strong parallels, and can be unified in a similar fashion;
there are several ways to do so. The central part of this paper,
\pagecref{sec:lcs}, defines an approach using categories and
commutative diagrams that is designed for easy exposition at the possible
expense of abstractness and generality. In particular, I have chosen to
assume the Axiom of Choice (AOC).

This paper treats atlases as objects of interest in their own right,
although it does not give them primacy. It introduces notions that
are convenient for use here and others that, while not used here, may be
useful for future work. It defines the new notions of
\hyperref[def:model]{model space}
\footnote{The phrase has been used before, but with a different
meaning.},
\hyperref[def:m-atlas]{m-atlas}
and of \hyperref[def:M-ATLmorph]{m-atlas morphism}.
Informally, a model space is a topological space with a category
specifying a family of open sets and functions satisfying specified
conditions.

Although this paper incidentally defines partial equivalents to
manifolds and fiber bundles using model spaces and model atlases, it
proposes the more general
\hyperref[def:LCS]{Local Coordinate Space (LCS)} in order to explicitly
reflect the role of the group in fiber bundles.

A local coordinate space (LCS) is a space (total space) with some
additional structure, including a coordinate model space and an atlas
whose transition functions are restricted to morphisms of the coordinate
model space; one can impose, e.g., differentiability restrictions, by
appropriate choice of the coordinate category. There is an equivalent
paradigm that avoids explicit mention of the total space by imposing
compatibility conditions on the transition functions, but that approach
is beyond the scope of this paper.

This paper defines functors among categories of atlases, categories of
model spaces, categories of local coordinate spaces, categories of
manifolds and categories of fiber bundles; it constructs more machinery
than is customary in order to facilitate the presentation of those
categories and functors.

\Crefrange{sec:conv}{sec:pre} present nomenclature and give basic
results.
\Cref{sec:m-charts} defines m-atlases, m-atlas morphisms and categories
of them; \cref{lem:ATLiscat} proves that the defined categories are
indeed categories.
\Cref{sec:lcs} defines local coordinate spaces and
categories of them; \pagecref{the:LCSiscat} proves that the defined
categories are indeed categories.

\Pagecref{Examples} gives some examples of structures that can be
represented as local coordinate spaces; \pagecref{sec:man} and
\pagecref{sec:bun} present two of the examples in detail, showing the
equivalence of manifolds and fiber bundles with special cases of local
coordinate spaces by explicitly exhibiting functors to and from local
coordinate spaces.
\begin{remark}
The unconventional definitions of manifold and fiber bundle are intended
to make their relationship to local coordinate spaces more natural.
\end{remark}

Most of the lemmata, theorems and corollaries in this paper should be
substantially identical to results that are familiar to the reader. What
is novel is the perspective and the material directly related to local
coordinate spaces. The presentation assumes only a basic knowledge of
Category Theory, such as may be found in the first chapter of
\cite{CftWM} or
{
  \showlabelsinline
  \cite{JoyCat}.
}

\subsection{New concepts and notation}
\label{sub:new}
This paper introduces a significant number of new concepts and some
modifications of the definitions for some conventional concepts. It also
introduces some notation of lesser importance.  The following are the
most important.

\begin{enumerate}
\item \hyperref[def:NCD]{Nearly commutative diagram (NCD)}, NCD at a point,
locally NCD and special cases with related nomenclature
\item \hyperref[def:model]{Model space} and related concepts
\item \hyperref[def:ModTop]{Model topology} and
\hyperref[def:M-para]{M-paracompactness}
\item \hyperref[sec:sig]{Signature},
\hyperref[def:Sigmacomm]{$\Sigma$-commutation} and related concepts
\item
\hyperref[def:LCS]{Local Coordinate Space (LCS)} and related concepts
\item
\hyperref[def:lin]{Linear space and related concepts}
\item
\hyperref[def:trivck]{Trivial $\Ck$ linear model space} and related
concepts
\item \hyperref[def:BunAtl]{$G$-$\rho$ bundle atlas}
\footnote{Similar to coordinate bundles}
and related concepts
\end{enumerate}

\section {Conventions}
\label{sec:conv}
A diagram arrow with an Equal-Tilde ($A \toiso_f B$) represents an
isomorphism. One with a hook ($A \underset{i}{\hookrightarrow} B$)
represents an inclusion map. One with a double arrowhead
($A \onto_\pi B$) represents a surjection.

All diagrams shown are commutative; none are exact.

Blackboard bold upper case will denote specific sets, e.g., the
Naturals.

Bold lower case italic letters will refer to sets, sequences and tuples
of functions, e.g.,
$\funcseqname{f} \defeq (\funcname{f}_1, \funcname{f}_2)$.

Bold lower case Latin letters will refer to sequence valued functions of
sequences and tuple valued functions of tuples, e.g., $\Range$ yields
the sequence of ranges of a sequence of functions.

Bold upper case calligraphic (script) letters will refer to
sequences of categories, e.g.,
$\catseqname{A} \defeq (\catname{A}_\alpha, \alpha \in \Alpha)$.

Bold upper case italic letters will refer to sequences or tuples, e.g.,
$\seqname{A}=(x,y,z)$, to sets of them, to sets of topological spaces or
to sets of open sets.

Fraktur will refer to topologies and to topology-valued functions, e.g.,
$\Topology$.

Functions have a range, domain and relation, not just a relation. Unless
otherwise stated, they asre assumed to be continuous.

Groups are assumed to be topological groups. The ambiguous notation
$x^{-1}$ will be used when it is obvious from context what the group
operation $\star$ and the group identity $\mathbf{1}_G$ are.

Lower case Greek letters other than $\pi$, $\rho$, $\sigma$, $\phi$ and
$\psi$ will refer to ordinals, possibly transfinite, and to formal
labels.  A letter with a Greek superscript and a letter with a Latin or
numeric superscript always refer to distinct variables.

Lower case $\pi$ will refer to a projection operator

Lower case $\rho$ will refer to a continuous effective group action,
i.e., a continuous representation of a group in a homeomorphism group.

Lower case $\sigma$ will refer to a sequence of ordinals, referred to as
a signature.

Lower case $\phi$ will refer to a coordinate function.

Lower case Latin letters will refer to
\begin{enumerate}
\item{elements of a set or sequence}
\item{functions}
\item{natural numbers}
\end{enumerate}

Upper case calligraphic (script) Latin letters will refer to categories
and functors.  Due to font limitations the special form $\mathcal{T}riv$
will be used instead of lower case calligraphic letters to refer to
constructed categories.

Upper case Greek letters other than $\Sigma$ may refer to
\begin{enumerate}
\item ordinal used as the limit of a sequence of consecutive ordinals,
e.g., $x_\alpha, \alpha \preceq \Alpha$
\item ordinal used as the order type of a sequence of consecutive
ordinals, e.g., $x_\alpha, \alpha \prec \Alpha$
\end{enumerate}

Upper case $\Sigma$ will refer to a sequence of signatures

Upper case Latin letters will refer to
\begin{enumerate}
\item Natural numbers
\item Topological spaces
\item Open sets
\item Elements of a sequence or tuple of functions, e.g.,
$\funcname{f}_E$ might be $\funcname{f}_0 \maps E_1 \to E_2$.
\end{enumerate}

Upright Latin letters will be used for long names.

The term $\Ck$ includes $\mathrm{C}^\infty$ (smooth) and
$\mathrm{C}^\omega$ (analytic).

This paper uses the term morphism in preference to arrow, but uses
the conventional $\Ar$.

The term sequence without an explicit reference to $\mathbb{N}$
will refer to a general ordinal sequence, possibly transfinite.

Sequence numbering, unlike tuple numbering, starts at 0 and the
exposition assumes a von Neumann definition of ordinals, so that
$\alpha \in \beta \equiv \alpha \prec \beta$.

Except where explicitly stated otherwise, all categories mentioned are
small categories with underlying sets, but the morphisms will often not
be set functions between the objects and there will not always be a
forgetful function to $\Set$ or $\Cat{Top}$. By abuse of language no
distinction will be made between a category $\catname{A}$ of topological
spaces and the concrete category $(\catname{A},\catname{U})$ over
$\Cat{Top}$.  Similarly, no distinction will be made among the object $U
\in \Ob(\catname{A})$, the topological space $\catname{U}(U)$ and the
underlying set.

When defining a category, the Ordered pair $(\seqname{O}, \seqname{M})$
refers to the smallest concrete category over $\Set$ or $\Cat{Top}$
whose objects are in $\seqname{O}$, whose morphisms from $o^1 \in
\seqname{O}$ to $o^2 \in \seqname{O}$ are functions $\funcname{f} \maps
o^1 \to o^2$ in $\seqname{M}$ and whose composition is function
composition.

When defining a category, the Ordered triple $(\seqname{O}, \seqname{M},
C)$ refers to the small category whose objects are in $\seqname{O}$,
whose morphisms are in $\seqname{M}$, whose $\Hom$ is

\begin{equation}
\Hom_{(\seqname{O}, \seqname{M}, C)}
  (o_1 \in \seqname{O}, o_2 \in \seqname{O})
\defeq
\set
  {
    (
      \funcseqname{f},
      o_1,
      o_2
    )
    \in \seqname{M}
  }
\end{equation}
and whose composition is C.

By abuse of language I may write
``$\catname{S}$'' for $\Ob(\catname{S})$,
``$A \in \catname{A}$'' for  $A \in\ \Ob(\catname{A})$,
``$A \subset \catname{A}$'' for $A \subset \Ob(\catname{A})$,
``$A \in \catname{A} \subset B \in \catname{B}$'' for
``the underlying set of $A$ is contained in the underlying set of $B$
and the inclusion $\funcname{i} \maps x \in A \hookrightarrow x \in B$
is a morphism'' and
``$\funcname{f} \maps A \to B$'' for
$\funcname{f} \in \Hom_{\catname{C}}(A,B)$, where $\catname{C}$ is
understood by context.

By abuse of language I shall use the same nomenclature for sequences and
tuples.

By abuse of language I shall use the $\times$ and $\bigtimes$ symbols
for both Cartesian products of sets and Cartesian products of
functions on those sets.

By abuse of language, and assuming AOC,  I shall refer to some sets as
ordinal sequences, e.g., ``$(C_\alpha, \alpha \in \Alpha)$'' for
``$\{C_\alpha \mid \alpha \in \Alpha\}$'', in cases where the order is
irrelevant.

By abuse of language, I may omit universal quantifiers in cases where
the intent is clear.

In some cases I define a notion similar to a conventional notion and
also need to refer to the conventional notion. In those cases I prefix
a letter or phrase to the term, e.g., m-paracompact versus paracompact.

\section {General notions}
\label{sec:notions}
This section describes nomenclature used throughout the paper. In
some cases this reflects new nomenclature or notions, in others it
simply makes a choice from among the various conventions in the
literature.

\begin{definition}[Operations on categories]
\label{def:catprop}
If $\catname{C}$ is a category then $x \objin \catname{C}$ iff $x$ is an
object of $\catname{C}$ and $y \arin \catname{C}$ iff $y$ is a morphism
of $\catname{C}$.

If $\catname{S}$ and $\catname{T}$ are categories then
$\catname{S} \subcat \catname{T}$ iff $S$ is a subcategory of
$\catname{T}$ and $\catname{S} \subcat[full-] \catname{T}$ iff $S$ is a
full subcategory of $\catname{T}$.

If $\catname{S}$ and $\catname{T}$ are categories then the category
union of $\catname{S}$ and $\catname{T}$, abbreviated
$\catname{S} \unioncat \catname{T}$, is the category whose objects are
in $\catname{S}$ or in $\catname{T}$ and whose morphisms are in
$\catname{S}$ or in $\catname{T}$.
\end{definition}

\begin{definition}[Identity]
$\Id_S$ is the identity function on the space $S$, $\Id_o$ is the
identity morphism for the object $o$\footnote{
The object is often expressed as a tuple, e.g.,
$\Id_{(\seqname{A}, \seqname{B})}$ is the identity morphism for the
object $(\seqname{A}, \seqname{B})$
},
$\Id_{U,V}$, for $U \subseteq V$, is the inclusion map,
$\Id_\catname{C}$ is the identity functor on the category $\catname{C}$.

$\ID_{\seqname{S}^i}$, $i=1,2$, is the sequence of identity functions
for the elements of the sequence
$\seqname{S}^i \defeq (\seqname{S}^1_\alpha,\ \alpha \prec \Alpha)$. Let
$\seqname{S}^1 \SUBSETEQ \seqname{S}^2$. Then
$\ID_{\seqname{S}^1,\seqname{S}^2}$ is the sequence of inclusion maps
$(\Id_{\seqname{S}^1_\alpha,\seqname{S}^2_\alpha}),\ \alpha \prec \Alpha$
for the elements of the sequences $\seqname{S}^i$.

The subscript may be omitted
when it is clear from context.
\end{definition}

\begin{definition}[Images]
$\funcname{f} [U] \defeq \set { {\funcname{f}(x)} }[x \in U]$ is the
image of $U$ under $\funcname{f}$ and
$\funcname{f}^{-1} [V] \defeq \set {x}[\funcname{f}(x) \in V]$ is the
inverse image of $V$ under $\funcname{f}$.

\begin{remark}
  This notation, adopted from \cite{GenTop}, avoids the ambiguity in
  the traditional $\funcname{f}(U)$ and $\funcname{f}^{-1}(V)$.
\end{remark}
\end{definition}

\begin{definition}[Projections]
\label{projections}
$\pi_\alpha$ is the projection function that maps a sequence into
element $\alpha$ of the sequence.  $\pi_i$ is also the projection
function that maps a tuple into element $i$ of the tuple.
\end{definition}

\begin{definition}[Topological category]
\label{def:topcat}
A topological category is a small subcategory of $\Cat{Top}$ or its
concrete category over $\Set$.

$\catname{T}$ is a full topological category iff it is a topological
category and whenever $U^i,V^i \objin \catname{T}$, $i=1,2$,
$V^i \subseteq U^i$,
$\funcname{f} \maps U^1 \to U^2 \arin \catname{T}$ and
$\funcname{f}[V^1] \subseteq V^2$ then \\
$\funcname{f} \maps V^1 \to V^2 \arin \catname{T}$.
\end{definition}

\begin{lemma}[Inclusions in topological categories are morphisms]
\label{lem:topInc}
Let $\catname{T}$ be a full topological category,
$S^i \objin \catname{T}$, $i=1,2$, and $S^1 \subseteq S^2$. Then
$\Id_{S^1,S^2}$ is a morphism of $\catname{T}$

\begin{proof}
$\Id_{S^2} \arin \catname{T}$, $S^1 \subseteq S^2$ by hypothesis and
$S^1 \subseteq S^2$, so $\Id_{S^1,S^2} \arin \catname{T}$ by
\cref{def:topcat}.
\end{proof}
\end{lemma}

\begin{definition}[Local morphisms]
\label{def:topLocal}
Let $\catname{T}^i$, $i=1,2$, be a full topological category and
$S^i \objin \catname{T}^i$, A continuous function
$\funcname{f} \maps S^1 \to S^2$ is locally a
$\catname{T}^1$-$\catname{T}^2$ morphism of $S^1$ to $S^2$ iff
$\catname{T}^1 \subcat[full-] \catname{T}^2$ and for every
$u \in S^1$ there is an open neighborhood $U_u$ for $u$ and an open
neighborhood $V_u$ for $v \defeq \funcname{f}(u)$ such that
$\funcname{f}[U_u] \subseteq V_u$ and $\funcname{f} \maps U_u \to V_u$
is a morphism of $\catname{T}^2$.

Let $\catname{T}$ be a full topological category and
$S^i \objin \catname{T}$, $i=1,2$. A continuous function
$\funcname{f} \maps S^1 \to S^2$ is locally a $\catname{T}$ morphism of
$S^1$ to $S^2$ iff it is locally a $\catname{T}$-$\catname{T}$ morphism
of $S^1$ to $S^2$.
\end{definition}

\begin{lemma}[Local morphisms]
\label{lem:topLocal}
Let $\catname{T}^i$, $i=1,2,3$, be a full topological category,
$\catname{T}^i \subcat[full-] \catname{T}^{i+1}$ and
$S^i \objin \catname{T}^i$,

If $\funcname{f}^i \maps S^i \to S^{i+1} \arin \catname{T}^{i+1}$ then
$\funcname{f}^i$ is locally a $\catname{T}^i$-$\catname{T}^{i+1}$
morphism of $S^i$ to $S^{i+1}$.

\begin{proof}
Let $u \in S^i$ and $v \defeq \funcname{f}^i(v) \in S^{i+1}$. $S^i$ is
an open for $u$, $S^{i+1}$ is an open neighborhood for $v$ and
$\funcname{f}^i \maps S^i \to S^{i+1} \arin \catname{T}^{i+1}$ by
hypothesis.
\end{proof}

If each $\funcname{f}^i \maps S^i \to S^{i+1}$, is locally a
$\catname{T}^i$-$\catname{T}^{i+1}$ morphism of $S^i$ to
$S^{i+1}$ then
$\funcname{f}^2 \compose \funcname{f}^1 \maps S^1 \to S^3$ is locally a
$\catname{T}^1$-$\catname{T}^3$ morphism of $S^1$ to $S^3$.

\begin{proof}
Since $\catname{T}^1 \subcat[full-] \catname{T}^2$ and
$\catname{T}^2 \subcat[full-] \catname{T}^3$,
$\catname{T}^1 \subcat[full-] \catname{T}^3$.
Let $u \in S^1$, $v \defeq \funcname{f}^1(u)$ and
$w \defeq \funcname{f}^2(v)$. There exist an open neighborhood $U_u$
for $u$, open neighborhoods $V_u$, $V'_v$ for $v$ and an open
neighborhood $W_v$ of $w$ such that
$\funcname{f}^1[U_u] \subseteq V_u$,
$\funcname{f}^1 \maps U_u \to V_u$ is a morphism of $\catname{T}^2$,
$\funcname{f}^2[V'_v] \subseteq W_v$ and
$\funcname{f}^2 \maps V'_v \to W_v$ is a morphism of $\catname{T}^3$.
Then $\hat{V_u} \defeq V_u \cap V'_v \neq \emptyset$, $\hat{V_u}$ is an
open neighborhood of $v$ and
$\hat{U_u} \defeq \funcname{f}^{i-1}_1[\hat{V_u}]$ is an open
neighborhood for $u$.  $\funcname{f}^1 \maps \hat{U_u} \to \hat{V_u}$
and $\funcname{f}^2 \maps \hat{V_u} \to W_v$ are morphisms of
$\catname{T}^3$ by \pagecref{def:topcat} and thus
$\funcname{f}^2 \compose \funcname{f}^1 \maps \hat{U_u} \to W_v$ is a
morphism of $\catname{T}^3$.
\end{proof}
\end{lemma}

\begin{corollary}[Local morphisms]
\label{cor:topLocal}
Let $\catname{T}^i$, $i=1,2$, be a full topological category,
$\catname{T}^i \subcat[full-] \catname{T}^{i+1}$,
$S^i \objin \catname{T}^i$ and $S^1 \subseteq S^2$. Then $\Id_{S^1,S^2}$
is locally a $\catname{T}^1$-$\catname{T}^2$ morphism of $S^1$ to $S^2$
and $\Id_{S^i}$ is locally a $\catname{T}$ morphism
of $S^i$ to $S^i$.

\begin{proof}
$S^1 \objin \catname{T}^2$ because $S^1 \objin \catname{T}^1$ and
$\catname{T}^1 \subcat \catname{T}^2$,
$S^2 \objin \catname{T}^2$ by hypothesis and
$S^1 \subseteq S^2$ by hypothesis, so
$\Id_{S^1,S^2} \arin \catname{T}^2$ by \cref{lem:topInc}.

$\Id_{S^i} \defeq \Id_{S^i,S^i}$.
\end{proof}
\end{corollary}

\begin{definition}[Sequence functions]
Let $\seqname{S} \defeq (s_\alpha, \alpha \prec \Alpha)$ be a sequence
of functions. Then
\begin{equation}
\Domain(\seqname{S}) \defeq \bigl ( \domain(s_\alpha), \alpha \prec \Alpha \bigr )
\end{equation}
\begin{equation}
\Range(\seqname{S}) \defeq \bigl ( \range(s_\alpha), \alpha \prec \Alpha \bigr )
\end{equation}

Let $\seqname{T} \defeq (t_\alpha, \alpha \prec \Alpha)$ be a sequence
of functions with $\Range(\seqname{S}) = \Domain(\seqname{T})$. Then
their composition is the sequence
$
  \seqname{T} \compose[()] \seqname{S} \defeq
  (t_\alpha \compose s_\alpha, \alpha \prec \Alpha)
$,

Let $\seqname{S} \defeq (s_\gamma, \gamma \preceq \Gamma)$, then these
functions extract information about the sequence:
\begin{equation}
\head(\seqname{S},\Omega) \defeq (s_\gamma, \gamma \prec \Omega)
\end{equation}
\begin{equation}
\head(\seqname{S}) \defeq \head(\seqname{S},\Gamma)
\end{equation}
\begin{equation}
\lengtho(\seqname{S}) \defeq \Gamma
\end{equation}
\begin{equation}
\tail(\seqname{S}) \defeq S_\Gamma
\end{equation}

Let $\seqname{S} \defeq (s_\gamma, \gamma \prec \Gamma)$, then
\begin{equation}
\length(\seqname{S}) \defeq \Gamma
\end{equation}

\begin{remark}
If $\lengtho(\seqname{S})$ is defined then
$\length(\seqname{S}) = \lengtho(\seqname{S}) + 1$.
$\lengtho(\seqname{S})$ is the ordinal type of
$\head(\seqname{S})$, not the ordinal type of $\seqname{S}$.
\end{remark}

Let $\catseqname{S} \defeq (\catname{S}_\alpha, \alpha \prec \Alpha)$
and $\catseqname{T} \defeq (\catname{T}_\alpha, \alpha \prec \Alpha)$ be
sequences of categories.  Then $\catseqname{S}$ is a subcategory
sequence of $\catseqname{T}$, abbreviated
$\catseqname{S} \SUBCAT \catseqname{T}$, iff every category in
$\catseqname{S}$ is a subcategory of the corresponding category in
$\catseqname{T}$, i.e.,
$
  \uquant%
    {\alpha \prec \Alpha}
    {\catname{S}_\alpha \subcat \catname{T}_\alpha}
$,
and $\catseqname{S}$ is a full subcategory sequence of $\catseqname{T}$,
abbreviated $\catseqname{S} \SUBCAT[full-] \catseqname{T}$, iff every
category in $\catseqname{S}$ is a full subcategory of the corresponding
category in $\catseqname{T}$, i.e.,
$
  \uquant%
    {\alpha \prec \Alpha}
    {\catname{S}_\alpha \subcat[full-] \catname{T}_\alpha}
$.

The category sequence union of $\catseqname{S}$ and $\catseqname{T}$,
abbreviated $\catseqname{S} \UNIONCAT \catseqname{T}$, is the sequence
of category unions of corresponding categories in $\catseqname{S}$ and
$\catseqname{T}$, i.e.,
$(\catseqname{S}_\alpha \unioncat \catseqname{T}_\alpha)$.

\end{definition}

\begin{lemma}[Sequence functions]
\label{lem:seqfunc}
Let $\funcseqname{f}^i \defeq (\funcname{f}^i_\alpha, \alpha \prec \Alpha)$,
$i=1,2,3$ be sequences of functions with
$\Domain(\funcseqname{f}^2) = \Range(\funcseqname{f}^1)$ and
$\Domain(\funcseqname{f}^3) = \Range(\funcseqname{f}^2)$. Then
$
(\funcseqname{f}^3 \compose[()] \funcseqname{f}^2) \compose[()] \funcseqname{f}^1 =
\funcseqname{f}^3 \compose[()] (\funcseqname{f}^2 \compose[()] \funcseqname{f}^1)
$.
\begin{proof}
\begin{equation}
\begin{split}
&
  (\funcseqname{f}^3 \compose[()] \funcseqname{f}^2) \compose[()] \funcseqname{f}^1 =
\\* &
  \bigl ( (\funcname{f}^3_\alpha \compose \funcname{f}^2_\alpha) \compose \funcname{f}^1_\alpha, \alpha \prec A \bigr ) =
\\* &
  \bigl ( \funcname{f}^3_\alpha \compose (\funcname{f}^2_\alpha \compose \funcname{f}^1_\alpha), \alpha \prec A \bigr ) =
\\* &
  \funcseqname{f}^3 \compose[()] (\funcseqname{f}^2 \compose[()] \funcseqname{f}^1)
\end{split}
\end{equation}
\end{proof}

Let $\funcseqname{f} \defeq (\funcname{f}_\alpha, \alpha \prec \Alpha)$
be a sequence of functions, $\seqname{D} = \Domain(\funcseqname{f})$ and
$\seqname{R} = \Range(\funcseqname{f})$. Then $ID_\seqname{R}$ is a left
$\compose[()]$ identity for $\funcseqname{f}$ and $ID_\seqname{D}$ is a
right $\compose[()]$ identity for $\funcseqname{f}$.

\begin{proof}
\begin{equation}
\begin{split}
& \ID_\seqname{R} \compose[()] \funcseqname{f} = \\*
& (\Id_{\range(\funcname{f}_\alpha)} \compose \funcname{f}_\alpha, \alpha \prec \Alpha) = \\*
& (\funcname{f}_\alpha, \alpha \prec \Alpha) = \\*
& \funcseqname{f}
\end{split}
\end{equation}
\begin{equation}
\begin{split}
& \funcseqname{f} \compose[()] \ID_\seqname{D} = \\*
& (\funcname{f}_\alpha \compose \Id_{\domain(\funcname{f}_\alpha)}, \alpha \prec \Alpha) = \\*
& (\funcname{f}_\alpha, \alpha \prec \Alpha) = \\*
& \funcseqname{f}
\end{split}
\end{equation}
\end{proof}
\end{lemma}

\begin{definition}[Tuple functions]
Let $\seqname{S} \defeq (s_n, n \in [1,N])$ be a tuple of functions. Then
\begin{equation}
\Domain(\seqname{S}) \defeq \bigl ( \domain(s_n), n \in [1,N] \bigr )
\end{equation}
\begin{equation}
\Range(\seqname{S}) \defeq \bigl ( \range(s_n), n \in [1,N] \bigr )
\end{equation}

Let $\seqname{T} \defeq (t_n, n \in [1,N])$ be a tuple of functions with
$\Range(\seqname{S})=\Domain(\seqname{T})$,
Then their composition is the tuple
$
  \seqname{T} \compose[()] \seqname{S} \defeq
  (t_n \compose s_n, n \in [1,N])
$

Let $\seqname{S} \defeq (s_m, m \in [1,M])$ and
$\seqname{T} \defeq (t_n, n \in [1,N])$ be tuples.
Then the following are tuple functioms
\begin{equation}
\head(S,I) \defeq (s_m, m \in [1,I])
\end{equation}
\begin{equation}
\head(S) \defeq \head(S,M-1)
\end{equation}
\begin{equation}
\tail(T,I) \defeq (t_n, n \in [I,N])
\end{equation}
\begin{equation}
\tail(T) \defeq t_N
\end{equation}
\begin{equation}
\join(S,T) \defeq (s_1, \dots, s_M, t_1, \dots, t_N)
\end{equation}
\end{definition}

\begin{lemma}[Tuple functions]
\label{lem:tupfunc}
Let $\funcseqname{f}^i \defeq (\funcname{f}^i_n, n \in [1,N])$,
$i=1,2,3$ be tuples of functions with
$\Domain(\funcseqname{f}^2) = \Range(\funcseqname{f}^1)$ and
$\Domain(\funcseqname{f}^3) = \Range(\funcseqname{f}^2)$. Then
$
(\funcseqname{f}^3 \compose[()] \funcseqname{f}^2) \compose[()] \funcseqname{f}^1 =
\funcseqname{f}^3 \compose[()] (\funcseqname{f}^2 \compose[()] \funcseqname{f}^1)
$.

\begin{proof}
\begin{equation}
\begin{split}
&
  (\funcseqname{f}^3 \compose[()] \funcseqname{f}^2) \compose[()] \funcseqname{f}^1 =
\\* &
  \bigl ( (\funcname{f}^3_n \compose \funcname{f}^2_n) \compose \funcname{f}^1_n, n \in [1,N] \bigr ) =
\\* &
  \bigl ( \funcname{f}^3_n \compose (\funcname{f}^2_n \compose \funcname{f}^1_n), n \in [1,N] \bigr ) =
\\* &
  \funcseqname{f}^3 \compose[()] (\funcseqname{f}^2 \compose[()] \funcseqname{f}^1)
\end{split}
\end{equation}
\end{proof}

Let $\funcseqname{f} \defeq (\funcname{f}_\alpha, \alpha \prec \Alpha)$
be a sequence of functions,
$\seqname{D} \defeq \Domain(\funcseqname{f})$ and \\*
$\seqname{R} \defeq \Range(\funcseqname{f})$. Then $\ID_\seqname{R}$ is a
left $\compose[()]$ identity for $\funcseqname{f}$ and $\ID_\seqname{D}$
is a right $\compose[()]$ identity for $\funcseqname{f}$.

\begin{proof}
\begin{equation}
\begin{split}
& \ID_\seqname{R} \compose[()] \funcseqname{f} = \\*
& (\Id_{\range(\funcname{f}_\alpha)} \compose \funcname{f}_\alpha, \alpha \prec \Alpha) = \\*
& (\funcname{f}_\alpha, \alpha \prec \Alpha) = \\*
& \funcseqname{f}
\end{split}
\end{equation}
\begin{equation}
\begin{split}
& \funcseqname{f} \compose[()] \ID_\seqname{D} = \\*
& (\funcname{f}_\alpha \compose \Id_{\domain(\funcname{f}_\alpha)}, \alpha \prec \Alpha) = \\*
& (\funcname{f}_\alpha, \alpha \prec \Alpha) = \\*
& \funcseqname{f}
\end{split}
\end{equation}
\end{proof}
\end{lemma}

\begin{definition}[Tuple composition for labeled morphisms]
\label{def:labcomp}
Let $\seqname{M}^i \defeq (\funcseqname{f}^i, o^i_1, o^i_2)$, $i=1,2$,
be tuples such that each $\funcseqname{f}^i$ is a sequence of functions
or each $\funcseqname{f}^i$ is a tuple of functions,
$\Range(\funcseqname{f}^1) = \Domain(\funcseqname{f}^2)$ and
$0^1_2=o^2_1$. Then

\begin{equation}
\seqname{M}^2 \compose[A] \seqname{M}^1 \defeq
\bigl (
  \funcseqname{f}^2 \compose[()] \funcseqname{f}^1,
  o^1_1,
  o^2_2
\bigr )
\end{equation}
\end{definition}

\begin{lemma}[Tuple composition for labeled morphisms]
\label{lem:atlcomp}
Let $\seqname{M}^i \defeq (\funcseqname{f}^i, o^i_1, o^i_2)$, $i=1,2,3$,
be tuples such that $\funcseqname{f}^i$ are sequences or tuples of
functions,
$\Range(\funcseqname{f}^i) = \Domain(\funcseqname{f}^{i+1})$
and $o^i_2 = o^{i+1}_1$, $i=1,2$.
Then
\begin{equation}
\seqname{M}^3 \compose[A] \bigl ( \seqname{M}^2 \compose[A] \seqname{M}^1 \bigr )
=
\bigl ( \seqname{M}^3 \compose[A] \seqname{M}^2 \bigr ) \compose[A] \seqname{M}^1
\end{equation}

\begin{proof}
From \fullcref{def:labcomp}, \\
\pagecref{lem:seqfunc} and
\fullcref{lem:tupfunc}, we have
\begin{equation}
\begin{split}
&
  \seqname{M}^3 \compose[A] (\seqname{M}^2 \compose[A] \seqname{M}^1)
=
\\* &
  \seqname{M}^3 \compose[A] (\funcseqname{f}^2 \compose[()] \funcseqname{f}^1, o^1_1, o^2_2) =
\\* &
  (\funcseqname{f}^3 \compose[()] \funcseqname{f}^2 \compose[()] \funcseqname{f}^1, o^1_1, o^3_2) =
\\* &
  (\funcseqname{f}^3 \compose[()] \funcseqname{f}^2, o^2_1, o^3_2) \compose[A] \seqname{M}^1 =
\\* &
  (\seqname{M}^3 \compose[A] \seqname{M}^2) \compose[A] \seqname{M}^1
\end{split}
\end{equation}
\end{proof}

Let $\seqname{D}^i \defeq \Domain(\funcseqname{f}^i)$ and
$\seqname{R}^i \defeq \Range(\funcseqname{f}^i)$. Then
$(\ID_{\seqname{R}^i}, o^i_2, o^i_2)$ is a left $\compose[A]$ identity for
$\seqname{M}^i$ and $(\ID_{\seqname{D}^i}, o^i_1, o^i_1)$ is a right
$\compose[A]$ identity for $\seqname{M}^i$.

\begin{proof}
\begin{equation}
\begin{split}
& (\ID_{\seqname{R}^i}, o^i_2, o^i_2) \compose[A] \seqname{M}^i = \\*
& (\ID_{\seqname{R}^i} \compose[()] \funcseqname{f}^i, o^i_1, o^i_2) = \\*
& (\funcseqname{f}^i, o^i_1, o^i_2) = \\*
& \seqname{M}^i
\end{split}
\end{equation}
\begin{equation}
\begin{split}
& \seqname{M}^i \compose[A] (\ID_{\seqname{D}^i}, o^i_1, o^i_1) = \\*
& (\funcseqname{f}^i \compose[()] \ID_{\seqname{D}^i}, o^i_1, o^i_2) = \\*
& (\funcseqname{f}^i, o^i_1, o^i_2) = \\*
& \seqname{M}^i
\end{split}
\end{equation}
\end{proof}
\end{lemma}

\begin{definition}[Cartesian product of sequence]
\label{def:Cart}
Let $\seqname{S}^i \defeq (S^i_\alpha, \alpha \prec \Alpha)$, $i=1,2$,
be a
sequence
and
$\funcseqname{f} \defeq (\funcname{f}_\alpha \maps S^1_\alpha \to
S^2_\alpha, \alpha \prec \Alpha)$
be a sequence of functions, then
$
  \bigtimes \seqname{S}^i \defeq
  \bigtimes_{\alpha \prec \Alpha} S^i_\alpha
$
is the generalized Cartesian product of the sequence $\seqname{S}^i$ and
$
  \bigtimes \funcseqname{f} \maps \seqname{S}^1 \to \seqname{S}^2
  \defeq
  \bigtimes_{\alpha \prec \Alpha}\funcname{f}_\alpha
$
is the generalized Cartesian product of the function sequence
$\funcseqname{f}$.
\end{definition}

\begin{definition}[underline]
Let $\seqname{S}^1 \defeq (S^1_\alpha, \alpha \preceq \Alpha)$,
$\seqname{S}^2 \defeq (S^2_\alpha, \alpha \preceq \Alpha)$ be
sequences and
$
  \funcseqname{f} \defeq
  (
    \funcname{f}_\alpha \maps S^1_\alpha \to S^2_\alpha,
    \alpha \preceq \Alpha
  )
$
be a sequence of functions, then
\begin{equation}
\underline{\funcname{f}} \maps head(S_1) \to head(S_2) \defeq
\bigtimes \head(\funcseqname{f}) =
\bigtimes_{\alpha \prec \Alpha} \funcname{f}_\alpha
\end{equation}
is the function mapping
$(s_\alpha \in S^1_\alpha, \alpha \prec \Alpha)$
into $(\funcname{f}_\alpha(s_\alpha), \alpha \prec \Alpha)$.
\end{definition}

\begin{definition}[Head and tail compositions]
Let \\
$\funcname{f}^1 \maps (S^1_\alpha, \alpha \prec \Alpha) \to S^1_\Alpha$,
$\funcname{f}^2 \maps (S^2_\alpha, \alpha \prec \Alpha) \to S^2_\Alpha$
and
$
  \funcseqname{g} \defeq
  (
    \funcname{g}_\alpha \maps S^1_\alpha \to S^2_\alpha,
    \alpha \preceq \Alpha
  )
$.
Then (see \cref{fig:fg,fig:gf})

\begin{subequations}
\begin{equation}
\funcname{f}^2 \composeh \funcseqname{g} \defeq
\funcname{f}^2 \compose \underline{\funcseqname{g}}
\end{equation}
\begin{equation}
\funcseqname{g} \composet \funcname{f}^1 \defeq
\tail(\funcseqname{g}) \compose \funcname{f}^1
\end{equation}
\end{subequations}
\end{definition}

\begin{figure}
\[ \bfig
\node s1(0,0)[{\head(S^1)}]
\node s2(0,-1000)[{\head(S^2)}]
\node s2a(1000,-1000)[{S^2_\Alpha}]
\arrow |l|[s1`s2;{\underline{\funcseqname{g}}}]
\arrow |r|[s1`s2a;{\funcseqname{f}^2 \composeh \funcseqname{g}}]
\arrow |b|[s2`s2a;{\funcseqname{f^2}}]
\efig \]
\caption{$\funcname{f}^2 \composeh \funcseqname{g}$}
\label{fig:fg}
\end{figure}

\begin{figure}
\[ \bfig
\node s1(0,0)[{\head(S^1)}]
\node s1a(0,-1000)[{S^1_\Alpha}]
\node s2a(1000,-1000)[{S^2_\Alpha}]
\arrow |l|[s1`s1a;{\funcname{f}^1}]
\arrow |r|[s1`s2a;{\funcseqname{g} \composet \funcseqname{f}^1}]
\arrow |b|[s1a`s2a;{\tail(\funcseqname{g}) = \funcseqname{g}_\Alpha}]
\efig \]
\caption{$\funcseqname{g} \composet \funcname{f}^1$}
\label{fig:gf}
\end{figure}

\begin{definition}[Topology functions]
Let $S$ be a topological space and $Y$ a subset. Then \\
\begin{enumerate}
\item $\Topology(S)$ is the topology of $S$.
\item
$\Topology(Y,S) \defeq \set{U \cap Y}[{U \in \Topology(S)}]$
is the relative topology of $Y$.
\item
$\Top(Y,S) \defeq \bigl ( Y, \Topology(Y,S) \bigr )$ is $Y$ with the
relative topology.
\item
$
  \op{S} \defeq
  \set
    {(U, \Topology(U,S))}%
    [{{ U \in \Topology(S) \setminus \{ \emptyset \} }}]
$
is the set of all non-null open subspaces of $S$.
\end{enumerate}

Let $\seqname{S}$ be a set of topological spaces. Then
$\op{\seqname{S}} \defeq \union [{S \in \seqname{S}}] { {\op{S}}}$ is the set
of open subspaces in $\seqname{S}$.

Let $S$ and $T$ be spaces, $S' \subseteq S$ and $T' \subseteq T$ be
subspaces and $\funcname{f} \maps S \to T$ a function such that
$\funcname{f}[S'] \subseteq T'$. Then $\funcname{f} \maps S' \to T'$,
also written $\funcname{f} \restriction_{S',T'}$, is
$\funcname{f} \restriction_{S'}$ considered as a function from $S'$ to
$T'$.

Let $\seqname{S}^i \defeq (S^i_\alpha, \alpha \preceq \Alpha)$, $i-1,2$,
be a sequence of spaces, $\seqname{S}^1 \SUBSETEQ \seqname{S}^2$ and \\
$\funcname{f}^2 \maps \head(\seqname{S}^2) \to \tail(\seqname{S}^2)$ a
function.
$
  \funcname{f}^2 \restriction_{\head(\seqname{S}^1)} \defeq
  \funcname{f}^2 \restriction_{\bigtimes \head(\seqname{S}^1)}
$.
If
$
  \funcname{f}^2 \restriction_{\head(\seqname{S}^1)}
    [\bigtimes \head(\seqname{S}^1)]
  \subseteq \tail(\seqname{S}^1)
$
then
$ \funcname{f}^2 \maps \head(\seqname{S}^1) \to \tail(\seqname{S}^1)$,
also written
$
  \funcname{f}^2 \restriction_%
    {\head(\seqname{S}^1),\tail(\seqname{S}^1)}
$,
is $\funcname{f}^2 \restriction_{\head(\seqname{S}^1)}$ considered as a
function from $\bigtimes \head(\seqname{S}^1)$ to $\tail(\seqname{S}^1)$.
\end{definition}

\begin{definition}[Truth space]
The truth set is $\truthset \defeq \set{\false, \true}$, where
$\false \defeq \emptyset$ and $\true \defeq \set{\emptyset}$.
$\truthtop \defeq \set {\emptyset,\truthset}$ and the truth space
$\truthspace \defeq (\truthset, \truthtop)$ is $\truthset$
with the indiscrete topology.
\end{definition}

\begin{definition}[Truth category]
The truth category is
$
   \truthcat \defeq \\
   \bigl (
     \truthspace,
     \{
       \false \to \false,
       \false \to \true,
       \true \to \true
     \}
   \bigr )
$.
The truth model space is
$\seqname{Truthspace} \defeq (\truthspace, \truthcat)$.
\end{definition}

\begin{definition}[Constraint functions]
A constraint function is a continuous function with range
$\truthspace$ or a model
function with range $\seqname{Truthspace}$.
\end{definition}

\begin{definition}[Sequence inclusion]
\label{def:seqin}
Let $\seqname{S} \defeq (S_\alpha, \alpha \prec \Alpha)$ and
$\seqname{T} \defeq (T_\alpha, \alpha \prec \Alpha)$ be sequences.
$\seqname{S} \seqin \seqname{T}$ iff
$\uquant{\alpha \prec \Alpha} {S_\alpha \in T_\alpha}$ or
$\uquant{\alpha \prec \Alpha} {S_\alpha \objin T_\alpha}$.
\end{definition}

\begin{lemma}[Sequence inclusion]
\label{lem:seqin}
Let $\seqname{S} \defeq (S_\alpha, \alpha \prec \Alpha)$
be a sequence,
$
  \catseqname{T}^i \defeq \\
  (\catname{T}^i_\alpha, \alpha \prec \Alpha)
$,
$i=1,2$, a sequence of
categories, $\catseqname{T}^1 \SUBCAT \catseqname{T}^2$ and
$\seqname{S} \seqin \catseqname{T}^1$ Then
$\seqname{S} \seqin \catseqname{T}^2$.

\begin{proof}
If
$\uquant{\alpha \prec \Alpha} {S_\alpha \objin \catname{T}^1_\alpha}$
then
$\uquant{\alpha \prec \Alpha} {S_\alpha \objin \catname{T}^2_\alpha}$.
\end{proof}
\end{lemma}

\section{Nearly commutative diagrams}
\label{sec:ncd}
The notion of commutative diagrams is very useful in, e.g., Algebraic
Topology. Often one encounters commutative diagrams in which two
outgoing terminal nodes can be connected by a bridging function such
that the resulting diagram is still commutative. This paper uses the
term nearly commutative to describe a restricted class of such diagrams.

\begin{figure}
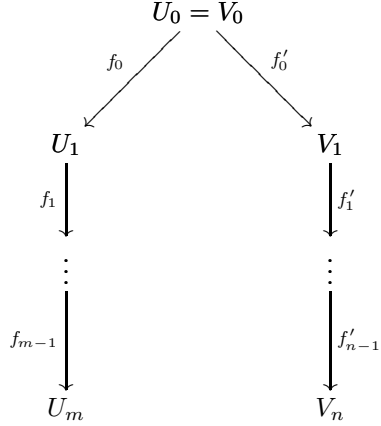

\[ \bfig
\node u0(0,0)[U_0=V_0]
\node u1(-500,-500)[U_1]
\node v1(500,-500)[V_1]
\node d1(-500,-1000)[\vdots]
\node d2(500,-1000)[\vdots]
\node um(-500,-1500)[U_m]
\node vn(500,-1500)[V_n]
\arrow |l|[u0`u1;\funcname{f}_0]
\arrow |r|[u0`v1;\funcname{f}'_0]
\arrow |l|[u1`d1;\funcname{f}_1]
\arrow |r|[v1`d2;\funcname{f}'_1]
\arrow |l|[d1`um;\funcname{f}_{m-1}]
\arrow |r|[d2`vn;\funcname{f}'_{n-1}]
\efig \]

\caption{Uncompleted nearly commutative diagram}
\label{fig:NCDb}
\end{figure}

Let $\catname{C}$ be a full topological category and $D$ a tree with two
branches, whose nodes are topological spaces $U_i$ and $V^j$ and whose
links are continuous functions $\funcname{f}_i \maps U_i \to U_{i+1}$
and $\funcname{f}'_j \maps U_j \to U_{j+1}$ between the sets:

\begin{align*}
D = \{
  &
    \funcname{f}_0 \maps U_0 = V_0 \to U_1,
    \dotsc,
    \funcname{f}_{m - 1} \maps U_{m - 1} \to U_m,
\\
  &
    \funcname{f}'_0 \maps  U_0 = V_0 \to V_1,
    \dotsc,
    \funcname{f}'_{m-1} \maps V_{m-1} \to V_n
\}
\end{align*}
with $U_0 = V_0$, $U_m \objin \catname{C}$ and
$V_n \objin \catname{C}$, as shown in \cref{fig:NCDb}.

\begin{figure}
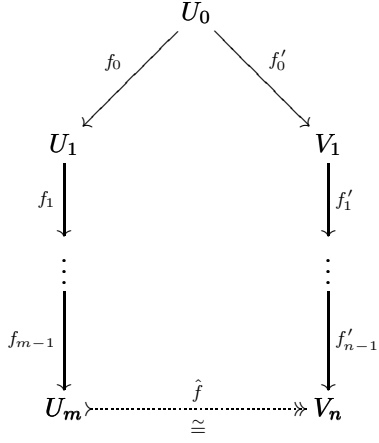

\[ \bfig
\node u0(0,0)[U_0]
\node u1(-500,-500)[U_1]
\node v1(500,-500)[V_1]
\node d1(-500,-1000)[\vdots]
\node d2(500,-1000)[\vdots]
\node um(-500,-1500)[U_m]
\node vn(500,-1500)[V_n]
\arrow |l|[u0`u1;\funcname{f}_0]
\arrow |r|[u0`v1;\funcname{f}'_0]
\arrow |l|[u1`d1;\funcname{f}_1]
\arrow |r|[v1`d2;\funcname{f}'_1]
\arrow |l|[d1`um;\funcname{f}_{m-1}]
\arrow |r|[d2`vn;\funcname{f}'_{n-1}]
\arrow |a|/>.>>/[um`vn;\hat{\funcname{f}}]
\arrow |b|//[um`vn;\iso]
\efig \]
\caption{Completed nearly commutative diagram}
\label{fig:NCDa}
\end{figure}

\begin{definition}[Nearly commutative diagrams in category $\catname{C}$]
\label{def:NCD}
$D$ is nearly commutative in category $\catname{C}$ iff the two final
nodes are in $\catname{C}$ and there is an isomorphism
$\hat{\funcname{f}} \maps U_m \toiso V_n$ in $\catname{C}$ making the
graph a commutative diagram, as shown in
\cref{fig:NCDa}.
\end{definition}

\begin{figure}
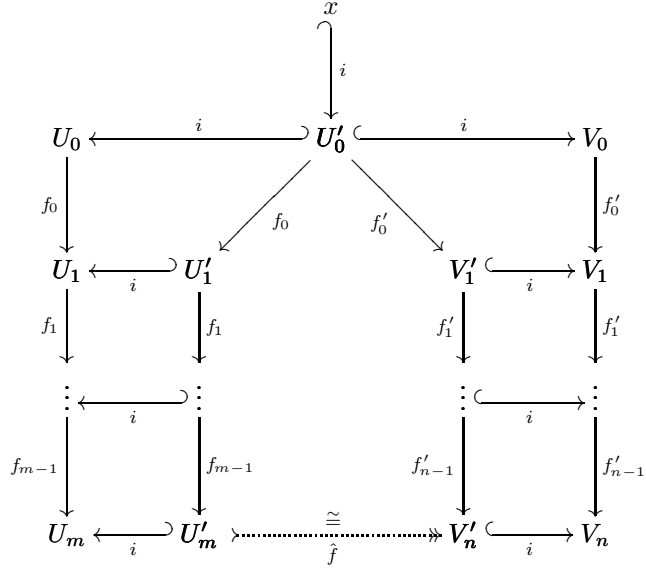

\[ \bfig
\node x(0,500)[x]
\node u0(-1000,0)[U_0]
\node up0(0,0)[U'_0]
\node v0(1000,0)[V_0]
\node u1(-1000,-500)[U_1]
\node up1(-500,-500)[U'_1]
\node vp1(500,-500)[V'_1]
\node v1(1000,-500)[V_1]
\node d1(-1000,-1000)[\vdots]
\node d2(-500,-1000)[\vdots]
\node d3(500,-1000)[\vdots]
\node d4(1000,-1000)[\vdots]
\node um(-1000,-1500)[U_m]
\node upm(-500,-1500)[U'_m]
\node vpn(500,-1500)[V'_n]
\node vn(1000,-1500)[V_n]
\arrow |r|/_{ (}->/[x`up0;i]
\arrow |a|/_{ (}->/[up0`u0;i]
\arrow |b|[up0`up1;\funcname{f}_0]
\arrow |b|[up0`vp1;\funcname{f}'_0]
\arrow |a|/^{ (}->/[up0`v0;i]
\arrow |l|[u0`u1;\funcname{f}_0]
\arrow |r|[v0`v1;\funcname{f}'_0]
\arrow |l|[u1`d1;\funcname{f}_1]
\arrow |b|/_{ (}->/[up1`u1;i]
\arrow |r|[up1`d2;\funcname{f}_1]
\arrow |b|/^{ (}->/[vp1`v1;i]
\arrow |l|[vp1`d3;\funcname{f}'_1]
\arrow |r|[v1`d4;\funcname{f}'_1]
\arrow |l|[d1`um;\funcname{f}_{m-1}]
\arrow |b|/_{ (}->/[d2`d1;i]
\arrow |r|[d2`upm;\funcname{f}_{m-1}]
\arrow |l|/^{ (}->/[d3`d4;i]
\arrow |l|[d3`vpn;\funcname{f}'_{n-1}]
\arrow |r|[d4`vn;\funcname{f}'_{n-1}]
\arrow |b|/_{ (}->/[upm`um;i]
\arrow |a|//[upm`vpn;\iso]
\arrow |b|/ >.>>/[upm`vpn;\hat{\funcname{f}}]
\arrow |b|/^{ (}->/[vpn`vn;i]
\efig \]
\caption{Local nearly commutative diagram}
\label{fig:NCDl}
\end{figure}

\begin{definition}[Nearly commutative diagrams in category $\catname{C}$ at a point]
Let $\catname{C}$ and $D$ be as above and $x$ be an element of the
initial node. $D$ is nearly commutative in $\catname{C}$ at $x$ iff
there are subobjects of the nodes such that the tree formed by replacing
the nodes is nearly commutative in $\catname{C}$ and $x$ is in the new
initial node, as shown in \pagecref{fig:NCDl}:
$x \in U'_0 = V'_0$, $U'_i \subseteq U_i$, $V'_j \subseteq V_j$,
$\funcname{f}_i \restriction_{U'_i} \maps U'_i \to U'_{i+1}$,
$\funcname{f}'_j \restriction_{V'_j} \maps V'_j \to V'_{j+1}$.
\end{definition}

\begin{definition}[Locally nearly commutative diagrams in category $\catname{C}$]
Let $\catname{C}$ and $D$ be as above. $D$ is locally nearly
commutative in $\catname{C}$ iff it is nearly commutative in
$\catname{C}$ at $x$ for every $x$ in the initial node.
\end{definition}

\begin{lemma}[Locally nearly commutative diagrams in category $\catname{C}$]
Let $\catname{C}$ and $D$ be as above. If $U_0 = \emptyset$ then $D$ is
locally nearly commutative in $\catname{C}$.

\begin{proof}
$D$ is vacuously locally nearly commutative at every $x \in U_0$ since
there is no such $x$.
\end{proof}

Let $D$ be locally nearly commutative in $\catname{C}$ and
let $\hat{U_0}=\hat{V_0} \subseteq U_0$ and $\hat{D}$ be $D$ with
$U_0=V_0$ replaced by $\hat{U_0}=\hat{V_0}$. Then $\hat{D}$ is locally
nearly commutative in $\catname{C}$.

\begin{proof}
If $x \in \hat{U_0}$ then $x \in U_0$ and hence $D$ is locally nearly
commutative in $\catname{C}$ at $x$.
Replacing $U'_0$ with $U'_0 \cap \hat{U_0}$ in the definition shows that
$\hat{D}$ is locally nearly commutative in $\catname{C}$ at $x$.
\end{proof}
\end{lemma}

\begin{remark}
It will often be clear from context what the relevant categories are.
This paper may use the term ``nearly commutative'' without explicitly
identifying the categories in which the modes are found.
\end{remark}

\section{Model spaces}
\label{sec:modeldef}
Let $S$ be a topological space. We need to formalize the notions of an
open cover by sets that are ``well behaved'' in some sense, e.g.,
convex, sufficiently small, and of "well behaved" functions among those
sets, e.g., preserving fibers, smooth. We do this by associating a
category of acceptable sets and functions.

\begin{remark}
Using pseudo-groups, as in \cite[p.~1]{FoundDiffGeo1}, would not allow
restricting model neighborhoods to, e.g., convex sets.
\end{remark}

\begin{definition}[Model spaces]
\label{def:model}
Let $S$ be a topological space and $\catname{S}$ a small category whose
objects are open subsets of $S$ and whose morphisms are continuous
functions. $\seqname{M} \defeq (S, \catname{S})$ is a model space for
$S$ iff

\begin{enumerate}
\item $\Ob(\catname{S})$ is an open cover for $S$. Note that it need
not be a basis for $S$.
\label{mod:cover}
\item $\Ob(\catname{S})$ is closed under finite intersections.
\label{mod:closed}
\item The morphisms of $\catname{S}$ are continuous functions in $S$.
\label{mod:continuous}
\item If $\funcname{f} \maps A \to B$ is a morphism,
$A' \in \Ob(\catname{S}) \subseteq A \in \Ob(\catname{S})$,
$B' \in \Ob(\catname{S}) \subseteq B \in \Ob(\catname{S})$
and $\funcname{f}[A'] \subseteq B'$ then
$\funcname{f} \restriction_{A'} \maps A' \to B'$ is a morphism.
\label{mod:restriction}
\item If $A' \objin \catname{S} \subseteq A \objin \catname{S}$ then
the inclusion map $\Id_{A',A} \maps A' \hookrightarrow A$ is a
morphism.
\label{mod:inclusion}
\begin{remark}
This is actually a consequence of \cref{mod:restriction}, but it is
convenient to give it here.
\end{remark}
\item Restricted sheaf condition:
informally, consistent morphisms can be glued together. Whenever
\begin{enumerate}
\item $U_\alpha$ and $V_\alpha$, $\alpha \prec \Alpha$,
are objects of $\catname{S}$.
\item $\funcname{f}_\alpha \maps U_\alpha \to V_\alpha$ are morphisms of
$\catname{S}$.
\item
$U \defeq \union[\alpha \prec \Alpha]{U_\alpha} \in \Ob(\catname{S})$,
\item
$V \defeq \union[\alpha \prec \Alpha]{V_\alpha} \in \Ob(\catname{S})$
\item
$\funcname{f} \maps U \to V$ is a continuous function and
for every $\alpha \prec \Alpha$,
$\funcname{f}$ agrees with $\funcname{f}_\alpha$ on $U_\alpha$
\end{enumerate}
then $\funcname{f}$ is a morphism of $\catname{S}$.
\label{mod:sheaf}
\end{enumerate}

$\Top(\seqname{M}) \defeq S = \pi_1(\seqname{M})$.

Let $U \objin \catname{S}$. Then $\Top(U,\seqname{M}) \defeq \Top(U,S)$
is $U$ with the relative topology.

$
  \Topcat(\seqname{M}) \defeq
  \bigl (
    \set
      {{\Top(U,S)}}%
      [{U \objin \catname{S}}],
    \Ar(\catname{S})
  \bigr )
$
is the topological category of $\seqname{M}$.

By abuse of language we write $U \subseteq \seqname{M}$.
\end{definition}

\begin{lemma}[The topological category of a model space is a full topological category]
\label{lem:modtopcat}
Let $\seqname{M} \defeq (S, \catname{S})$ be a model space for $S$. Then
$\Topcat(\seqname{M})$ is a full topological category.
\begin{proof}
$\Topcat(\seqname{M})$ is a small subcategory of $\Cat{Top}$ by
construction. $\Topcat(\seqname{M})$ is a full topological category by
\cref{mod:restriction} of
{
  \showlabelsinline
  \cref{def:model}.
}
\end{proof}
\end{lemma}

\begin{definition}[Model neighborhoods]
Let $\seqname{S} \defeq (S, \catname{S})$ be a model space for $S$.
Then the objects of $\catname{S}$ are model neighborhoods of
$\seqname{S}$. If $u \in U \objin \catname{S}$ then $U$ is a model
neighborhood for $u$. If $U_i \objin \catname{S}$, $i=1,2$, and
$U_1 \subseteq U_2$ then $U_1$ is a model subneighborhood of $U_2$.
\end{definition}

\begin{definition}[Model subspaces]
$\seqname{M} \defeq (X, \catname{X})$ is a model subspace of
$\seqname{N} \defeq (Y, \catname{Y})$, abbreviated
$\seqname{M} \submod \seqname{N}$, iff

\begin{enumerate}
\item $X$ is a subspace of $Y$
\item $\catname{X} \subcat[full-] \catname{Y}$
\item Every object of $\catname{Y}$ contained in $X$ is an object of
$\catname{X}$
\item When $\funcname{f} \maps o_1 \objin \catname{X} \to o_2 \objin \catname{X}$
is a morphism of $\catname{Y}$ then $\funcname{f} \maps o_1 \to o_2$ is a
morphism of $\catname{X}$.
\end{enumerate}

Let $\seqname{M} \defeq (X, \catname{X})$ be a model space and $Y$ a
model neighborhood of $\seqname{M}$. Then $\Mod(Y,\seqname{M})$, the
relative model space of $Y$, is $(Y, \catname{Y})$, where $\catname{Y}$
is the full subcategory of $\catname{X}$ containing all model
subneighborhoods of $Y$.
\end{definition}

\subsection{M-nearly commutative diagrams}
\label{sub:m-ncd}
Let $\seqname{M} \defeq (M,\catname{M})$ be a model space and $D$ a tree
with two branches, whose nodes are topological spaces $U_i$ and $V^j$
and whose links are continuous functions
$\funcname{f}_i \maps U_i \to U_{i+1}$ and
$\funcname{f}'_j \maps U_j \to U_{j+1}$ between the spaces:

\begin{align*}
D = \{
  &
    \funcname{f}_0 \maps U_0 = V_0 \to U_1,
    \dotsc,
    \funcname{f}_{m - 1} \maps U_{m - 1} \to U_m,
\\
  &
    \funcname{f}'_0 \maps  U_0 = V_0 \to V_1,
    \dotsc,
    \funcname{f}'_{m-1} \maps V_{m-1} \to V_n
\}
\end{align*}
with $U_0 = V_0$, $U_m \objin \catname{M}$ and
$V_n \objin \catname{M}$, as shown in \pagecref{fig:NCDb}.

\begin{definition}[M-nearly commutative diagrams]
\label{def:M-NCD}
$D$ is M-nearly commutative in model space $\seqname{M}$ iff $D$ is
nearly commutative in category $\Topcat(\seqname{M})$.
\end{definition}

\begin{definition}[M-nearly commutative diagrams at a point]
Let $\seqname{M}$ and $D$ be as above and $x$ be an element of the
initial node. $D$ is M-nearly commutative in $\seqname{M}$ at $x$ iff
$D$ is nearly commutative in category $\Topcat(\seqname{M})$ at $x$.
\end{definition}

\begin{definition}[M-locally nearly commutative diagrams]
Let $\catname{C}$ and $D$ be as above. $D$ is M-locally nearly
commutative in $\seqname{M}$ iff $D$ is nearly commutative in
category $\Topcat(\seqname{M})$ at every point.
\end{definition}

\subsection{Trivial model spaces}
\label{sub:triv}
Informally, a trivial model space of a specific type is one that does
not restrict the potential objects and morphisms of its type.

\begin{definition}[Trivial model spaces]
\label{def:trivmod}
Let $S$ be a topological space and $\catname{S}$ the category of all
continuous functions between open sets of $S$. Then
$\Triv{S} \defeq (S, \catname{S})$ is the trivial model space
of $S$ and $\Triv{S}$ is a trivial model space.

Let $\seqname{S}$ be a set of topological spaces. Then
$
  \Triv{\seqname{S}}
  \defeq
  \set
    {\Triv{S'}}%
    [S' \in \seqname{S}]
$
is the set of all trivial model spaces in $\seqname{S}$ and
$\Trivcat{\seqname{S}}$, the category of all continuous functions
between objects of $\Triv{\seqname{S}}$, is the trivial model category
of $\seqname{S}$.
\end{definition}

\begin{lemma}[The trivial model space of $S$ is a model space]
\label{lem:trivmod}
Let $S$ be a topological space. Then  $\Triv{S}$ is a model space.
\begin{proof}
Let $\catname{S} \defeq \pi_2 \bigl ( \Triv{S} \bigr )$,
$f \maps A \to B$ be a morphism of $\catname{S}$,
$A' \in \Ob(\catname{S}) \subseteq A$ and
$B' \in \Ob(\catname{S}) \subseteq B$. Then the definition of continuity and $\Triv{S}$ imply each of the following.
\begin{enumerate}
\item $\Ob(\catname{S})=\Topology(S)$ and thus is an open cover of $S$.
\item $\Ob(\catname{S})=\Topology(S)$ and thus is closed under finite intersections.
\item All morphisms in $\Ar(\catname{S})$ are continuous.
\item Since $f \maps A \to B$ is a morphism of $\catname{S}$, $f$ is continuous.
If $A' \in \Ob(\catname{S})$ then $A'$ is open.
If $B' \in \Ob(\catname{S})$ then $B'$ is open.
If $f \maps A \to B$ is continuous then $f \restriction_{A'} \maps A' \to B$
is continuous.
If $f[A'] \subseteq B'$ then $f \restriction_{A'} \maps A' \to B'$ is
well defined and continuous, hence a morphism.
\item If $A' \in \Ob(\catname{S}) \subseteq A \in \Ob(\catname{S})$ then
\begin{enumerate}
\item $A$ and $A'$ are open.
\item
The inclusion map $\Id_{A'} \maps A' \hookrightarrow A$ is continuous.
\item The inclusion map $\funcname{i} \maps A' \hookrightarrow A$ is a
morphism $\catname{S}$ by the definition of $\Triv{S}$.
\end{enumerate}
\end{enumerate}
\end{proof}
\end{lemma}

\subsection{Minimal model spaces}
\label{sub:min}
Even if a set of open sets fails one of \cref{mod:cover,mod:closed} or a
set of functions among them fails one of
\crefrange{mod:continuous}{mod:sheaf} in the definition of a model
space, there is a minimal model space containing them.

\begin{definition}[Minimal model spaces]
\label{def:minmod}
Let $S$ be a topological space, $\seqname{O}$ a set of open sets in $S$,
$\funcseqname{f}$ a set of continuous functions between elements of
$\seqname{O}$ and $\catname{S}$ the smallest concrete category  over
$\Top$ having all sets in $\seqname{O}$ as objects, having all functions
in $\funcseqname{f}$ as morphisms and satisfying
\crefrange{mod:continuous}{mod:sheaf} of \pagecref{def:model}.

Then
\begin{equation}
\minimal{\Mod}(S, \seqname{O}, \funcseqname{f}) \defeq
(\union{O}, \catname{S})
\end{equation}
is the minimal model space of $S$ with neighborhoods $\seqname{O}$ and
neighborhood mappings $\funcseqname{f}$.
\begin{remark}
The trivial model space $\Triv{S}$ is a special case.
\end{remark}
\end{definition}

\begin{lemma}[Minimal model spaces are model spaces]
\label{lem:minmod}
Let $S$ be a topological space, $\seqname{O}$ a set of open sets in $S$
and $\funcseqname{f}$ a set of continuous functions between elements of
$\seqname{O}$. Then
$(C,\catname{C}) \defeq \minimal{\Mod}(S, \seqname{O}, \funcseqname{f})$
is a model space.

\begin{proof}
$(C,\catname{C})$ satisfies the conditions of \pagecref{def:model}
\begin{enumerate}
\item Finite intersections of open sets are open and
$C = \union[U \in \seqname{O}]{U}$ by construction
\item $\Ob(\catname{C})$ is closed under finite intersections by construction
\item Compositions of continuous functions, inclusion maps and
restrictions of continuous functions are continuous
\item Restrictions of morphisms are morphisms by construction
\item Inclusion maps are morphisms by construction
\item The restricted sheaf condition holds by construction
\end{enumerate}
\end{proof}
\end{lemma}

\subsection{M-paracompact model spaces}
\label{sub:m-para}
Paracompactness is an important property for topological spaces because
of partitions of unity. There is an analogous property for model
spaces.

\begin{definition}[Model topology]
\label{def:ModTop}
Let $\seqname{M} \defeq (S, \catname{S})$ be a model space for $S$.
Then the model topology $\topname{M}^*$ for $M$ is the topology
generated by $\Ob(\catname{S})$.
\begin{remark}
$\topname{M}^*$ is not guarantied to be T0 even if $S$ is T4. However,
$\topname{M}^*$ may be normal or regular even if $S$ is not.
\end{remark}
\end{definition}

\begin{definition}[m-paracompactness]
\label{def:M-para}
A model space $\seqname{M} \defeq (S, \catname{S})$ is \\
m-paracompact iff $\topname{M}^*$ is regular and every cover of $S$ by
model neighborhoods has a locally finite refinement by model neighborhoods.
\begin{remark}
This is a stronger condition than merely requiring $\topname{M}^*$ to
be paracompact.
\end{remark}
\end{definition}

\begin{theorem}[m-paracompactness and paracompactness]
If $\seqname{M} = (S, \catname{S})$ is m-paracompact and the model
neighborhoods form a basis for $S$ then $S$ is paracompact.
\begin{proof}
Let $R$ be an open cover of $S$.  Since the model neighborhoods form a
basis, every set in $R$ is a union of model neighborhoods and thus there
is a refinement $R_1$ by model neighborhoods.  Since by hypothesis
$\seqname{M}$ is m-paracompact, $R_1$ has a locally finite refinement
$R_2$ by model neighborhoods.  Since model neighborhoods are open sets,
$R_2$ is a locally finite refinement of $R$ in the conventional sense.
\end{proof}
\end{theorem}

\subsection{Model functions and model categories}
\label{sub:modcat}
It is convenient to have a notion of mappings between model spaces that
are well behaved in some sense, e.g., fiber preserving; using that
notion it is then possible to group model spaces into categories.

\begin{definition}[Model functions]
Let $\seqname{M}^i \defeq (S^i, \catname{S}^i)$, $i=1,2$,
be a model space and
$\funcname{f} \maps S^1 \to S^2$ be a continuous function.
$\funcname{f}$ is a model function iff the inverse images of model
neighborhoods are model neighborhoods and the images of model
neighborhoods are contained in model neighborhoods.

\begin{equation}
\uquant%
  {V \objin \catname{S}^2}
  {\funcname{f}^{-1}[V] \objin \catname{S}^1}
\end{equation}

\begin{equation}
\uquant%
  {U \objin \catname{S}^1}
  {
    \equant {V \objin \catname{S}^2} {\funcname{f}[U] \subseteq V}
  }
\end{equation}

By abuse of language we write
$\funcname{f} \maps \seqname{M}^1 \to \seqname{M}^2$ both for
$\funcname{f}$ considered as a model function and for $\funcname{f}$
considered as a continuous function.

Let $\seqname{M}'^i \defeq (S'^i, \catname{S}'^i)$
be model spaces,
$\seqname{M}'^i \submod \seqname{M}^i$ and
$\funcname{f}[S'^1] \subseteq S'^2$. Then by abuse of language we write
$\funcname{f} \maps \seqname{M}'^1 \to \seqname{M}'^2$ for
$\funcname{f} \restriction_{S'^1,S'^2}$.
\end{definition}

\begin{lemma}[Composition of model functions]
\label{lem:modcomp}
Let $\seqname{M}_i \defeq (S_i,\catname{S}_i)$, $i=1,2,3$, be model
spaces and $\funcname{f}_1 \maps \seqname{M}_1 \to \seqname{M}_2$,
$\funcname{f}_2 \maps \seqname{M}_2 \to \seqname{M}_3$ model functions.
Then $\funcname{f}_2 \compose \funcname{f}_1$ is a model function.
\begin{proof}
Let $U_i \objin \catname{S}_i,$ $i=1,2,3.$
Since $\funcname{f}_2$ is a model function,
$\funcname{f}^{-1}_2[U_3] \objin \catname{S}_2$.
Since $\funcname{f}_1$ is a model function,
$
  (\funcname{f}_2 \compose \funcname{f}_1)^{-1}[U_3]
  =
  \funcname{f}^{-1}_1[\funcname{f}^{-1}_2[U_3]]
  \objin
  \catname{S}_1
$.
Since $\funcname{f}_1$ is a model functions, $\funcname{f}_1[U_1]$
is contained in a model neighborhood $V_2$.
Since $\funcname{f}_2$ is a model functions, $\funcname{f}_2[V_2]$
is contained in a model neighborhood $V_3$.
Then $(\funcname{f}_2 \compose \funcname{f}_1)[U_1] \subseteq V_3$.
\end{proof}
\end{lemma}

\begin{definition}[Model homeomorphisms]
Let $\seqname{M}^1 \defeq (S^1, \catname{S}^1)$ and \\
$\seqname{M}^2 \defeq (S^2, \catname{S}^2)$ be model spaces and
$\funcname{f} \maps S^1 \to S^2$ be a model function. $\funcname{f}$ is
a model homeomorphism iff it is also invertible and its inverse is a
model function.
\end{definition}

\begin{definition}[Model categories]
\label{def:modcat}
A category $\catname{M}$ is a model category iff

\begin{enumerate}
\item the objects of $\catname{M}$ are model spaces
\item the morphisms of $\catname{M}$ are model functions.
\item composition is functional composition.
\item Every model subspace of an object in $\catname{M}$ is in
$\catname{M}$ and the inclusion map is a morphism.
\item
\label{modcat:rest}
If $\seqname{S}^i \defeq ({S}^i, \catname{S}^i)$, $i=1,2$, are
subspaces of $\seqname{S} \defeq (S, \catname{S}) \objin \catname{M}$
and $\funcname{f} \maps S^1 \to S^2$ is a morphism of $\catname{S}$ then
$\funcname{f} \maps \seqname{S}^1 \to \seqname{S}^2$ is a morphism of
$\catname{M}$.
\end{enumerate}
\end{definition}

\begin{definition}[Trivial model categories]
Let $\seqname{M}$ be a set of model spaces. Then
$\Triv{\Mod}(\seqname{M})$ is the category of all model functions
between model subspaces of model spaces in $\seqname{M}$
\end{definition}

\begin{definition}[Local m-morphisms]
\label{def:modlocal}
Let $\catname{M}^i$, $i=1,2$, be model categories and
$\seqname{S}^i \defeq (S^i, \catname{S}^i) \objin \catname{M}^i$.

A continuous function $\funcname{f} \maps S^1 \to S^2$ is locally an
m-morphism of $\seqname{S}^1$ to $\seqname{S}^2$ iff
$\seqname{S}^1 \submod \seqname{S}^2$ and for every $u \in S^i$ there is
a model neighborhood $U_u$ for $u$ and a model neighborhood $V_u$ for
$v \defeq \funcname{f}(u)$ such that $\funcname{f}[U_u] \subseteq V_u$
and $\funcname{f} \maps U_u \to V_u$ is a morphism of $\catname{S}^2$.

A continuous function $\funcname{f} \maps S^i \to S^i$ is locally an
m-morphism of $\seqname{S}^i$ iff it is locally an m-morphism of
$\seqname{S}^i$ to $\seqname{S}^i$.

A continuous function $\funcname{f} \maps S^1 \to S^2$ is locally an
$\catname{M}^1$-$\catname{M}^2$ m-morphism of $\seqname{S}^1$ to
$\seqname{S}^2$ iff $\catname{M}^1 \subcat[full-] \catname{M}^2$ and for
every $u \in S^1$ there is a model neighborhood $U_u$ for $u$ and a
model neighborhood $V_u$ for $v \defeq \funcname{f}(u)$ such that
$\funcname{f}[U_u] \subseteq V_u$ and $\funcname{f} \maps U_u \to V_u$
is a morphism of $\seqname{M}^2$.
\end{definition}

\begin{lemma}[Local m-morphisms]
\label{lem:modlocal}
Let $\catname{M}^i$, $i=1,2,3$, be model categories and
$\seqname{S}^i \defeq (S^i, \catname{S}^i) \objin \catname{M}^i$.

If $\funcname{f}^i_j \maps S^i \to S^i$, $j=1,2$, is locally an
m-morphism of $\seqname{S}^i$ then $\funcname{f}^i_2 \compose
\funcname{f}^i_1 \maps S^i \to S^i$ is locally an m-morphism of
$\seqname{S}^i$.

\begin{proof}
Let $u \in S^i$, $v \defeq \funcname{f}^i_1(u)$ and
$w \defeq \funcname{f}^i_2(v)$. There exist a model neighborhood $U_u$
for $u$, model neighborhoods $V_u$, $V'_u$ for $v$ and a model
neighborhood $W_v$ of $w$ such that
$\funcname{f}^i_1[U_u] \subseteq V'_u$,
$\funcname{f}^i_2[V_v] \subseteq W_v$,
$\funcname{f}^i_1 \maps U_u \to V'_u$ is a morphism of $\catname{S}^i$
and $\funcname{f}^i_2 \maps V_v \to W_v$ is a morphism of
$\catname{S}^i$.  Then $\hat{V_u} \defeq V_v \cap V'_u \neq \emptyset$,
$\hat{V_u}$ is a model neighborhood of $v$ and
$\hat{U_u} \defeq \funcname{f}^{i-1}_1[\hat{V_u}]$ is a model
neighborhood for $u$.  $\funcname{f}^i_1 \maps \hat{U_u} \to \hat{V_u}$
and $\funcname{f}^i_2 \maps \hat{V_u} \to W_v$ are morphisms of
$\catname{S}^i$ by
{
  \showlabelsinline
  \cref{mod:restriction}
}
of \pagecref{def:model} and thus
$\funcname{f}^i_2 \compose \funcname{f}^i_1 \maps \hat{U_u} \to W_v$ is a
morphism of $\catname{S}^i$.
\end{proof}

If $S$ is a model neighborhood of $\seqname{S} \defeq (S,\catname{S})$
then every function $\funcname{f} \maps S \to S$ that is locally an
m-morphism of $\seqname{S}$ is a morphism of $\catname{S}$.

\begin{proof}
For every $u\in S$ there is a model neighborhood $U_u$ for $u$ and a
model neighborhood $V_u$ for $v \defeq \funcname{f}(u)$ such that
$\funcname{f}[U_u] \subseteq V_u$ and $\funcname{f} \maps U_u \to V_u$
is a morphism of $\catname{S}$. $\Id_{V_u,S} \arin \catname{S}$ so
$\Id_{V_u,S} \compose (\funcname{f} \maps U_u \to V_u)$ is a morphism of
$\catname{S}$. Since $\union[{u \in S}]{U_u} = S$ and
$\union[{u \in S}]{S}$ are model neighborhoods, the result follows by
{
  \showlabelsinline
  \cref{mod:sheaf}
}
of \pagecref{def:model}.
\end{proof}

If each $\funcname{f}^i \maps \seqname{S}^i \to \seqname{S}^{i+1}$ is
locally an $\catname{M}^i$-$\catname{M}^{i+1}$ m-morphism of
$\seqname{S}^1$ to $\seqname{S}^2$ then
$\funcname{f}^2 \compose \funcname{f}^1 \maps S^1 \to S^3$ is locally an
$\catname{M}^1$-$\catname{M}^3$ m-morphism of $\seqname{S}^1$ to
$\seqname{S}^3$.

\begin{proof}
Since $\catname{M}^1 \subcat[full-] \catname{M}^2$ and
$\catname{M}^2 \subcat[full-] \catname{M}^3$,
$\catname{M}^1 \subcat[full-] \catname{M}^3$.
Let $u \in S^1$, $v \defeq \funcname{f}^1(u)$ and
$w \defeq \funcname{f}^2(v)$. There exist a model neighborhood $U_u$
for $u$, model neighborhoods $V_v$, $V'_u$ for $v$ and a model
neighborhood $W_v$ of $w$ such that
$\funcname{f}^1[U_u] \subseteq V'_u$,
$\funcname{f}^1 \maps U_u \to V'_u$ is a morphism of $\catname{M}^2$,
$\funcname{f}^2[V_v] \subseteq W_v$ and
$\funcname{f}^2 \maps V_v \to W_v$ is a morphism of $\catname{M}^3$.
Then $\hat{V_u} \defeq V_v \cap V'_u \neq \emptyset$, $\hat{V_u}$ is a
model neighborhood of $v$ and
$\hat{U_u} \defeq \funcname{f}^{i-1}_1[\hat{V_u}]$ is a model
neighborhood for $u$.  $\funcname{f}^1 \maps \hat{U_u} \to \hat{V_u}$
and $\funcname{f}^2 \maps \hat{V_u} \to W_v$ are morphisms of
$\catname{M}^3$ by
{
  \showlabelsinline
  \cref{modcat:rest}
}
of \pagecref{def:modcat} and thus
$\funcname{f}^2 \compose \funcname{f}^1 \maps \hat{U_u} \to W_v$ is a
morphism of $\catname{M}^3$.
\end{proof}

If $\seqname{S}^i \submod \seqname{S}^{i+1}$ then $\Id_{S^i,S^{i+1}}$ is
locally an m-morphism of $\seqname{S}^i$.  If each ${S}^i$ is a model
neighborhood of $\seqname{S}^i$ then $\Id_{S^i,S^{i+1}}$ is a
morphism of $\catname{S}^{i+1}$.

\begin{proof}
Let $u \in S^1$ and $U_u$ a model neighborhood of $\seqname{S}^i$ for
$u$. Since $\seqname{S}^i \submod \seqname{S}^{i+1}$, $U_u$ is also a
model neighborhood of $\seqname{S}^{i+1}$ for $u$, and hence
$\Id_{U_u} \arin \catname{S}^{i+1}$.

If each ${S}^i$ is a model neighborhood of $\seqname{S}^i$ then each
$\Id_{S^i}$ is a morphism of $\catname{S}^i$. Since
$\seqname{S}^i \submod \seqname{S}^{i+1}$,
${S}^i \objin \catname{S}^{i+1}$ and
$\Id_{U_u,S^{i+1}}$ is a morphism of $\catname{S}^{i+1}$ by
{
  \showlabelsinline
  \cref{modcat:rest}
}
of \pagecref{def:modcat}.
Since $\union[{u \in S^1}]{U_u} = S^1$ and
$\union[{u \in S^1}]{S}^2$ are model neighborhoods, the result follows
by
{
  \showlabelsinline
  \cref{mod:sheaf}
}
of \pagecref{def:model}.
\end{proof}

If $\catname{M}^i \subcat \catname{M}^{i+1}$ and
$\seqname{S}^i \submod \seqname{S}^{i+1}$ then $\Id_{S^i,S^{i+1}}$ is a
morphism of $\catname{M}^{i+1}$.

\begin{proof}
$\seqname{S}^i \objin \catname{M}^{i+1}$ because
$\catname{M}^i \subcat \catname{M}^{i+1}$.
$\Id_{S^{i+1}} \arin \catname{M}^{i+1}$. Then
$\Id_{S^i,S^{i+1}} \arin \catname{M}^{i+1}$ by
{
  \showlabelsinline
  \cref{modcat:rest}
}
of \pagecref{def:modcat}
\end{proof}
\end{lemma}

\begin{corollary}[Local m-morphisms]
\label{cor:modlocal}
Let $\catname{M}^i$, $i=1,2$, be model categories and
$\seqname{S}^i \defeq (S^i, \catname{S}^i) \objin \catname{M}^i$.

If $S^i$ is a model neighborhood of $\seqname{S}^i$ then
$\Id_{\seqname{S}^i}$ is a morphism of $\catname{S}^i$.

\begin{proof}
$\seqname{S}^i \submod \seqname{S}^i$.
\end{proof}
\end{corollary}

\subsection{Spaces and proper functions}
\label{sub:spaces}
Several of the following definitions involve spaces of a restrictive
character and specific types of mappings among them. Except where
otherwise qualified, the word \textit{space} will have this restricted
meaning.

\begin{definition}
A space is a topological space, a model space or either with an
additional associated structure.
\end{definition}

\begin{definition}
Let $S^1$, $S^2$ be spaces and $\funcname{f} \maps S^1 \to S^2$ a
continuous function; $\funcname{f}$ need not preserve any associated
algebraic structure\footnote{
  However, in practice commutation relations will often enforce the
  preservation of algebraic structures.
}.
$\funcname{f}$ is a proper\footnote{
  In the spirit of \cite[footnote, p.~112] {GenTop}
  \begin{quote}
    This nomenclature is an excellent example of the time-honored
    custom of referring to a problem we cannot handle as abnormal,
    irregular, improper, degenerate, inadmissible, and otherwise
    undesirable.
  \end{quote}
}
function iff

\begin{enumerate}
\item $S^1$ and $S^2$ are both $\truthspace$ or both
$\seqname{Truthspace}$, and $\funcname{f}(\true)=\true$.

\item $S^1$, $S^2$ are topological spaces other than
$\truthspace$.

\item $S^1$ is a topological space other than
$\truthspace$, $S^2$ is a model space other than
$\seqname{Truthspace}$ and the images of open sets are contained in
model neighborhoods.

\item $S^1$ is a model space other than $\seqname{Truthspace}$, $S^2$ is
a topological spaces other than $\truthspace$ and the
inverse images of open sets are model neighborhoods.

\item $S^1$, $S^2$ are both model spaces other than
$\seqname{Truthspace}$ and $\funcname{f}$ is a model function.
\end{enumerate}
\end{definition}

\begin{definition}
\label {def:singcat}
Let $S$ be a topological space. Then the singleton category of $S$,
abbreviated $\singcat{S}$, is the category whose sole object is $S$ and
whose morphisms are all of the continuous functions from $S$ to itself.

Let $\seqname{S}$ be a model space. Then the singleton category of
$\seqname{S}$, abbreviated $\singcat{\seqname{S}}$, is the category whose
sole object is $\seqname{S}$ and whose morphisms are all of the model
functions from $\seqname{S}$ to itself that are locally morphisms of
$\seqname{S}$.

Let $\seqname{S} \defeq (S^\alpha, \alpha \prec \Alpha)$ be a sequence
of spaces. Then the singleton category sequence of $\seqname{S}$,
abbreviated $\Singcat{\seqname{S}}$, is
$(\singcat{S^\alpha}, \alpha \prec \Alpha)$.

Let $\seqname{S}$ be a set of spaces. Then the singleton category of
$\seqname{S}$, abbreviated $\Singcat{\seqname{S}}$, is
$\union[S \in \seqname{S}]{\singcat{S}}$.

\begin{remark}
By abuse of language the notation $\Singcat{\catseqname{S}}$ will be used
to name sequences of categories constructed with this and similar
functions.
\end{remark}
\end{definition}

%

\section{Signatures}
\label{sec:sig}
Given a sequences of spaces, we need a way to characterise a function
that takes arguments in those spaces and has a value in them. For this
purpose we use a sequence of ordinals where the last ordinal in the
sequence identifies the space containing the function's value.

\begin{definition}[Signature]
Let $\seqname{S} \defeq (S_\alpha, \alpha \prec \Alpha)$ be a sequence
and $\sigma \defeq (\alpha_\beta \prec \Alpha, \beta \preceq \Beta)$ be
a sequence of ordinals in $\Alpha$. Then $\sigma$ is a signature over
$\seqname{S}$.
\end{definition}

\begin{definition}[$\seqname{S}$-signature of function]
Let $\seqname{S} \defeq (S_\alpha, \alpha \prec \Alpha)$ be a
sequence,
$\sigma \defeq (\alpha_\beta \prec \Alpha, \beta \preceq \Beta)$ be a
signature over $\seqname{S}$ and
$\funcname{f} \maps (S_{\alpha_\beta}, \beta \prec \Beta) \to
S_{\alpha_\Beta}$ continuous.
Then $\sigma$ is the $\seqname{S}$-signature for $\funcname{f}$ and
$\funcname{f}$ has $\seqname{S}$-signature $\sigma$.

\begin{remark}
If $n \in \mathbb{N}$ and $\sigma \defeq (\alpha, \beta \preceq n)$
then $\funcname{f}$ is an n-ary operation over $S_\alpha$.
\end{remark}
$\funcname{f}$ is proper iff either
\begin{enumerate}
\item
$S_{\alpha_\beta}$, $\beta \preceq \Beta)$ are all topological spaces.
\item $S_{\alpha_\beta}$, $\beta \preceq \Beta)$ are all model spaces.
\end{enumerate}
\end{definition}

\begin{definition}[$\catseqname{S}$-signature of function]
Let $\seqname{S} \defeq (S_\alpha, \alpha \prec \Alpha)$ be a
sequence,
$\catseqname{S} \defeq (\catname{S}_\alpha, \alpha \prec \Alpha)$ be a
sequence of categories with
$
  \uquant%
    {\alpha \prec \Alpha}
    {S_\alpha \objin \catname{S}_\alpha}
$,
$\sigma \defeq (\alpha_\beta \prec \Alpha, \beta \preceq \Beta)$ be a
signature over $\catseqname{S}$ and
$\funcname{f} \maps (S_{\alpha_\beta}, \beta \prec \Beta) \to
S_{\alpha_\Beta}$ continuous.
Then $\sigma$ is the $\catseqname{S}$-signature for $\funcname{f}$ and
$\funcname{f}$ has $\catseqname{S}$-signature $\sigma$.
\end{definition}

\begin{definition}[$\seqname{S}$-signature of function sequence]
Let \\
$\Sigma \defeq (\sigma_\gamma, \gamma \prec \Gamma) \defeq
\bigl ( (\sigma_{\gamma,{\alpha_\beta}}, \alpha_\beta \prec \Alpha,
\beta \preceq \Beta_\gamma), \gamma \prec \Gamma \bigr )$
be a sequence of signatures over
$\seqname{S} \defeq (S_\alpha, \alpha \prec \Alpha)$,
$\funcseqname{F}=(\funcname{F}_\gamma, \gamma \prec \Gamma)$ a
sequence of functions where $\funcname{F}_\gamma$ has
$\seqname{S}$-signature $\sigma_\gamma$, then $\Sigma$ is the
$\seqname{S}$-signature for $\funcseqname{F}$ and
$\funcseqname{F}$ has $\seqname{S}$-signature $\Sigma$.
\end{definition}

\begin{definition}[$\catseqname{S}$-signature of function sequence]
Let \\
$\Sigma \defeq (\sigma_\gamma, \gamma \prec \Gamma) \defeq
\bigl ( (\sigma_{\gamma,{\alpha_\beta}},
\alpha_\beta \prec \Alpha,
\beta \preceq \Beta_\gamma), \gamma \prec \Gamma \bigr )$
be a sequence of signatures over
$\catseqname{S} \defeq (\catname{S}_\alpha, \alpha \prec \Alpha)$,
$\funcseqname{F}=(\funcname{F}_\gamma, \gamma \prec \Gamma)$
a sequence of functions where $\funcname{F}_\gamma$ has
$\catseqname{S}$-signature $\sigma_\gamma$,
then $\Sigma$ is the
$\catseqname{S}$-signature for $\funcseqname{F}$ and
$\funcseqname{F}$ has $\catseqname{S}$-signature $\Sigma$.
\end{definition}

\subsection{Functions associated with signatures and function sequences}
\label{sub:sig}
When characterizing functions with signatures, several auxillary
functions are helpful.

\begin{definition}[Signature operators]
\label{def:sigop}
Let $\catseqname{S}^i \defeq (\catname{S}^i_\alpha, \alpha \prec \Alpha)$, $i=1,2$,
be sequences of categories with $\seqname{S}^i \defeq (S^i_\alpha, \alpha \prec
\Alpha)$ sequences of sets in those categories, $\funcname{f}^i,i=1,2$ functions with
$\catseqname{S}^i$-signatures $\sigma = (\alpha_\beta \prec \Alpha, \beta \preceq \Beta)$,
$\funcseqname{g} \defeq (\funcname{g}_\alpha \maps S^1_\alpha
\to S^2_\alpha, \alpha \prec \Alpha)$
a sequence of functions then
\begin{equation}
\langle \sigma | \seqname{S}^i \rangle
\defeq \overset{\sigma}{\seqname{S}^i} \defeq (S^i_{\alpha_\beta}, \beta \preceq \Beta)
\end{equation}
selects sets from $\seqname{S}^i$ according to signature $\sigma$.
\begin{equation}
\langle \sigma | \funcseqname{g} \rangle
\defeq \overset{\sigma}{\funcseqname{g}} \defeq (\funcname{g}_{\alpha_\beta}, \beta \preceq \Beta)
\end{equation}
selects functions from $\funcseqname{g}$ according to signature $\sigma$.
\begin{equation}
\overset{\sigma}{\underline{\funcseqname{g}}} \maps
\bigtimes_{\beta \prec \Beta} S^1_{\alpha_\beta} \to
\bigtimes_{\beta \prec \Beta} S^2_{\alpha_\beta}
= \underline{\langle \sigma | \funcseqname{g} \rangle} =
\bigtimes_{\beta \prec \Beta} \funcname{g}_{\alpha_\beta}
\end{equation}
is the product of the functions selected by all but the last ordinal in signature $\sigma$
\begin{equation}
\funcseqname{g} \composet[\sigma] \funcname{f}^1 \defeq
\tail(\langle \sigma | \funcseqname{g} \rangle) \compose \funcname{f}^1
\end{equation}
composes the last function selected from $\funcseqname{g}$ by signature
$\sigma$ with $\funcname{f}^1$.
\begin{equation}
\funcname{f}^2 \composeh[\sigma] \funcseqname{g} \defeq \funcname{f}^2 \compose
\underline{\overset{\sigma}{\funcseqname{g}}}
\end{equation}
composes $\funcname{f}$ with all but the last function in
$\funcseqname{g}$ selected by signature $\sigma$.
\end{definition}

\begin{figure}
\[ \bfig
\node s1ab(0,0)[\underline{\overset{\sigma}{\seqname{S}^1}}]
\node s2ab(0,-1000)[\underline{\overset{\sigma}{\seqname{S}^2}}]
\node s2aB(1000,-1000)[\tail(\overset{\sigma}{\seqname{S}^2})]
\arrow |l|[s1ab`s2ab;%
  {%
    \underline{\overset{\sigma}{\funcseqname{g}}} \defeq%
    \bigtimes \head(\overset{\sigma}{\funcseqname{g}})
  }]
\arrow |r|[s1ab`s2aB;{{\funcname{f}_2 \composeh[\sigma] \funcseqname{g}}}]
\arrow |r|[s2ab`s2aB;\funcname{f}_2]
\efig \]
\caption{$\funcname{f}_2 \composeh[\sigma] \funcseqname{g}$}
\label{fig:fgsigma}
\end{figure}

\begin{figure}
\[ \bfig
\node s1ab(0,0)[\underline{\overset{\sigma}{\seqname{S}^1}}]
\node s1aB(0,-1000)[\tail(\overset{\sigma}{\seqname{S}^1})]
\node s2aB(1000,-1000)[\tail(\overset{\sigma}{\seqname{S}^2})]
\arrow |l|[s1ab`s1aB;\funcname{f}_1]
\arrow |r|[s1ab`s2aB;{{\funcseqname{g} \composet[\sigma] \funcname{f}_1}}]
\arrow |b|[s1aB`s2aB;\tail(\overset{\sigma}{\funcseqname{g}})]
\efig \]
\caption{$\funcseqname{g} \composet[\sigma] \funcname{f}_1$}
\label{fig:gsigmaf}
\end{figure}

\subsection{Commutative diagrams with signatures}
\label{sub:cdsig}
The conventional usage of the word commutative diagram is for functions
of a single arguement. It is convenient to extend the nomenclature to
functions of multiple arguments and to sequences of such functions. This
is convenient, e.g., as a means of specifying that functions preserve
algebraic structures.

\begin{definition}[$\Sigma$-commutation]
\label{def:Sigmacomm}
Let
$\catseqname{S}^i \defeq (\catname{S}^i_\alpha, \alpha \prec \Alpha)$,
$i=1,2$, be a sequence of categories, with
$
  \seqname{S}^i \defeq (S^i_\alpha, \alpha \prec \Alpha)
  \seqin \catseqname{S}^i
$
a sequence of spaces,
$\funcseqname{F}^i \defeq (\funcname{F}^i_\gamma, \gamma \prec \Gamma)$
be a sequence of functions with
$\catseqname{S}^i$-signature
$
  \Sigma \defeq (\sigma_{\gamma \prec \Gamma}) \defeq \\
    \bigl (
      (
        \sigma_{\gamma,\alpha_\beta}, \alpha_\beta \prec \Alpha,
        \beta \preceq \Beta_\gamma
      ),
      \gamma \prec \Gamma
    \bigr )
$,
$
  \funcseqname{f} \defeq
  (
    \funcname{f}_\alpha \maps S^1_\alpha \to S^2_\alpha,
    \alpha \prec \Alpha
  )
$
a sequence of functions satisfying
$
  \funcseqname{f} \composet[\sigma_\gamma] \funcname{F}^1_\gamma
  = \funcname{F}^2_\gamma \composeh[\sigma_\gamma] \funcseqname{f},
  \gamma \prec \Gamma
$,
then
$\funcseqname{f}$ $\Sigma$-commutes with the function sequences
$\funcseqname{F}^1, \funcseqname{F}^2$.
\end{definition}

\begin{remark}
If all of the $\sigma_{\gamma,\alpha_\beta}$ are equal to the same
ordinal $\hat{\alpha}$ and each $\Beta_\gamma$ is finite then
$\set{{(S^i_{\hat{\alpha}},\funcname{F}^i_\gamma)}}[{\gamma \prec \Gamma}]$ is
an $\Omega$-algebra and $\funcseqname{f}$ is an $\Omega$-homomorphism.
\end{remark}

\begin{figure}
\[ \bfig
\node s1b(0,0)[\underline{\overset{\sigma_\gamma}{\seqname{S}^1_\gamma}}]
\node s1B(1000,0)[\tail \left ( \overset{\sigma_\gamma}{\seqname{S}^1_\gamma} \right )]
\node s2b(0,-1000)[\underline{\overset{\sigma_\gamma}{\seqname{S}^2_\gamma}}]
\node s2B(1000,-1000)[\tail \left ( \overset{\sigma_\gamma}{\seqname{S}^2_\gamma} \right )]
\arrow |a|[s1b`s1B;F^1_\gamma]
\arrow |l|[s1b`s2b;\underline{\overset{\sigma_\gamma}{\funcseqname{f}}} ]
\arrow |r|[s1B`s2B;\tail \left ( \overset{\sigma_\gamma}{\funcseqname{f}} \right )]
\arrow |a|[s2b`s2B;F^2_\gamma]
\efig \]
\caption{$\funcseqname{f}$ $\Sigma$-commutes with the function sequences $\funcseqname{F}^1, \funcseqname{F}^2$}
\label{fig:fF1F2}
\end{figure}

\begin{lemma}[$\Sigma$-commutation]
\label{lem:Sigmacomm}
Let
$
  \catseqname{S}^i \defeq (\catname{S}^i_\alpha, \alpha \prec \Alpha)
$,
$i \in [1,3]$, be a sequence of categories, with
$
  \seqname{S}^i \defeq (S^i_\alpha, \alpha \prec \Alpha)
  \seqin \catseqname{S}^i
$
a sequence of spaces and
$\funcseqname{F}^i \defeq (\funcname{F}^i_\gamma, \gamma \prec \Gamma)$
be a sequence of functions with $\catseqname{S}^i$-signature
$
  \Sigma \defeq (\sigma_\gamma, \gamma \prec \Gamma) \defeq
  \bigl (
    (
      \sigma_{\gamma,\alpha_\beta}, \alpha_\beta \prec \Alpha,
      \beta \preceq \Beta_\gamma
    ),
    \gamma \prec \Gamma
  \bigr )
$.

Let
$
  \funcseqname{f}^i \defeq
  (
    \funcname{f}^i_\alpha \maps S^i_\alpha \to S^{i+1}_\alpha,
    \alpha \prec \Alpha
  )
$,
$i=1,2$, $\Sigma$-commute with $\funcseqname{F}^i,\funcseqname{F}^{i+1}$.
Then $\funcseqname{f}^2 \compose[()] \funcseqname{f}^1$ $\Sigma$-commutes
with $\funcseqname{F}^1,\funcseqname{F}^3$.

\begin{figure}
\[ \bfig
\node s1b(0,2000)[\underline{\overset{\sigma_\gamma}{\seqname{S}^1_\gamma}}]
\node s1B(1000,2000)[\tail \left ( \overset{\sigma_\gamma}{\seqname{S}^1_\gamma} \right )]
\node s2b(0,1000)[\underline{\overset{\sigma_\gamma}{\seqname{S}^2_\gamma}}]
\node s2B(1000,1000)[\tail \left ( \overset{\sigma_\gamma}{\seqname{S}^2_\gamma} \right )]
\node s3b(0,0)[\underline{\overset{\sigma_\gamma}{\seqname{S}^3_\gamma}}]
\node s3B(1000,0)[\tail \left ( \overset{\sigma_\gamma}{\seqname{S}^3_\gamma} \right )]
\arrow |a|[s1b`s1B;F^1_\gamma]
\arrow |l|[s1b`s2b;\underline{\overset{\sigma_\gamma}{\funcseqname{f}^1}} ]
\arrow |r|[s1B`s2B;\tail \left ( \overset{\sigma_\gamma}{\funcseqname{f}^1} \right )]
\arrow |a|[s2b`s2B;F^2_\gamma]
\arrow |l|[s2b`s3b;\underline{\overset{\sigma_\gamma}{\funcseqname{f}^2}} ]
\arrow |r|[s2B`s3B;\tail \left ( \overset{\sigma_\gamma}{\funcseqname{f}^2} \right )]
\arrow |a|[s3b`s3B;F^3_\gamma]
\efig \]
\caption{$\funcseqname{f}^2 \compose[()] \funcseqname{f}^1$ $\Sigma$-commute with $\funcseqname{F}^1,\funcseqname{F}^3$}
\label{fig:fF1F2F3}
\end{figure}

\begin{proof}
The result follows by diagram chasing in \cref{fig:fF1F2F3}.
\end{proof}

Let
$\seqname{S}^1 \SUBSETEQ \seqname{S}^2$
and
$
  \funcname{F}^1_\gamma =
  \funcname{F}^2_\gamma \maps
    \bigtimes \underline{\overset{\sigma_\gamma}{\seqname{S}^1}} \to
    \tail \left ( \overset{\sigma_\gamma}{\seqname{S}^1} \right )
$.
Then $\ID_{\seqname{S}^1,\seqname{S}^2}$ $\Sigma$-commutes with
$\funcseqname{F}^1,\funcseqname{F}^2$.

\begin{proof}
The result follows from \cref{def:sigop,def:Sigmacomm} on
\cpageref{def:sigop,def:Sigmacomm}: \\
$
    \ID_{\seqname{S}^1,\seqname{S}^2} \composet[\sigma_\gamma]
    \funcname{F}^1_\gamma
  =
    \Id_%
      {
        S^1_{\gamma,\alpha_{\Beta_\gamma}},
        S^2_{\gamma,\alpha_{\Beta_\gamma}}
      }
    \compose
    \funcname{F}^1_\gamma
  = \funcname{F}^2_\gamma \restriction_%
      {
        \bigtimes \underline{\overset{\sigma_\gamma}{\seqname{S}^1}},
      }
  = \funcname{F}^2_\gamma \compose
    \bigtimes_{\beta \prec \Beta_\gamma}
      \Id_{S^1_{\alpha_\beta}, S^2_{\alpha_\beta}}
  = \funcname{F}^2_\gamma \composeh[\sigma_\gamma]
    \ID_{\seqname{S}^1,\seqname{S}^2}
$
\end{proof}
\end{lemma}

\begin{corollary}[$\Sigma$-commutation]
\label{cor:Sigmacomm}
Let
$
  \catseqname{S} \defeq (\catname{S}_\alpha, \alpha \prec \Alpha)
$,
be a sequence of categories, with
$
  \seqname{S} \defeq (S_\alpha, \alpha \prec \Alpha)
  \seqin \catseqname{S}
$
a sequence of spaces and
$\funcseqname{F} \defeq (\funcname{F}_\gamma, \gamma \prec \Gamma)$
be a sequence of functions with $\catseqname{S}$-signature
$
  \Sigma \defeq (\sigma_\gamma, \gamma \prec \Gamma) \defeq
  \bigl (
    (
      \sigma_{\gamma,\alpha_\beta}, \alpha_\beta \prec \Alpha,
      \beta \preceq \Beta_\gamma
    ),
    \gamma \prec \Gamma
  \bigr )
$.
Then $\ID_{\seqname{S}}$ $\Sigma$-commutes with
$\funcseqname{F},\funcseqname{F}$.

\begin{proof}
$\seqname{S} \SUBSETEQ \seqname{S}$,
$
  \funcname{F}_\gamma =
  \funcname{F}_\gamma \maps
    \bigtimes \underline{\overset{\sigma_\gamma}{\seqname{S}}} \to
    \tail \left ( \overset{\sigma_\gamma}{\seqname{S}} \right )
$
and $\ID_{\seqname{S}} = \ID_{\seqname{S},\seqname{S}}$.
\end{proof}
\end{corollary}

\section{Prestructures}
\label{sec:pre}
Prestructures and prestructure morphisms are generalizations of
$\Omega$-algebras and $\Omega$-homomorphisms, and are notational
conveniences to simplify imposing commutation relations on multiple
unrelated functions of multiple variables.

\begin{definition}[Prestructures]
\label{def:pre}
Let $\catseqname{S} \defeq (\catname{S}_\alpha, \alpha \prec \Alpha)$ be
a sequence of categories,
$
  \seqname{S} \defeq (S_\alpha, \alpha \prec \Alpha)
$
a sequence of spaces,
$\funcseqname{F}=(F_\gamma, \gamma \prec \Gamma)$
a sequence of continuous functions and
$
  \Sigma \defeq
  (\sigma_\gamma, \gamma \prec \Gamma) \defeq
  \bigl (
    (
      \sigma_{\gamma.{\alpha_\beta}},
      \beta \preceq \Beta_\gamma
    ),
    \gamma \prec \Gamma
  \bigr )
$
a sequence of signatures.
Then $\seqname{P} \defeq (\catseqname{S}, \seqname{S}, \Sigma, F)$ is a
$\catseqname{S}-\Sigma$ prestructure iff
\begin{enumerate}
\item $\catname{S}_\alpha$ is either a topological category or a model category.
\item $\seqname{S} \seqin \catseqname{S}$.
\item $\funcseqname{F}$ has $\catseqname{S}$-signature $\Sigma$.
\end{enumerate}
\end{definition}

\begin{lemma}[Prestructures]
\label{lem:pre}
Let $\catseqname{S}^i$, $i=1,2$, be a sequence of
categories, \\
$\catseqname{S}^1 \SUBCAT \catseqname{S}^2$,
$\seqname{S} \seqin \catseqname{S}^1$,
$\funcseqname{F}$ a sequence of continuous functions with
$\catseqname{S}^1$-signatures $\Sigma$ and
$\seqname{P}^1 \defeq (\catseqname{S}^1, \seqname{S}, \Sigma, F)$ a
$\catseqname{S}^1-\Sigma$ prestructure.  Then
$\seqname{P}^2 \defeq (\catseqname{S}^2, \seqname{S}, \Sigma, F)$ is a
$\catseqname{S}^2-\Sigma$ prestructure.

\begin{proof}
$\seqname{S} \seqin \catseqname{S}^2$ by \pagecref{lem:seqin}. The other
conditions do not depend on the category sequence.
\end{proof}

Let $\catseqname{S}^i$, $i=1,2$, be a sequence of categories,
$\catseqname{S}^1 \SUBCAT \catseqname{S}^2$,
$\seqname{S}^i \seqin \catseqname{S}^i$,
$\seqname{S}^1 \SUBSETEQ \seqname{S}^2$,
$
  \funcseqname{F}^2 \defeq \\
  (\funcname{F}^2_\gamma, \gamma \prec \Gamma)
$
a sequence of continuous functions with $\catseqname{S}^2$-signatures
$
  \Sigma \defeq \\
  (\sigma_\gamma, \gamma \prec \Gamma) \defeq
  \bigl (
    (
      \sigma_{\gamma.{\alpha_\beta}},
      \beta \preceq \Beta_\gamma
    ),
    \gamma \prec \Gamma
  \bigr )
$,
$
  \funcname{F}^2_\gamma
    \left [
      \bigtimes \underline{\overset{\sigma_\gamma}{\seqname{S}^1}}
    \right ]
  \subseteq
  \tail \left ( \overset{\sigma_\gamma}{\seqname{S}^1} \right )
$,
$
  \funcseqname{F}^1 \defeq \\
  \left (
    \funcname{F}^1_\gamma \defeq
        \funcname{F}^2_\gamma \maps
          \bigtimes \underline{\overset{\sigma_\gamma}{\seqname{S}^1}}
          \to
          \tail \left ( \overset{\sigma_\gamma}{\seqname{S}^1} \right ),
        \gamma  \prec \Gamma
  \right )
$
and
$
  \seqname{P}^2 \defeq
  (\catseqname{S}^2, \seqname{S}^2, \Sigma, \funcseqname{F}^2)
$
a $\catseqname{S}^2-\Sigma$ prestructure.
Then
$
  \seqname{P}^1 \defeq
  (\catseqname{S}^1, \seqname{S}^1, \Sigma, \funcseqname{F}^1)
$
is a $\catseqname{S}^1-\Sigma$ prestructure.

\begin{proof}
$\seqname{S}^1 \seqin \catseqname{S}^1$ by hypothesis.
$\funcseqname{F}^1$ is continuous and has $\catseqname{S}^1$-signatures
$\Sigma$ by construction.
\end{proof}
\end{lemma}

\begin{definition}[Morphisms of prestructures]
\label{def:premorph}
Let $\catseqname{S}^i=(\catname{S}^i_\alpha, \alpha \prec \Alpha)$,
$i=1,2$, be a sequence of categories,
$
  \seqname{S}^i \defeq (S^i_\alpha \alpha \prec \Alpha)
  \seqin \catseqname{S}^i
$
be a sequence of spaces,
$\funcseqname{F}^i=(\funcname{F}^i_\gamma, \gamma \prec \Gamma)$
be a sequence of continuous functions with
$\catseqname{S}^i$-signatures $\Sigma$,
$
\funcseqname{f} \defeq \\
  (
    \funcname{f}_\alpha \maps S^1_\alpha \to S^2_\alpha,
    \alpha \prec \Alpha
  )
$
be a sequence of proper functions and
$
  \seqname{P}^i \defeq
  (\catseqname{S}^i, \seqname{S}^i, \Sigma, \funcseqname{F}^i)
$
be a $\catseqname{S}^i-\Sigma$ prestructure. $\funcseqname{f}$ is a
(strict) prestructure morphism of $\seqname{P}^1$ to $\seqname{P}^2$ iff
it $\Sigma$-commutes with $\funcseqname{F}^1,\funcseqname{F}^2$, as
shown in \pagecref{fig:fF1F2}.

It is a semi-strict prestructure morphism iff
$\catseqname{S}^1 \SUBCAT[full-] \catseqname{S}^2$ and for each $\alpha$

\begin{enumerate}
\item Each $\catname{S}^i_\alpha$, $i=1,2$, is a model category:
$\funcname{f}_\alpha$ is locally an
$\catname{S}^1_\alpha$-$\catname{S}^2_\alpha$ m-morphism of $S^1_\alpha$
to $S^2_\alpha$.
\item Each $\catname{S}^i_\alpha$, $i=1,2$, is a topological category:
$\funcname{f}_\alpha$ is locally a
$\catname{S}^1_\alpha$-$\catname{S}^2_\alpha$ morphism of $S^1_\alpha$
to $S^2_\alpha$.
\item Otherwise
$\funcname{f}_\alpha$ is a morphism of $\catname{S}^2_\alpha$.
\end{enumerate}

It is a strict prestructure morphism iff
$\catseqname{S}^1 \SUBCAT[full-] \catseqname{S}^2$ and each function
$\funcname{f}_\alpha$ is a morphism of $\catname{S}^2_\alpha$.

\begin{remark}
It is not sufficient to require that each $\funcname{f}_\alpha$ be a
morphism of $\catname{S}^1_\alpha \unioncat \catname{S}^2_\alpha$
because that would not ensure that the composition of strict
prestructure morphisms is strict.
\end{remark}

If $\seqname{S}^1 \SUBSETEQ \seqname{S}^2$ and
$
  \funcname{F}^1_\gamma =
  \funcname{F}^2_\gamma \maps
    \bigtimes \underline{\overset{\sigma_\gamma}{\seqname{S}^1}} \to
    \tail \left ( \overset{\sigma_\gamma}{\seqname{S}^1} \right )
$
then
$
  \Id_{\seqname{P}^1,\seqname{P}^2} \defeq
  \ID_{\seqname{S}^1,\seqname{S}^2}
$.

$
  \Id_{\seqname{P}^i} \defeq
  \ID_{\seqname{S}^i} = \Id_{\seqname{P}^i,\seqname{P}^i}
$.
\end{definition}

\begin{lemma}[Prestructure morphisms]
\label{lem:premorph}
Let $\catseqname{S}^i=(\catname{S}^i_\alpha, \alpha \prec \Alpha)$,
$i \in [1,4]$, be a sequence of categories,
$
  \seqname{S}^i \defeq (S^i_\alpha, \alpha \prec \Alpha)
  \seqin \catseqname{S}^i
$
a sequence of spaces,
$\funcseqname{F}^i=(\funcname{F}^i_\gamma, \gamma \prec \Gamma)$,
be a sequence of continuous functions with
$\seqname{S}^i$-signatures $\Sigma$ and
$
  \seqname{P}^i \defeq
  (\catseqname{S}^i, \seqname{S}^i, \Sigma, \funcseqname{F}^i)
$
be a $\catseqname{S}^i-\Sigma$ prestructure.

If $\seqname{S}^i \SUBSETEQ \seqname{S}^{i+1}$ and each
$
  \funcname{F}^i_\gamma =
  \funcname{F}^{i+1}_\gamma \maps
    \bigtimes \underline{\overset{\sigma_\gamma}{\seqname{S}^i}} \to
    \tail \left ( \overset{\sigma_\gamma}{\seqname{S}^i} \right )
$
then $\Id_{\seqname{P}^i,\seqname{P}^{i+1}}$ is a prestructure morphism from
$\seqname{P}^i$ to $\seqname{P}^{i+1}$.

$\ID_{\seqname{S}^i,\seqname{S}^{i+1}}$ is a semi-strict prestructure
morphism from $\seqname{P}^i$ to $\seqname{P}^{i+1}$ iff
$\catseqname{S}^i \SUBCAT[full-] \catseqname{S}^{i+1}$.

$\ID_{\seqname{S}^i,\seqname{S}^{i+1}}$ is a strict prestructure
morphism from $\seqname{P}^i$ to $\seqname{P}^{i+1}$ iff each \\
$\Id_{S^1_\alpha,S^{i+1}_\alpha} \arin \catname{S}^{i+1}_\alpha$ and
$\catseqname{S}^i \SUBCAT[full-] \catseqname{S}^{i+1}$.

\begin{proof}
$\ID_{\seqname{S}^i,\seqname{S}^{i+1}}$ $\Sigma$-commutes with
$\funcseqname{F}^i,\funcseqname{F}^{i+1}$ by \pagecref{lem:Sigmacomm}.

The hypotheses for strictness are precisely those of
\fullcref{def:premorph}.
\end{proof}

Let
$
  \funcseqname{f}^i \defeq
  (
    \funcname{f}^i_\alpha \maps S^i_\alpha \to S^{i+1}_\alpha,
    \alpha \prec \Alpha
  )
$,
$i=1,2,3$, be a sequence of proper functions,

If $\catseqname{S}'^i=(\catname{S}'^i_\alpha, \alpha \prec \Alpha)$,
$i \in [1,4]$, is a sequence of categories,
$\catseqname{S}^i \SUBCAT \catname{S}'^i$ and
$
  \seqname{P}'^i \defeq
  (\catseqname{S}'^i, \seqname{S}^i, \Sigma, \funcseqname{F}^i)
$,
then $\funcseqname{f}^i$ is a prestructure morphism of $\seqname{P}'^i$
to $\seqname{P}'^{i+1}$ iff $\funcseqname{f}^i$ is a
prestructure morphism of $\seqname{P}^i$ to
$\seqname{P}^{i+1}$.

If $\catseqname{S}^i \SUBCAT[full-] \catname{S}'^i$ and
$\funcseqname{f}^i$ is a semi-strict (strict) prestructure morphism of
$\seqname{P}^i$ to $\seqname{P}^{i+1}$ then $\funcseqname{f}^i$ is a
semi-strict (strict) prestructure morphism of $\seqname{P}'^i$ to
$\seqname{P}'^{i+1}$.

\begin{proof}
The commutation relations do not depend on the choice of categories.

If $\catseqname{S}^i \SUBCAT[full-] \catname{S}'^i$ and
$\funcseqname{f}^i$ is a semi-strict prestructure morphism of
$\seqname{P}^i$ to $\seqname{P}^{i+1}$ then for each $\alpha$, if

\begin{enumerate}
\item Each $\catname{S}^i_\alpha$, $j=i,i+1$, is a model category:
$\funcname{f}_\alpha$ is locally an
$\catname{S}^i_\alpha$-$\catname{S}^{i+i}_\alpha$ m-morphism of
$S^i_\alpha$ to $S^{i+i}_\alpha$ and thus locally an
$\catname{S}'^i_\alpha$-$\catname{S}'^{i+i}_\alpha$ m-morphism of
$S^i_\alpha$ to $S^{i+i}_\alpha$.

\item Each $\catname{S}^i_\alpha$, $j=i,i+1$, is a topological category:
$\funcname{f}_\alpha$ is locally a
$\catname{S}^i_\alpha$-$\catname{S}^{i+i}_\alpha$ morphism of
$S^i_\alpha$ to $S^{i+i}_\alpha$ and thus locally a
$\catname{S}'^i_\alpha$-$\catname{S}'^{i+i}_\alpha$ morphism of
$S^i_\alpha$ to $S^{i+i}_\alpha$.

\item Otherwise
$\funcname{f}_\alpha$ is a morphism of $\catname{S}^{i+i}_\alpha$ and
thus a morphism of $\catname{S}'^{i+i}_\alpha$.
\end{enumerate}

If $\catseqname{S}^i \SUBCAT[full-] \catname{S}'^i$ and
$\funcseqname{f}^i$ is a strict prestructure morphism of $\seqname{P}^i$
to $\seqname{P}^{i+1}$ then each
$
  \funcseqname{f}^i_\alpha \arin \catseqname{S}^{i+1}_\alpha
  \subcat[full-] \catname{S}'^{i+1}_\alpha
$,
$\alpha \prec \Alpha$, and thus $\funcseqname{f}^i$ is a strict
prestructure morphism of $\seqname{P}'^i$ to $\seqname{P}'^{i+1}$.
\end{proof}

$\funcseqname{f}^{i+1} \compose[()] \funcseqname{f}^i$, $i \in [1,3]$,
is a prestructure morphism and
$
    \funcseqname{f}^3 \compose[()]
    \bigl ( \funcseqname{f}^2 \compose[()] \funcseqname{f}^1 \bigr )
  = \bigr ( \funcseqname{f}^3 \compose[()] \funcseqname{f}^2 \bigr )
    \compose[()] \funcseqname{f}^1
$.

\begin{proof}
If
$\funcseqname{f}^1$ $\Sigma$-commutes with $\funcseqname{F}^1,\funcseqname{F}^2$ and
$\funcseqname{f}^2$ $\Sigma$-commutes with $\funcseqname{F}^2,\funcseqname{F}^3$
then by \pagecref{lem:Sigmacomm}
$\funcseqname{f}^2 \compose[()] \funcseqname{f}^1$ $\Sigma$-commutes with $\funcseqname{F}^1,\funcseqname{F}^3$.
Thus $\funcseqname{f}^2 \compose[()] \funcseqname{f}^1$ is a morphism.

Associativity is just \pagecref{lem:seqfunc}.
\end{proof}

$\Id_{\seqname{P}^i}$ is an identity morphism of $\seqname{P}^i$.

\begin{proof}
$\funcseqname{f}^i \compose[()] \ID_{\seqname{S}^i} = \funcseqname{f}^i$
and
$\ID_{\seqname{S}^{i+1}} \compose[()] \funcseqname{f}^i = \funcseqname{f}^i$
\end{proof}

Let each $\funcseqname{f}^i$ be a semi-strict prestructure morphism.
Then $\funcseqname{f}^2 \compose[()] \funcseqname{f}^1$ is a semi-strict
prestructure morphism.

\begin{proof}
If $\catseqname{S}^1 \SUBCAT[full-] \catseqname{S}^2$ and
$\catseqname{S}^2 \SUBCAT[full-] \catseqname{S}^3$ then
$\catseqname{S}^1 \SUBCAT[full-] \catseqname{S}^3$.
For each $\alpha \prec \Alpha$:

\begin{enumerate}
\item Each $\catname{S}^i_\alpha$, $i=1,2,3$, is a model category:
\newline
Each $\funcname{f}^i_\alpha$ is locally an
$\catname{S}^1_\alpha$-$\catname{S}^{i+1}_\alpha$ m-morphism of
$S^i_\alpha$ to $S^{i+1}_\alpha$.
Then $\funcname{f}^2_\alpha \compose \funcname{f}^1_\alpha$ is locally
an $\catname{S}^i_\alpha$-$\catname{S}^{i+2}_\alpha$ m-morphism of
$S^i_\alpha$ to $S^{i+2}_\alpha$ by \pagecref{lem:modlocal}.

\item Each $\catname{S}^i_\alpha$, $i=1,2,3$, is a topological category:
\newline
Each $\funcname{f}^i_\alpha$ is locally a
$\catname{S}^i_\alpha$-$\catname{S}^{i+1}_\alpha$ morphism of
$S^i_\alpha$ to $S^{i+1}_\alpha$.
Then $\funcname{f}^2_\alpha \compose \funcname{f}^1_\alpha$ is locally
an $\catname{S}^i_\alpha$-$\catname{S}^{i+2}_\alpha$ morphism of
$S^i_\alpha$ to $S^{i+2}_\alpha$ by \pagecref{lem:topLocal}.

\item Otherwise
\newline
Each $\funcname{f}^i_\alpha$ is a morphism of $\catname{S}^{i+1}_\alpha$.
Then $\funcname{f}^2_\alpha \compose \funcname{f}^1_\alpha$ is a
morphism of $\catname{S}^{i+2}_\alpha$.
\end{enumerate}
\end{proof}

Let each $\funcseqname{f}^i$ be a strict prestructure morphism. Then
$\funcseqname{f}^2 \compose[()] \funcseqname{f}^1$ is a strict
prestructure morphism.

\begin{proof}
If $\catseqname{S}^1 \SUBCAT[full-] \catseqname{S}^2$ and
$\catseqname{S}^2 \SUBCAT[full-] \catseqname{S}^3$ then
$\catseqname{S}^1 \SUBCAT[full-] \catseqname{S}^3$.

If $\funcname{f}^1_\alpha$ is a morphism of
$\catname{S}^2_\alpha \subcat[full-] \catname{S}^3_\alpha$ and
$\funcname{f}^2_\alpha$ is a morphism of $\catname{S}^3_\alpha$ then
$\funcname{f}^2_\alpha \compose \funcname{f}^1_\alpha$ is a morphism of
$\catname{S}^3_\alpha$.
\end{proof}
\end{lemma}

\begin{corollary}[Prestructure morphisms]
\label{cor:premorph}
Let $\catseqname{S}=(\catname{S}_\alpha, \alpha \prec \Alpha)$
be a sequence of categories,
$
  \seqname{S} \defeq (S_\alpha \alpha \prec \Alpha)
  \seqin \catseqname{S}
$
be a sequence of spaces,
$\funcseqname{F}=(\funcname{F}_\gamma, \gamma \prec \Gamma)$
be a sequence of continuous functions with $\catseqname{S}$-signatures
$\Sigma$ and
$
  \seqname{P} \defeq
  (\catseqname{S}, \seqname{S}, \Sigma, \funcseqname{F})
$
be a $\catseqname{S}-\Sigma$ prestructure. Then $\Id_{\seqname{P}}$ is a
strict prestructure morphism from $\seqname{P}$ to $\seqname{P}$.

\begin{proof}
$\seqname{S} \SUBSETEQ \seqname{S}$ and
$
  \funcname{F}_\gamma =
  \funcname{F}_\gamma \restriction_%
    {
      \bigtimes \underline{\overset{\sigma_\gamma}{\seqname{S}}},
      \tail \left ( \overset{\sigma_\gamma}{\seqname{S}} \right )
    }
$.

Each
$\Id_{S_\alpha,S_\alpha} = \Id_{S_\alpha} \arin \catname{S}_\alpha$ and
$\catseqname{S} \SUBCAT[full-] \catseqname{S}$.
\end{proof}
\end{corollary}

\section{M-charts and m-atlases}
\label{sec:m-charts}
The conventional definitions of fiber bundles and manifolds use the
language of charts, atlases and transition functions, with slight
technical differences among them. This paper continues that usage, but
modifies the definitions of charts in order to make them fit more
natuarally into the context of local coordinate spaces. It adds a
prefix, e.g., $\Ck$, in order to avoid confusion with the conventional
definitions. An m-atlas based on a topological spaces is closer to the
conventional definitions of an atlas for a manifold while an m-atlas
based on a model space is suitable for defining both manifolds and fiber
bundles. Although simple manifolds and fiber bundles could both be
defined directly in terms of maximal m-atlases, this paper has a
different perspective, and uses the m-atlases as part of the more
general Local Coordinate Space (LCS), presented in
\pagecref{sec:lcs}.

This section defines M-atlases, categories of M-atlases and functors,
and proves some basic results.

\subsection{M-charts}
\label{sub:m-charts}
\begin{definition}[M-charts]
Let $\seqname{E} \defeq (E, \catname{E})$ and $\catname{C} \defeq (C,\catname{C})$ be model spaces. An m-chart
$(U, V, \phi)$ of
$\seqname{E}$ in the coordinate space $\catname{C}$ consists of
\begin{enumerate}
\item A model neighborhood $U\objin \catname{E}$, known as a coordinate
patch
\item A model neighborhood $V\objin \catname{C}$
\item A model homeomorphism $\phi \maps U \toiso V$, known as a
coordinate function
\end{enumerate}
\begin{remark}
I consider it clearer to explicate the range, rather than the
conventional usage of specifying only the domain and function or the
minimalist usage of specifying only the function.
\end{remark}

Let $E$ be a topological space and  $\catname{C} \defeq (C,\catname{C})$ a model space. An
m-chart $(U, V, \phi)$ of $E$ in the coordinate space $\catname{C}$
consists of
\begin{enumerate}
\item An open set $U$ of $E$, known as a coordinate patch
\item A model neighborhood $V \objin \catname{C}$
\item A homeomorphism $\phi \maps U \toiso V$, known as a coordinate
function, that maps open sets into model neighborhoods
\end{enumerate}
\end{definition}

\begin{lemma}[M-charts]
Let $\seqname{E} \defeq (E, \catname{E})$ and
$\seqname{C} \defeq (C,\catname{C})$ be model spaces.

Let $(U, V, \phi)$ be an m-chart of $\seqname{E}$ in the coordinate
space $\seqname{C}$ such that every open subset of $V$ is a model
neighborhood of $\seqname{C}$. Then $(U, V, \phi)$ is an m-chart of $E$
in the coordinate space $\seqname{C}$.

\begin{proof}
$U$, $V$ and $\phi$ satisfy the conditions of the definition:
\begin{enumerate}
\item Since $U$ is a model neighborhood, it is open.
\item $V$ is a model neighborhood of $\seqname{C}$,
\item Since $\phi \maps U \toiso V$ is a model homeomorphism, it is a
homeomorphism. If $U' \subseteq U$ is open then $\phi[U'] \subseteq V$
is open and thus a model neighborhood of $(C,\catname{C})$ by
hypothesis. Thus $\phi$ is a homeomorphism that maps open sets into
model neighborhood of $(C,\catname{C})$.

\end{enumerate}
\end{proof}

Let $(U, V, \phi)$ be an m-chart of $E$ in the coordinate space
$\seqname{C}$. $(U, V, \phi)$ is an m-chart of $\seqname{E}$ in the
coordinate space $\seqname{C}$ iff every open subset of $U$ is a model
neighborhood of $\seqname{E}$.

\begin{proof}
If $(U, V, \phi)$ is an m-chart of $E$ in the coordinate space
$\seqname{C}$ and also an m-chart of $\seqname{E}$ in the coordinate
space $\seqname{C}$ then let $U' \subseteq U$ be open. Since
$(U, V, \phi)$ is an m-chart of $E$ in the coordinate space
$\seqname{C}$, $\phi[U']$ is a model neighborhood of $\seqname{C}$.
Since $(U, V, \phi)$ is an m-chart of $\seqname{E}$ in the coordinate
space $\seqname{C}$, $\phi$ is a model function and thus
$\phi^{-1}[\phi[U']]$ is a model neighborhood of $\seqname{E}$.

If $(U, V, \phi)$ is an m-chart of $E$ in the coordinate space
$\seqname{C}$ and every open $U'\subseteq U$ is a model neighborhood of
$\seqname{E}$ then
\begin{enumerate}
\item Since $U \subseteq U$ is open, $U$ is a model neighborhood of $\seqname{E}$.
\item $V$ is a model neighborhood of $\catname{C}$.
\item $\phi$ is a model homeomorphism: If $U' \subseteq U$ is a model
neighborhood of $\seqname{E}$ then $U'$ is open and thus $\phi[U']$ is a
model neighborhood of $\seqname{C}$. If $V' \subseteq V$ is a model
neighborhood of $\seqname{C}$ then $\phi^{-1}[V']$ is open and thus
a model neighborhood of $\seqname{E}$.
\end{enumerate}

By definition, every open set of $E$ is a model neighborhood of $\Triv{E}$.
\end{proof}
\end{lemma}

\begin{corollary}[M-charts of $\Triv{E}$]
\label{cor:triv}
Let $E$ be a topological space. $(U,V,\phi)$ is an m-chart of $E$ in
the coordinate space $\seqname{C}$ iff $(U,V,\phi)$ is an m-chart of
$\Triv{E}$ in the coordinate space $\seqname{C}$.

\begin{proof}
Let $(U, V, \phi)$ be an m-chart of $E$ in the coordinate space
$\seqname{C}$. Every open subset of $E$ is a model neighborhood of
$\Triv{E}$ by definition, hence $(U, V, \phi)$ is an m-chart of
$\Triv{E}$ in the coordinate space $\seqname{C}$.

Let $(U, V, \phi)$ be an m-chart of $\Triv{E}$ in the coordinate space
$\seqname{C}$ and let $V' \subseteq V$ be open.
$\phi^{-1}[V'] \subseteq U$ is open, hence by definition a model
neighborhood of $\Triv{E}$. Then $\phi[\phi^{-1}[V']] = V'$ is a model
neighborhood of $\seqname{C}$.
\end{proof}
\end{corollary}

\begin{definition}[Subcharts]
Let $(U, V, \phi)$ be an m-chart of $\seqname{E} \defeq (E, \catname{E})$ in the coordinate
space $(C,\catname{C})$ and $U' \subseteq U$ a model neighborhood of
$(E,\catname{E})$. Then \\
$(U', V', \phi') \defeq (U', \phi[U'], \phi \restriction_{U',V'})$
is a subchart of $(U, V, \phi)$.

Let $(U, V, \phi)$ be an m-chart of $E$ in the coordinate
space $(C,\catname{C})$ and $U' \subseteq U$ an open set of
$E$. Then
$(U', V', \phi') \defeq (U', \phi[U'], \phi \restriction_{U',V'})$
is a subchart of $(U, V, \phi)$.
\end{definition}

\begin{lemma}[Subcharts]
Let $\seqname{E} \defeq (E, \catname{E})$ and
$\seqname{C} \defeq (C,\catname{C})$ be model spaces.

Let $(U, V, \phi)$ be an m-chart of $\seqname{E}$ in the coordinate
space $\seqname{C}$ and $(U', V', \phi')$ a subchart of $(U, V, \phi)$.
Then $(U', V', \phi')$ is an m-chart of $\seqname{E}$ in the coordinate
space $\seqname{C}$.

\begin{proof}
$U'$, $V'$ and $\phi'$ satisfy the conditions of the definition of an m-chart:
\begin{enumerate}
\item $U'$ is a model neighborhood of $\seqname{E}$ by the definition of
subchart.
\item Since $U'$ is a model neighborhood and $\phi$ is a model
homeomorphism, $V'=\phi[U']$ is a model neighborhood.
\item Since $\phi$ and $\phi^{-1}$ are model functions, so are
$\phi' = \phi \restriction_{U',V'}$ and
$\phi'^{-1} = \phi^{-1} \restriction_{V',U'}$ Thus $\phi'$ is a model
homeomorphism.
\end{enumerate}
\end{proof}

Let $(U, V, \phi)$ be an m-chart of $E$ in the coordinate space
$\seqname{C}$ and $(U', V', \phi')$ a subchart of $(U, V, \phi)$. Then
$(U', V', \phi')$ is an m-chart of $E$ in the coordinate space
$\seqname{C}$.

\begin{proof}
$U'$, $V'$ and $\phi'$ satisfy the conditions of the definition of an m-chart:
\begin{enumerate}
\item $U'$ is open by the definition of subchart.
\item Since $\phi$ is a homeomorphism and $U'$ is open, $V'=\phi[U']$ is open.
\item Since $\phi$ is a homeomorphism that maps open sets into model
neighborhoods, $\phi \restriction_{U',V'}$ is a homeomorphism that maps
open sets into model neighborhoods.
\end{enumerate}
\end{proof}
\end{lemma}

\begin{definition}[M-compatibility]
Let $(U, V, \phi)$ and  $(U',V', \phi')$ be m-charts of $(E, \catname{E})$
in the coordinate space $(C,\catname{C})$. Then $(U, V, \phi)$ is
compatible with $(U', V', \phi')$ iff either

\begin{enumerate}
\item $U$ and  $U'$ are disjoint
\item The transition function
$\funcname{t}=\phi' \compose \phi^{-1} \restriction_{\phi[U \cap U']}$
is an isomorphism of $\catname C$
\end{enumerate}

Let $(U, V, \phi)$ and  $(U',V', \phi')$ be m-charts of $E$
in the coordinate space $(C,\catname{C})$. Then $(U, V, \phi)$ is
m-compatible with $(U', V', \phi')$ in the coordinate space
$(C,\catname{C})$ iff either

\begin{enumerate}
\item $U$ and  $U'$ are disjoint
\item The transition function
$\funcname{t}=\phi' \compose \phi^{-1} \restriction_{\phi[U \cap U']}$
is an isomorphism of $\catname C$
\end{enumerate}
\end{definition}

\begin{lemma}[Symmetry of m-compatibility]
Let $(U, V, \phi)$ and  $(U',V', \phi')$ be m-charts of $(E, \catname{E})$
in the coordinate space $(C,\catname{C})$. Then $(U, V, \phi)$ is
m-compatible with $(U', V', \phi')$ in the coordinate space
$(C,\catname{C})$ iff $(U', V', \phi')$ is m-compatible with
$(U,V,\phi)$ in the coordinate space $(C,\catname{C})$.

Let $(U, V, \phi)$ and  $(U',V', \phi')$ be m-charts of $E$
in the coordinate space $(C,\catname{C})$. Then $(U, V, \phi)$ is
m-compatible with $(U', V', \phi')$ iff $(U', V', \phi')$ is m-compatible
with $(U,V,\phi)$.

\begin{proof}
The same proof applies to both cases.
It suffices to prove the implication in only one direction.
\begin{enumerate}
\item $U \cap U' = U' \cap U$.
\item Since the transition function
$\funcname{t}=\phi' \compose \phi^{-1} \restriction_{\phi[U \cap U']}$
is an isomorphism of $\catname C$, so is
$\funcname{t}^{-1} = \phi \compose \phi'^{-1} \restriction_{\phi'[U \cap U']}$.
\end{enumerate}
\end{proof}
\end{lemma}

\begin{lemma}[M-compatibility of subcharts]
Let $(U_i, V_i, \phi_i)$, $i=1,2$, be m-charts of $(E, \catname{E})$ in
the coordinate space $(C,\catname{C})$, $(U'_i, V'_i, \phi'_i)$ be
subcharts and $(U_1, V_1, \phi_1)$ be m-compatible with
$(U_2, V_2, \phi_2)$.  Then $(U'_1, V'_1, \phi'_1)$ is m-compatible with
$(U'_2, V'_2, \phi'_2)$.

\begin{proof}
Since subcharts are charts, $U'_i$ and $V'_i$ are model neighborhoods.
If $U_1 \cap U_2 = \emptyset$ then $U'_1 \cap U'_2 = \emptyset$. If
$U'_1 \cap U'_2 = \emptyset$ then $(U'_1, V'_1, \phi'_1)$ is
m-compatible with $(U'_2, V'_2, \phi'_2)$.  Otherwise, the transition
function
$
  t^1_2 \defeq
  \phi_2 \compose \phi^{-1}_1 \restriction_{\phi_1[U_1 \cap U_2]}
$
is a model homeomorphism and hence
$
  t^1_2 \restriction_{\phi_1[U'_1 \cap U'_2]} \maps
    \phi_1[U'_1 \cap U'_2] \toiso
    \phi_2[U'_1 \cap U'_2]
$
is a model homeomorphism.
\end{proof}

Let $(U_i, V_i, \phi_i)$, $i=1,2$, be m-charts of $E$ in the coordinate
space $(C,\catname{C})$ and $(U'_i, V'_i, \phi'_i)$ subcharts. Then
$(U'_i, V'_i, \phi'_i)$ is m-compatible with $(U_i,V_i,\phi_i)$.

\begin{proof}
Since subcharts are charts, $U'_i$ and $V'_i$ are open. If
$U_1 \cap U_2 = \emptyset$ then $U'_1 \cap U'_2 = \emptyset$. If
$U'_1 \cap U'_2 = \emptyset$ then $(U'_1, V'_1, \phi'_1)$ is
m-compatible with $(U'_2, V'_2, \phi'_2)$.  Otherwise, the transition
function
$
  t^1_2 \defeq
  \phi_2 \compose \phi^{-1}_1 \restriction_{\phi_1[U_1 \cap U_2]}
$
is a model homeomorphism and hence
$
  t^1_2 \restriction_{\phi_1[U'_1 \cap U'_2]} \maps
  \phi_1[U'_1 \cap U'_2] \toiso
  \phi_2[U'_1 \cap U'_2]
$
is a model homeomorphism.
\end{proof}
\end{lemma}

\begin{corollary}[M-compatibility with subcharts]
Let $(U, V, \phi)$ be an m-chart of $(E, \catname{E})$ in
the coordinate space $(C,\catname{C})$ and $(U', V', \phi')$ a
subchart. Then $(U', V', \phi')$ is m-compatible with
$(U,V,\phi)$.

\begin{proof}
$(U, V, \phi)$ is m-compatible with itself and is a subchart of
itself,
\end{proof}

Let $(U, V, \phi)$ be an m-chart of $E$ in the coordinate
space $(C,\catname{C})$ and $(U', V', \phi')$ a subchart. Then
$(U', V', \phi')$ is m-compatible with $(U,V,\phi)$.

\begin{proof}
$(U, V, \phi)$ is m-compatible with itself and is a subchart of
itself,
\end{proof}
\end{corollary}

\begin{definition}[Covering by m-charts]
Let $\seqname{A}$ be a set of charts of the topological space $E$ in
the coordinate space
$\seqname{C}=(C,\catname{C})$.
$\seqname{A}$ covers $E$ iff $\pi_1[\seqname{A}]$ covers $E$.

Let $\seqname{A}$ be a set of charts of the model space
$\seqname{E}=(E, \catname{E})$ in the coordinate space
$\seqname{C}=(C,\catname{C})$.  $\seqname{A}$ covers $\seqname{E}$ iff
$\pi_1[\seqname{A}]$ covers $E$.
\end{definition}

\subsection{M-atlases}
\label{sub:m-atlas}
A set of charts can be atlases for different coordinate model spaces
even if it is for the same total model space. In order to aggregate
atlases into categories, there must be a way to distinguish them.
Including the two\footnote{
The spaces are redundant, but convenient.}
spaces in the definitions of the categories serves the purpose.

\begin{definition}[M-atlases]
\label{def:m-atlas}
Let $\seqname{A}$ be a set of mutually m-compatible m-charts of
$\seqname{E}=(E, \catname{E})$ in the coordinate space
$\seqname{C}=(C,\catname{C})$. Then $\seqname{A}$ is an m-atlas of
$\seqname{E}$ in the coordinate space $\seqname{C}$, abbreviated
$\isAtl_\Ob(\seqname{A}, \seqname{E}, \seqname{C})$, iff $\seqname{A}$
covers $E$. $\seqname{A}$ is a full atlas of $\seqname{E}$ in the
coordinate space $\seqname{C}$, abbreviated
$\full{\isAtl_\Ob}(\seqname{A}, E, \seqname{C})$, iff $\seqname{A}$
covers $E$ and $\pi_2[\seqname{A}]$ covers $C$. The triple
$(\seqname{A}, \seqname{E}, \seqname{C})$ refers to $\seqname{A}$
considered as an m-atlas of $\seqname{E}$ in the coordinate space
$\seqname{C}$.

Let $\seqname{E}=(E, \catname{E})$ and $\seqname{C}=(C,\catname{C})$ be
model spaces. Then
\begin{equation}
\Atl_\Ob(\seqname{E}, \seqname{C}) \defeq
  \set
    {{(\seqname{A}, \seqname{E}, \seqname{C})}}%
    [{{\isAtl_\Ob(\seqname{A}, \seqname{E}, \seqname{C})}}]
\end{equation}
\begin{equation}
\full{\Atl_\Ob}(\seqname{E}, \seqname{C}) \defeq
  \set
    {{(\seqname{A}, \seqname{E}, \seqname{C})}}%
    [{{\full{\isAtl_\Ob}(\seqname{A}, \seqname{E}, \seqname{C})}}]
\end{equation}

Let $E$ be a topological space, $\seqname{C}=(C,\catname{C})$ a model
space and $\seqname{A}$ be a set of mutually m-compatible m-charts of
$E$ in the coordinate space $\seqname{C}$. Then $\seqname{A}$ is an
m-atlas of $E$ in the coordinate space $\seqname{C}$, abbreviated as
$\isAtl_\Ob(\seqname{A}, E, \seqname{C})$, iff $\pi_1[\seqname{A}]$
covers $E$. $\seqname{A}$ is a full atlas of $E$ in the coordinate space
$\seqname{C}$, abbreviated $\full{\isAtl_\Ob}(\seqname{A}, E,
\seqname{C})$, iff $\pi_1[\seqname{A}]$ covers covers $E$ and
$\pi_2[\seqname{A}]$ covers $C$. The triple
$(\seqname{A}, E, \seqname{C})$ refers to $\seqname{A}$ considered as an
m-atlas of $E$ in the coordinate space $\seqname{C}$.

Let $E$ be a topological space and $\seqname{C}=(C,\catname{C})$ a
model space. Then
\begin{equation}
\Atl_\Ob(E, \seqname{C}) \defeq
  \set
    {{(\seqname{A}, E, \seqname{C})}}%
    [{{\isAtl_\Ob(\seqname{A}, E, \seqname{C})}}]
\end{equation}
\begin{equation}
\full{\Atl_\Ob}(E, \seqname{C}) \defeq
  \set
    {{(\seqname{A}, E, \seqname{C})}}%
    [{{\full{\isAtl_\Ob}(\seqname{A}, E, \seqname{C})}}]
\end{equation}
\end{definition}

By abuse of language we write $U \in \seqname{A}$ for
$U \in \pi_1[\seqname{A}]$.

\begin{definition}[M-compatibility with atlases]
An m-chart $(U, V, \phi)$ of $(E, \catname{E})$ is m-compatible with an
m-atlas $\seqname{A}$ iff it is m-compatible with every chart in the m-atlas.

An m-chart $(U, V, \phi)$ of $E$ is m-compatible with an
m-atlas $\seqname{A}$ iff it is m-compatible with every chart in the m-atlas.
\end{definition}

\begin{lemma}[M-compatibility of subcharts with atlases]
\label{lem:M-Compat}
Let $\seqname{A}$ be an m-atlas of $(E, \catname{E})$ in the coordinate
space $(C,\catname{C})$ and $\seqname{C}_1 = (U_1, V_1, \phi_1)$ an
m-chart in $\seqname{A}$. Then any subchart of $\seqname{C}_1$ is
m-compatible with $\seqname{A}$.

Let $\seqname{A}$ be an m-atlas of $E$ in the coordinate
space $(C,\catname{C})$ and $\seqname{C}_1 = (U_1, V_1, \phi_1)$ an
m-chart in $\seqname{A}$. Then any subchart of $\seqname{C}_1$ is
m-compatible with $\seqname{A}$.

\begin{proof}
The same proof applies in both cases.
Let $\seqname{C}' = (U', V', \phi')$ be a subchart of $\seqname{C}_1$ and
$\seqname{C}_2 = (U_2, V_2, \phi_2)$ another chart in $\seqname{A}$.
\begin{enumerate}
\item If $U_1 \cap U_2 = \emptyset$, then $U' \cap U_2 = \emptyset$
\item If $U' \cap U_2 = \emptyset$ then $\seqname{C}'$ is m-compatible
with $\seqname{C}_2$
\item Otherwise the transition function
$
  t^1_2 \defeq
  \phi_2 \compose \phi^{-1}_1 \restriction_{\phi_1[U_1 \cap U_2]}
$
is an isomorphism of $(C,\catname{C})$. Since $\phi_1[U_1 \cap U_2]$ and
$V'$ are model neighborhoods, so is $\phi_1[U' \cap U_2]$ and thus
$t^1_2 \restriction_{\phi_1[U' \cap U_2]}$ is an isomorphism of
$(C,\catname{C})$.
\end{enumerate}
\end{proof}
\end{lemma}

\begin{lemma}[Extensions of m-atlases]
\label{m-atl:extensions}
Let $\seqname{E}=(E,\catname{E})$ and $\seqname{C}=(C,\catname{C})$ be
model spaces, $\seqname{A}$ an m-atlas of $\seqname{E}$ in the
coordinate space $\seqname{C}$, and $(U,V,\phi)$, $(U',V',\phi')$
m-charts of $\seqname{E}$ in the coordinate space $\seqname{C}$
m-compatible with $\seqname{A}$ in the coordinate space $\seqname{C}$.
Then $(U,V,\phi)$ is m-compatible with $(U',V',\phi')$ in the coordinate
space $\seqname{C}$.

\begin{proof}
If $U \cap U' = \emptyset$ then $(U,V,\phi)$ is m-compatible with
$(U',V',\phi')$. Otherwise, $\phi$ is a model homeomorphism,
$U \cap U'$ is a model neighborhood of $\seqname{E}$,
$\phi[U \cap U']$ is a model neighborhood of
$\seqname{C}$, $\phi'[U \cap U']$ is a model neighborhood of
$\seqname{C}$ and
$\phi' \compose \phi^{-1} \maps \phi[U \cap U'] \toiso \phi'[U \cap U']$
is a homeomorphism.  It remains to show that
$\phi' \compose \phi^{-1} \restriction_{\phi[U \cap U']}$ is a model
homeomorphism.  Let
$(U_\alpha,V_\alpha,\phi_\alpha)$, $\alpha \prec \Alpha$, be charts in
$\seqname{A}$ such that
$U \cap U' \subseteq \union[\alpha \prec \Alpha]{U_\alpha}$ and
$U \cap U' \cap U_\alpha \neq \emptyset$, $\alpha \prec \Alpha$.  Since
the charts are m-compatible with $(U_\alpha,V_\alpha,\phi_\alpha)$,
$
  \phi' \compose \phi^{-1}_\alpha
  \restriction_{\phi_\alpha[U \cap U' \cap U_\alpha]}
$
and
$
  \phi_\alpha \compose \phi^{-1}
  \restriction_{\phi[U \cap U' \cap U_\alpha]}
$
are model homeomorphisms and thus
$
  \phi' \compose \phi^{-1} \restriction_{\phi[U \cap U' \cap U_\alpha]} =
  \phi' \compose \phi^{-1}_\alpha \compose \phi_\alpha \compose \phi^{-1}
  \restriction_{\phi[U \cap U' \cap U_\alpha]}
$
is a model homeomorphism. Then by
{
  \showlabelsinline
  \cref{mod:sheaf}
}
of \pagecref{def:model}\!, $\phi' \compose \phi^{-1}$ is a model
homeomorphism.
\end{proof}
\end{lemma}

\begin{definition}[Maximal m-atlases]
\label{def:maxM}
Let $\seqname{E}=(E, \catname{E})$ and $\seqname{C}=(C,\catname{C})$ be
model spaces and $\seqname{A}$ an m-atlas of $\seqname{E}$ in the
coordinate space $\seqname{C}$. $\seqname{A}$ is a maximal m-atlas of
$\seqname{E}$ in the coordinate space $\seqname{C}$, abbreviated
$\maximal{\isAtl_\Ob}(\seqname{A}, \seqname{E}, \seqname{C})$, iff
$\seqname{A}$ cannot be extended by adding an additional m-compatible
chart. $\seqname{A}$ is a semi-maximal m-atlas of $\seqname{E}$ in the
coordinate space $\seqname{C}$, abbreviated
$\maximal[S-]{\isAtl_\Ob}(\seqname{A}, \seqname{E}, \seqname{C})$, iff
whenever
$(U,V,\phi) \in \seqname{A}$, $U' \subseteq U, V' \subseteq V$ and
$V'' \objin \catname{C}$ are model neighborhoods, $\phi[U'] = V'$ and
$\phi' \maps V' \toiso V''$ is an isomorphism of $\catname{C}$ then
$(U', V'', \phi' \compose \phi \restriction_{U'}) \in \seqname{A}$.

Let $E$ be a topological space, $\seqname{C}=(C,\catname{C})$ be a model
space and $\seqname{A}$ be an m-atlas of $E$ in the coordinate space
$\seqname{C}$. $\seqname{A}$ is a maximal m-atlas of $E$ in the
coordinate space $\seqname{C}$, abbreviated
$\maximal{\isAtl_\Ob}(\seqname{A}, E, \seqname{C})$, iff $\seqname{A}$
is an m-atlas that cannot be extended by adding an additional
m-compatible chart.  $\seqname{A}$ is a semi-maximal m-atlas of $E$ in
the coordinate space $\seqname{C}$, abbreviated
$\maximal[S-]{\isAtl_\Ob}(\seqname{A}, E, \seqname{C})$, iff whenever
$(U,V,\phi) \in \seqname{A}$, $U' \subseteq U, V' \subseteq V$ and
$V'' \objin \catname{C}$ are model neighborhoods, $\phi[U'] = V'$ and
$\phi' \maps V' \toiso V''$ is an isomorphism of $\catname{C}$ then
$(U', V'', \phi' \compose \phi \restriction_{U'}) \in \seqname{A}$.

\begin{remark}
There is no sheaf condition;\footnote{
However, note \cref{mod:sheaf} (restricted sheaf condition) of
\pagecref{def:model}.}
the union of a set of coordinate patches in the maximal atlas whose
coordinate functions match on the intersections need not be a coordinate
patch in the atlas.
\end{remark}

Let $(E, \catname{E})$ and $(C,\catname{C})$ be model spaces. Then
\begin{equation}
  \maxfull{\isAtl_\Ob}(\seqname{A}, \seqname{E}, C) \defeq
    \full{\isAtl_\Ob}(\seqname{A}, \seqname{E}, C) \land
    \maximal{\isAtl_\Ob}(\seqname{A}, \seqname{E}, C)
\end{equation}
\begin{equation}
  \maxfull{\isAtl_\Ob}(\seqname{A}, E, C) \defeq
    \full{\isAtl_\Ob}(\seqname{A}, E, C) \land
    \maximal{\isAtl_\Ob}(\seqname{A}, E, C)
\end{equation}
\begin{equation}
  \maxfull[S-]{\isAtl_\Ob}(\seqname{A}, \seqname{E}, C) \defeq
    \full{\isAtl_\Ob}(\seqname{A}, \seqname{E}, C) \land
    \maximal[S-]{\isAtl_\Ob}(\seqname{A}, \seqname{E}, C)
\end{equation}
\begin{equation}
  \maxfull[S-]{\isAtl_\Ob}(\seqname{A}, E, C) \defeq
    \full{\isAtl_\Ob}(\seqname{A}, E, C) \land
    \maximal[S-]{\isAtl_\Ob}(\seqname{A}, E, C)
\end{equation}
\begin{equation}
\maximal{\Atl_\Ob} \bigl ( (E,\catname{E}), (C,\catname{C}) \bigr )
  \defeq
  \set
    {{\bigl ( \seqname{A}, (E,\catname{E}), (C,\catname{C}) \bigr ) }}%
    [{{\maximal{\isAtl_\Ob}\bigl ( \seqname{A}, E, (C,\catname{C}) \bigr )}}]
\end{equation}
\begin{equation}
\maximal[S-]{\Atl_\Ob} \bigl ( (E,\catname{E}), (C,\catname{C}) \bigr )
  \defeq
  \set
    {{\bigl ( \seqname{A}, (E,\catname{E}), (C,\catname{C}) \bigr ) }}%
    [{{\maximal[S-]{\isAtl_\Ob}\bigl ( \seqname{A}, E, (C,\catname{C}) \bigr )}}]
\end{equation}
\begin{equation}
\maxfull{\Atl_\Ob} \bigl ( (E,\catname{E}), (C,\catname{C}) \bigr )
  \defeq
  \set
    {{\bigl ( \seqname{A}, (E,\catname{E}), (C,\catname{C}) \bigr ) }}%
    [{{\maxfull{\isAtl_\Ob}\bigl ( \seqname{A}, E, (C,\catname{C}) \bigr )}}]
\end{equation}
\begin{equation}
\maxfull[S-]{\Atl_\Ob} \bigl ( (E,\catname{E}), (C,\catname{C}) \bigr )
  \defeq
  \set
    {{\bigl ( \seqname{A}, (E,\catname{E}), (C,\catname{C}) \bigr ) }}%
    [{{\maxfull[S-]{\isAtl_\Ob}\bigl ( \seqname{A}, E, (C,\catname{C}) \bigr )}}]
\end{equation}

Let $E$ be a topological space and $(C,\catname{C})$ be a model space. Then
\begin{equation}
\maximal{\Atl_\Ob} \bigl ( E, (C,\catname{C}) \bigr ) \defeq
  \set
    {{ \bigl ( \seqname{A}, E, (C,\catname{C}) \bigr ) }}%
    [{{\maximal{\isAtl_\Ob} \bigl ( E, (C,\catname{C}) \bigr ) }}]
\end{equation}
\begin{equation}
\maximal[S-]{\Atl_\Ob} \bigl ( E, (C,\catname{C}) \bigr ) \defeq
  \set
    {{ \bigl ( \seqname{A}, E, (C,\catname{C}) \bigr ) }}%
    [{{\maximal[S-]{\isAtl_\Ob} \bigl ( E, (C,\catname{C}) \bigr ) }}]
\end{equation}
\begin{equation}
\maxfull{\Atl_\Ob} \bigl ( E, (C,\catname{C}) \bigr ) \defeq
  \set
    {{ \bigl ( \seqname{A}, E, (C,\catname{C}) \bigr ) }}%
    [{{\maxfull{\isAtl_\Ob} \bigl ( E, (C,\catname{C}) \bigr ) }}]
\end{equation}
\begin{equation}
\maxfull[S-]{\Atl_\Ob} \bigl ( E, (C,\catname{C}) \bigr ) \defeq
  \set
    {{ \bigl ( \seqname{A}, E, (C,\catname{C}) \bigr ) }}%
    [{{\maxfull[S-]{\isAtl_\Ob} \bigl ( E, (C,\catname{C}) \bigr ) }}]
\end{equation}

Let $\seqname{E}$ and $\seqname{C}$ be sets of model spaces. Then
\begin{multline}
\maximal{\Atl_\Ob}(\seqname{E}, \seqname{C}) \defeq
\set
{{
   \bigl ( \seqname{A}, (E,\catname{E}), (C,\catname{C}) \bigr )
}}%
[
  {(E,\catname{E}) \in \seqname{E}},
  {(C,\catname{C}) \in \seqname{C}},
  {\maximal{\isAtl_\Ob} \bigl ( \seqname{A}, (E,\catname{E}), (C,\catname{C}) \bigr ) }
]*
\end{multline}
\begin{multline}
\maximal[S-]{\Atl_\Ob}(\seqname{E}, \seqname{C}) \defeq
\set
{{
   \bigl ( \seqname{A}, (E,\catname{E}), (C,\catname{C}) \bigr )
}}%
[
  {(E,\catname{E}) \in \seqname{E}},
  {(C,\catname{C}) \in \seqname{C}},
  {\maximal[S-]{\isAtl_\Ob} \bigl ( \seqname{A}, (E,\catname{E}), (C,\catname{C}) \bigr ) }
]*
\end{multline}
\begin{multline}
\maxfull{\Atl_\Ob}(\seqname{E}, \seqname{C}) \defeq
\set
{{
   \bigl ( \seqname{A}, (E,\catname{E}), (C,\catname{C}) \bigr )
}}%
[
  {(E,\catname{E}) \in \seqname{E}},
  {(C,\catname{C}) \in \seqname{C}},
  {\maxfull{\isAtl_\Ob} \bigl ( \seqname{A}, (E,\catname{E}), (C,\catname{C}) \bigr ) }
]*
\end{multline}
\begin{multline}
\maxfull[S-]{\Atl_\Ob}(\seqname{E}, \seqname{C}) \defeq
\set
{{
   \bigl ( \seqname{A}, (E,\catname{E}), (C,\catname{C}) \bigr )
}}%
[
  {(E,\catname{E}) \in \seqname{E}},
  {(C,\catname{C}) \in \seqname{C}},
  {\maxfull[S-]{\isAtl_\Ob} \bigl ( \seqname{A}, (E,\catname{E}), (C,\catname{C}) \bigr ) }
]*
\end{multline}

Let $\seqname{E}$ be a set of topological spaces and $\seqname{C}$ a
set of model spaces. Then
\begin{multline}
\maximal{\Atl_\Ob}(\seqname{E}, \seqname{C}) \defeq
\set
{{
   \bigl ( \seqname{A}, E, (C,\catname{C}) \bigr )
}}%
[
  E \in \seqname{E},
  {(C,\catname{C}) \in \seqname{C}},
  {\maximal{\isAtl_\Ob} \bigl ( \seqname{A}, E, (C,\catname{C}) \bigr )  }
]*
\end{multline}
\begin{multline}
\maximal[S-]{\Atl_\Ob}(\seqname{E}, \seqname{C}) \defeq
\set
{{
   \bigl ( \seqname{A}, E, (C,\catname{C}) \bigr )
}}%
[
  E \in \seqname{E},
  {(C,\catname{C}) \in \seqname{C}},
  {\maximal[S-]{\isAtl_\Ob} \bigl ( \seqname{A}, E, (C,\catname{C}) \bigr ) }
]*
\end{multline}
\begin{multline}
\maxfull[S-]{\Atl_\Ob}(\seqname{E}, \seqname{C}) \defeq
\set
{{
   \bigl ( \seqname{A}, (E,\catname{E}), (C,\catname{C}) \bigr )
}}%
[
  {(E,\catname{E}) \in \seqname{E}},
  {(C,\catname{C}) \in \seqname{C}},
  {\maxfull[S-]{\isAtl_\Ob} \bigl ( \seqname{A}, (E,\catname{E}), (C,\catname{C}) \bigr ) }
]*
\end{multline}
\end{definition}

\begin{lemma}[Maximal m-atlases are semi-maximal m-atlases]
Let $\seqname{E}=(E, \catname{E})$ and $\seqname{C}=(C,\catname{C})$ be
model spaces and $\seqname{A}$ a maximal m-atlas of $\seqname{E}$ in the
coordinate space $\seqname{C}$. Then $\seqname{A}$ is a semi-maximal
m-atlas of $\seqname{E}$ in the coordinate space $\seqname{C}$.

\begin{proof}
Let $(U,V,\phi) \in \seqname{A}$, $U' \subseteq U, V' \subseteq V$
and $V'' \objin \catname{C}$ be model neighborhoods, $\phi[U'] = V'$ and
$\phi' \maps V' \toiso V''$ be an isomorphism of $\catname{C}$.
$(U', V', \phi)$ is a subchart of $(U,V,\phi)$ and by
\pagecref{lem:M-Compat} is m-compatible with the charts of
$\seqname{A}$.  Since $\phi'$ is a model homeorphism,
$(U', V'', \phi' \compose \phi)$ is m-compatible with the charts of
$\seqname{A}$.  Since $\seqname{A}$ is maximal,
$(U', V'', \phi' \compose \phi)$ is a chart of $\seqname{A}$.
\end{proof}

Let $E$ be a topological space, $\seqname{C}=(C,\catname{C})$ a model
space and $\seqname{A}$ a maximal m-atlas of $E$ in the coordinate space
$\seqname{C}$. Then $\seqname{A}$ is a semi-maximal m-atlas of $E$ in
the coordinate space $\seqname{C}$.

\begin{proof}
Let $(U,V,\phi) \in \seqname{A}$, $U' \subseteq U$ open, $V' \subseteq V$
and $V'' \objin \catname{C}$ be model neighborhoods, $\phi[U'] = V'$ and
$\phi' \maps V' \toiso V''$ be an isomorphism of $\catname{C}$.
$(U', V', \phi)$ is a subchart of $(U,V,\phi)$ and by
\pagecref{lem:M-Compat} is m-compatible with the charts of
$\seqname{A}$.  Since $\phi'$ is a model homeorphism,
$(U', V'', \phi' \compose \phi)$ is m-compatible with the charts of
$\seqname{A}$.  Since $\seqname{A}$ is maximal,
$(U', V'', \phi' \compose \phi)$ is a chart of $\seqname{A}$.
\end{proof}
\end{lemma}

\begin{theorem}[Existence and uniqueness of maximal m-atlases]
Let $\seqname{A}$ be an m-atlas of $\seqname{E} \defeq (E, \catname{E})$
in the coordinate space $\seqname{C} \defeq (C,\catname{C})$. Then there
exists a unique maximal m-atlas
$\maximal{\Atlas}(\seqname{A}, \seqname{E}, \seqname{C})$ of $\seqname{E}$
in the coordinate space $\seqname{C}$ m-compatible with $\seqname{A}$.

\begin{proof}
Let $\seqname{P}$ be the set of all m-atlases of $\seqname{E}$ in the
coordinate space $\seqname{C}$ containing $\seqname{A}$ and m-compatible
in the coordinate space $\seqname{C}$ with all of the m-charts in
$\seqname{A}$. Let $\maximal{\seqname{P}}$ be a maximal chain of
$\seqname{A}$. Then $A'=\union{\maximal{\seqname{P}}}$ is a maximal
m-atlas of $\seqname{E}$ in the coordinate space $\seqname{C}$
m-compatible with $\seqname{A}$.
Uniqueness follows from \pagecref{m-atl:extensions}.
\end{proof}

Let $\seqname{A}$ be an m-atlas of $E$
in the coordinate space $\seqname{C} \defeq (C,\catname{C})$. Then there
exists a unique maximal m-atlas
$\maximal{\Atlas}(\seqname{A}, E, \seqname{C})$ of $E$
in the coordinate space $\seqname{C}$ m-compatible with $\seqname{A}$.

\begin{proof}
Since $\seqname{A}$ is an m-atlas of $E$ in the coordinate space
$\seqname{C} \defeq (C,\catname{C})$, then by \pagecref{cor:triv}
$\seqname{A}$ is an m-atlas of $\Triv{E}$ in the coordinate space
$\seqname{C}$. Let $\seqname{A}'$ be a maximal m-atlas of $\Triv{E}$ in
the coordinate space $\seqname{C}$. Then $\seqname{A}$ is an m-atlas of
$E$ in the coordinate space $\seqname{C}$.
\end{proof}

\end{theorem}

\subsection{M-atlas morphisms and functors}
\label{sec:M-ATLmorph}
\begin{definition}[M-atlas morphisms]
\label{def:M-ATLmorph}
Let $\catname{E}^i$, $\catname{C}^i$, $i=1,2$ be model categories,
$\seqname{E}^i \in \Ob(\catname{E}^i)$,
$\seqname{C}^i \in \Ob(\catname{C}^i)$ and $\seqname{A}^i$ be m-atlases
of $\seqname{E}^i$ in the coordinate spaces $\seqname{C}^i$.
$
  \funcseqname{f} \defeq
  (
    \funcname{f}_0 \maps \seqname{E}^1 \to \seqname{E}^2,
    \funcname{f}_1 \maps \seqname{C}^1 \to \seqname{C}^2
  )
$
is a (strict) $\seqname{E}^1$-$\seqname{E}^2$ m-atlas morphism of
$\seqname{A}^1$ to $\seqname{A}^2$ in the coordinate spaces
$\seqname{C}^1$, $\seqname{C}^2$, abbreviated as
$\isAtl_\Ar
  (
    \seqname{A}^1,
    \seqname{E}^1,
    \seqname{C}^1,
    \seqname{A}^2,
    \seqname{E}^2,
    \seqname{C}^2
    \funcname{f}_0,
    \funcname{f}_1
  )
$,
and a (strict) $\catname{E}^1$-$\catname{E}^2$ m-atlas morphism of
$\seqname{A}^1$ to $\seqname{A}^2$ in the coordinate model categories
$\catname{C}^1$, $\catname{C}^2$, abbreviated as
$
  \isAtl_\Ar
  (
    \seqname{A}^1,
    \catname{E}^1,
    \catname{C}^1,
    \seqname{A}^2,
    \catname{E}^2,
    \catname{C}^2
    \funcname{f}_0,
    \funcname{f}_1
  )
$,
iff for any
$(U^1, V^1, \phi^1 \maps U^1 \toiso V^1) \in \seqname{\seqname{A}}^1$,
$(U^2 ,V^2, \phi^2 \maps U^2 \toiso V^2) \in \seqname{\seqname{A}}^2$,
the diagram
$
D \defeq
\bigl (
  \{I \defeq U^1 \cap \funcname{f}_0^{-1}[U^2], V^1, E^2, U^2, V^2 \}
$,
$
  \{ \funcname{f}_0, \phi^2, \phi^1, \funcname{f}_1 \}
\bigr )
$
is M-locally nearly commutative in $\seqname{C}^2$, i.e., for any
$x \in I$ there are objects $U'^1 \subseteq I$, $V'^1 \subseteq V^1$,
$U'^2 \subseteq U^2$, $V'^2 \subseteq V^2$, $\hat{V}'^2 \subseteq C^2$
and an isomorphism $\hat{\funcname{f}} \maps V'^2 \toiso \hat{V}'^2$
such that \crefrange{eq:xin}{eq:fphigeqgphi} below hold.
The triple
$
  \bigl (
    \funcseqname{f},
    (\seqname{A}^1, \seqname{E}^1, \seqname{C}^1),
    (\seqname{A}^2, \seqname{E}^2, \seqname{C}^2)
  \bigr )
$
will refer to $\funcseqname{f}$ considered as an
$\seqname{E}^1$-$\seqname{E}^2$ m-atlas morphism of $\seqname{A}^1$ to
$\seqname{A}^2$ in the coordinate spaces
$\seqname{C}^1$, $\seqname{C}^2$.

\begin{figure}
\[ \bfig
\node uint(0,0)[{I = U^1 \cap \funcname{f}^{-1}_0[U^2]}]
\node u2(2000,0)[U^2]
\node v1(0,-1000)[V^1]
\node c21(1000,-1000)[C^2]
\node v2(2000,-1000)[V^2]
\arrow |a|[uint`u2;\funcname{f}_0]
\arrow |r|/>->/[uint`v1;\phi^1]
\arrow |r|[u2`v2;\phi^2]
\arrow |l|//[u2`v2;\iso]
\arrow |a|[v1`c21;\funcname{f}_1]
\efig \]
\caption{Uncompleted atlas morphism}
\label{fig:AtlM1}
\end{figure}

\begin{figure}
\[ \bfig
\node x(500,0)[x]
\node uint(-300,-500)[{I=U^1 \cap \funcname{f}^{-1}_0[U^2]}]
\node up1(500,-500)[U'^1]
\node up2(1500,-500)[U'^2]
\node u2(2000,-500)[U^2]
\node v1(-300,-1000)[V^1]
\node vp1(500,-1000)[V'^1]
\node vhp2(1000,-1000)[\hat{V}'^2]
\node vp2(1500,-1000)[V'^2]
\node v2(2000,-1000)[V^2]
\arrow |r|/^{ (}->/[x`up1;\funcname{i}]
\arrow |r|/>->/[uint`v1;\phi^1]
\arrow |a|/_{ (}->/[up1`uint;\funcname{i}]
\arrow |a|[up1`up2;\funcname{f}_0]
\arrow |r|[up1`vp1;{\phi^1}]
\arrow |l|//[up1`vp1;\iso]
\arrow |a|/^{ (}->/[up2`u2;i]
\arrow |r|[up2`vp2;\phi^2]
\arrow |l|//[up2`vp2;\iso]
\arrow |r|[u2`v2;\phi^2]
\arrow |l|//[u2`v2;\iso]
\arrow |a|/_{ (}->/[vp1`v1;\funcname{i}]
\arrow |a|[vp1`vhp2;\funcname{f}_1]
\arrow |a|/ >.>>/[vp2`vhp2;\hat{\funcname{f}}]
\arrow |b|//[vp2`vhp2;\iso]
\arrow |a|/^{ (}->/[vp2`v2;\funcname{i}]
\efig \]
\caption{Completed atlas morphism}
\label{fig:AtlM2}
\end{figure}

It is a semi-strict $\seqname{E}^1$-$\seqname{E}^2$ m-atlas morphism of
$\seqname{A}^1$ to $\seqname{A}^2$ in the coordinate spaces
$\seqname{C}^1$, $\seqname{C}^2$, abbreviated as
$
  \strict[semi-]{\isAtl_\Ar}
  (
    \seqname{A}^1,
    \seqname{E}^1,
    \seqname{C}^1,
    \seqname{A}^2,
    \seqname{E}^2,
    \seqname{C}^2,
    \funcname{f}_0,
    \funcname{f}_1
  )
$,
iff $\seqname{E}^1\submod \seqname{E}^2$,
$\seqname{C}^1\submod \seqname{C}^2$,
$\funcname{f}_0$ is locally an m-morphism of $\seqname{E}^2$ and
$\funcname{f}_1$ is locally an m-morphism of $\seqname{C}^2$.

It is a semi-strict $\seqname{E}^1$-$\seqname{E}^2$ m-atlas morphism of
$\seqname{A}^1$ to $\seqname{A}^2$ in the coordinate space
$\seqname{C}^1$, abbreviated as
$
  \strict[semi-]{\isAtl_\Ar}
  (
    \seqname{A}^1,
    \seqname{E}^1,
    \seqname{C}^1,
    \seqname{A}^2,
    \seqname{E}^2,
    \funcname{f}_0,
    \funcname{f}_1
  )
$,
iff it is a semi-strict $\seqname{E}^1$-$\seqname{E}^2$ m-atlas morphism
of $\seqname{A}^1$ to $\seqname{A}^2$ in the coordinate spaces
$\seqname{C}^1$, $\seqname{C}^1$.

It is a semi-strict $\seqname{E}^1$ m-atlas morphism of $\seqname{A}^1$
to $\seqname{A}^2$ in the coordinate space $\seqname{C}^1$, abbreviated
as
$
  \strict[semi-]{\isAtl_\Ar}
  (
    \seqname{A}^1,
    \seqname{E}^1,
    \seqname{C}^1,
    \seqname{A}^2,
    \funcname{f}_0,
    \funcname{f}_1
  )
$,
iff it is a semi-strict $\seqname{E}^1$-$\seqname{E}^1$ m-atlas morphism
of $\seqname{A}^1$ to $\seqname{A}^2$ in the coordinate spaces
$\seqname{C}^1$, $\seqname{C}^1$.

It is a semi-strict
$\seqname{E}^1$-$\seqname{E}^2$-$\catname{C}^1$-$\catname{C}^2$ m-atlas
morphism of $\seqname{A}^1$ to $\seqname{A}^2$ in the coordinate spaces
$\seqname{C}^1$, $\seqname{C}^2$, abbreviated as
$
  \strict[semi-]{\isAtl_\Ar}
  (
    \seqname{A}^1,
    \seqname{E}^1,
    \seqname{C}^1,
    \seqname{A}^2,
    \seqname{E}^2,
    \seqname{C}^2,
    \funcname{f}_0,
    \funcname{f}_1,
    \catname{C}^1,
    \catname{C}^2
  )
$,
iff $\seqname{E}^1\submod \seqname{E}^2$,
$\catname{C}^1 \subcat[full-] \catname{C}^2$,
$\funcname{f}_0$ is locally an m-morphism of $\seqname{E}^1$ to
$\seqname{E}^2$ and
$\funcname{f}_1$ is locally a $\catname{C}^1$-$\catname{C}^2$ morphism of
$\seqname{C}^1$ to $\seqname{C}^2$.

It is a semi-strict $\seqname{E}^1$-$\catname{C}^1$-$\catname{C}^2$
m-atlas morphism of $\seqname{A}^1$ to $\seqname{A}^2$ in the coordinate
spaces $\seqname{C}^1$, $\seqname{C}^2$, abbreviated as
$
  \strict[semi-]{\isAtl_\Ar}
  (
    \seqname{A}^1,
    \seqname{E}^1,
    \seqname{C}^1,
    \seqname{A}^2,
    \seqname{C}^2,
    \funcname{f}_0,
    \funcname{f}_1,
    \catname{C}^1,
    \catname{C}^2
  )
$,
iff it is a semi-strict
$\seqname{E}^1$-$\seqname{E}^1$-$\catname{C}^1$-$\catname{C}^2$ m-atlas
morphism of $\seqname{A}^1$ to $\seqname{A}^2$ in the coordinate spaces
$\seqname{C}^1$, $\seqname{C}^2$.

It is a semi-strict
$\catname{E}^1$-$\catname{E}^2$-$\seqname{E}^1$-$\seqname{E}^2$ m-atlas
morphism of $\seqname{A}^1$ to $\seqname{A}^2$ in the coordinate spaces
$\seqname{C}^1$, $\seqname{C}^2$, abbreviated as \\
$
  \strict[semi-]{\isAtl_\Ar}
  (
    \seqname{A}^1,
    \seqname{E}^1,
    \seqname{C}^1,
    \seqname{A}^2,
    \seqname{E}^2,
    \seqname{C}^2,
    \funcname{f}_0,
    \funcname{f}_1,
    \catname{E}^1,
    \catname{E}^2
  )
$,
iff $\catname{E}^1 \subcat[full-] \catname{E}^2$,
$\seqname{C}^1 \submod \seqname{C}^2,$ $\funcname{f}_0$ is locally
an $\catname{E}^1$-$\catname{E}^2$ morphism of $\seqname{E}^1$ to
$\seqname{E}^2$ and $\funcname{f}_1$ is locally an m-morphism of
$\seqname{C}^1$ to $\seqname{C}^2$.

It is a semi-strict
$\catname{E}^1$-$\catname{E}^2$-$\seqname{E}^1$-$\seqname{E}^2$ m-atlas
morphism of $\seqname{A}^1$ to $\seqname{A}^2$ in the coordinate space
$\seqname{C}^1$, abbreviated as
$
  \strict[semi-]{\isAtl_\Ar}
  (
    \seqname{A}^1,
    \seqname{E}^1,
    \seqname{C}^1,
    \seqname{A}^2,
    \seqname{E}^2,
    \funcname{f}_0,
    \funcname{f}_1,
    \catname{E}^1,
    \catname{E}^2
  )
$,
iff it is a semi-strict
$\catname{E}^1$-$\catname{E}^2$-$\seqname{E}^1$-$\seqname{E}^2$ m-atlas
morphism of $\seqname{A}^1$ to $\seqname{A}^2$ in the coordinate spaces
$\seqname{C}^1$, $\seqname{C}^1$.

It is a semi-strict
$\catname{E}^1$-$\catname{E}^2$-%
$\seqname{E}^1$-$\seqname{E}^2$-%
$\catname{C}^1$-$\catname{C}^2$
m-atlas morphism of $\seqname{A}^1$ to $\seqname{A}^2$ in the coordinate
spaces $\seqname{C}^1$, $\seqname{C}^2$, abbreviated as \
$
  \strict[semi-]{\isAtl_\Ar}
  (
    \seqname{A}^1,
    \seqname{E}^1,
    \seqname{C}^1,
    \seqname{A}^2,
    \seqname{E}^2,
    \seqname{C}^2,
    \funcname{f}_0,
    \funcname{f}_1,
    \catname{E}^1,
    \catname{E}^2,
    \catname{C}^1
    \catname{C}^2
  )
$,
iff $\catname{E}^1 \subcat[full-] \catname{E}^2$,
$\catname{C}^1 \subcat[full-] \catname{C}^2$, $\funcname{f}_0$ is locally an $\catname{E}^1$-$\catname{E}^2$
morphism of $\seqname{E}^1$ to $\seqname{E}^2$
and $\funcname{f}_1$ is locally a $\catname{C}^1$-$\catname{C}^2$ morphism of
$\seqname{C}^1$ to $\seqname{C}^2$.

It is a strict $\seqname{E}^1$-$\seqname{E}^2$ m-atlas morphism of
$\seqname{A}^1$ to $\seqname{A}^2$ in the coordinate spaces
$\seqname{C}^1$, $\seqname{C}^2$, abbreviated as
$
  \strict{\isAtl_\Ar}
  (
    \seqname{A}^1,
    \seqname{E}^1,
    \seqname{C}^1,
    \seqname{A}^2,
    \seqname{E}^2,
    \seqname{C}^2,
    \funcname{f}_0,
    \funcname{f}_1
  )
$,
iff $\seqname{E}^1\submod \seqname{E}^2$,
$\seqname{C}^1\submod \seqname{C}^2$,
$\funcname{f}_0$ is a morphism of $\seqname{E}^2$ and
$\funcname{f}_1$ is a morphism of $\seqname{C}^2$.

It is a strict $\seqname{E}^1$ m-atlas morphism of $\seqname{A}^1$ to
$\seqname{A}^2$ in the coordinate space $\seqname{C}^1$, abbreviated as
$
  \strict{\isAtl_\Ar}
  (
    \seqname{A}^1,
    \seqname{E}^1,
    \seqname{C}^1,
    \seqname{A}^2,
    \funcname{f}_0,
    \funcname{f}_1
  )
$,
iff
it is a strict $\seqname{E}^1$-$\seqname{E}^1$ m-atlas morphism of
$\seqname{A}^1$ to $\seqname{A}^2$ in the coordinate spaces
$\seqname{C}^1$, $\seqname{C}^1$.

It is a strict
$\seqname{E}^1$-$\seqname{E}^2$-$\catname{C}^1$-$\catname{C}^2$ m-atlas
morphism of $\seqname{A}^1$ to $\seqname{A}^2$ in the coordinate spaces
$\seqname{C}^1$, $\seqname{C}^2$, abbreviated as
$
  \strict{\isAtl_\Ar}
  (
    \seqname{A}^1,
    \seqname{E}^1,
    \seqname{C}^1,
    \seqname{A}^2,
    \seqname{E}^2,
    \seqname{C}^2,
    \funcname{f}_0,
    \funcname{f}_1,
    \catname{C}^1,
    \catname{C}^2
  )
$,
iff
$\seqname{E}^1\submod \seqname{E}^2$,
$\catname{C}^1 \subcat[full-] \catname{C}^2$,
$\funcname{f}_0$ is a morphism of $\seqname{E}^2$ and
$\funcname{f}_1$ is a morphism of $ \catname{C}^2$.

It is a strict $\seqname{E}^1$-$\catname{C}^1$-$\catname{C}^2$ m-atlas
morphism of $\seqname{A}^1$ to $\seqname{A}^2$ in the coordinate spaces
$\seqname{C}^1$, $\seqname{C}^2$, abbreviated as

$
  \strict{\isAtl_\Ar}
  (
    \seqname{A}^1,
    \seqname{E}^1,
    \seqname{C}^1,
    \seqname{A}^2,
    \seqname{C}^2,
    \funcname{f}_0,
    \funcname{f}_1,
    \catname{C}^1,
    \catname{C}^2
  )
$,
iff it is a strict
$\seqname{E}^1$-$\seqname{E}^1$-$\catname{C}^1$-$\catname{C}^2$ m-atlas
morphism of $\seqname{A}^1$ to $\seqname{A}^2$ in the coordinate spaces
$\seqname{C}^1$, $\seqname{C}^2$.

It is a strict
$\catname{E}^1$-$\catname{E}^2$-$\seqname{E}^1$-$\seqname{E}^2$ m-atlas
morphism of $\seqname{A}^1$ to $\seqname{A}^2$ in the coordinate spaces
$\seqname{C}^1$, $\seqname{C}^2$, abbreviated as
$
  \strict{\isAtl_\Ar}
  (
    \seqname{A}^1,
    \seqname{E}^1,
    \seqname{C}^1,
    \seqname{A}^2,
    \seqname{E}^2,
    \seqname{C}^2,
    \funcname{f}_0,
    \funcname{f}_1,
    \catname{E}^1,
    \catname{E}^2
  )
$,
iff $\catname{E}^1 \subcat[full-] \catname{E}^2$,
$\seqname{C}^1 \submod \seqname{C}^2,$ $\funcname{f}_0$ is a morphism of
$\catname{E}^2$ and $\funcname{f}_1$ is a morphism of $\seqname{C}^2$.

It is a strict
$\catname{E}^1$-$\catname{E}^2$-$\seqname{E}^1$-$\seqname{E}^2$ m-atlas
morphism of $\seqname{A}^1$ to $\seqname{A}^2$ in the coordinate space
$\seqname{C}^1$, abbreviated as
$
  \strict{\isAtl_\Ar}
  (
    \seqname{A}^1,
    \seqname{E}^1,
    \seqname{C}^1,
    \seqname{A}^2,
    \seqname{E}^2,
    \funcname{f}_0,
    \funcname{f}_1,
    \catname{E}^1,
    \catname{E}^2
  )
$,
iff it is a strict
$\catname{E}^1$-$\catname{E}^2$-$\seqname{E}^1$-$\seqname{E}^2$ m-atlas
morphism of $\seqname{A}^1$ to $\seqname{A}^2$ in the coordinate spaces
$\seqname{C}^1$, $\seqname{C}^1$.

It is a strict
$\catname{E}^1$-$\catname{E}^2$-%
$\seqname{E}^1$-$\seqname{E}^2$-%
$\catname{C}^1$-$\catname{C}^2$
m-atlas morphism of $\seqname{A}^1$ to $\seqname{A}^2$ in the coordinate
spaces $\seqname{C}^1$, $\seqname{C}^2$, abbreviated as
$
  \strict{\isAtl_\Ar}
  (
    \seqname{A}^1,
    \seqname{E}^1,
    \seqname{C}^1,
    \seqname{A}^2,
    \seqname{E}^2,
    \seqname{C}^2,
    \funcname{f}_0,
    \funcname{f}_1,
    \catname{E}^1,
    \catname{E}^2,
    \catname{C}^1
    \catname{C}^2
  )
$,
iff $\catname{E}^1 \subcat[full-] \catname{E}^2$,
$\catname{C}^1 \subcat[full-] \catname{C}^2$, $\funcname{f}_0$ is a
morphism of $\catname{E}^2$ and $\funcname{f}_1$ is a morphism of
$\catname{C}^2$.

\begin{equation}
x \in U'^1
\label{eq:xin}
\end{equation}
\begin{equation}
\funcname{f}_1 \compose \phi^1(x) \in \hat{V}'^2
\label{eq:fphixin}
\end{equation}
\begin{equation}
\funcname{f}_0[U'^1] \subseteq U'^2
\end{equation}
\begin{equation}
\phi^1[U'^1] \subseteq V'^1
\end{equation}
\begin{equation}
\funcname{f}_1[V'^1] \subseteq \hat{V}'^2
\end{equation}
\begin{equation}
\phi^2[U'^2] \subseteq V'^2
\end{equation}
\begin{equation}
\hat{\funcname{f}} \compose \phi^2 \compose \funcname{f}_0 =
\funcname{f}_1 \compose \phi^1
\label{eq:fphigeqgphi}
\end{equation}

The identity morphism of $(\seqname{A}^i, \seqname{E}^i, \seqname{C}^i)$
is
\begin{equation}
\Id_{(\seqname{A}^i, \seqname{E}^i, \seqname{C}^i)} \defeq
\bigl (
  (\Id_{\seqname{E}^i}, \Id_{\seqname{C}^i}),
  (\seqname{A}^i, \seqname{E}^i, \seqname{C}^i),
  (\seqname{A}^i, \seqname{E}^i, \seqname{C}^i)
\bigr )
\end{equation}

Let $\catname{C}^i$, $i=1,2$, be model categories,
$\seqname{C}^i \objin \catname{C}^i$, $E^i$ be topological spaces and
$\seqname{A}^i$ be m-atlases of $E^i$ in the coordinate spaces
$\seqname{C}^i$.
A pair of functions
$
  \funcseqname{f} \defeq \\
  (
    \funcname{f}_0 \maps E^1 \to E^2,
    \funcname{f}_1 \maps \seqname{C}^1 \to \seqname{C}^2
  )
$,
is a (strict) $E^1$-$E^2$ m-atlas morphism of $\seqname{A}^1$ to
$\seqname{A}^2$ in the coordinate spaces $\seqname{C}^i$, abbreviated as
$
 \isAtl_\Ar(\seqname{A}^1, E^1, \seqname{C}^1$,
            $\seqname{A}^2$, $E^2$, $\seqname{C}^2$,
            $\funcname{f}_0$, $\funcname{f}_1)$,
and a (strict)
$E^1$-$E^2$ m-atlas morphism of $\seqname{A}^1$ to $\seqname{A}^2$ in
the coordinate model categories $\catname{C}^i$, abbreviated as
$\isAtl_\Ar(\seqname{A}^1, E^1, \catname{C}^1)$,
            $\seqname{A}^2$, $E^2$, $\catname{C}^2$,
iff it is a (strict) $\Triv{E^1}$-$\Triv{E^2}$ m-atlas morphism of
$\seqname{A}^1$ to $\seqname{A}^2$ in the coordinate model categories
$\catname{C}^1$, $\catname{C}^2$.
The triple
$
  \bigl (
    \funcseqname{f},
    (\seqname{A}^1, E^1, \seqname{C}^1),
    (\seqname{A}^2, E^2, \seqname{C}^2)
  \bigr )
$
will refer to $\funcseqname{f}$ considered as an $E^1$-$E^2$ m-atlas
morphism of $\seqname{A}^1$ to $\seqname{A}^2$ in the coordinate spaces
$\seqname{C}^1$, $\seqname{C}^2$.

It is a semi-strict $E^1$-$E^2$ m-atlas morphism of $\seqname{A}^1$ to
$\seqname{A}^2$ in the coordinate spaces $\seqname{C}^1$,
$\seqname{C}^2$, abbreviated as
$
  \strict[semi-]{\isAtl_\Ar}
  (
    \seqname{A}^1,
    E^1,
    \seqname{C}^1,
    \seqname{A}^2,
    E^2,
    \seqname{C}^2,
    \funcname{f}_0,
    \funcname{f}_1
  )
$,
iff $\seqname{C}^1 \submod \seqname{C}^2$ and $\funcname{f}_1$ is
locally an m-morphism of $\seqname{C}^2$.

It is a semi-strict $E^1$-$E^2$ m-atlas morphism of $\seqname{A}^1$ to
$\seqname{A}^2$ in the coordinate space $\seqname{C}^1$, abbreviated as
$
  \strict[semi-]{\isAtl_\Ar}
  (
    \seqname{A}^1,
    E^1,
    \seqname{C}^1,
    \seqname{A}^2,
    E^2,
    \funcname{f}_0,
    \funcname{f}_1
  )
$,
iff it is a semi-strict $E^1$-$E^2$ m-atlas morphism of $\seqname{A}^1$
to $\seqname{A}^2$ in the coordinate spaces $\seqname{C}^1$,
$\seqname{C}^1$.

It is a semi-strict $E^1$-$E^2$-$\catname{C}^1$-$\catname{C}^2$ m-atlas
morphism of $\seqname{A}^1$ to $\seqname{A}^2$ in the coordinate spaces
$\seqname{C}^1$, $\seqname{C}^2$, abbreviated as
$
  \strict[semi-]{\isAtl_\Ar}
  (
    \seqname{A}^1,
    E^1,
    \seqname{C}^1,
    \seqname{A}^2,
    E^2,
    \seqname{C}^2,
    \funcname{f}_0,
    \funcname{f}_1,
    \catname{C}^1,
    \catname{C}^2
  )
$,
iff $\catname{C}^1 \subcat[full-] \catname{C}^2$ and $\funcname{f}_1$ is a
morphism of $\catname{C}^2$.

It is a semi-strict $E^1$-$E^2$-$\catname{C}^1$ m-atlas
morphism of $\seqname{A}^1$ to $\seqname{A}^2$ in the coordinate spaces
$\seqname{C}^1$, abbreviated as
$
  \strict[semi-]{\isAtl_\Ar}
  (
    \seqname{A}^1,
    E^1,
    \seqname{C}^1,
    \seqname{A}^2,
    E^2,
    \funcname{f}_0,
    \funcname{f}_1,
    \catname{C}^1
  )
$,
iff
it is a semi-strict $E^1$-$E^2$-$\catname{C}^1$-$\catname{C}^1$ m-atlas
morphism of $\seqname{A}^1$ to $\seqname{A}^2$ in the coordinate spaces
$\seqname{C}^1$, $\seqname{C}^1$.

It is a strict $E^1$-$E^2$ m-atlas morphism of $\seqname{A}^1$ to
$\seqname{A}^2$ in the coordinate spaces $\seqname{C}^1$,
$\seqname{C}^2$, abbreviated as
$
  \strict{\isAtl_\Ar}
  (
    \seqname{A}^1,
    E^1,
    \seqname{C}^1,
    \seqname{A}^2,
    E^2,
    \seqname{C}^2,
    \funcname{f}_0,
    \funcname{f}_1
  )
$,
iff $\seqname{C}^1\submod \seqname{C}^2$ and $\funcname{f}_1$ is a
morphism of $\seqname{C}^2$.

It is a strict $E^1$-$E^2$ m-atlas morphism of $\seqname{A}^1$ to
$\seqname{A}^2$ in the coordinate space $\seqname{C}^1$, abbreviated as

$
  \strict{\isAtl_\Ar}
  (
    \seqname{A}^1,
    E^1,
    \seqname{C}^1,
    \seqname{A}^2,
    E^2,
    \funcname{f}_0,
    \funcname{f}_1
  )
$,
iff it is a strict $E^1$-$E^2$ m-atlas morphism of $\seqname{A}^1$ to
$\seqname{A}^2$ in the coordinate spaces $\seqname{C}^1$,
$\seqname{C}^1$.

It is a strict $E^1$-$E^2$-$\catname{C}^1$-$\catname{C}^2$ m-atlas
morphism of $\seqname{A}^1$ to $\seqname{A}^2$ in the coordinate spaces
$\seqname{C}^1$, $\seqname{C}^2$, abbreviated as
$
  \strict{\isAtl_\Ar}
  (
    \seqname{A}^1,
    E^1,
    \seqname{C}^1,
    \seqname{A}^2,
    E^2,
    \seqname{C}^2,
    \funcname{f}_0,
    \funcname{f}_1,
    \catname{C}^1,
    \catname{C}^2
  )
$,
iff $\catname{C}^1 \subcat[full-] \catname{C}^2$ and $\funcname{f}_1$ is a
morphism of $\catname{C}^2$.

It is a strict $E^1$-$E^2$-$\catname{C}^1$-$\catname{C}^2$ m-atlas
morphism of $\seqname{A}^1$ to $\seqname{A}^2$ in the coordinate space
$\seqname{C}^1$, abbreviated as
$
  \strict{\isAtl_\Ar}
  (
    \seqname{A}^1,
    E^1,
    \seqname{C}^1,
    \seqname{A}^2,
    E^2,
    \funcname{f}_0,
    \funcname{f}_1,
    \catname{C}^1
  )
$,
iff it is a strict $E^1$-$E^2$-$\catname{C}^1$-$\catname{C}^1$ m-atlas
morphism of $\seqname{A}^1$ to $\seqname{A}^2$ in the coordinate spaces
$\seqname{C}^1$, $\seqname{C}^1$.

The identity morphism of $(\seqname{A}^i, E^i, \seqname{C}^i)$ is
\begin{equation}
\Id_{(\seqname{A}^i, E^i, \seqname{C}^i)} \defeq
\bigl (
  (\Id_{E^i}, \Id_{\seqname{C}^i}),
  (\seqname{A}^i, E^i, \seqname{C}^i),
  (\seqname{A}^i, E^i, \seqname{C}^i)
\bigr )
\end{equation}
\end{definition}

\begin{lemma}[M-atlas morphisms]
\label{lem:M-ATLmorph}
Let $\catname{E}^i$, $\catname{C}^i$,
$\catname{E}'^i$, $\catname{C}'^i$, $i=1,2,3$, be model categories,
$\catname{E}^i \subcat[full-] \catname{E}'^i$,
$\catname{C}^i \subcat[full-] \catname{C}'^i$.
$\seqname{E}^i \objin \catname{E}^i$,
$\seqname{C}^i \objin \catname{C}^i$, $\seqname{A}^i$ an m-atlas of
$\seqname{E}^i$ in the coordinate spaces $\seqname{C}^i$,
$\seqname{A}^2$ semi-maximal and
$
  \funcseqname{f}^i \defeq
  (
    \funcname{f}^i_0 \maps \seqname{E}^i \to \seqname{E}^{i+1},
    \funcname{f}^i_1 \maps \seqname{C}^i \to \seqname{C}^{i+1}
  )
$
a (strict) $\seqname{E}^i$-$\seqname{E}^{i+1}$ m-atlas morphism of
$\seqname{A}^i$ to $\seqname{A}^{i+1}$ in the coordinate spaces
$\seqname{C}^i$, $\seqname{C}^{i+1}$.

If $\catname{E}^i \subcat[full-] \catname{E}^{i+1}$,
$\catname{C}^i \subcat[full-] \catname{C}^{i+1}$ and
$\seqname{A}^i = \seqname{A}^{i+1}$ then
$
\ID_%
  {(\seqname{E}^i,\seqname{C}^i),(\seqname{E}^{i+1},\seqname{C}^{i+1})}
$
is a strict
$\catname{E}^i$-$\catname{E}^{i+1}$-%
$\seqname{E}^i$-$\seqname{E}^{i+1}$-%
$\catname{C}^i$-$\catname{C}^{i+1}$ m-atlas morphism of
$\seqname{A}^i$ to $\seqname{A}^{i+1}$ in the coordinate spaces
$\seqname{C}^i$, $\seqname{C}^{i+1}$.

\begin{proof}
A commutative diagram is M-locally commutative.
$
  \Id_{\seqname{E}^i} \arin \catname{E}^i
    \subcat[full-] \catname{E}^{i+1}
$ and
$
  \Id_{\seqname{C}^i} \arin \catname{C}^i
    \subcat[full-] \catname{C}^{i+1}
$.
\end{proof}

If $\funcseqname{f}^i$ is a semi-strict (strict)
$\catname{E}^i$-$\catname{E}^{i+1}$-$\seqname{E}^i$-$\seqname{E}^{i+1}$-$\catname{C}^i$-$\catname{C}^{i+1}$
m-atlas morphism of $\seqname{A}^i$ to $\seqname{A}^{i+1}$ in the
coordinate spaces $\seqname{C}^i$, $\seqname{C}^{i+1}$ then
$\funcseqname{f}^i$ is a semi-strict (strict)
$\catname{E}'^i$-$\catname{E}'^{i+1}$-$\seqname{E}^i$-$\seqname{E}^{i+1}$-$\catname{C}'^i$-$\catname{C}'^{i+1}$
m-atlas morphism of $\seqname{A}^i$ to $\seqname{A}^{i+1}$ in the
coordinate spaces $\seqname{C}^i$, $\seqname{C}^{i+1}$.

\begin{proof}
If $\funcseqname{f}^i_0$ is locally a morphism of $\catname{E}^i$ then
$\funcseqname{f}^i_0$ is locally a morphism of $\catname{E}'^i$. If
$\funcseqname{f}^i_1$ is locally a morphism of $\catname{C}^i$ then
$\funcseqname{f}^i_1$ is locally a morphism of $\catname{c}'^i$. If
$\funcseqname{f}^i_0 \arin \catname{E}^i$ then $\funcseqname{f}^i_0
\arin \catname{E}'^i$.  If $\funcseqname{f}^i_1 \arin \catname{C}^i$
then $\funcseqname{f}^i_1 \arin \catname{C}'^i$.
\end{proof}

$
  \funcseqname{f}^2 \compose[()] \funcseqname{f}^1 =
  (
    \funcname{f}^2_0 \compose \funcname{f}^1_0,
    \funcname{f}^2_1 \compose \funcname{f}^1_1
  )
$
is a (strict) $\seqname{E}^1$-$\seqname{E}^3$ m-atlas morphism of
$\seqname{A}^1$ to $\seqname{A}^3$ in the coordinate spaces
$\seqname{C}^1$, $\seqname{C}^3$.

Let $\catname{C}^i$, $i=1,2,3$, be a model category,
$E^i$ be a topological space,
$\seqname{C}^i \objin \catname{C}^i$, $\seqname{A}^i$ an m-atlas
of $E^i$ in the coordinate spaces $\seqname{C}^i$,
$\seqname{A}^2$ semi-maximal and
$
  \funcseqname{f}^i \defeq \\
  (
    \funcname{f}^i_0 \maps E^i \to E^{i+1},
    \funcname{f}^i_1 \maps \seqname{C}^i \to \seqname{C}^{i+1}
  )
$
a (strict) $E^i$-$E^{i+1}$ m-atlas morphism of
$\seqname{A}^i$ to $\seqname{A}^{i+1}$ in the coordinate spaces
$\seqname{C}^i$, $\seqname{C}^{i+1}$. Then
$
  \funcseqname{f}^2 \compose[()] \funcseqname{f}^1 =
  (
    \funcname{f}^2_0 \compose \funcname{f}^1_0,
    \funcname{f}^2_1 \compose \funcname{f}^1_1
  )
$
is a (strict) $E^1$-$E^3$ m-atlas morphism of
$\seqname{A}^1$ to $\seqname{A}^3$ in the coordinate spaces
$\seqname{C}^1$, $\seqname{C}^3$.

\begin{proof}
\begin{figure}
\[ \bfig
\node uint13(0,2500)[I^{1,3}]
\node uint12(500,2500)[I^{1,2}]
\node v1(1000,2500)[V^1]
\node c21(1000,2000)[C^2_1]
\node uint23(0,1500)[I^{2,3}]
\node u2(500,1500)[U^2]
\node v2(1000,1500)[V^2]
\arrow |r|[uint13`uint23;\funcname{f}^1_0]
\arrow |a|/^{ (}->/[uint13`uint12;\funcname{i}]
\arrow |r|[uint12`u2;\funcname{f}^1_0]
\arrow |a|/>->>/[uint12`v1;\phi^1]
\arrow |b|//[uint12`v1;\iso]
\arrow |a|/^{ (}->/[uint23`u2;\funcname{i}]
\arrow |a|/>->>/[u2`v2;\phi^2]
\arrow |b|//[u2`v2;\iso]
\arrow |r|[v1`c21;\funcname{f}^1_1]
\arrow |a/>->>/|[u2`v2;\phi^2]
\arrow |b|//[u2`v2;\iso]
\efig \]
\[ \bfig
\node uint23(0,1500)[I^{2,3}]
\node u2(500,1500)[]
\node v2(1000,1500)[V^2]
\node c31(1000,1000)[C^3_1]
\node u3l(0,500)[U^3]
\node u3(500,500)[U^3]
\node v3(1000,500)[V^3]
\arrow |r|[uint23`u3l;\funcname{f}^2_0]
\arrow |a|/^{ (}->/[uint23`u2;\funcname{i}]
\arrow |a|/>->>/[u2`v2;\phi^2]
\arrow |b|//[u2`v2;\iso]
\arrow |r|[u2`u3;\funcname{f}^2_0]
\arrow |a/>->>/|[u2`v2;\phi^2]
\arrow |b|//[u2`v2;\iso]
\arrow |a|/^{ (}->/[u3l`u3;\funcname{i}]
\arrow |a|/>->>/[u3`v3;\phi^2]
\arrow |b|//[u3`v3;\iso]
\arrow |r|[v2`c31;\funcname{f}^2_1]
\efig \]
\caption{Uncompleted atlas morphisms}
\label{fig:AtlM3}
\end{figure}

If each $\funcseqname{f}^i$ is a (strict)
$\seqname{E}^i$-$\seqname{E}^{i+1}$ m-atlas morphism of $\seqname{A}^i$
to $\seqname{A}^{i+1}$ in the coordinate spaces $\seqname{C}^i$,
$\seqname{C}^{i+1}$, then $\funcname{f}^2_0 \compose \funcname{f}^1_0$
and $\funcname{f}^2_1 \compose \funcname{f}^1_1$ are model functions by
\pagecref{lem:modcomp}\!. If each $\funcseqname{f}^i$ is a (strict)
$E^i$-$E^{i+1}$ m-atlas morphism of $\seqname{A}^i$ to
$\seqname{A}^{i+1}$ in the coordinate spaces $\seqname{C}^i$,
$\seqname{C}^{i+1}$, then $\funcname{f}^2_0 \compose \funcname{f}^1_0$
is continuous and $\funcname{f}^2_1 \compose \funcname{f}^1_1$ is a
model function by \cref{lem:modcomp}\!.

Let $(U^i, V^i, \phi^i)$ be charts in $\seqname{A}^i$,
$
  I^{1,2} \defeq
  U^1 \cap {\funcname{f}^1_0}^{-1}[U^2]
$
$
  I^{1,3} \defeq
  U^1 \cap (\funcname{f}^2_0 \compose \funcname{f}^1_0)^{-1}[U^3]
$
and
$
  I^{2,3} \defeq
  U^2 \cap {\funcname{f}^2_0}^{-1}[U^3]
$.
$
  I^{1,3} \subseteq U^1 \cap { \funcname{f}^1_0}^{-1}[U^2]
$.
If $I^{1,3} = \emptyset$ then
$
  D^{1,3} \defeq \\
    \bigl (
      \{I^{1,3}, V^1, E^3, U^3, V^3 \},
      \{ \funcname{f}_0, \phi^2, \phi^1, \funcname{f}_1 \}
    \bigr )
$
is vacuously M-locally nearly commutative.
Otherwise, since
$
  (
    \funcname{f}^1_0 \maps \seqname{E}^1 \to \seqname{E}^2,
    \funcname{f}^1_1 \maps \seqname{C}^1 \to \seqname{C}^2),
    \seqname{\seqname{A}}^1,
    \seqname{\seqname{A}}^2
  )
$
is a $\seqname{E}^1$-$\seqname{E}^2$ m-atlas morphism of
$\seqname{C}^1$, $\seqname{C}^2$,
$D \defeq \bigl ( \{I^{1,2}, V^1, E^2, U^2, V^2 \}$,
$\{ \funcname{f}_0, \phi^2, \phi^1, \funcname{f}_1 \} \bigr )$,
\fullcref{fig:AtlM3}\!,
is M-locally nearly commutative in $\seqname{C}^2$ and for any
$x \in I^{1,3} \subseteq I^{1,2}$
there are objects $U'^1 \subseteq I^{1,2}$,
$V'^1 \subseteq V^1$,
$U'^2 \subseteq U^2$, $V'^2 \subseteq V^2$,
$\hat{V}'^2 \subseteq C^2$ and an isomorphism
$\hat{\funcname{f}} \maps V'^2 \toiso \hat{V}'^2$ such that
{
  \showlabelsinline
  \crefrange{eq:xin}{eq:fphigeqgphi}
}
in \fullcref{def:M-ATLmorph} hold.

\begin{figure}
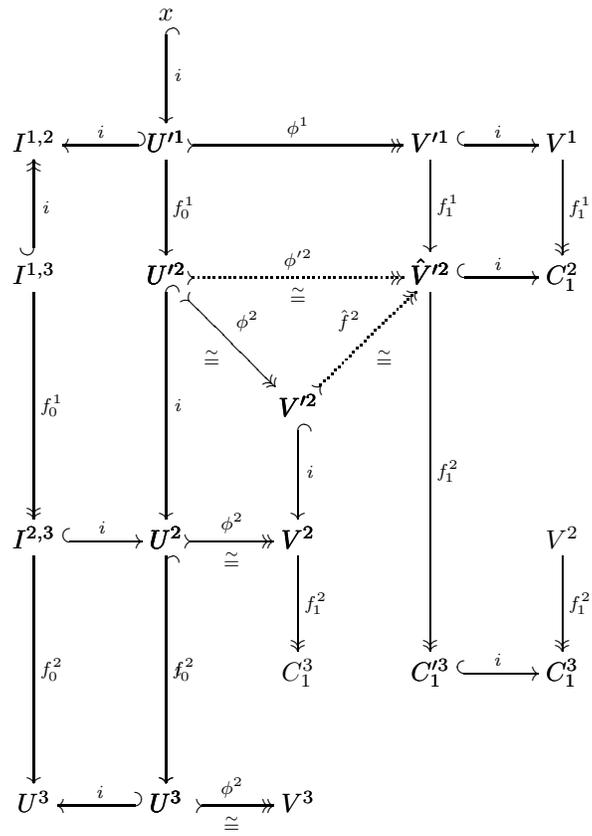

\[ \bfig
\node x(500,3500)[x]
\node uint12(0,3000)[I^{1,2}]
\node up1(500,3000)[U'^1]
\node vp1(1500,3000)[V'^1]
\node v1(2000,3000)[V^1]
\node uint13(0,2500)[I^{1,3}]
\node vhp2(1500,2500)[\hat{V}'^2]
\node c21(2000,2500)[C^2_1]
\node up2(500,2500)[U'^2]
\node vp2(1000,2000)[V'^2]
\node v2r(2000,1500)[V^2]
\node uint23(0,1500)[I^{2,3}]
\node u2(500,1500)[U^2]
\node v2(1000,1500)[V^2]
\node v2r(2000,1500)[V^2]
\node c31(1000,1000)[C^3_1]
\node cp31(1500,1000)[C'^3_1]
\node c31r(2000,1000)[C^3_1]
\node u3l(0,500)[U^3]
\node u3(500,500)[U^3]
\node v3(1000,500)[V^3]
\node v3r(2000,500)[V^3]
\arrow |r|/^{ (}->/[x`up1;\funcname{i}]
\arrow |r|/^{ (}->>/[uint13`uint12;\funcname{i}]
\arrow |a|/>->>/[up1`vp1;\phi^1]
\arrow |r|[up1`up2;\funcname{f}^1_0]
\arrow |r|[vp1`vhp2;\funcname{f}^1_1]
\arrow |r|/->>/[v1`c21;\funcname{f}^1_1]
\arrow |a|/_(->/[up1`uint12;\funcname{i}]
\arrow |a|/^{ (}->/[vp1`v1;\funcname{i}]
\arrow |r|/->>/[vhp2`cp31;\funcname{f}^2_1]
\arrow |a|/^{ (}->/[vhp2`c21;\funcname{i}]
\arrow |a|/>.>>/[vp2`vhp2;\hat{\funcname{f}}^2]
\arrow |b|//[vp2`vhp2;\iso]
\arrow |a|/^{ (}->/[uint23`u2;\funcname{i}]
\arrow |a|/>.>>/[up2`vhp2;\phi'^2]
\arrow |b|//[up2`vhp2;\iso]
\arrow |a|/ >->>/[up2`vp2;\phi^2]
\arrow |b|//[up2`vp2;\iso]
\arrow |r|/^(->/[up2`u2;\funcname{i}]
\arrow |r|[uint23`u3l;\funcname{f}^2_0]
\arrow |r|/->>/[uint13`uint23;\funcname{f}^1_0]
\arrow |a|/>->>/[u2`v2;\phi^2]
\arrow |b|//[u2`v2;\iso]
\arrow |r|/^{ (}->/[u2`u3;\funcname{i}]
\arrow |r|[u2`u3;\funcname{f}^2_0]
\arrow |r|/^{ (}->/[vp2`v2;\funcname{i}]
\arrow |a|/ >->>/[u3`v3;\phi^2]
\arrow |b|//[u3`v3;\iso]
\arrow |r|/->>/[v2`c31;\funcname{f}^2_1]
\arrow |r|/->>/[v2r`c31r;\funcname{f}^2_1]
\arrow |a|/^{ (}->/[cp31`c31r;\funcname{i}]
\arrow |a|/_{ (}->/[u3`u3l;\funcname{i}]
\efig \]
\caption{Partially completed atlas morphisms}
\label{fig:AtlM4}
\end{figure}

Let
$y = \funcname{f}^1_0(x)$, $I'^{2,3} = U'^2 \cap I^{2,3}$ and
$
  \phi'^2 \maps U'^2 \toiso \hat{V}'^2 \defeq
  \hat{\funcname{f}^2} \compose \phi^2
$.
Since $\seqname{A}^2$ is semi-maximal,
$\bigl ( U'^2, \hat{V}'^2, \phi'^2 \bigr )$
is a chart of $\seqname{A}^2$.

Since
$
  (
    \funcname{f}^2_0 \maps \seqname{E}^2 \to \seqname{E}^3,
    \funcname{f}^2_1 \maps \seqname{C}^2 \to \seqname{C}^3),
    \seqname{\seqname{A}}^2,
    \seqname{\seqname{A}}^3
  )
$
is a $\seqname{E}^2$-$\seqname{E}^3$ m-atlas morphism of
$\seqname{C}^2$, $\seqname{C}^3$, the diagram
$
D \defeq
\bigl (
  \{ I'^{2,3}, \hat{V}'^2, E^3, U^3, V^3 \}
$,
$
  \{ \funcname{f}^2_0, \phi^3, \phi'^2, \funcname{f}^2_1 \}
\bigr )
$
is M-locally nearly commutative in $\seqname{C}^3$ and hence there
are objects $U''^2 \subseteq I'^{2,3}$, $V''^2 \subseteq \hat{V}'^2$,
$U''^3 \subseteq U^3$, $V''^3 \subseteq V^3$,
$\hat{V}''^3 \subseteq C^3$ and an isomorphism
$\hat{\funcname{f}}^3 \maps V''^3 \toiso \hat{V}''^3$ such that
\crefrange{eq:yin2}{eq:fphigeqgphi2} below hold.

\begin{equation}
y = \funcname{f}^1_0(x) \in U''^2
\label{eq:yin2}
\end{equation}
\begin{equation}
\funcname{f}^2_1 \compose \phi'^2(y) \in \hat{V}''^3
\end{equation}
\begin{equation}
\funcname{f}^2_0[U''^2] \subseteq U''^3
\end{equation}
\begin{equation}
\phi'^2[U''^2] \subseteq V''^2
\end{equation}
\begin{equation}
\funcname{f}^2_1[V''^2] \subseteq \hat{V}''^3
\end{equation}
\begin{equation}
\phi^3[U''^3] \subseteq V''^3
\end{equation}
\begin{equation}
\hat{\funcname{f}}^3 \compose \phi^3 \compose \funcname{f}^2_0 =
\funcname{f}^2_1 \compose \phi'^2
\label{eq:fphigeqgphi2}
\end{equation}

Let
$
  I''^{1,3} \defeq
  U'^1 \cap (\funcname{f}^2_0 \compose \funcname{f}^1_0)^{-1}[U''^3]
$
and
$
  V''^1 \defeq
  V'^1 \cap (\funcname{f}^3_0 \compose \funcname{f}^2_0)^{-1}[V''^3]
$
Then \\
$\phi^1 \maps U'''^1 \toiso V''^1$,
$\hat{\funcname{f}}^2 \maps \hat{V}'^2 \toiso V'''^2$,
$\hat{\funcname{f}}^3 \maps \hat{V}''^3 \toiso V''^2$,
$\phi'^2 \maps U'''^2 \toiso V'''^2$ and
$\phi^3 \maps U'''^3 \toiso V'''^3$ are isomorphisms and

\begin{figure}
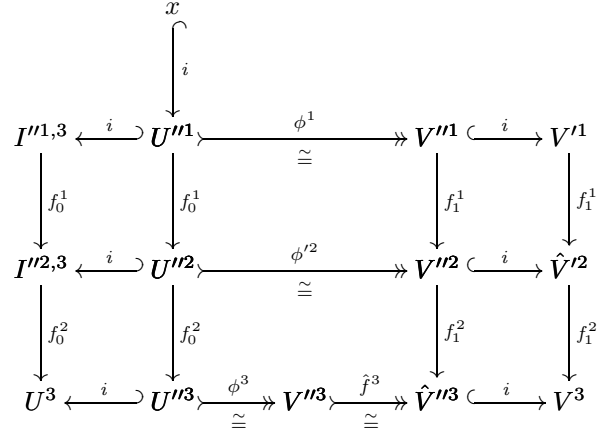

\[ \bfig
\node x(500,2000)[x]
\node uppint13(0,1500)[I''^{1,3}]
\node upp1(500,1500)[U''^1]
\node vpp1(1500,1500)[V''^1]
\node vp1(2000,1500)[V'^1]
\node uppint23(0,1000)[I''^{2,3}]
\node upp2(500,1000)[U''^2]
\node vpp2(1500,1000)[V''^2]
\node vhp2(2000,1000)[\hat{V}'^2]
\node u3(0,500)[U^3]
\node upp3(500,500)[U''^3]
\node vpp3(1000,500)[V''^3]
\node vhpp3(1500,500)[\hat{V}''^3]
\node v3(2000,500)[V^3]
\arrow |r|/^{ (}->/[x`upp1;\funcname{i}]
\arrow |r|[uppint13`uppint23;\funcname{f}^1_0]
\arrow |a|/_{ (}->/[upp1`uppint13;\funcname{i}]
\arrow |a|/>->>/[upp1`vpp1;\phi^1]
\arrow |b|//[upp1`vpp1;\iso]
\arrow |r|[upp1`upp2;\funcname{f}^1_0]
\arrow |r|[vpp1`vpp2;\funcname{f}^1_1]
\arrow |r|[vp1`vhp2;\funcname{f}^1_1]
\arrow |a|/^{ (}->/[vpp1`vp1;\funcname{i}]
\arrow |r|[uppint23`u3;\funcname{f}^2_0]
\arrow |a|/_{ (}->/[upp2`uppint23;\funcname{i}]
\arrow |r|[upp2`upp3;\funcname{f}^2_0]
\arrow |a|/>->>/[upp2`vpp2;\phi'^2]
\arrow |b|//[upp2`vpp2;\iso]
\arrow |r|[vpp2`vhpp3;\funcname{f}^2_1]
\arrow |r|[vhp2`v3;\funcname{f}^2_1]
\arrow |a|/^{ (}->/[vpp2`vhp2;\funcname{i}]
\arrow |a|/_{ (}->/[upp3`u3;\funcname{i}]
\arrow |a|/>->>/[upp3`vpp3;\phi^3]
\arrow |b|//[upp3`vpp3;\iso]
\arrow |a|/>->>/[vpp3`vhpp3;\hat{\funcname{f}}^3]
\arrow |b|//[vpp3`vhpp3;\iso]
\arrow |a|/^{ (}->/[vhpp3`v3;\funcname{i}]
\efig \]
\caption{Completed atlas morphisms}
\label{fig:AtlM5}
\end{figure}

\begin{equation}
x \in U'''^1
\label{eq:xin3}
\end{equation}
\begin{equation}
\funcname{f}^2_1 \compose \funcname{f}^1_1 \compose \phi^1(x) \in \hat{V}''^3
\end{equation}
\begin{equation}
\funcname{f}^2_0 \compose \funcname{f}^1_0[U''^1] \subseteq U''^3
\end{equation}
\begin{equation}
\phi^1[U''^1] \subseteq V''^1
\end{equation}
\begin{equation}
\funcname{f}^2_1 \compose \funcname{f}^1_1[V''^1] \subseteq \hat{V}''^3
\end{equation}
\begin{equation}
\phi^3[U''^3] \subseteq V''^3
\end{equation}
\begin{equation}
\begin{split}
& \hat{\funcname{f}}^3 \compose \phi^3 \compose \funcname{f}^2_0 \compose \funcname{f}^1_0 = \\*
& \funcname{f}^2_1 \compose \phi'^2 \compose \funcname{f}^1_0 =                              \\*
& \funcname{f}^2_1 \compose \funcname{f}^1_1 \compose \phi^1
\end{split}
\label{eq:fphigeqgphi3}
\end{equation}

\Cref{def:M-ATLmorph} has several definitions for a semi-strict and
strict m-atlas morphism, differing in the restrictions on
$\funcseqname{f}^i_j$.

If each $\seqname{E}^i \submod \seqname{E}^{i+1}$ and each
$\funcseqname{f}^i_0$ is locally a morphism of $\seqname{E}^i$ to
$\seqname{E}^{i+1}$, then by \pagecref{lem:modlocal}
$\funcseqname{f}^{i+1}_0 \compose \funcseqname{f}^i_0$ is locally a
morphism of $\seqname{E}^i$ to
of $\seqname{E}^{i+3}$.

If each $\catname{E}^i \subcat[full-] \catname{E}^{i+1}$ and each
$\funcseqname{f}^i_0$ is a morphism of $\catname{E}^{i+1}$, then
$\funcseqname{f}^{i+1}_0 \compose \funcseqname{f}^i_0$ is a morphism of
$\catname{E}^{i+2}$.

If each $\funcseqname{f}^i_0$ is continuous then
$\funcseqname{f}^{i+1}_0 \compose \funcseqname{f}^i_0$ is continuous.

If each $\seqname{C}^i \submod \seqname{C}^{i+1}$ and each
$\funcseqname{f}^i_1$ is locally a morphism of $\seqname{C}^i$ to
$\seqname{C}^{i+1}$, then by \cref{lem:modlocal}
$\funcseqname{f}^{i+1}_1 \compose \funcseqname{f}^i_1$ is locally a
morphism of $\seqname{C}^i$ to $\seqname{C}^{i+2}$.

If each $\catname{C}^i \subcat[full-] \catname{C}^{i+1}$ and each
$\funcseqname{f}^i_1$ is a morphism of $\catname{C}^i$, then
$\funcseqname{f}^{i+1}_1 \compose \funcseqname{f}^i_1$ is a
morphism of $\catname{C}^{i+2}$.
\end{proof}

The identity morphism of $(\seqname{A}^i, \seqname{E}^i, \seqname{C}^i)$
is an identity morphism.
\begin{proof}
The result follows from \pagecref{lem:atlcomp}.
\end{proof}
\end{lemma}

\begin{corollary}[M-atlas morphisms]
\label{cor:M-ATLmorph}
Let $\catname{C}^i$, $i=1,2,3$, be model categories,
$\seqname{C}^i \objin \catname{C}^i$, $E^i$ be topological spaces,
$\seqname{A}^i$ be m-atlases of $E^i$ in the coordinate spaces
$\seqname{C}^i$,
$\seqname{A}^2$ semi-maximal and
$
  \funcseqname{f}^i \defeq
    (
      \funcname{f}^i_0 \maps E^i \to E^{i+1},
      \funcname{f}^i_1 \maps \seqname{C}^i \to \seqname{C}^{i+1}
    \ )
$
be (strict) $E^i$-$E^{i+1}$ m-atlas morphisms of
$\seqname{A}^i$ to $\seqname{A}^{i+1}$ in the coordinate spaces
$\seqname{C}^i$, $\seqname{C}^{i+1}$. Then
$
  \funcseqname{f}^{i+1} \compose[()] \funcseqname{f}^i =
    (
      \funcname{f}^2_0 \compose \funcname{f}^1_0,
      \funcname{f}^2_1 \compose \funcname{f}^1_1
    )
$
is a (strict) $E^1$-$E^3$ m-atlas morphisms of
$\seqname{A}^1$ to $\seqname{A}^3$ in the coordinate spaces
$\seqname{C}^1$, $\seqname{C}^3$.

\begin{proof}
Since
$
    (
      \funcname{f}^i_0 \maps E^i \to E^{i+1},
      \funcname{f}^i_1 \maps \seqname{C}^i \to \seqname{C}^{i+1}
    )
$
is a (strict) $\Triv{E}^i$-$\Triv{E}^{i+1}$ m-atlas morphisms of
$\seqname{A}^i$ to $\seqname{A}^{i+1}$ in the coordinate spaces
$\seqname{C}^i$, $\seqname{C}^{i+1}$, then
$
    (
      \funcname{f}^2_0 \compose \funcname{f}^1_0,
      \funcname{f}^2_1 \compose \funcname{f}^1_1
    )
$
is a (strict) $\Triv{E}^1$-$\Triv{E}^3$ m-atlas morphisms of
$\seqname{A}^1$ to $\seqname{A}^3$ in the coordinate spaces
$\seqname{C}^1$, $\seqname{C}^3$.
\end{proof}

Let $\catname{E}^i$, $\catname{C}^i$, $i=1,2,3$ be model categories,
$\seqname{E}^i \in \Ob(\catname{E}^i)$,
$\seqname{C}^i \in \Ob(\catname{C}^i)$ $\seqname{A}^i$ m-atlases
of $\seqname{E}^i$ in the coordinate spaces $\seqname{C}^i$,
$\seqname{A}^2$ maximal and
$
  \funcseqname{f}^i \defeq \\
  (
    \funcname{f}^i_0 \maps \seqname{E}^i \to \seqname{E}^{i+1},
    \funcname{f}^i_1 \maps \seqname{C}^i \to \seqname{C}^{i+1}
  )
$
(strict) $\seqname{E}^i$-$\seqname{E}^{i+1}$ m-atlas morphisms of
$\seqname{A}^i$ to $\seqname{A}^{i+1}$ in the coordinate spaces
$\seqname{C}^i$, $\seqname{C}^{i+1}$. Then
$
  \funcseqname{f}^{i+1} \compose[()] \funcseqname{f}^i =
  (
    \funcname{f}^2_0 \compose \funcname{f}^1_0,
    \funcname{f}^2_1 \compose \funcname{f}^1_1
  )
$
is a (strict) $\seqname{E}^1$-$\seqname{E}^3$ m-atlas morphism of
$\seqname{A}^1$ to $\seqname{A}^3$ in the coordinate spaces
$\seqname{C}^1$, $\seqname{C}^3$.

\begin{proof}
Since $\seqname{A}^2$ is maximal it is semi-maximal.
\end{proof}

Let $\catname{C}^i$, $i=1,2,3$, be model categories,
$\seqname{C}^i \objin \catname{C}^i$, $E^i$ be topological spaces,
$\seqname{A}^i$ be m-atlases of $E^i$ in the coordinate spaces
$\seqname{C}^i$,
$\seqname{A}^2$ maximal and
$
  \funcseqname{f}^i \defeq \\
  (
    \funcname{f}^i_0 \maps E^i \to E^{i+1},
    \funcname{f}^i_1 \maps \seqname{C}^i \to \seqname{C}^{i+1}
  )
$
be (strict) $E^i$-$E^{i+1}$ m-atlas morphisms of
$\seqname{A}^i$ to $\seqname{A}^{i+1}$ in the coordinate spaces
$\seqname{C}^i$, $\seqname{C}^{i+1}$. Then
$
  \funcseqname{f}^{i+1} \compose[()] \funcseqname{f}^i =
  (
    \funcname{f}^2_0 \compose \funcname{f}^1_0,
    \funcname{f}^2_1 \compose \funcname{f}^1_1
  )
$
is a (strict) $E^1$-$E^3$ m-atlas morphisms of
$\seqname{A}^1$ to $\seqname{A}^3$ in the coordinate spaces
$\seqname{C}^1$, $\seqname{C}^3$.

\begin{proof}
Since $\seqname{A}^2$ is maximal it is semi-maximal.
\end{proof}
\end{corollary}

\begin{definition}[Sets of m-atlas morphisms]
Let $\seqname{E}^i$ and $\seqname{C}^i$. $i=1,2$, be model spaces.
Then
\begin{multline}
\Atl_\Ar(\seqname{E}^1, \seqname{C}^1, \seqname{E}^2, \seqname{C}^2)
\defeq \\
  \set
  {{
    \bigl (
      (
        \funcname{f}_0,
        \funcname{f}_1
      ),
      (\seqname{A}^1, \seqname{E}^1, \seqname{C}^1),
      (\seqname{A}^2, \seqname{E}^2, \seqname{C}^2)
    \bigr )
  }}%
  [
    {
      \isAtl_\Ar
      (
        \seqname{A}^1, \seqname{E}^1, \seqname{C}^1,
        \seqname{A}^2, \seqname{E}^2, \seqname{C}^2,
        \funcname{f}_0, \funcname{f}_1
      )
    }
  ]
\end{multline}

Let $E^i$, $i=1,2$ be topological spaces and
$\seqname{C}^1=(C^i, \catname{C}^i)$, $i=1,2$ be model spaces.
Then
\begin{multline}
\Atl_\Ar(E^1, \seqname{C}^1, E^2, \seqname{C}^2) \defeq \\
  \set
  {{
    \bigl (
      (
        \funcname{f}_0,
        \funcname{f}_1
      ),
      (\seqname{A}^1, \seqname{E}^1, \seqname{C}^1),
      (\seqname{A}^2, \seqname{E}^2, \seqname{C}^2)
    \bigr )
  }}%
  [
    {
      \isAtl_\Ar
      (
        \seqname{A}^1, E^1, \seqname{C}^1,
        \seqname{A}^2, E^2, \seqname{C}^2,
        \funcname{f}_0, \funcname{f}_1
      )
    }
  ]
\end{multline}
\end{definition}

\begin{definition}[Category $\Atl(\seqname{E}, \seqname{C})$]
Let $\seqname{E}$ and $\seqname{C}$ be sets of model spaces.
Let $P \defeq \seqname{E} \times \seqname{C}$.
Then
\begin{equation}
  \Atl_\Ob(\seqname{E}, \seqname{C})
  \defeq
  \union
    [
      {(\seqname{E}^\mu, \seqname{C}^\mu) \in \seqname{P}},
    ]
    {{
      \Atl_\Ob(\seqname{E}^\mu, \seqname{C}^\mu)
    }}
\end{equation}
\begin{equation}
  \Atl_\Ar(\seqname{E}, \seqname{C})
  \defeq
  \union
    [
      {(\seqname{E}^\mu, \seqname{C}^\mu) \in \seqname{P}},
      {(\seqname{E}^\nu, \seqname{C}^\nu) \in \seqname{P}}
    ]
    {{
      \Atl_\Ar
      (
        \seqname{E}^\mu,
        \seqname{C}^\mu,
        \seqname{E}^\nu,
        \seqname{C}^\nu
      )
    }}
\end{equation}
\begin{equation}
  \Atl(\seqname{E}, \seqname{C})
  \defeq
  \bigl (
    \Atl_\Ob(\seqname{E}, \seqname{C}),
    \Atl_\Ar(\seqname{E}, \seqname{C},
    \compose[A]
  \bigr )
\end{equation}
\end{definition}

\begin{lemma}[$\Atl(\seqname{E}, \seqname{C})$ is a category]
\label{lem:ATLiscat}
Let $\seqname{E}$ and $\seqname{C}$ be sets of model spaces. Then
$\Atl(\seqname{E}, \seqname{C})$ is a category and the identity
morphism for an object $(\seqname{A}^i, \seqname{E}^i, \seqname{C}^i)$
of $\Atl_\Ob(\seqname{E}, \seqname{C})$ is
$\Id_{(\seqname{A}^i, \seqname{E}^i, \seqname{C}^i}$.

\begin{proof}
Let $(\seqname{A}^i,\seqname{E}^i,\seqname{C}^i)$, $i=1,2,3$
be objects of $\Atl(\seqname{E}, \seqname{C})$ and
let \\
$
  \funcname{m}^i \defeq
  \bigl (
    (\funcname{f}_0^i,\funcname{f}_1^i),
    (\seqname{A}^i,\seqname{E}^i,\seqname{C}^i),
    (\seqname{A}^{i+1},\seqname{E}^{i+1},\seqname{C}^{i+1})
  \bigr )
$
be morphisms of $\Atl(\seqname{E}, \seqname{C})$. Then
\begin{enumerate}
\item Composition: \newline
$
  \bigl (
    \funcseqname{f}^2 \compose[()] \funcseqname{f}^1 =
    (
      \funcname{f}_0^2 \compose \funcname{f}_0^1,
      \funcname{f}_1^2 \compose \funcname{f}_1^1
    ),
    (\seqname{A}^1,\seqname{E}^1,\seqname{C}^1),
    (\seqname{A}^3,\seqname{E}^3,\seqname{C}^3)
  \bigr )
$
is a morphism of $\Atl(\seqname{E}, \seqname{C})$ by \pagecref{lem:M-ATLmorph}.
\item Associativity: \newline
Composition is associative by \pagecref{lem:atlcomp}.
\item Identity: \newline
$\Id_{(\seqname{A}^i, \seqname{E}^i, \seqname{C}^i)}$ is an identity
morphism by \cref{lem:atlcomp}.
\end{enumerate}
\end{proof}
\end{lemma}

\begin{definition}[Category $\Atl(\catname{E}, \catname{C})$]
Let $\catname{E}$ and $\catname{C}$ be model categories.
Let $\seqname{P} \defeq Ob(\catname{E}) \times \Ob(\catname{C})$.
Then
\begin{equation}
  \Atl_\Ob(\catname{E}, \catname{C})
  \defeq
  \union
    [
        {\seqname{E} \objin \catname{E}},
        {\seqname{C} \objin \catname{C}}
    ]
    {{
      \Atl_\Ob(\seqname{E}, \seqname{C})
    }}
\end{equation}
\begin{multline}
\Atl_\Ar(\catname{E}, \catname{C})
\defeq
  \set
  {{
    \bigl (
      (
        \funcname{f}_0,
        \funcname{f}_1,
      ),
      (\seqname{A}^1, \seqname{E}^1, \seqname{C}^1),
      (\seqname{A}^2, \seqname{E}^2, \seqname{C}^2)
    \bigr )
  }}%
  [
    {
      (\seqname{E}^i, \seqname{C}^i) \in \seqname{P}
    },
    {
      \strict{\isAtl_\Ar}
      (
        \seqname{A}^1,
        \seqname{E}^1,
        \seqname{C}^1,
        \seqname{A}^2,
        \seqname{E}^2,
        \seqname{C}^2,
        \funcname{f}_0,
        \funcname{f}_1,
        \catname{E},
        \catname{C}
      )
    }
  ]*
\end{multline}
\begin{equation}
  \Atl(\catname{E}, \catname{C})
  \defeq
  \bigl (
    \Atl_\Ob(\catname{E}, \catname{C}),
    \Atl_\Ar(\catname{E}, \catname{C},
    \compose[A]
  \bigr )
\end{equation}

Let $(\seqname{A}^1, \seqname{E}^1, \seqname{C}^1) \in \Atl_\Ob(\catname{E}, \catname{C})$.
\begin{equation}
\Id_{(\seqname{A}^1, \seqname{E}^1, \seqname{C}^1)} \defeq
\bigl (
  (\Id_{\seqname{E}^1}, \Id_{\seqname{C}^1}),
  (\seqname{A}^1, \seqname{E}^1, \seqname{C}^1),
  (\seqname{A}^1, \seqname{E}^1, \seqname{C}^1)
\bigr )
\end{equation}
\end{definition}

\begin{lemma}[$\Atl(\catname{E}, \catname{C})$ is a category]
Let $\catname{E}$ and $\catname{C}$ be model categories. Then
$\Atl(\catname{E}, \catname{C})$ is a category and the identity morphism
for an object $(\seqname{A}^i, \seqname{E}^i, \seqname{C}^i)$ of
$\Atl_\Ob(\catname{E}, \catname{C})$ is
$\Id_{(\seqname{A}^i, \seqname{E}^i, \seqname{C}^i}$.

\begin{proof}
Let $(\seqname{A}^i,\seqname{E}^i,\seqname{C}^i)$, $i=1,2,3$
be objects of $\Atl(\catname{E}, \catname{C})$ and
let \\
$
  \funcname{m}^i \defeq
  \bigl (
    (\funcname{f}_0^i,\funcname{f}_1^i),
    (\seqname{A}^i,\seqname{E}^i,\seqname{C}^i),
    (\seqname{A}^{i+1},\seqname{E}^{i+1},\seqname{C}^{i+1})
  \bigr )
$
be morphisms of $\Atl(\catname{E}, \catname{C})$. Then
\begin{enumerate}
\item Composition: \newline
$\funcname{f}_0^2 \compose \funcname{f}_0^1$ is a morphism of $\catname{E}$,
$\funcname{f}_1^2 \compose \funcname{f}_1^1$ is a morphism of $\catname{C}$ and
$
  \funcseqname{f}^2 \compose[()] \funcseqname{f}^1 =
  (
    \funcname{f}_0^2 \compose \funcname{f}_0^1,
    \funcname{f}_1^2 \compose \funcname{f}_1^1
  )
$
is an $E^1$-$E^2$ m-atlas morphism of $\seqname{A}^1$ to $\seqname{A}^2$
in the coordinate spaces $\seqname{C}^1$, $\seqname{C}^2$ by
\pagecref{lem:M-ATLmorph}.
\item Associativity: \newline
Composition is associative by \pagecref{lem:atlcomp}.
\item Identity: \newline
$\Id_{(\seqname{A}^i, \seqname{E}^i, \seqname{C}^i)}$ is an identity
morphism by \cref{lem:atlcomp}.
\end{enumerate}
\end{proof}
\end{lemma}

\begin{definition}[Category $\Atl(\seqname{E}, \seqname{C})$ of topological spaces]
Let $\seqname{E}$ be a set of topological spaces and $\seqname{C}$ a
set of model spaces.
Let $\seqname{P} \defeq \seqname{E} \times \seqname{C}$.
Then
\begin{equation}
  \Atl_\Ob(\seqname{E}, \seqname{C})
  \defeq
  \union
    [
        {E^1 \in \seqname{E}},
        {\seqname{C}^1 \in \seqname{C}}
    ]
    {{
      \Atl_\Ob(E^1, \seqname{C}^1)
    }}
\end{equation}
\begin{equation}
  \Atl_\Ar(\seqname{E}, \seqname{C})
  \defeq
  \union
    [
        {(E^1, \seqname{C}^1) \in \seqname{P}},
        {(E^2, \seqname{C}^2) \in \seqname{P}}
    ]
    {{
      \Atl_{\Ar}
      (
        E^1,
        \seqname{C}^1,
        E^3,
        \seqname{C}^3
      )
    }}
\end{equation}
\begin{equation}
  \Atl(\seqname{E}, \seqname{C})
  \defeq
  \bigl (
    \Atl_\Ob(\seqname{E}, \seqname{C}),
    \Atl_\Ar(\seqname{E}, \seqname{C}),
    \compose[A]
  \bigr )
\end{equation}

Let $(\seqname{A}^1, \seqname{E}^1, \seqname{C}^1) \in \Atl_\Ob(\seqname{E}, \seqname{C})$.
\begin{equation}
\Id_{(\seqname{A}^1, \seqname{E}^1, \seqname{C}^1)} \defeq
\bigl (
  (\Id_{\seqname{E}^1}, \Id_{\seqname{C}^1}),
  (\seqname{A}^1, \seqname{E}^1, \seqname{C}^1),
  (\seqname{A}^1, \seqname{E}^1, \seqname{C}^1)
\bigr )
\end{equation}
\end{definition}

\begin{lemma}[$\Atl(\seqname{E}, \seqname{C})$ of spaces is a category]
Let $\seqname{E}$ be a set of topological spaces and $\seqname{C}$ a
set of model spaces. Then $\Atl(\seqname{E}, \seqname{C})$ is a category
and the identity morphism for an object
$(\seqname{A}^1, E^1, \seqname{C}^1)$ is
$\Id_{(\seqname{A}^1, E^1, \seqname{C}^1)}$.

\begin{proof}
Let $(\seqname{A}^i,E^i,\seqname{C}^i)$, $i=1,2,3$
be objects of $\Atl(E, \seqname{C})$. and let
$
  \funcname{m}^i \defeq \\
  \bigl (
    (\funcname{f}_0^i,\funcname{f}_1^i),
    (\seqname{A}^i,E^i,\seqname{C}^i),
    (\seqname{A}^{i+1},E^{i+1},\seqname{C}^{i+1})
  \bigr )
$
be morphisms of $\Atl(\seqname{E}, \seqname{C})$. Then
\begin{enumerate}
\item Composition: \newline
$
  \bigl (
    (
      \funcname{f}_0^2 \compose \funcname{f}_0^1,
      \funcname{f}_1^2 \compose \funcname{f}_1^1
    ),
    (\seqname{A}^1,E^1,\seqname{C}^1),
    (\seqname{A}^3,E^3,\seqname{C}^3)
  \bigr )
$
is a morphism of $\Atl(\seqname{E}, \seqname{C})$ by \pagecref{cor:M-ATLmorph}.
\item Associativity: \newline
Composition is associative by \pagecref{lem:atlcomp}.
\item Identity: \newline
$\Id_{(\seqname{A}^i, \seqname{E}^i, \seqname{C}^i)}$ is an identity
morphism by \cref{lem:atlcomp}.
\end{enumerate}
\end{proof}
\end{lemma}

\begin{definition}
Let $\seqname{E}$ be a set of topological spaces and $\seqname{C}$ a
set of model spaces. Then
\begin{equation}
\Triv{\seqname{E}} \defeq \set {\Triv{E^\mu}}[E^\mu \in \seqname{E}]
\end{equation}
\begin{equation}
\Triv{\Atl}(\seqname{E}, \seqname{C}) \defeq
\Atl \Big ( \Triv{\seqname{E}}, \seqname{C} \Big )
\end{equation}
\end{definition}

\begin{definition}[$\Functor^\Top$ on objects]
\label{def:FtopOb}
Let $E$ be a topological space, $\seqname{C}$ a model space
and $\seqname{A}$ an m-atlas of $E$ in the coordinate space $\seqname{C}$.
Then
\begin{equation}
\Functor^\Top (\seqname{A}, E, \seqname{C})
\defeq
\Big ( \seqname{A}, \Triv{E}, \seqname{C} \Big )
\end{equation}
\end{definition}

\begin{definition}[$\Functor^\Top$ on morphisms]
\label{def:FtopAr}
Let $\catname{C}^i$, $i=1,2$, be model categories,
$\seqname{C}^i \objin \catname{C}^i$, $E^i$ topological spaces,
$\seqname{A}^i$ m-atlases of $E^i$ in the coordinate spaces
$\seqname{C}^i$, $(\funcname{f}_0 \maps E^1 \to E^2$ and
$\funcname{f}_1 \maps \seqname{C}^1 \to \seqname{C}^2)$. Then
\begin{multline}
\Functor^\Top
  \bigl (
    (
      \funcname{f}_0 \maps E^1 \to E^2,
      \funcname{f}_1 \maps \seqname{C}^1 \to \seqname{C}^2,
    ),
    (\seqname{A}^1, E^1, \seqname{C}^1),
    (\seqname{A}^2, E^2, \seqname{C}^2)
  \bigr )
\defeq \\
  \biggl (
    \Bigl (
      \funcname{f}_0 \maps \Triv{E^1} \to \Triv{E^2},
      \funcname{f}_1 \maps \seqname{C}^1 \to \seqname{C}^2,
     \Bigr ) ,
    \Bigl ( \seqname{A}^1, \Triv{E^1}, \seqname{C}^1 \Bigr ) ,
    \Bigl ( \seqname{A}^2, \Triv{E^2}, \seqname{C}^2 \Bigr )
  \biggr )
\end{multline}
\end{definition}

\begin{theorem}[$\Functor^\Top$ is a functor]
\label{the:Ftop}
Let $\seqname{E}$ be a set of topological spaces and $\seqname{C}$ a
set of model spaces.
$\Functor^\Top$ is a functor from
$\Atl(\seqname{E}, \seqname{C})$
to
$\Atl(\Triv{\seqname{E}}, \seqname{C})$

\begin{proof}
Let $\seqname{C}^1=(\seqname{A}^1, E^1, \seqname{C}^1)$
be an object of $\Atl(\seqname{E}, \seqname{C})$ and \\
$
  \bigl (
    (\Id_{E^1}, \Id_{\seqname{C}^1}),
    (\seqname{A}^1, E^1, \seqname{C}^1),
    (\seqname{A}^1, E^1, \seqname{C}^1)
  \bigr )
$
be an identity morphism of \\
$(\seqname{A}^1, E^1, \seqname{C}^1)$.
Then
$
  \bigl (
    (\Id_{\Triv{E}^1}, \Id_{\seqname{C}^1}),
    (\seqname{A}^1, \Triv{E}^1, \seqname{C}^1),
    (\seqname{A}^1, \Triv{E}^1, \seqname{C}^1)
  \bigr )
$
is an identity morphism of
$(\seqname{A}^1, \Triv{E}^1, \seqname{C}^1)$.

Let
$
  m^1 \defeq
  \bigl (
    (f_0^1, f_1^1),
    (\seqname{A}^1, E^1, \seqname{C}^1),
    (\seqname{A}^2, E^2, \seqname{C}^2)
  \bigr )
$
and \\
$
  m_2 \defeq
  \bigl (
    (f_0^2, f_1^2),
    (\seqname{A}^2, E^2, \seqname{C}^2),
    (\seqname{A}^3, E^3, \seqname{C}^3)
  \bigr )
$
be morphisms of $\Atl(\seqname{E}, \seqname{C})$.
\begin{equation}
\begin{split}
  &
  \Functor^\Top (m2 \compose[A] m^1) \defeq
\\*
  &
  \Functor^\Top
    \bigl (
      (f_0^2 \compose f_0^1, f_1^2 \compose f_1^1),
      (\seqname{A}^1, E^1, \seqname{C}^1),
      (\seqname{A}^3, E^3, \seqname{C}^3)
    \bigr ) =
\\*
  &
    \bigl (
      (f_0^2 \compose f_0^1, f_1^2 \compose f_1^1),
      (\seqname{A}^1, \Triv{E^1}, \seqname{C}^1),
      (\seqname{A}^3, \Triv{E^3}, \seqname{C}^3)
    \bigr ) =
\\*
  &
    \bigl (
      (f_0^2, f_1^2),
      (\seqname{A}^1, \Triv{E^1}, \seqname{C}^1),
      (\seqname{A}^2, \Triv{E^2}, \seqname{C}^2)
    \bigr )
    \compose[A]
    \bigl (
      (f_0^1, f_1^1),
      (\seqname{A}^2, \Triv{E^2}, \seqname{C}^2),
      (\seqname{A}^3, \Triv{E^3}, \seqname{C}^3)
    \bigr ) =
\\*
  &
  \Functor^\Top
    \bigl (
      (f_0^2, f_1^2),
      (\seqname{A}^2, E^2, \seqname{C}^2),
      (\seqname{A}^3, E^3, \seqname{C}^3)
    \bigr )
  \compose[A]
\\*
  &
    \Functor^\Top
    \bigl (
      (f_0^1, f_1^1),
      (\seqname{A}^1, E^1, \seqname{C}^1),
      (\seqname{A}^2, E^2, \seqname{C}^2)
    \bigr ) =
\\*
  &
  \Functor^\Top (m2)
  \compose[A]
  \Functor^\Top (m1)
  \end{split}
\end{equation}
\end{proof}
\end{theorem}

\section{Local Coordinate Spaces}
\label{sec:lcs}

This section of the paper defines local coordinate spaces, morphisms
among them and categories of them.

\begin{definition}[Local $*-\Sigma$ coordinate spaces]
\label{def:LCS}
Let $\catseqname{M} \defeq (\catname{M}_\alpha, \alpha \prec \Alpha)$ be
a sequence of categories,
$\seqname{M} \defeq (M_\alpha, \alpha \prec \Alpha)$ be a
sequence of spaces and
$\funcseqname{F} \defeq (\funcname{F}_\gamma, \gamma \prec \Gamma)$ a
sequence of functions.
$
  \seqname{L} \defeq
  (\catseqname{M}, \seqname{M}, \seqname{A}, \funcseqname{F}, \Sigma)
$
is a local $\seqname{M}-\Sigma$ coordinate space, abbreviated
$\isLCS_\Ob \seqname{L}$, a local $\catseqname{M}-\Sigma$ coordinate
space, abbreviated $\isLCS_\Ob \seqname{L}$, and a local
$\Alpha-\Sigma$ coordinate space, abbreviated
$\isLCS_\Ob(\seqname{L},\Alpha)$, iff

\begin{enumerate}
\item $\catname{M}_1$ is a model category
\item $(\catseqname{M}, \seqname{M}, \Sigma, \funcseqname{F})$ is a
$\catseqname{M}$-$\Sigma$ prestructure
\item $\seqname{A}$ is a maximal m-atlas of $M_0$ in $M_1$.
\item $\forall \gamma \prec \Gamma$, if
$\funcname{F}_\gamma \maps (M_{\sigma_{\gamma,\alpha_\beta}}, \beta
\prec \Beta_\gamma) \to \truthset$ is a constraint function then \\
$\funcname{F}_\gamma
  [
    \bigtimes_{\beta \prec \Beta_\gamma}
    M_{\sigma_{\gamma,\alpha_\beta}}
  ]
= \{ \true \}$, i.e., \\
$\uquant
  {
    {
      (
        s_\beta \in M_{\sigma_{\gamma,\alpha_\beta}},
        \beta \prec \Beta_\gamma
      )
    }
  }
  {{\funcname{F}_\gamma(s_\beta, \beta \prec \Beta_\gamma) = \true}}
  $.
\end{enumerate}

Define the total space $(E, \catname{E})$ to be $M_0$, the Coordinate
space $(C,\catname{C})$ to be $M_1$, the adjunct spaces, if any, to be
$M_\alpha$, $\alpha >1$, and the adjunct functions if any, to be
$F_\alpha$, $\alpha \prec \Alpha$.
\end{definition}

\begin{remark}
There are alternative approaches that are beyond the scope of this
paper, e.g.,
\begin{enumerate}
\item Requiring the atlas to be full would eliminate certain pathologies.
\item A more complex definition would more easily accommodate structures
with more than one atlas, e.g.,
\begin{enumerate}
\item A fiber bundle with a differential structure
\item Multiple fiber bundles on the same base space.
\end{enumerate}
\end{enumerate}
\end{remark}

\begin{lemma}[Local $*-\Sigma$ coordinate spaces]
\label{lem:LCS}
let $\catseqname{M}^i$, $i=1,2$, be a sequence of categories,
$\catseqname{M}^2_j$, $j=0,1$, model categories,
$\catseqname{M}^1 \SUBCAT[full-] \catseqname{M}^2$ and
$
  \seqname{L}^1 \defeq \\
  (\catseqname{M}^1, \seqname{M}, \seqname{A}, \funcseqname{F}, \Sigma)
$
a local $\catseqname{M}^1-\Sigma$ coordinate space, Then
$
  \seqname{L}^2 \defeq
  (\catseqname{M}^2, \seqname{M}, \seqname{A}, \funcseqname{F}, \Sigma)
$
is a local $\catseqname{M}^2-\Sigma$ coordinate space,

\begin{proof}
$\seqname{L}^2$ satisfies the conditions of \fullcref{def:LCS}:
\begin{enumerate}
\item $\catname{M}^2_0$ and $\catname{M}^2_1$ are model categories by
hypothesis.
\item $(\catseqname{M}^2, \seqname{M}, \Sigma, \funcseqname{F})$ is a
$\catseqname{M}^2$-$\Sigma$ prestructure by \pagecref{lem:pre}.
\item $\seqname{A}$ is a maximal m-atlas of $M_0$ in $M_1$ by hypothesis.
\item $\forall \gamma \prec \Gamma$, if
$\funcname{F}_\gamma \maps (M_{\sigma_{\gamma,\alpha_\beta}}, \beta
\prec \Beta_\gamma) \to \truthset$ is a constraint function then \\
$
  \funcname{F}_\gamma
    [
      \bigtimes_{\beta \prec \Beta_\gamma}
      M_{\sigma_{\gamma,\alpha_\beta}}
    ]
  = \{ \true \}
$,
i.e., \\
$
  \uquant
  {
    {
      (
        s_\beta \in M_{\sigma_{\gamma,\alpha_\beta}},
        \beta \prec \Beta_\gamma
      )
    }
  }
  {{\funcname{F}_\gamma(s_\beta, \beta \prec \Beta_\gamma) = \true}}
$.
\end{enumerate}
\end{proof}
\end{lemma}

\begin{definition}[$\LCS_\Ob$]
\begin{equation}
  \LCS_\Ob(\catseqname{M}, \Sigma) \defeq
  \set
    {{
      \seqname{L} \defeq
      (\catseqname{M}, \seqname{M}, \seqname{A}, \funcseqname{F}, \Sigma)
    }}%
    [
      {{
        \seqname{M} \seqin \catseqname{M} \land \isLCS_\Ob \; \seqname{L}
      }}
    ]
\end{equation}
\end{definition}

A variety of conditions can be imposed on the transition functions
$\funcname{t}_\beta^\alpha=\phi_\beta \compose \phi_\alpha^{-1}$ by
appropriate choice of category.

\begin{remark}
While some of cases in \pagecref{Examples} require constraint functions,
the specific examples worked out in \cref{sec:man,sec:bun} below
do not.
\end{remark}

\subsection {Morphisms of local coordinate spaces}
\label{sub:lcsmorph}
Informally, a morphism between two local coordinate spaces is a sequence
of functions that is compatible with the functions and m-atlases of the
two local coordinate spaces. The definition implicitly uses commutation
relations, which allows a variety of properties to fall out
automatically, e.g., preservation of fibers, being a homomorphism.

\begin{definition}[morphisms of local coordinate spaces]
\label{def:LCSmorph}
Let
$\catseqname{M}^i \defeq (\catname{M}^i_\alpha, \alpha \prec \Alpha)$,
$i=1,2$, be a sequence of categories,
$
  \seqname{M}^i \defeq (\seqname{M}^i_\alpha, \alpha \prec \Alpha)
  \seqin \catseqname{M}^i
$
a sequence of spaces,
$
  \seqname{L}^i \defeq
  (
    \catseqname{M}^i,
    \seqname{M}^i,
    \seqname{A}^i,
    \funcseqname{F}^i,
    \Sigma
  )
$,
a local $\seqname{M}^i$-$\Sigma$ coordinate space and
$
  \funcseqname{f} \defeq \\
    (
      \funcname{f}_\alpha \maps
        \seqname{M}^1_\alpha \to \seqname{M}^2_\alpha,
      \alpha \prec \Alpha
    )
$.

$\funcseqname{f}$ is a morphism from $\seqname{L}^1$ to $\seqname{L}^2$,
abbreviated \\
$
  \isLCS_\Ar
    (
      \catseqname{M}^1,
      \seqname{M}^1,
      \seqname{A}^1,
      \funcseqname{F}^1,
      \Sigma,
      \catseqname{M}^2,
      \seqname{M}^2,
      \seqname{A}^2,
      \funcseqname{F}^2,
      \funcseqname{f}
    )
$,
iff

\begin{enumerate}
\item $\funcseqname{f}$ is a morphism from the prestructure
$\seqname{P}^1 \defeq (\catseqname{M}^1, \seqname{M}^1, \Sigma, \funcseqname{F}^1)$
to the prestructure
$\seqname{P}^2 \defeq (\catseqname{M}^2, \seqname{M}^2, \Sigma, \funcseqname{F}^2)$.
\item $(\funcname{f}_0, \funcname{f}_1)$ is a
$\catname{M}^1_0$-$\catname{M}^2_0$ m-atlas morphism from
$\seqname{A}^1$ to $\seqname{A}^2$ in the coordinate model categories
$\catname{M}^1_1$-$\catname{M}^2_1$.
\end{enumerate}

$\funcseqname{f}$ is a semi-strict morphism from $\seqname{L}^1$ to
$\seqname{L}^2$, abbreviated \\
$
  \strict[semi-]{\isLCS_\Ar}
    (
      \catseqname{M}^1,
      \seqname{M}^1,
      \seqname{A}^1,
      \funcseqname{F}^1,
      \Sigma,
      \catseqname{M}^2,
      \seqname{M}^2,
      \seqname{A}^2,
      \funcseqname{F}^2,
      \funcseqname{f}
    )
$,
iff

\begin{enumerate}
\item $\funcseqname{f}$ is a semi-strict morphism from the prestructure
$\seqname{P}^1 \defeq (\catseqname{M}^1, \seqname{M}^1, \Sigma, \funcseqname{F}^1)$
to the prestructure
$\seqname{P}^2 \defeq (\catseqname{M}^2, \seqname{M}^2, \Sigma, \funcseqname{F}^2)$.
\item $(\funcname{f}_0, \funcname{f}_1)$ is a semi-strict
$\catname{M}^1_0$-$\catname{M}^2_0$ m-atlas morphism from
$\seqname{A}^1$ to $\seqname{A}^2$ in the coordinate model categories
$\catname{M}^1_1$-$\catname{M}^2_1$.
\end{enumerate}

$\funcseqname{f}$ is a strict morphism from
$
  (
    \catseqname{M}^1,
    \seqname{M}^1,
    \seqname{A}^1,
    \funcseqname{F}^1,
    \Sigma
  )
$ to
$
  (
    \catseqname{M}^2,
    \seqname{M}^2,
    \seqname{A}^2,
    \funcseqname{F}^2,
    \Sigma
  )
$,
abbreviated
$
  \strict{\isLCS_\Ar}
  (
    \catseqname{M}^1,
    \seqname{M}^1,
    \seqname{A}^1,
    \funcseqname{F}^1,
    \Sigma,
    \catseqname{M}^2,
    \seqname{M}^2,
    \seqname{A}^2,
    \funcseqname{F}^2,
    \funcseqname{f}
  )
$,
iff

\begin{enumerate}
\item $\funcseqname{f}$ is a strict morphism from the prestructure
$\seqname{P}^1 \defeq (\catseqname{M}^1, \seqname{M}^1, \Sigma, \funcseqname{F}^1)$
to the prestructure
$\seqname{P}^2 \defeq (\catseqname{M}^2, \seqname{M}^2, \Sigma, \funcseqname{F}^2)$.
\item $(\funcname{f}_0, \funcname{f}_1)$ is a strict
$\catname{M}^1_0$-$\catname{M}^2_0$-$\seqname{M}^1_0$-$\seqname{M}^2_0$%
-$\catname{M}^1_1$-$\catname{M}^2_1$
m-atlas morphism from $\seqname{A}^1$ to $\seqname{A}^2$ in the
cordinate spaces $\seqname{M}^1_1$, $\seqname{M}^2_1$.
\end{enumerate}

\end{definition}

\begin{lemma}[morphisms of local coordinate spaces]
\label{lem:LCSmorph}
Let
$\catseqname{M}^i \defeq (\catname{M}^i_\alpha, \alpha \prec \Alpha)$,
$i\in [1,4]$, be a sequence of categories,
$
  \seqname{M}^i \defeq (\seqname{M}^i_\alpha, \alpha \prec \Alpha)
  \seqin \catseqname{M}^i
$
a sequence of spaces,
$
  \seqname{P}^i \defeq
  (
    \catseqname{M}^i,
    \seqname{M}^i,
    \funcseqname{F}^i,
    \Sigma
  )
$,
$
  \seqname{L}^i \defeq
  (
    \catseqname{M}^i,
    \seqname{M}^i,
    \seqname{A}^i,
    \funcseqname{F}^i,
    \Sigma
  )
$
a local $\seqname{M}^i$-$\Sigma$ coordinate space
and
$
  \funcseqname{f}^i \defeq
    (
      \funcname{f}^i_\alpha \maps
        \seqname{M}^i_\alpha \to \seqname{M}^{i+1}_\alpha,
      \alpha \prec \Alpha
    )
$,
$i=1,2,3$,
a (strict) morphism from $\seqname{L}^i$ to $\seqname{L}^{i+1}$.

If $\catseqname{M}^i \SUBCAT[full-] \catseqname{M}^{i+1}$,
$\seqname{S}^i \SUBSETEQ \seqname{S}^{i+1}$,
$\seqname{A}^i = \seqname{A}^{i+1}$ and each
$
  \funcname{F}^i_\gamma = \\
  \funcname{F}^{i+1}_\gamma \maps
    \bigtimes \underline{\overset{\sigma_\gamma}{\seqname{S}^i}} \to
    \tail \left ( \overset{\sigma_\gamma}{\seqname{S}^i} \right )
$
then
$\Id_{\seqname{M}^i,\seqname{M}^{i+1}}$ is a strict morphism from
$\seqname{L}^i$ to $\seqname{L}^{i+1}$.

\begin{proof}
$\Id_{\seqname{P}^i,\seqname{P}^{i+1}}$ is a strict morphism from
$\seqname{P}^i$ to $\seqname{P}^{i+1}$
by \pagecref{lem:premorph}

$
\ID_%
  {
    (\seqname{M}^i_0,\seqname{M}^i_1),
    (\seqname{M}^{i+1}_0,\seqname{M}^{i+1}_1)
  }
$
is a strict
$\catname{M}^i_0$-$\catname{M}^{i+1}_0$-%
$\seqname{M}^i_0$-$\seqname{M}^{i+1}_0$-%
$\catname{M}^i_1$-$\catname{M}^{i+1}_1$ m-atlas morphism of
$\seqname{A}^i$ to $\seqname{A}^{i+1}$ in the coordinate spaces
$\seqname{M}^i_1$, $\seqname{M}^{i+1}_1$.
by \pagecref{lem:M-ATLmorph}.
\end{proof}

Let
$\catseqname{M}'^i \defeq (\catname{M}'^i_\alpha, \alpha \prec \Alpha)$,
$i\in [1,4]$, be a sequence of categories and
$
  \seqname{L}'^i \defeq \\
  (
    \catseqname{M}'^i,
    \seqname{M}^i,
    \seqname{A}^i,
    \funcseqname{F}^i,
    \Sigma
  )
$
be a local coordinate space. Then $\funcseqname{f}^i$ is a morphism from
$\seqname{L}'^i$ to $\seqname{L}'^{i+1}$

\begin{proof}
The commutation relations do not depend on the categories.
\end{proof}

Let
$\catseqname{M}'^i \defeq (\catname{M}'^i_\alpha, \alpha \prec \Alpha)$,
$i\in [1,4]$, be a sequence of categories,
$\catname{M}'^i_0$ and $\catname{M}'^i_1$ be model categories,
$\catseqname{M}^i \SUBCAT[full-] \catseqname{M}'^i$,
$
  \seqname{P}'^i \defeq
  (
    \catseqname{M}'^i,
    \seqname{M}^i,
    \funcseqname{F}^i,
    \Sigma
  )
$ and
$
  \seqname{L}'^i \defeq \\
  (
    \catseqname{M}'^i,
    \seqname{M}^i,
    \seqname{A}^i,
    \funcseqname{F}^i,
    \Sigma
  )
$.
Then
\begin{enumerate}
\item $\seqname{L}'^i$ is a local $\seqname{M}^i$-$\Sigma$ coordinate
space and $\funcseqname{f}^i$ is a (strict) morphism from
$\seqname{L}'^i$ to $\seqname{L}'^{i+1}$.

\item If $\funcseqname{f}^i$ is a semi-strict morphism from
$\seqname{L}^i$ to $\seqname{L}^{i+1}$ then $\funcseqname{f}^i$ is
a semi-strict morphism from $\seqname{L}'^i$ to $\seqname{L}'^{i+1}$.

\item If $\funcseqname{f}^i$ is a strict morphism from $\seqname{L}^i$
to $\seqname{L}^{i+1}$ then $\funcseqname{f}^i$ is a strict
morphism from $\seqname{L}'^i$ to $\seqname{L}'^{i+1}$.
\end{enumerate}

\begin{proof}
Each $\seqname{L}'^i$ is a local $\seqname{M}^i$-$\Sigma$
coordinate space:

\begin{enumerate}
\item $\catseqname{M}'^i_0$ and $\catseqname{M}'^i_1$ are model
categories by hypothesis.

\item $\seqname{P}'^i$ is a $\catseqname{M}'^i$-$\Sigma$ prestructure
by \pagecref{lem:pre}.

\item $\seqname{A}^i$ is a maximal m-atlas of $\seqname{M}^i_0$ in
$\seqname{M}^i_1$ by hypothesis.

\item All constraint functions evaluate to $\true$ by hypothesis.
\end{enumerate}

Each $\funcseqname{f}^i$ is a (strict) morphism from $\seqname{L}'^i$ to
$\seqname{L}'^{i+1}$:

\begin{enumerate}
\item $\funcseqname{f}^i$ is a morphism from the prestructure
$\seqname{P}'^i$ to the prestructure $\seqname{P}'^{i+1}$ by
\pagecref{lem:premorph}.

If $\funcseqname{f}^i$ is a semi-strict (strict) morphism from the
prestructure $\seqname{P}^i$ to the prestructure $\seqname{P}^{i+1}$
then $\funcseqname{f}^i$ is a semi-strict (strict) morphism from the
prestructure $\seqname{P}'^i$ to the prestructure $\seqname{P}'^{i+1}$
by \cref{lem:premorph}.

\item $(\funcname{f}^i_0, \funcname{f}^i_1)$ is a
$\seqname{M}^i_0$-$\seqname{M}^{i+1}_0$ m-atlas morphism from
$\seqname{A}^i$ to $\seqname{A}^{i+1}$ in the coordinate model
spaces $\seqname{M}^i_i$, $\seqname{M}^{i+1}_i$ by hypothesis.

If $(\funcname{f}^i_0, \funcname{f}^i_1)$ is a semi-strict (strict)
$\catname{M}^i_0$-$\catname{M}^{i+1}_0$-%
$\seqname{M}^i_0$-$\seqname{M}^{i+1}_0$-
$\catname{M}^i_1$-$\catname{M}^{i+1}_1$ m-atlas morphism from
$\seqname{A}^i$ to $\seqname{A}^{i+1}$ in the coordinate model spaces
$\seqname{M}^i_i$, $\seqname{M}^{i+1}_i$ then
$(\funcname{f}^i_0, \funcname{f}^i_1)$ is a semi-strict (strict)
$\catname{M}'^i_0$-$\catname{M}'^{i+1}_0$-%
$\seqname{M}^i_0$-$\seqname{M}^{i+1}_0$-%
$\catname{M}'^i_1$-$\catname{M}'^{i+1}_1$ m-atlas morphism from
$\seqname{A}^i$ to $\seqname{A}^{i+1}$ in the coordinate model spaces
$\seqname{M}^i_i$, $\seqname{M}^{i+1}_i$ by \pagecref{lem:M-ATLmorph}.
\end{enumerate}
\end{proof}

$\funcseqname{f}^{i+1} \compose \funcseqname{f}^i$ is
a (strict) morphism from $\seqname{L}^i$ to $\seqname{L}^{i+2}$:

\begin{proof}

\begin{enumerate}
\item If each $\funcseqname{f}^j$ is a (strict) morphism from the
prestructure $\seqname{P}^j$ to the prestructure $\seqname{P}^{j+1}$
then $\funcseqname{f}^{i+1} \compose \funcseqname{f}^i$ is a (strict)
morphism from the prestructure $\seqname{P}^i$ to the prestructure
$\seqname{P}^{i+2}$ by \pagecref{lem:premorph}.

If each $\funcseqname{f}^j$ is a semi-strict (strict) morphism from the
prestructure $\seqname{P}^j$ to the prestructure $\seqname{P}^{j+1}$
then $\funcseqname{f}^{i+1} \compose \funcseqname{f}^i$ is a semi-strict
(strict) morphism from the prestructure $\seqname{P}^i$ to the
prestructure $\seqname{P}^{i+2}$ by \cref{lem:premorph}.

\item
$
  (
    \funcname{f}^{i+1}_0 \compose \funcname{f}^i_0,
    \funcname{f}^{i+1}_1 \compose \funcname{f}^i_1
  )
$
is a semi-strict (strict)
$\catname{M}^i_0$-$\catname{M}^{i+1}_0$-%
$\seqname{M}^i_0$-$\seqname{M}^{i+1}_0$-%
$\catname{M}^i_1$-$\catname{M}^{i+1}_1$ m-atlas morphism from
$\seqname{A}^i$ to $\seqname{A}^{i+2}$ in the coordinate model spaces
$\seqname{M}^i_i$, $\seqname{M}^{i+2}_i$ by \pagecref{lem:M-ATLmorph}.
If
$
  (
    \funcname{f}^{i+1}_0 \compose \funcname{f}^i_0,
    \funcname{f}^{i+1}_1 \compose \funcname{f}^i_1
  )
$
is a semi-strict (strict)
$\catname{M}^i_0$-$\catname{M}^{i+1}_0$-%
$\seqname{M}^i_0$-$\seqname{M}^{i+1}_0$-%
$\catname{M}^i_1$-$\catname{M}^{i+1}_1$
m-atlas morphism from $\seqname{A}^i$ to $\seqname{A}^{i+2}$ in the
coordinate model spaces $\seqname{M}^i_i$, $\seqname{M}^{i+2}_i$ then
$
  (
    \funcname{f}^{i+1}_0 \compose \funcname{f}^i_0,
    \funcname{f}^{i+1}_1 \compose \funcname{f}^i_1
  )
$
is a semi-strict (strict)
$\catname{M}^i_0$-$\catname{M}^{i+1}_0$-%
$\seqname{M}^i_0$-$\seqname{M}^{i+1}_0$-%
$\catname{M}^i_1$-$\catname{M}^{i+1}_1$
m-atlas morphism from $\seqname{A}^i$ to $\seqname{A}^{i+2}$ in the
coordinate model spaces $\seqname{M}^i_i$, $\seqname{M}^{i+2}_i$ by
\cref{lem:M-ATLmorph}.
\end{enumerate}
\end{proof}
\end{lemma}

\begin{definition}[$\LCS_\Ar$]
\label{def:LCSmorphAr}
Let
$\catseqname{M}^i \defeq (\catname{M}^i_\alpha, \alpha \prec \Alpha)$,
$i=1,2$, be sequences of categories and
$\Sigma \defeq (\sigma_\gamma, \gamma \prec \Gamma) \defeq
\bigl ( (\sigma_{\gamma,{\alpha_\beta}}, \alpha_\beta \prec \Alpha,
\beta \preceq \Beta_\gamma), \gamma \prec \Gamma \bigr )$
a signature over $\catseqname{M}^i$.

\begin{multline}
\LCS_\Ar(\catseqname{M}^1, \Sigma, \catseqname{M}^2) \defeq
\set
  {{
    \bigl (
      \funcseqname{f},
      (
        \catseqname{M}^1,
        \seqname{M}^1,
        \seqname{A}^1,
        \funcseqname{F}^1,
        \Sigma
      ),
      (
        \catseqname{M}^2,
        \seqname{M}^2,
        \seqname{A}^2,
        \funcseqname{F}^2,
        \Sigma
      )
    \bigr )
  }}%
  [
    {
      \seqname{M}^i \seqin \catseqname{M}^i
    },
    {
    \isLCS_\Ar
      (
        \catseqname{M}^1,
        \seqname{M}^1,
        \seqname{A}^1,
        \funcseqname{F}^1,
        \Sigma,
        \catseqname{M}^2,
        \seqname{M}^2,
        \seqname{A}^2,
        \funcseqname{F}^2,
        \funcseqname{f}
      )
    }
  ]*
\end{multline}
\begin{equation}
\LCS
  (\catseqname{M}^i, \Sigma)
\defeq
  \bigl (
    \LCS_\Ob
      (
        \catseqname{M}^i, \Sigma
      ),
    \LCS_\Ar
      (
        \catseqname{M}^i,
        \Sigma,
        \catseqname{M}^i
      ),
    \compose[A]
  \bigr )
\end{equation}

If $\catseqname{M}^1 \SUBCAT[full-] \catseqname{M}^2$ then

\begin{multline}
\strict[semi-]{\LCS_\Ar}(\catseqname{M}^1, \Sigma, \catseqname{M}^2)
\defeq
\set
  {{
    \bigl (
      \funcseqname{f}.
      (
        \catseqname{M}^1,
        \seqname{M}^1,
        \seqname{A}^1,
        \funcseqname{F}^1,
        \Sigma
      ),
      (
        \catseqname{M}^2,
        \seqname{M}^2,
        \seqname{A}^2,
        \funcseqname{F}^2,
        \Sigma
      )
    \bigr )
  }}%
  [
    {
      \seqname{M}^i \seqin \catseqname{M}^i
    },
    {
      \strict[semi-]{\isLCS_\Ar}
        (
          \catseqname{M}^1,
          \seqname{M}^1,
          \seqname{A}^1,
          \funcseqname{F}^1,
          \Sigma,
          \seqname{M}^2,
          \seqname{A}^2,
          \funcseqname{F}^2,
          \funcseqname{f}
        )
    }
  ]*
\end{multline}

\begin{multline}
\strict{\LCS_\Ar}(\catseqname{M}^1, \Sigma, \catseqname{M}^2) \defeq
\set
  {{
    \bigl (
      \funcseqname{f}.
      (
        \catseqname{M}^1,
        \seqname{M}^1,
        \seqname{A}^1,
        \funcseqname{F}^1,
        \Sigma
      ),
      (
        \catseqname{M}^2,
        \seqname{M}^2,
        \seqname{A}^2,
        \funcseqname{F}^2,
        \Sigma
      )
    \bigr )
  }}%
  [
    {
      \seqname{M}^i \seqin \catseqname{M}^i
    },
    {
      \strict{\isLCS_\Ar}
        (
          \catseqname{M}^1,
          \seqname{M}^1,
          \seqname{A}^1,
          \funcseqname{F}^1,
          \Sigma,
          \seqname{M}^2,
          \seqname{A}^2,
          \funcseqname{F}^2,
          \funcseqname{f}
        )
    }
  ]*
\end{multline}

\begin{equation}
\strict[semi-]{\LCS}(\catseqname{M}^1, \Sigma) \defeq
\Bigl (
  \LCS_\Ob(\catseqname{M}^1, \Sigma),
  \strict[semi-]{\LCS_\Ar}(\catseqname{M}^1, \Sigma,\catseqname{M}^1),
  \compose[A]
\Bigr )
\end{equation}
\begin{equation}
\strict{\LCS}(\catseqname{M}^1, \Sigma) \defeq
\Bigl (
  \LCS_\Ob(\catseqname{M}^1, \Sigma),
  \strict{\LCS_\Ar}(\catseqname{M}^1, \Sigma,\catseqname{M}^1),
  \compose[A]
\Bigr )
\end{equation}

Let
$
  \seqname{M}^i \defeq (\seqname{M}^i_\alpha, \alpha \prec \Alpha)
  \seqin \catseqname{M}^i
$,
$i=1,2$,
be a sequence of spaces,
$
  \seqname{P}^i \defeq
  (
    \catseqname{M}^i,
    \seqname{M}^i,
    \funcseqname{F}^i,
    \Sigma
  )
$
and
$
  \seqname{L}^i \defeq
  (
    \catseqname{M}^i,
    \seqname{M}^i,
    \seqname{A}^i,
    \funcseqname{F}^i,
    \Sigma
  )
$
be a local $\seqname{M}^i$-$\Sigma$ coordinate space.

If $\catseqname{M}^1 \SUBCAT[full-] \catseqname{M}^2$ then
$\seqname{S}^1 \SUBSETEQ \seqname{S}^2$,
$\seqname{A}^1 = \seqname{A}^2$ and each
$
  \funcname{F}^1_\gamma =
  \funcname{F}^2_\gamma \maps
    \bigtimes \underline{\overset{\sigma_\gamma}{\seqname{S}^1}} \to
    \tail \left ( \overset{\sigma_\gamma}{\seqname{S}^1} \right )
$
then

\begin{multline}
\Id_{\seqname{L}^1,\seqname{L}^2} \defeq
\bigl (
  \ID_{\seqname{M}^1,\seqname{M}^2},
  \seqname{L}^1,
  \seqname{L}^2
\bigr )
\end{multline}
\begin{multline}
\Id_{\seqname{L}^i} \defeq \Id_{\seqname{L}^i,\seqname{L}^i}
\end{multline}

This nomenclature is justified below.
\end{definition}

\begin{theorem}[$\LCS(\catseqname{M}, \Sigma)$ is a category]
\label{the:LCSiscat}
Let $\catseqname{M} \defeq (\catname{M}_\alpha, \alpha \prec \Alpha)$
be a sequence of categories and
$
\Sigma \defeq
(\sigma_\gamma, \gamma \prec \Gamma) \defeq
\bigl (
  (
    \sigma_{\gamma,{\alpha_\beta}}, \beta \preceq \Beta_\gamma
  ),
\gamma \prec \Gamma
\bigr )
$ \\
a signature over $\catseqname{M}$. Then
$\LCS(\catseqname{M}, \Sigma)$,
$\strict[semi-]{\LCS}(\catseqname{M}, \Sigma)$ and \\
$\strict{\LCS}(\catseqname{M}, \Sigma)$ are categories.

Let $\seqname{L}^i \objin \LCS(\catseqname{M}, \Sigma)$.
Then $\Id_{\seqname{L}^i}$ is the identity morphism of $\seqname{L}^i$:

\begin{proof}
Let
$
  \seqname{L}^i
  \defeq
  (
    \seqname{M}^i,
    \seqname{A}^i,
    \funcseqname{F}^i,
    \Sigma
  )
$
be objects of $\LCS(\catseqname{M}, \Sigma)$, let
$
  \seqname{P}^i \defeq \\
  (
    \catseqname{M},
    \seqname{M}^i,
    \Sigma,
    \funcseqname{F}^i
  )
$
and let
$
  m^i \defeq
  \bigl (
    \funcseqname{f}^i \defeq (\funcname{f}^i_\alpha \maps
    \seqname{M}^i_\alpha \to \seqname{M}^{i+1}_\alpha,
    \alpha \prec \Alpha),
    L^i,
    L^{i+1},
  \bigr )
$
be morphisms.
\begin{enumerate}
\item Composition: \newline
$\funcseqname{f}^{i+1} \compose[()] \funcseqname{f}^i$ is a
prestructure morphism from $\seqname{M}^i$ to $\seqname{M}^{i+2}$
by \pagecref{lem:premorph}
and
$
\bigl (
  \funcseqname{f}^{i+1}_0 \compose \funcseqname{f}^i_0,
  \funcseqname{f}^{i+1}_1 \compose \funcseqname{f}^i_1
\bigr )
$
is an m-atlas morphism from $\seqname{A}^i$ to $\seqname{A}^{i+1}$ by
\pagecref{lem:M-ATLmorph}.

If each $\funcseqname{f}^i$ is a strict prestructure morphism from
$\seqname{P}^i$ to $\seqname{P}^{i+1}$ then
$\funcseqname{f}^{i+1} \compose[()] \funcseqname{f}^i$ is a strict
prestructure morphism from $\seqname{M}^i$ to $\seqname{M}^{i+2}$
by \cref{lem:premorph}

If each $(\funcseqname{f}^i_0, \funcseqname{f}^i_1)$ is a strict (semi-strict)
$\catname{M}_0$-$\catname{M}_0$-%
$\seqname{M}^i_0$-$\seqname{M}^{i+1}_0$-%
$\catname{M}_1$-$\catname{M}_1$
m-atlas morphism from $\seqname{A}^i$ to $\seqname{A}^{i+1}$ in the
coordinate spaces $\seqname{M}^i_1$, $\seqname{M}^{i+1}_1$ then
$
\bigl (
  \funcseqname{f}^{i+1}_0 \compose \funcseqname{f}^i_0,
  \funcseqname{f}^{i+1}_1 \compose \funcseqname{f}^i_1
\bigr )
$
is a strict (semi-strict)
$\catname{M}_0$-$\catname{M}_0$-%
$\seqname{M}^i_0$-$\seqname{M}^{i+1}_0$-%
$\catname{M}_1$-$\catname{M}_1$
m-atlas morphism from $\seqname{A}^i$ to $\seqname{A}^{i+1}$ in the
coordinate spaces $\seqname{M}^i_1$, $\seqname{M}^{i+1}_1$ by
\cref{lem:M-ATLmorph}.

\item Associativity: \newline
Composition is associative by \pagecref{lem:atlcomp}.

\item Identity: \newline
$\Id_{\seqname{L}^1}$ is the identity morphism of $\seqname{L}^1$
by \cref{lem:atlcomp}.

\end{enumerate}
\end{proof}
\end{theorem}

\subsection{Examples}
\label{sub:ex}
\label{Examples}
Functors can be constructed among special cases of local coordinate
spaces and many mathematical structures, although this paper only gives
details for two of them.

\begin{example}[Manifolds]
Choosing the coordinate category as open subsets of a Banach space or
more generally a Fr\'echet space, yields a manifold;
the definitions and results of
\pagecref{sec:man}\!, also cover manifolds with boundaries and manifolds
with tails.
\end{example}

\begin{example}[Fiber bundles]
Choosing the coordinate category as the space of products of open sets
in a topological space (base) with a fixed fiber and the morphisms as
fiber-preserving maps yields a fiber bundle; \pagecref{sec:bun} covers
the more general case of a restricted group of transition functions on
fibers.
\end{example}

\begin{example}[Lie groups]
Let $G$ be a topological group with group operation $\star$ and identity $\mathbf{1}_G$,
$C$ a linear space
\footnote{See \pagecref{def:lin} and \pagecref{def:trivck}.} and
$\seqname{A}$ a maximal $\Ck$-atlas\footnote{See \pagecref{def:Ck-ATLmax}.}
of $G$ in the coordinate space $C$. Define
$\hat{\star} \maps G \times G \to G$,
$\isCk_{\star,\seqname{A}} \maps G \times G \to \truthspace$ and
$\seqname{L}$ by:

\begin{equation}
 \hat{\star}(g_1,g_2)  \defeq  g_1 \star g_2^{-1}
\end{equation}
\begin{equation}
  \isCk_{\star,\mathbf{1}_G,\seqname{A}}(g_1,g_2) \defeq
  \begin{cases}
    \true  & \text{$\hat{\star}$ is $\Ck$ at $(g_1,g_2)$} \\*
    \false & \text{$\hat{\star}$ is not $\Ck$ at $(g_1,g_2)$}
  \end{cases}
\end{equation}
\begin{equation}
  \seqname{M} \defeq \Bigl ( \Triv{G}, \Triv[\Ck-]{C}, \seqname{Truthspace} \Bigr )
\end{equation}
\begin{equation}
  \catseqname{M} \defeq \Singcat{\seqname{M}}
\end{equation}
\begin{equation}
  \funcseqname{F} \defeq (\star, \isCk_{\star,\seqname{A}})
\end{equation}
\begin{equation}
  \Sigma \defeq
  \bigl (
    (0,0,0),
    (0,0,2)
  \bigr )
\end{equation}
\begin{equation}
  \seqname{L} \defeq
  (
    \catseqname{M},
    \seqname{M},
    \seqname{A},
    \funcseqname{F},
    \Sigma
  )
\end{equation}
Then $\seqname{L}$
is a local coordinate space iff $\hat{\star}$ is $\Ck$; the constrain
function $\isCk_{\star,\seqname{A}}$ expresses this condition.
For $k = \infty$, that LCS is equivalent to a Lie group.
\end{example}

\begin{example}[Fiber bundle with global sections]
Let $\seqname{B} \defeq (E, X, Y, \pi, G, \rho, \seqname{A})$ be a fiber
bundle\footnote{See \pagecref{def:Bun}.}, $\funcname{s} \maps X \to E$
and
\begin{equation}
  \mathrm{isSection}_{E,X,\pi,s} (x) \defeq
  \begin{cases}
    \true  & \pi \compose \funcname{s} (x) = x \\*
    \false & \pi \compose \funcname{s} (x) \ne x
  \end{cases}
\end{equation}
a constraint function. Then a simple modification of
\pagecref{def:BunLCS} gives a local coordinate space equivalent to
$\seqname{B}$ with the global section $\funcname{s}$.
\end{example}

\begin{example}[Minkowsky manifolds with foliations]
Similarly to the other examples, a local coordinate space can represent
a manifold, a global section of the tensor bundle, a constraint function
for the signature, a global time coordinate, constraint functions
enforcing diferentiability and a constraint function enforcing a foliation.
\end{example}

\section{Equivalence of manifolds}
\label{sec:man}
For a manifold\footnote{
  The literature has several different definitions of a manifold.
  This paper uses one chosen for ease of exposition.
}
the coordinate category is open subsets of a Banach space or more
generally a Fr\'echet space, with an appropriate choice of morphisms.
Choosing a separating hyperplane and half space, with open sets in the
chosen half space, allows manifolds with boundary. Similarly, choosing a
ball with tails allows a manifold with tails.

For differentiable
manifolds the coordinate category is similar, but the morphisms are
limited to those sufficiently differentiable, in order to impose a diferentiability
constraint on the transition functions
$\funcname{t}^\alpha_\beta=\phi_\beta \compose \phi_\alpha^{-1}$.

This section defines $\Ck$-atlases, $\Ck$-manifolds, local coordinate
spaces equivalent to $\Ck$-manifolds, categories of them and functors,
and gives some basic results.

\subsection {Linear spaces and linear model spaces}
\label{sub:lin}
\begin{definition}[Linear spaces]
\label{def:lin}
Let $S$ be a locally arcwise connected topological subspace, with
non-void interior, of a real (complex) Banach or Fr\'echet space. Then
$S$ is a real (complex) linear space.
\begin{remark}
This paper uses the term \textit{ball} to refer to balls of the
underlying space but uses the terms \textit{open set} and
\textit{neighborhoods} to refer to the relative topology.
\end{remark}

Let $\catname{S}$ be a small category whose objects are real (complex)
linear spaces and whose morphisms are $\Ck$ functions. Then
$\catname{S}$ is a $\Ck$ linear category.

Let $S$ be a real (complex) linear space. Then $\trivcat[\Ck-]{S}$, the
category of all $\Ck$ functions between open subspaces of $S$, is the
trivial $\Ck$ linear category of $S$.
\end{definition}

\begin{lemma}[Open subsets of linear spaces are locally arcwise connected]
Let $U$ be an open subset of the real (complex) linear space $S$. Then
$U$ is locally arcwise connected.

\begin{proof}
An open subset of a locally connected space is locally connected.
\end{proof}
\end{lemma}

\begin{definition}[Linear model spaces]
Let $S$ be a real (complex) linear space and
$\seqname{S} \defeq (S, \catname{S})$ a model space for $S$.
Then $\seqname{S}$ is a real (complex) linear model space.

Let $\seqname{S} \defeq (S, \catname{S})$ be a real (complex) linear
model space such that every morphism of $\catname{S}$ is a $\Ck$
function. Then $\seqname{S}$ is a real (complex) $\Ck$ linear model
space.

Let $\catname{S}$ be a small category whose objects are real (complex)
$\Ck$ linear model spaces and whose morphisms are $\Ck$ model functions.
Then $\catname{S}$ is a $\Ck$ linear model category.
\end{definition}

\begin{definition}[Trivial $\mathrm{C^k}$ linear model spaces]
\label{def:trivck}
Let $S$ be a real (complex) linear space and $\catname{S}$ the category
of all $\Ck$ functions between open sets of $S$. Then
$\Triv[\Ck-]{S} \defeq (S, \catname{S})$ is the trivial $\Ck$ linear
model space of $S$ and $\Triv[\Ck-]{S}$ is a real (complex) trivial
$\Ck$ linear model space.

Let $\seqname{S}$ be a set of real (complex) linear spaces.

The category of open trivial $\Ck$ model spaces in $\seqname{S}$,
abbreviated $\optriv[\Ck-]{\seqname{S}}$, is the category whose objects
are $\set {\Triv[\Ck-]{U}}[U \in \op{\seqname{S}}]$, the trivial $\Ck$
linear model spaces of non-null open sets of spaces in $\seqname{S}$,
and whose morphisms are all the $\Ck$ functions among them.

The set of trivial $\Ck$ linear model spaces of $\seqname{S}$ is
$
  \Triv[\Ck-]{\seqname{S}}
  \defeq
  \set
    {\Triv[\Ck-]{S'}}%
    [S' \in \seqname{S}]
$.

The category of trivial $\Ck$ linear model spaces of $\seqname{S}$,
abbreviated $\Trivcat[\Ck-]{\seqname{S}}$, is the category whose
objects are $\Triv[\Ck-]{\seqname{S}}$ and whose morphisms are all the
$\Ck$ functions among them.
\end{definition}

\begin{lemma}[The trivial $\mathrm{C^k}$ linear model space of $S$ is a linear model space]
Let $S$ be a real (complex) linear space. Then
$\Triv[\Ck-]{S}=(S,\catname{S})$ is a linear $\Ck$ model space.

\begin{proof}
$\Triv[\Ck-]{S}$ satisfies the conditions in \pagecref{def:model}:
\begin{enumerate}
\item $\Ob(\catname{S})$ is an open cover for $S$.
\item $\Ob(\catname{S})$ is closed under finite intersections.
\item The morphisms of $\catname{S}$ are $\Ck$, hence continuous.
\item If $f \maps A \to B$ is a morphism,
$A' \in \Ob(\catname{S}) \subseteq A \in \Ob(\catname{S})$,
$B' \in \Ob(\catname{S}) \subseteq B \in \Ob(\catname{S})$
and $f[A'] \subseteq B'$ then $f \restriction_{A'} \maps A' \to B'$ is
continuous and thus a morphism.
\item If $A' \in \Ob(\catname{S}) \subseteq A \in \Ob(\catname{S})$ then
the inclusion map $\funcname{i} \maps A' \hookrightarrow A$ is
$\Ck$ and thus a morphism.
\item Restricted sheaf condition: Whenever
\begin{enumerate}
\item $U_\alpha$ and $V_\alpha$, $\alpha \prec \Alpha$,
are objects of $\catname{S}$.
\item $\funcname{f}_\alpha \maps U_\alpha \to V_\alpha$ are morphisms of
$\catname{S}$.
\item
$U \defeq \union[\alpha \prec \Alpha]{U_\alpha} \in \Ob(\catname{S})$,
\item
$V \defeq \union[\alpha \prec \Alpha]{V_\alpha} \in \Ob(\catname{S})$
\item
$\funcname{f} \maps U \to V$ is a continuous function and
for every $\alpha \prec \Alpha$,
$\funcname{f}$ agrees with $\funcname{f}_\alpha$ on $U_\alpha$
\end{enumerate}
then $\funcname{f}$ is $\Ck$ and thus a morphism of $\catname{S}$.
\end{enumerate}
\end{proof}
\end{lemma}

\begin{definition}[ $\mathrm{C^k}$ singleton categories]
\label {def:singCk}
Let $\seqname{C}$ be a $\Ck$ linear model space. Then the $\Ck$
singleton category of $\seqname{C}$, abbreviated
$\singcat[\Ck-]{\seqname{C}}$, is the category whose sole object is
$\seqname{C}$ and whose morphisms are all of the $\Ck$ model functions
from $\seqname{C}$ to itself.
\end{definition}

\subsection{$\Ck$-nearly commutative diagrams}
\label{sub:ckncd}
Let $C$ be a linear space, $\catname{C} \defeq \trivcat[\Ck]{S}$ and
$D$ a tree with two branches, whose nodes are topological spaces $U_i$
and $V^j$ and whose links are continuous functions
$\funcname{f}_i \maps U_i \to U_{i+1}$ and
$\funcname{f}'_j \maps U_j \to U_{j+1}$ between the spaces:

\begin{align*}
D = \{
  &
    \funcname{f}_0 \maps U_0 = V_0 \to U_1,
    \dotsc,
    \funcname{f}_{m - 1} \maps U_{m - 1} \to U_m,
\\
  &
    \funcname{f}'_0 \maps  U_0 = V_0 \to V_1,
    \dotsc,
    \funcname{f}'_{m-1} \maps V_{m-1} \to V_n
\}
\end{align*}
with $U_0 = V_0$, $U_m \subseteq C$ and $V_n \subseteq C$, as shown in
\pagecref{fig:NCDb}.

\begin{definition}[$\mathrm{C^k}$-nearly commutative diagrams]
$D$ is $\Ck$-nearly commutative in linear space $C$ iff $D$ is nearly
commutative in category $\catname{C}$.
\end{definition}

\begin{definition}[$\mathrm{C^k}$-nearly commutative diagrams at a point]
\label{def:CkNCD}
Let $C$, $\catname{C}$ and $D$ be as above and $x$ be an element of the
initial node. $D$ is $\Ck$-nearly commutative in $C$ at $x$
iff $D$ is nearly commutative in $\catname{C}$ at $x$.
\end{definition}

\begin{definition}[$\mathrm{C^k}$-locally nearly commutative diagrams]
Let $C$, $\catname{C}$ and $D$ be as above. $D$ is $\Ck$-locally nearly
commutative in $C$ iff $D$ is locally nearly commutative in
$\catname{C}$.
\end{definition}

\subsection{$\Ck$ charts}
\label{sub:ckcharts}
\begin{definition}[$\mathrm{C^k}$ charts]
\label{def:Ck chart}
Let $C$ be a linear space and $E$ a topological space. A $\Ck$\footnote{
With $k \in \mathbb{N} \cup \{ \infty, \omega \}$.}
chart $(U, V, \phi)$ of $E$ in the coordinate space $C$ consists of

\begin{enumerate}
\item An open subset $U \subseteq E$, known as a coordinate patch
\item An open subset $V \subseteq C$
\item A homeomorphism $\phi \maps U \toiso V$, known as a
coordinate function
\end{enumerate}
\begin{remark}
I consider it clearer to explicate the range, rather than the
conventional usage of specifying only the domain and function or the
minimalist usage of specifying only the function.
\end{remark}
\end{definition}

\begin{definition}[$\mathrm{C^k}$ subcharts]
Let $(U, V, \phi)$ be a $\Ck$ chart of $E$ in the coordinate space $C$
and $U' \subseteq U$ open. Then
$(U', V', \phi') \defeq (U', \phi[U'], \phi \restriction_{U',V'})$
is a subchart of $(U, V, \phi)$.

By abuse of language we will write $(U', V', \phi)$ for
$(U', V', \phi')$.
\end{definition}

\begin{lemma}[$\mathrm{C^k}$ subcharts]
Let $(U, V, \phi)$ be a $\Ck$ chart of $E$ in the coordinate space $C$
and $(U', V', \phi')$ a subchart of $(U, V, \phi)$. Then \\
$(U', V', \phi')$ is a $\Ck$ chart of $E$ in the coordinate space $C$.
\begin{proof}
$(U',V',\phi')$ satisfies the conditions of \cref{def:Ck chart}
\begin{enumerate}
\item
$U'$ is open by hypothesis.
\item
$\phi$ is a homeomorphism, so $V'=\phi[U']$ is also open.
\item
$\phi$ is a homeomorphism, so
$\phi \restriction_{U'} \maps U' \toiso \phi[U']$ is also.
\end{enumerate}
\end{proof}
\end{lemma}

\begin{definition}[$\mathrm{C^k}$ compatibility]
Let $(U, V, \phi)$ and  $(U',V', \phi')$ be $\Ck$ charts of $E$ in the
coordinate space $C$. Then $(U, V, \phi)$ is $\Ck$ compatible with
$(U', V', \phi')$ iff either
\begin{enumerate}
\item $U$ and  $U'$ are disjoint
\item The transition function $\funcname{t}=\phi' \compose \phi^{-1} \restriction_{\phi[U \cap U']}$
is a $\Ck$ diffeomorphism.
\end{enumerate}
\end{definition}

\begin{lemma}[Symmetry of $\mathrm{C^k}$ compatibility]
Let $(U, V, \phi)$ and  $(U',V', \phi')$ be $\Ck$ charts of $E$ in the
coordinate space $C$.
Then $(U, V, \phi)$ is $\Ck$ compatible with $(U', V', \phi')$ iff
$(U', V', \phi')$ is $\Ck$ compatible with $(U,V,\phi)$.

\begin{proof}
It suffices to prove the implication in only one direction.
\begin{enumerate}
\item $U \cap U' = U' \cap U$.
\item Since the transition function
$\funcname{t}=\phi' \compose \phi^{-1} \restriction_{\phi[U \cap U']}$
is a  $\Ck$ diffeomorphism of $C$, so is
$\funcname{t}^{-1} = \phi \compose \phi'^{-1} \restriction_{\phi'[U \cap U']}$.
\end{enumerate}
\end{proof}
\end{lemma}

\begin{lemma}[$\mathrm{C^k}$ compatibility of subcharts]
Let $(U_i, V_i, \phi_i)$, $i=1,2$, be $\Ck$ charts of $E$ in the
coordinate space $C$, $(U'_i, V'_i, \phi'_i)$ be subcharts and
$(U_1, V_1, \phi_1)$ be $\Ck$ compatible with $(U_2, V_2, \phi_2)$. Then
$(U'_1, V'_1, \phi'_1)$ is $\Ck$ compatible with
$(U'_2, V'_2, \phi'_2)$.

\begin{proof}
If $U_1 \cap U_2 = \emptyset$ then $U'_1 \cap U'_2 = \emptyset$. If
$U'_1 \cap U'_2 = \emptyset$ then $(U'_1, V'_1, \phi'_1)$ is $\Ck$
compatible with $(U'_2, V'_2, \phi'_2)$. Otherwise, the transition
function
$
  t^1_2 \defeq
  \phi_2 \compose \phi^{-1}_1 \restriction_{\phi_1[U_1 \cap U_2]}
$
is a $\Ck$ diffeomorphism and hence
$
  t^1_2 \restriction_{\phi_1[U'_1 \cap U'_2]} \maps
  \phi_1[U'_1 \cap U'_2] \toiso
  \phi_2[U'_1 \cap U'_2]
$
is a $\Ck$ diffeomorphism.
\end{proof}
\end{lemma}

\begin{corollary}[$\mathrm{C^k}$ compatibility with subcharts]
Let $(U, V, \phi)$ be a $\Ck$ charts of $E$ in
the coordinate space $C$ and $(U', V', \phi')$ a
subchart. Then $(U', V', \phi')$ is $\Ck$ compatible with
$(U,V,\phi)$.

\begin{proof}
$(U, V, \phi)$ is $\Ck$ compatible with itself and is a subchart of
itself,
\end{proof}
\end{corollary}

\begin{definition}[Covering by $\mathrm{C^k}$ charts]
Let $\seqname{A}$ be a set of charts of $E$ in the coordinate space $C$.
$\seqname{A}$ covers $E$ iff $\pi_1[\seqname{A}]$ covers $E$, i.e.,
$E = \union{{\pi_1[\seqname{A}]}}$.
\end{definition}

\subsection{$\Ck$-atlases}
\label{sub:ckatlases}
A set of charts can be atlases for different coordinate spaces even if
it is for the same total space.  In order to aggregate them into
categories, there must be a way to distinguish them. Including the
two\footnote{The total space is redundant, but convenient.}
spaces in the definitions of the categories serves the purpose.

\begin{definition}[$\mathrm{C^k}$-atlases]
Let $\seqname{A}$ be a set of mutually $\Ck$ compatible charts of $E$ in
the coordinate space $C$. $\seqname{A}$ is a $\Ck$-atlas of $E$ in
the coordinate space $C$, abbreviated
$\isAtl^\Ck_\Ob(\seqname{A}, E, C)$, iff
\begin{enumerate}
\item $\seqname{A}$ covers $E$
\item There is at least one chart $(U,V,\phi) \in \seqname{A}$
where $V$ contains a ball of the underlying Banach or Fr\'echet space.
\end{enumerate}

$\seqname{A}$ is a full $\Ck$-atlas of $E$ in
the coordinate space $C$, abbreviated
$\full{\isAtl^\Ck_\Ob}(\seqname{A}, E, C)$, iff
\begin{enumerate}
\item $\pi_1[\seqname{A}]$ covers $E$
\item $\pi_2[\seqname{A}]$ covers $C$.
\item There is at least one chart $(U,V,\phi) \in \seqname{A}$
where $V$ contains a ball of the underlying Banach or Fr\'echet space.
\end{enumerate}

By abuse of language we write $U \in \seqname{A}$ for
$U \in \pi_1[\seqname{A}]$.

Let $E$ be a topological spaces and $C$ a linear space. Then
\begin{equation}
\Atl^\Ck_\Ob(E,C) \defeq
\set
{{
  (\seqname{A}, E, C)
}}%
[
  {{
    \isAtl^\Ck_\Ob(\seqname{A}, E, C}
  )}
]
\end{equation}
\begin{equation}
\full{\Atl^\Ck_\Ob}(E,C) \defeq
\set
{{
  (\seqname{A}, E, C)
}}%
[
  {{
    \full{\isAtl^\Ck_\Ob}(\seqname{A}, E, C}
  )}
]
\end{equation}

Let $\seqname{E}$ be a set of topological spaces and $\seqname{C}$ a
set of linear spaces. Then
\begin{equation}
\Atl^\Ck_\Ob(\seqname{E}, \seqname{C}) \defeq
  \union
    [
        {E_\mu \in \seqname{E}},
        {C_\mu \in \seqname{C}}
    ]
    {{
      \Atl_\Ob(E_\mu, C_\mu)
    }}
\end{equation}
\begin{equation}
\full{\Atl^\Ck_\Ob}(\seqname{E}, C) \defeq
\set
{{
  (\seqname{A}, E, C)
}}%
[
  {{
    \full{\isAtl^\Ck_\Ob}(\seqname{A}, E, C}
  )}
]
\end{equation}
\end{definition}

\begin{definition}[Compatibility of charts with $\mathrm{C^k}$-atlases]
A chart $(U, V, \phi)$ of $E$ in the coordinate space $C$ is $\Ck$
compatible with a $\Ck$-atlas $\seqname{A}$ iff it is $\Ck$ compatible with
every chart in the atlas.
\end{definition}

\begin{lemma}[Compatibility of subcharts with $\mathrm{C^k}$-atlases]
\label{lem:CkCompat}
Let $\seqname{A}$ be a $\Ck$-atlas of $E$ in the coordinate space $C$
and $\seqname{C}_1 = (U_1, V_1, \phi_1)$ a $\Ck$ chart in
$\seqname{A}$. Then any subchart of $\seqname{C}_1$ is $\Ck$ compatible
with $\seqname{A}$.

\begin{proof}
Let $\seqname{C}' = (U', V', \phi')$ be a subchart of $\seqname{C}_1$ and
$\seqname{C}_2 = (U_2, V_2, \phi_2)$ another chart in $\seqname{A}$.
\begin{enumerate}
\item If $U_1 \cap U_2 = \emptyset$, then $U' \cap U_2 = \emptyset$.
\item If $U' \cap U_2 = \emptyset$ then $\seqname{C}'$ is $\Ck$
compatible with $\seqname{C}_2$.
\item Otherwise the transition function
$t^1_2 \defeq \phi_2 \compose \phi^{-1}_1 \restriction_{\phi_1[U_1 \cap U_2]}$
is a $\Ck$ diffeomorphism and thus
$t^1_2 \restriction_{\phi_1[U' \cap U_2]}$ is a $\Ck$ diffeomorphism.
\end{enumerate}
\end{proof}
\end{lemma}

\begin{lemma}[Extensions of $\mathrm{C^k}$-atlases]
\label{Ck-atl:extensions}
Let $\seqname{A}$ be a $\Ck$ atlas of $\seqname{E}$ in the coordinate
space $C$ and $(U_i,V_i,\phi_i)$, $i=1,2$ be $\Ck$ charts of $E$ in
the coordinate space $C$ $\Ck$ compatible with $\seqname{A}$ in the
coordinate space $C$. Then $(U_1,V_1,\phi_1)$ is $\Ck$ compatible with
$(U_2,V_2,\phi_2)$ in the coordinate space $C$.

\begin{proof}
If $U_1 \cap U_2 = \emptyset$ then $(U_1,V_1,\phi_1)$ is
$\Ck$ compatible with $(U_2,V_2,\phi_2)$. Otherwise,
$\phi_2 \compose \phi^{-1}_1 \restriction_{\phi_1[U_1 \cap U_2]} \maps
\phi_1[U_1 \cap U_2] \toiso \phi_2[U_1 \cap U_2]$
is a homeomorphism.  It remains to show that
$\phi_2 \compose \phi^{-1}_1 \restriction_{\phi_1[U_1 \cap U_2]}$
is a $\Ck$ diffeomorphism.
Let $(U'_\alpha,V'_\alpha,\phi'_\alpha)$, $\alpha \prec \Alpha$, be
charts in $\seqname{A}$ such that
$U_1 \cap U_2 \subseteq \union[\alpha \prec \Alpha]{U'_\alpha}$ and
$U_1 \cap U_2 \cap U'_\alpha \neq \emptyset$, $\alpha \prec \Alpha$.
Since the charts are $\Ck$ compatible with $(U'_\alpha,V'_\alpha,\phi'_\alpha)$,
$\phi_2 \compose \phi'^{-1}_\alpha \restriction_{U_1 \cap U_2 \cap U'_\alpha}$
and
$\phi'_\alpha \compose \phi^{-1}_1 \restriction_{U_1 \cap U_2 \cap U'_\alpha}$
are $\Ck$ diffeomorphisms and thus
$\phi_2 \compose \phi^{-1}_1 = \phi_2 \compose \phi'^{-1}_\alpha \compose \phi'_\alpha \compose \phi^{-1}_1$
is a $\Ck$ diffeomorphism.
\end{proof}
\end{lemma}

\begin{definition}[Maximal $\mathrm{C^k}$-atlases]
\label{def:Ck-ATLmax}
Let $E$ be a topological spaces and $C$ a linear space. Then
$\seqname{A}$ is a maximal $\Ck$-atlas of $E$ in the coordinate space
$C$, abbreviated $\maximal{\isAtl^\Ck_\Ob}(\seqname{A}, E,C)$, iff
$\seqname{A}$ is a $\Ck$-atlas that cannot be extended by adding an
additional $\Ck$ compatible chart.
$\seqname{A}$ is a semi-maximal $\Ck$-atlas of $E$ in
the coordinate space $\seqname{C}$, abbreviated
$\maximal[S-]{\isAtl^\Ck_\Ob}(\seqname{A}, E, \seqname{C})$, iff whenever
$(U,V,\phi) \in \seqname{A}$, $U' \subseteq U, V' \subseteq V$ and
$V'' \subseteq C$ are open, $\phi[U'] = V'$ and
$\phi' \maps V' \toiso V''$ is a $\Ck$ diffeomorphism then
$(U', V'', \phi' \compose \phi) \in \seqname{A}$.

\begin{equation}
  \maxfull{\isAtl^\Ck_\Ob}(\seqname{A}, E, C) \defeq
    \full{\isAtl^\Ck_\Ob}(\seqname{A}, E, C) \land
    \maximal{\isAtl^\Ck_\Ob}(\seqname{A}, E, C)
\end{equation}
\begin{equation}
  \maxfull[S-]{\isAtl^\Ck_\Ob}(\seqname{A}, E, C) \defeq
    \full{\isAtl^\Ck_\Ob}(\seqname{A}, E, C) \land
    \maximal[S-]{\isAtl^\Ck_\Ob}(\seqname{A}, E, C)
\end{equation}

\begin{equation}
\maximal{\Atl^\Ck_\Ob}(E,C) \defeq
  \set
    {{(\seqname{A}, E, C)}}%
    [{{\maximal{\isAtl^\Ck_\Ob}(\seqname{A}, E, C)}}]
\end{equation}
\begin{equation}
\maximal[S-]{\Atl^\Ck_\Ob}(E,C) \defeq
  \set
    {{(\seqname{A}, E, C)}}%
    [{{\maximal[S-]{\isAtl^\Ck_\Ob}(\seqname{A}, E, C)}}]
\end{equation}
\begin{equation}
\maxfull[S-]{\Atl^\Ck_\Ob}(E,C) \defeq
  \set
    {{(\seqname{A}, E, C)}}%
    [{{\maxfull[S-]{\isAtl^\Ck_\Ob}(\seqname{A}, E, C)}}]
\end{equation}

Let $\seqname{E}$ be a set of topological spaces and $\seqname{C}$ a
set of linear spaces. Then
\begin{equation}
\maximal{\Atl^\Ck_\Ob}(\seqname{E}, \seqname{C}) \defeq
\set
{{
  (\seqname{A}, E \in \seqname{E}, C \in \seqname{C})
}}%
[
  {{
    \maximal{\isAtl^\Ck_\Ob}(\seqname{A}, E, C}
  )}
]
\end{equation}
\begin{equation}
\maximal[S-]{\Atl^\Ck_\Ob}(\seqname{E}, \seqname{C}) \defeq
\set
{{
  (\seqname{A}, E \in \seqname{E}, C \in \seqname{C})
}}%
[
  {{
    \maximal[S-]{\isAtl^\Ck_\Ob}(\seqname{A}, E, C}
  )}
]
\end{equation}
\begin{equation}
\maxfull[S-]{\Atl^\Ck_\Ob}(\seqname{E}, \seqname{C}) \defeq
\set
{{
  (\seqname{A}, E \in \seqname{E}, C \in \seqname{C})
}}%
[
  {{
    \maxfull[S-]{\isAtl^\Ck_\Ob}(\seqname{A}, E, C}
  )}
]
\end{equation}
\end{definition}

\begin{lemma}[Maximal $\Ck$-atlases are semi-maximal $\Ck$-atlases]
Let $E$ be a topological space, $C$ a $\Ck$-linear space and
$\seqname{A}$ a maximal $\Ck$-atlas of $E$ in the coordinate space $C$.
Then $\seqname{A}$ is a semi-maximal $\Ck$-atlas of $E$ in the
coordinate space $C$.

\begin{proof}
Let $(U,V,\phi) \in \seqname{A}$, $U' \subseteq U, V' \subseteq V$ and
$V'' \subseteq C$ be open, $\phi[U'] = V'$ and
$\phi' \maps V' \toiso V''$ be a $\Ck$ diffeomorphism.  $(U', V', \phi)$
is a subchart of $(U,V,\phi)$ and by \pagecref{lem:CkCompat} is
$\Ck$ compatible with the charts of $\seqname{A}$. Since $\phi'$ is a
$\Ck$ diffeomorphism, $(U', V'', \phi' \compose \phi)$ is $\Ck$
compatible with the charts of $\seqname{A}$.  Since $\seqname{A}$ is
maximal, $(U', V'', \phi' \compose \phi)$ is a chart of $\seqname{A}$.
\end{proof}
\end{lemma}

\begin{theorem}[Existence and uniqueness of maximal $\mathrm{C^k}$-atlases]
Let $\seqname{A}$ be a $\Ck$-atlas of $E$ in the coordinate space $C$.
Then there exists a unique maximal $\Ck$-atlas
$\maximal{\Atlas^\Ck}(\seqname{A}, E, C)$ of $E$ in the coordinate space $C$
compatible with $\seqname{A}$.

\begin{proof}
Let $\seqname{P}$ be the set of all $\Ck$-atlases of $\seqname{E}$ in
the coordinate space $C$ containing $\seqname{A}$ and $\Ck$ compatible
in the coordinate space $\seqname{C}$ with all of the $\Ck$ charts in
$\seqname{A}$. Let $\maximal{\seqname{P}}$ be a maximal chain of
$\seqname{A}$. Then $A'=\union{\maximal{\seqname{P}}}$ is a maximal
$\Ck$ atlas of $\seqname{E}$ in the coordinate space $\seqname{C}$
$\Ck$ compatible with $\seqname{A}$.
Uniqueness follows from \pagecref{Ck-atl:extensions}.
\end{proof}
\end{theorem}

\subsection{$\Ck$-atlas morphisms and functors}
\label{sub:ckmorph}
This section defines categories of $\Ck$-atlases
$\bigl ( \Atl^\Ck(\seqname{E}, \seqname{C})$,
$\full{\Atl^\Ck}(\seqname{E}, \seqname{C})$ and
$\maximal{\Atl^\Ck}(\seqname{E}, \seqname{C}) \bigr )$, constructs
functors \\
$\bigl ( \Functor^\Ck_{\Ck,\M}\bigr )$
between them and categories of m-atlases
$\bigl ( \Atl(\Triv{\seqname{E}}$, $\Triv[\Ck-]{\seqname{C}})\bigr )$
and constructs inverse functors $\bigl (\Functor^\Ck_{\M,\Ck}\bigr )$.

\begin{definition}[$\mathrm{C^k}$-atlas morphisms]
\label{def:Ck-ATLmorph}
Let $E^i$, $i=1,2$, be topological spaces, $C^i$ linear spaces and
$\seqname{A}^i$ $\Ck$-atlases of $E^i$ in the coordinate spaces $C^i$. A
pair\footnote{
  The conventional definition uses only the first of the two
  functions and a slightly different compatibility condition.
}
of functions $(\funcname{f}_0, \funcname{f}_1)$ is a (full)
$E^1$-$E^2$ $\Ck$-morphism of $\seqname{A}^1$ to $\seqname{A}^2$ in the
coordinate spaces $C^1$, $C^2$, abbreviated as
$\isAtl^\Ck_\Ar(\seqname{A}^1, E^1, C^1$,
                  $\seqname{A}^2$, $E^2$, $C^2$,
                  $\funcname{f}_0$, $\funcname{f}_1)$,
iff
\begin{enumerate}
\item $\funcname{f}_0 \maps E^1 \to E^2$ is a continuous function.
\item $\funcname{f}_1 \maps C^1 \to C^2$ is a $\Ck$ function.
\item for any
$(U^1, V^1, \phi^1 \maps U^1 \toiso V^1) \in \seqname{A}^1$,
$(U^2, V^2, \phi^2 \maps U^2 \toiso V^2) \in \seqname{A}^2$,
the diagram
$D \defeq (\{I \defeq U^1 \cap \funcname{f}_0^{-1}[U^2], V^1, E^2, U^2, V^2 \}$,
$\{ \funcname{f}_0, \phi^2, \phi^1, \funcname{f}_1 \})$
is $\Ck$-locally nearly commutative in $\seqname{C}^2$, i.e., for any
$x \in I$
there
are open sets $U'^1 \subseteq I$,
$V'^1 \subseteq V^1$,
$U'^2 \subseteq U^2$, $V'^2 \subseteq V^2$,
$\hat{V}'^2 \subseteq C^2$ and a $\Ck$ diffeomorphism
$\hat{\funcname{f}} \maps \hat{V}'^2 \toiso V'^2$ such that
\crefrange{eq:xin}{eq:fphigeqgphi} on \cpagerefrange{eq:xin}{eq:fphigeqgphi}
in \pagecref{def:M-ATLmorph} hold.
\end{enumerate}

$(\funcname{f}_0, \funcname{f}_1)$ is also a full
$E^1$-$E^2$ $\Ck$-morphism of $\seqname{A}^1$ to $\seqname{A}^2$ in the
coordinate spaces $C^1$, $C^2$, abbreviated as
$
  \full{\isAtl^\Ck_\Ar}
    (
      \seqname{A}^1, E^1, C^1,
      \seqname{A}^2, E^2, C^2,
      \funcname{f}_0, \funcname{f}_1
    )
$,
iff
$
  \union{{ \pi_2[\seqname{A}^1] }} = C^1 \land
  \union{{ \pi_2[\seqname{A}^2] }} = C^2
$.

$(\funcname{f}_0, \funcname{f}_1)$ is also a maximal
$E^1$-$E^2$ $\Ck$-morphism of $\seqname{A}^1$ to $\seqname{A}^2$ in the
coordinate spaces $C^1$, $C^2$, abbreviated as
$
  \maximal{\isAtl^\Ck_\Ar}
    (
      \seqname{A}^1, E^1, C^1,
      \seqname{A}^2, E^2, C^2,
      \funcname{f}_0, \funcname{f}_1
    )
$ iff \\
$
  \maximal{\isAtl^\Ck_\Ob}
    (
      \seqname{A}^1, E^1, C^1
    )
$ and
$
  \maximal{\isAtl^\Ck_\Ob}
    (
     \seqname{A}^2, E^2, C^2
    )
$.

$(\funcname{f}_0, \funcname{f}_1)$ is also a semi-maximal
$E^1$-$E^2$ $\Ck$-morphism of $\seqname{A}^1$ to $\seqname{A}^2$ in the
coordinate spaces $C^1$, $C^2$, abbreviated as
$
  \maximal[S-]{\isAtl^\Ck_\Ar}
    (
      \seqname{A}^1, E^1, C^1,
      \seqname{A}^2, E^2, C^2,
      \funcname{f}_0, \funcname{f}_1
    )
$ iff \\
$
  \maximal[S-]{\isAtl^\Ck_\Ob}
    (
      \seqname{A}^1, E^1, C^1
    )
$ and
$
  \maximal[S-]{\isAtl^\Ck_\Ob}
    (
     \seqname{A}^2, E^2, C^2
    )
$.

$(\funcname{f}_0, \funcname{f}_1)$ is also a full maximal
$E^1$-$E^2$ $\Ck$-morphism of $\seqname{A}^1$ to $\seqname{A}^2$ in the
coordinate spaces $C^1$, $C^2$, abbreviated as
$
  \maxfull{\isAtl^\Ck_\Ar}
    (
      \seqname{A}^1, E^1, C^1,
      \seqname{A}^2, E^2, C^2,
      \funcname{f}_0, \funcname{f}_1
    )
$ iff \\
$
  \maximal{\isAtl^\Ck_\Ob}
    (
      \seqname{A}^1, E^1, C^1
    )
$,
$
  \maximal{\isAtl^\Ck_\Ob}
    (
     \seqname{A}^2, E^2, C^2
    )
$
and
$
  \union{{ \pi_2[\seqname{A}^1] }} = C^1 \land
  \union{{ \pi_2[\seqname{A}^2] }} = C^2
$.

$(\funcname{f}_0, \funcname{f}_1)$ is also a full semi-maximal
$E^1$-$E^2$ $\Ck$-morphism of $\seqname{A}^1$ to $\seqname{A}^2$ in the
coordinate spaces $C^1$, $C^2$, abbreviated as
$
  \maxfull[S-]{\isAtl^\Ck_\Ar}
    (
      \seqname{A}^1, E^1, C^1,
      \seqname{A}^2, E^2, C^2,
      \funcname{f}_0, \funcname{f}_1
    )
$ iff \\
$
  \maximal[S-]{\isAtl^\Ck_\Ob}
    (
      \seqname{A}^1, E^1, C^1
    )
$,
$
  \maximal[S-]{\isAtl^\Ck_\Ob}
    (
     \seqname{A}^2, E^2, C^2
    )
$
and
$
  \union{{ \pi_2[\seqname{A}^1] }} = C^1 \land
  \union{{ \pi_2[\seqname{A}^2] }} = C^2
$.

The identity morphism of $(\seqname{A}^i, E^i, C^i)$ is
\begin{equation}
\Id_{(\seqname{A}^i, E^i, C^i)} \defeq
\bigl (
  (\Id_{E^i}, \Id_{C^i}),
  (\seqname{A}^i, E^i, C^i),
  (\seqname{A}^i, E^i, C^i)
\bigr )
\end{equation}
This nomenclature will be justified later.

Let $E^i$, $i=1,2$, be topological spaces and $C^i$ be linear
spaces. Then
\begin{multline}
\Atl^\Ck_\Ar(E^1, C^1, E^2, C^2) \defeq
  \set
  {{
    \bigl (
      (
        \funcname{f}_0,
        \funcname{f}_1
      ),
      (\seqname{A}^1, E^1, C^1),
      (\seqname{A}^2, E^2, C^2)
    \bigr )
  }}%
  [{{
    \isAtl^\Ck_\Ar
    (
      \seqname{A}^1,
      E^1,
      C^1,
      \seqname{A}^2,
      E^2,
      C^2,
      \funcname{f}_0,
      \funcname{f}_1
    )
  }}]*
\end{multline}
\begin{multline}
\full{\Atl^\Ck_\Ar}(E^1, C^1, E^2, C^2) \defeq
  \set
  {{
    \bigl (
      (
        \funcname{f}_0,
        \funcname{f}_1
      ),
      (\seqname{A}^1, E^1, C^1),
      (\seqname{A}^2, E^2, C^2)
    \bigr )
  }}%
  [{{
    \full{\isAtl^\Ck_\Ar}
    (
      \seqname{A}^1,
      E^1,
      C^1,
      \seqname{A}^2,
      E^2,
      C^2,
      \funcname{f}_0,
      \funcname{f}_1
    )
  }}]*
\end{multline}
\begin{multline}
\maximal{\Atl^\Ck_\Ar}(E^1, C^1, E^2, C^2) \defeq
  \set
  {{
    \bigl (
      (
        \funcname{f}_0,
        \funcname{f}_1
      ),
      (\seqname{A}^1, E^1, C^1),
      (\seqname{A}^2, E^2, C^2)
    \bigr )
  }}%
  [{{
    \maximal{\isAtl^\Ck_\Ar}
    (
      \seqname{A}^1,
      E^1,
      C^1,
      \seqname{A}^2,
      E^2,
      C^2,
      \funcname{f}_0,
      \funcname{f}_1
    )
  }}]*
\end{multline}
\begin{multline}
\maximal[S-]{\Atl^\Ck_\Ar}(E^1, C^1, E^2, C^2) \defeq
  \set
  {{
    \bigl (
      (
        \funcname{f}_0,
        \funcname{f}_1
      ),
      (\seqname{A}^1, E^1, C^1),
      (\seqname{A}^2, E^2, C^2)
    \bigr )
  }}%
  [{{
    \maximal[S-]{\isAtl^\Ck_\Ar}
    (
      \seqname{A}^1,
      E^1,
      C^1,
      \seqname{A}^2,
      E^2,
      C^2,
      \funcname{f}_0,
      \funcname{f}_1
    )
  }}]*
\end{multline}
\begin{multline}
\maxfull{\Atl^\Ck_\Ar}(E^1, C^1, E^2, C^2) \defeq
  \set
  {{
    \bigl (
      (
        \funcname{f}_0,
        \funcname{f}_1
      ),
      (\seqname{A}^1, E^1, C^1),
      (\seqname{A}^2, E^2, C^2)
    \bigr )
  }}%
  [{{
    \maxfull{\isAtl^\Ck_\Ar}
    (
      \seqname{A}^1,
      E^1,
      C^1,
      \seqname{A}^2,
      E^2,
      C^2,
      \funcname{f}_0,
      \funcname{f}_1
    )
  }}]*
\end{multline}
\begin{multline}
\maxfull[S-]{\Atl^\Ck_\Ar}(E^1, C^1, E^2, C^2) \defeq
  \set
  {{
    \bigl (
      (
        \funcname{f}_0,
        \funcname{f}_1
      ),
      (\seqname{A}^1, E^1, C^1),
      (\seqname{A}^2, E^2, C^2)
    \bigr )
  }}%
  [{{
    \maxfull[S-]{\isAtl^\Ck_\Ar}
    (
      \seqname{A}^1,
      E^1,
      C^1,
      \seqname{A}^2,
      E^2,
      C^2,
      \funcname{f}_0,
      \funcname{f}_1
    )
  }}]*
\end{multline}
\end{definition}

\begin{lemma}[$\mathrm{C^k}$-atlas morphisms]
\label{lem:Ck-ATLmorph}
Let $\seqname{E}$ be a set of topological spaces, $\seqname{C}$ a set of
$\Ck$ linear spaces, $E^i \in \seqname{E}$, $i=1,2$, $C^i \in
\seqname{C}$ and $\seqname{A}^i$ $\Ck$-atlases of $E^i$ in the
coordinate spaces $C^i$. A pair of functions
$\funcseqname{f} \defeq (\funcname{f}_0, \funcname{f}_1)$ is an
$E^1$-$E^2$-$\Ck$-morphism of $\seqname{A}^1$ to $\seqname{A}^2$ in the
coordinate spaces $C^1$, $C^2$ iff $\funcseqname{f}$ is a strict
$\Trivcat{\seqname{E}}$-%
$E^1$-$E^2$-%
$\Trivcat[\Ck-]{\seqname{C}}$-%
$\Trivcat[\Ck-]{\seqname{C}}$
morphism of $\seqname{A}^1$ to $\seqname{A}^2$ in the coordinate spaces
$\Triv[\Ck-]{\seqname{C}^1}$, $\Triv[\Ck-]{\seqname{C}^2}$.

\begin{proof}
The model neighborhoods of $\Triv[\Ck-]{\seqname{C}^i}$ are the open sets
of $\seqname{C}^i$ and the morphisms of
$\Trivcat[\Ck-]{\seqname{C}}$ are the $\Ck$ functions between
spaces in $\seqname{C}$.
\end{proof}
\end{lemma}

\begin{corollary}[$\mathrm{C^k}$-atlas morphisms]
\label{cor:Ck-ATLmorph}
Let $E^i$, $i=1,2,3$, be topological spaces, $C^i$ linear spaces,
$\seqname{A}^i$ $\Ck$-atlases of $E^i$ in the coordinate spaces $C^i$ and
$(\funcname{f}^i_0, \funcname{f}^i_1)$ $E^i$-$E^{i+1}$
$\Ck$-morphisms of $\seqname{A}^i$ to $\seqname{A}^{i+1}$ in the
coordinate spaces $C^i$, $C^{i+1}$. Then
$
  (
    \funcname{f}^2_0 \compose \funcname{f}^1_0,
    \funcname{f}^2_1 \compose \funcname{f}^1_1
  )
$
is a $E^1$-$E^3$ $\Ck$-morphism of $\seqname{A}^1$ to $\seqname{A}^3$ in
the coordinate spaces $C^1$, $C^3$.
\begin{proof}
The result follows from \pagecref{lem:M-ATLmorph}.
\end{proof}
\end{corollary}

\begin{definition}[Categories of $\mathrm{\Ck}$ atlases]
\label{def:Ck-ATLcat}
Let $\seqname{E}$ be a set of topological spaces and $\seqname{C}$ a
set of linear spaces.
Let $P \defeq \seqname{E} \times \seqname{C}$. Then
\begin{equation}
\Atl^\Ck_\Ar(\seqname{E}, \seqname{C})
  \defeq
  \union
    [
        {(E^\mu, C^\mu) \in \seqname{P}},
        {(E^\nu, C^\nu) \in \seqname{P}}
    ]
    {{
      \Atl^\Ck_\Ar
      (
        E^\mu,
        C^\mu,
        E^\nu,
        C^\nu
      )
    }}
\end{equation}
\begin{equation}
  \Atl^\Ck(\seqname{E}, \seqname{C})
  \defeq
  \Bigl (
    \Atl^\Ck_\Ob(\seqname{E}, \seqname{C}),
    \Atl^\Ck_\Ar(\seqname{E}, \seqname{C}),
    \compose[A]
  \Bigr )
\end{equation}
\begin{equation}
\full{\Atl^\Ck_\Ar}(\seqname{E}, \seqname{C})
  \defeq
  \union
    [
        {(E^\mu, C^\mu) \in \seqname{P}},
        {(E^\nu, C^\nu) \in \seqname{P}}
    ]
    {{
      \full{\Atl^\Ck_\Ar}
      (
        E^\mu,
        C^\mu,
        E^\nu,
        C^\nu
      )
    }}
\end{equation}
\begin{equation}
  \full{\Atl^\Ck}(\seqname{E}, \seqname{C})
  \defeq
  \Bigl (
    \full{\Atl^\Ck_\Ob}(\seqname{E}, \seqname{C}),
    \full{\Atl^\Ck_\Ar}(\seqname{E}, \seqname{C}),
    \compose[A]
  \Bigr )
\end{equation}
\begin{equation}
\maximal{\Atl^\Ck_\Ar}(\seqname{E}, \seqname{C})
  \defeq
  \union
    [
        {(E^\mu, C^\mu) \in \seqname{P}},
        {(E^\nu, C^\nu) \in \seqname{P}}
    ]
    {{
      \maximal{\Atl^\Ck_\Ar}
      (
        E^\mu,
        C^\mu,
        E^\nu,
        C^\nu
      )
    }}
\end{equation}
\begin{equation}
  \maximal{\Atl^\Ck}(\seqname{E}, \seqname{C})
  \defeq
  \Bigl (
    \maximal{\Atl^\Ck_\Ob}(\seqname{E}, \seqname{C}),
    \maximal{\Atl^\Ck_\Ar}(\seqname{E}, \seqname{C}),
    \compose[A]
  \Bigr )
\end{equation}
\begin{equation}
\maximal[S-]{\Atl^\Ck_\Ar}(\seqname{E}, \seqname{C})
  \defeq
  \union
    [
        {(E^\mu, C^\mu) \in \seqname{P}},
        {(E^\nu, C^\nu) \in \seqname{P}}
    ]
    {{
      \maximal[S-]{\Atl^\Ck_\Ar}
      (
        E^\mu,
        C^\mu,
        E^\nu,
        C^\nu
      )
    }}
\end{equation}
\begin{equation}
  \maximal[S-]{\Atl^\Ck}(\seqname{E}, \seqname{C})
  \defeq
  \Bigl (
    \maximal[S-]{\Atl^\Ck_\Ob}(\seqname{E}, \seqname{C}),
    \maximal[S-]{\Atl^\Ck_\Ar}(\seqname{E}, \seqname{C}),
    \compose[A]
  \Bigr )
\end{equation}
\begin{equation}
\maxfull{\Atl^\Ck_\Ar}(\seqname{E}, \seqname{C})
  \defeq
  \union
    [
        {(E^\mu, C^\mu) \in \seqname{P}},
        {(E^\nu, C^\nu) \in \seqname{P}}
    ]
    {{
      \maxfull{\Atl^\Ck_\Ar}
      (
        E^\mu,
        C^\mu,
        E^\nu,
        C^\nu
      )
    }}
\end{equation}
\begin{equation}
  \maxfull{\Atl^\Ck}(\seqname{E}, \seqname{C})
  \defeq
  \Bigl (
    \maxfull{\Atl^\Ck_\Ob}(\seqname{E}, \seqname{C}),
    \maxfull{\Atl^\Ck_\Ar}(\seqname{E}, \seqname{C}),
    \compose[A]
  \Bigr )
\end{equation}
\begin{equation}
\maxfull[S-]{\Atl^\Ck_\Ar}(\seqname{E}, \seqname{C})
  \defeq
  \union
    [
        {(E^\mu, C^\mu) \in \seqname{P}},
        {(E^\nu, C^\nu) \in \seqname{P}}
    ]
    {{
      \maxfull[S-]{\Atl^\Ck_\Ar}
      (
        E^\mu,
        C^\mu,
        E^\nu,
        C^\nu
      )
    }}
\end{equation}
\begin{equation}
  \maxfull[S-]{\Atl^\Ck}(\seqname{E}, \seqname{C})
  \defeq
  \Bigl (
    \maxfull[S-]{\Atl^\Ck_\Ob}(\seqname{E}, \seqname{C}),
    \maxfull[S-]{\Atl^\Ck_\Ar}(\seqname{E}, \seqname{C}),
    \compose[A]
  \Bigr )
\end{equation}
\end{definition}

\begin{lemma}[$\Atl^\Ck(\seqname{E}, \seqname{C})$ is a category]
\label{lem:Ck-ATLiscat}
Let $\seqname{E}$ be a set of topological spaces and $\seqname{C}$ a set
of linear spaces. Then $\Atl^\Ck(\seqname{E}, \seqname{C})$,
$\full{\Atl^\Ck_\Ar}(\seqname{E}, \seqname{C})$ and
$\maximal{\Atl^\Ck}(\seqname{E}, \seqname{C})$ are categories.

Let
$
  (\seqname{A}^i, E^i, C^i)
  \in
  \Atl^\Ck_\Ob(\seqname{E}, \seqname{C})
$.
Then $\Id_{(\seqname{A}^i, E^i, C^i)}$ is the
identity morphism for $(\seqname{A}^i, E^i, C^i)$.

\begin{proof}
Let $(\seqname{A}^i,\seqname{E}^i,\seqname{C}^i)$, $i=1,2,3$
be objects of $\Atl^\Ck(\seqname{E}, \seqname{C})$ and
let
$
  \funcname{m}^i \defeq \\
  \bigl (
    (\funcname{f}_0^i,\funcname{f}_1^i),
    (\seqname{A}^i,\seqname{E}^i,\seqname{C}^i),
    (\seqname{A}^{i+1},\seqname{E}^{i+1},\seqname{C}^{i+1})
  \bigr )
$
be morphisms of $\Atl^\Ck(\seqname{E}, \seqname{C})$. \\
Then
\begin{enumerate}
\item Composition: \newline
$
  \bigl (
    (
      \funcname{f}_0^2 \compose \funcname{f}_0^1,
      \funcname{f}_1^2 \compose \funcname{f}_1^1
    ),
    (\seqname{A}^1,\seqname{E}^1,\seqname{C}^1),
    (\seqname{A}^3,\seqname{E}^3,\seqname{C}^3)
  \bigr )
$
is a morphism of $\Atl^\Ck(\seqname{E}, \seqname{C})$ by \pagecref{cor:Ck-ATLmorph}.
\item Associativity: \newline
Composition is associative by \pagecref{lem:atlcomp}.
\item Identity: \newline
$\Id_{(\seqname{A}^i, \seqname{E}^i, \seqname{C}^i)}$ is an identity
morphism by \cref{lem:atlcomp}.
\end{enumerate}

The proofs for $\full{\Atl^\Ck}(\seqname{E}, \seqname{C})$ and
$\maximal{\Atl^\Ck}(\seqname{E}, \seqname{C})$ are the same.
\end{proof}
\end{lemma}

\begin{definition}[Functors from $\mathrm{C^k}$ atlases to m-atlases]
\label{def:FuncCktom}
Let $E^i$, $i=1,2$, be topological spaces, $C^i$ linear spaces,
$\seqname{A}^i$ $\Ck$-atlases of $E^i$ in the coordinate space $C^i$,
$\funcname{f}_0 \maps E^1 \to E^2$ continuous and
$\funcname{f}_1 \maps C^1 \to C^2$ $\Ck$. Then
\begin{equation}
\Functor^\Ck_{\Ck,\M} (\seqname{A}^i, E^i, C^i)
\defeq
  \bigl (
     \seqname{A}^i,
     \Triv{E}^i,
     \Triv[\Ck-]{C^i}
   \bigr )
\end{equation}
\begin{multline}
\Functor^\Ck_{\Ck,\M}
  \bigl (
    (
      \funcname{f}_0,
      \funcname{f}_1
    ),
    (\seqname{A}^1, E^1, C^1),
    (\seqname{A}^2, E^2, C^2)
  \bigr )
\defeq \\
  \biggl (
    \Bigl (
      \funcname{f}_0,
      \funcname{f}_1
     \Bigr ) ,
    \Bigl ( \seqname{A}^1, \Triv{E^1}, \Triv[\Ck-]{C^1} \Bigr ) ,
    \Bigl ( \seqname{A}^2, \Triv{E^2}, \Triv[\Ck-]{C^2} \Bigr )
  \biggr )
\end{multline}
\end{definition}

\begin{theorem}[Functors from $\mathrm{C^k}$ atlases to m-atlases]
Let $\seqname{E}$ be a set of topological spaces and $\seqname{C}$ a
set of linear spaces. Then
$\Functor^\Ck_{\Ck,\M}$ is a functor from
$\Atl^\Ck(\seqname{E}, \seqname{C})$
to
$\Atl(\Triv{\seqname{E}}, \Triv[\Ck-]{\seqname{C}})$

\begin{proof}
Let $\seqname{o}^i \defeq (\seqname{A}^i, E^i, C^i)$,
$i \in [1,3]$, be objects of $\Atl^\Ck(\seqname{E}, \seqname{C})$ and \\
$
\seqname{m}^i \defeq
  \bigl (
    (\funcseqname{f}^i_0, \funcseqname{f}^i_1).
    \seqname{o}^i,
    \seqname{o}^{i+1}
  \bigr )
$, $i=1,2$,
be morphisms from $\seqname{o}^i$ to $\seqname{o}^{i+1}$.

$\Functor^\Ck_{\Ck,\M} (\seqname{m}^1)$ is a morphism from
$\Functor^\Ck_{\Ck,\M} \seqname{o}^1$ to $\Functor^\Ck_{\Ck,\M} \seqname{o}^2$:

\begin{equation}
\begin{split}
  &
\Functor^\Ck_{\Ck,\M} (\seqname{m}^1) =
\\*
  &
\Functor^\Ck_{\Ck,\M}
  \bigl (
    (
      \funcseqname{f}^1_0,
      \funcseqname{f}^1_1
    ),
    (\seqname{A}^1, E^1, C^1),
    (\seqname{A}^2, E^2, C^2)
  \bigr )
=
\\*
  &
  \bigl (
    (
      \funcseqname{f}^1_0,
      \funcseqname{f}^1_1
    ),
    (\seqname{A}^1, \Triv{E}^1, \Triv[\Ck-]{C^1}),
    (\seqname{A}^32 \Triv{E}^2, \Triv[\Ck-]{C^2})
  \bigr )
\end{split}
\end{equation}

$\Functor^\Ck_{\Ck,\M}$ maps identity functions to identity functions:

\begin{equation}
\begin{split}
  &
\Functor^\Ck_{\Ck,\M} \Id_{(\seqname{A}^i, E^i, C^i)}
=
\\*
  &
\Functor^\Ck_{\Ck,\M}
  \bigl (
    (\Id_{E^1}, \Id_{C^1}),
    (\seqname{A}^i, E^i, C^i),
    (\seqname{A}^i, E^i, C^i)
  \bigr )
=
\\*
  &
  \bigl (
    (\Id_{\Triv{E}^i}, \Id_{C^i}),
    (\seqname{A}^i, \Triv{E}^i, \Triv[\Ck-]{C^i}),
    (\seqname{A}^i, \Triv{E}^i, \Triv[\Ck-]{C^i})
  \bigr )
=
\\*
  &
  \bigl (
    (\Id_{\Triv{E}^i}, \Id_{C^i}),
    \Functor^\Ck_{\Ck,\M} (\seqname{A}^i, E^i, C^i),
    \Functor^\Ck_{\Ck,\M} (\seqname{A}^i, E^i, C^i)
  \bigr )
=
\\*
  &
\Id_{\Functor^\Ck_{\Ck,\M} (\seqname{A}^i, E^i, C^i)}
\end{split}
\end{equation}

$
  \Functor^\Ck_{\Ck,\M} (m^2) \compose[A] \Functor^\Ck_{\Ck,\M} (m^1)
  =
  \Functor^\Ck_{\Ck,\M} (m^2 \compose[A] m^1)
$:

\begin{enumerate}
\item
$
  m^2 \compose[A] m^1 =
  \bigl (
    (f_0^2 \compose f_0^1, f_1^2 \compose f_1^1),
    (\seqname{A}^1, E^1, C^1),
    (\seqname{A}^3, E^3, C^3)
  \bigr )
  $
\item
$
  \Functor^\Ck_{\Ck,\M}
    \bigl (
       (\seqname{A}^i, E^i, C^i)
    \bigr )
  =
  (\seqname{A}^i, \Triv{E}^i, \Triv[\Ck-]{C^i})
$
\item
$
  \Functor^\Ck_{\Ck,\M} (m^i)
    =
    \bigl (
      (f_0^i, f_1^i),
      (\seqname{A}^i, \Triv{E}^i, \Triv[\Ck-]{C^i}),
      (\seqname{A}^{i+1}, \Triv{E}^{i+1}, \Triv[\Ck-]{C^{i+1}})
    \bigr )
$
\item
$
  \Functor^\Ck_{\Ck,\M} (m^2) \compose[A] \Functor^\Ck_{\Ck,\M} (m^1)
  = \\
    \bigl (
      (f_0^2 \compose f_0^1, f_1^2 \compose f_1^2),
      (\seqname{A}^1, \Triv{E}^1, \Triv[\Ck-]{C^1}),
      (\seqname{A}^3, \Triv{E}^3, \Triv[\Ck-]{C^3})
    \bigr )
$
\item
$
  \Functor^\Ck_{\Ck,\M} (m^2 \compose m^1)
  = \\
    \bigl (
      (f_0^2 \compose f_0^1, f_1^2 \compose f_1^2),
      (\seqname{A}^1, \Triv{E}^1, \Triv[\Ck-]{C^1}),
      (\seqname{A}^3, \Triv{E}^3, \Triv[\Ck-]{C^3})
    \bigr )
$
\end{enumerate}
\end{proof}
\end{theorem}

\begin{definition}[Functors from m-atlases to $\mathrm{C^k}$ atlases]
Let $\seqname{E}^i$, $i=1,2$, be model spaces, $\seqname{C}^i$ be linear
model spaces, $\seqname{A}^i$ a maximal m-atlas of $\seqname{E}^i$ in
the coordinate space $\seqname{C}^i$ and $(\funcname{f}_0,
\seqname{f}_1)$ an $\seqname{E}^1$-$\seqname{E}^2$ m-atlas morphism of
$\seqname{A}^1$ to $\seqname{A}^2$ in the coordinate spaces
$\seqname{C}^1$, $\seqname{C}^2$, Then
\begin{equation}
\Functor^\Ck_{\M,\Ck} (\seqname{A}^i, \seqname{E}^i, \seqname{C}^i)
\defeq
  \bigl (
     \seqname{A}^i,
     \pi_1(\seqname{E}^i),
     \pi_1(\seqname{C}^i)
   \bigr )
\end{equation}
\begin{multline}
\Functor^\Ck_{\M,\Ck}
  \bigl (
    (
      \funcname{f}_0,
      \funcname{f}_1
    ),
    (\seqname{A}^1, \seqname{E}^1, \seqname{C}^1),
    (\seqname{A}^2, \seqname{E}^2, \seqname{C}^2)
  \bigr )
\defeq \\
  \Bigl (
    \bigl (
      \funcname{f}_0,
      \funcname{f}_1
    \bigr ),
    \bigl ( \seqname{A}^1, \pi_1(\seqname{E}^1), \pi_1(\seqname{C}^1) \bigr ),
    \bigl ( \seqname{A}^2, \pi_1(\seqname{E}^2), \pi_1(\seqname{C}^2) \bigr )
  \Bigr )
\end{multline}
\end{definition}

\begin{theorem}[Functors from m-atlases to $\mathrm{C^k}$ atlases]
Let $\seqname{E}$ be a set of model spaces and $\seqname{C}$ a
set of linear model spaces. Then
$\Functor^\Ck_{\M,\Ck}$ is a functor from
$\Atl(\seqname{E}, \seqname{C})$
to
$\Atl^\Ck \bigl (\pi_1[\seqname{E}], \pi_1[\seqname{C}] \bigr )$.

\begin{proof}
Let
$
\seqname{o}^i \defeq
\bigl (
  \seqname{A}^i, (E^i,\catname{E}^i), (C^i,\catname{C}^i)
\bigr )
$,
$i \in [1,3]$, be objects of $\Atl^\Ck(\seqname{E}, \seqname{C})$ and
$
\seqname{m}^i \defeq
  \bigl (
    (\funcseqname{f}^i_0, \funcseqname{f}^i_1).
    \seqname{o}^i,
    \seqname{o}^{i+1}
  \bigr )
$, $i=1,2$,
be morphisms from $\seqname{o}^i$ to $\seqname{o}^{i+1}$.

$\Functor^\Ck_{\M,\Ck} (\seqname{m}^1)$ is a morphism from
$\Functor^\Ck_{\M,\Ck} \seqname{o}^1$ to $\Functor^\Ck_{\M,\Ck} \seqname{o}^2$:

\begin{equation}
\begin{split}
  &
\Functor^\Ck_{\M,\Ck} (\seqname{m}^1) =
\\*
  &
\Functor^\Ck_{\M,\Ck}
  \bigl (
    (
      \funcseqname{f}^1_0,
      \funcseqname{f}^1_1
    ),
    (\seqname{A}^1, \seqname{E}^1, \seqname{C}^1),
    (\seqname{A}^2, \seqname{E}^2, \seqname{C}^2)
  \bigr )
=
\\*
  &
  \Bigl (
    (
      \funcseqname{f}^1_0,
      \funcseqname{f}^1_1
    ),
    \bigl (
      \seqname{A}^1, \pi_1(\seqname{E}^1), \pi_1(\seqname{C}^1)
    \bigr ),
    \bigl (
      \seqname{A}^2, \pi_1(\seqname{E}^2), \pi_1(\seqname{C}^2)
    \bigr )
  \Bigr )
\end{split}
\end{equation}

$\Functor^\Ck_{\M,\Ck}$ maps identity functions to identity functions:

\begin{equation}
\begin{split}
  &
\Functor^\Ck_{\M,\Ck} \Id_{(\seqname{A}^i, \seqname{E}^i, \seqname{C}^i)}
=
\\*
  &
\Functor^\Ck_{\M,\Ck}
  \bigl (
    (\Id_{\seqname{E}^i}, \Id_{\seqname{C}^i}),
    (\seqname{A}^i, \seqname{E}^i, \seqname{C}^i),
    (\seqname{A}^i, \seqname{E}^i, \seqname{C}^i)
  \bigr )
=
\\*
  &
  \bigl (
    (\Id_{\pi_1(\seqname{E}^i)}, \Id_{\pi_1(\seqname{C}^i)}),
    (\seqname{A}^i, \pi_1(\seqname{E}^i), \pi_1(\seqname{C}^i)
    ),
    (\seqname{A}^i, \pi_1(\seqname{E}^i), \pi_1(\seqname{C}^i)
    )
  \bigr )
=
\\*
  &
  \bigl (
    (\Id_{E^i}, \Id_{C^i}),
    \Functor^\Ck_{\M,\Ck} (\seqname{A}^i, \seqname{E}^i, \seqname{C}^i),
    \Functor^\Ck_{\M,\Ck} (\seqname{A}^i, \seqname{E}^i, \seqname{C}^i)
  \bigr )
=
\\*
  &
\Id_{\Functor^\Ck_{\M,\Ck} (\seqname{A}^i, \seqname{E}^i, \seqname{C}^i)}
\end{split}
\end{equation}

$
  \Functor^\Ck_{\M,\Ck} (m^2) \compose[A] \Functor^\Ck_{\M,\Ck} (m^1)
  =
  \Functor^\Ck_{\M,\Ck} (m^2 \compose[A] m^1)
$:

\begin{enumerate}
\item
$
  m^2 \compose[A] m^1 =
  \bigl (
    (f_0^2 \compose f_0^1, f_1^2 \compose f_1^1),
    (\seqname{A}^1, \seqname{E}^1, \seqname{C}^1),
    (\seqname{A}^3, \seqname{E}^3, \seqname{C}^3)
  \bigr )
$
\item
$
  \Functor^\Ck_{\M,\Ck}
    \bigl (
       (\seqname{A}^i, \seqname{E}^i, \seqname{C}^i)
    \bigr )
  =
  \bigl (
    \seqname{A}^i,
    \pi_1(\seqname{E}^i),
    \pi_1(\seqname{C}^i)
  \bigr )
$
\item
$
  \Functor^\Ck_{\M,\Ck} (m^i) = \\
    \Bigl (
      (f_0^i, f_1^i),
      \bigl (
        \seqname{A}^i,
        \pi_1(\seqname{E}^i),
        \pi_1(\seqname{C}^i)
      \bigr ),
      \bigl (
        \seqname{A}^{i+1},
        \pi_1(\seqname{E}^{i+1}),
        \pi_1(\seqname{C}^{i+1})
      \bigr )
    \Bigr )
$
\item
$
  \Functor^\Ck_{\M,\Ck} (m^2) \compose[A] \Functor^\Ck_{\M,\Ck} (m^1)
  = \\
    \Bigl (
      (f_0^2 \compose f_0^1, f_1^2 \compose f_1^2),
      \bigl (
        \seqname{A}^1,
        \pi_1(\seqname{E}^1),
        \pi_1(\seqname{C}^1)
      \bigr ),
      \bigl (
        \seqname{A}^3,
        \pi_1(\seqname{E}^3),
        \pi_1(\seqname{C}^3)
      \bigr )
    \Bigr )
$
\item
$
  \Functor^\Ck_{\M,\Ck} (m^2 \compose m^1) = \\
    \Bigl (
      (f_0^2 \compose f_0^1, f_1^2 \compose f_1^2),
      \bigl (
        \seqname{A}^1,
        \pi_1(\seqname{E}^1),
        \pi_1(\seqname{C}^1)
      \bigr ),
      \bigl (
        \seqname{A}^3,
        \pi_1(\seqname{E}^3),
        \pi_1(\seqname{C}^3)
      \bigr )
    \Bigr )
$
\end{enumerate}

\end{proof}
\end{theorem}

\subsection{Associated model spaces and functors}
\label{sub:assmodC}
\begin{definition}[Coordinate model spaces associated with $\mathrm{C^k}$-atlases]
\label{def:ismodC}
Let $\seqname{A}^i$, $i=1,2$, be a $\Ck$-atlas of $E^i$ in
the coordinate space $C^i$, $\funcname{f}_0 \maps E^1 \to E^2$ a
continuous function and $\funcname{f}_1 \maps C^1 \to C^2$ a $\Ck$
function. Then

\begin{multline}
\Functor^{\minimal{\Ck}}_2 (\seqname{A}^i, E^i, C^i)
\defeq \\
\minimal{\Mod}
  \left (
    C^i,
    \pi_2[\seqname{A}^i],
    \set
      {\phi' \compose \phi^{-1}}%
      [
        \equant
          {
            {(U,V,\phi) \in \seqname{A}^i},
            {(U',V',\phi') \in \seqname{A}^i}
          }
          {
            U \cap U' \ne \emptyset
          }
      ]
  \right )
\end{multline}
\begin{multline}
\Functor^{\minimal{\Ck}}_2
  \bigl (
    (
      \funcname{f}_0,
      \funcname{f}_1
    ),
    (\seqname{A}^1, E^1, C^1),
    (\seqname{A}^2, E^2, C^2)
  \bigr )
\defeq \\
  \funcname{f}_1
  \maps
  \Functor^{\minimal{\Ck}}_2 (\seqname{A}^1, E^1, C^1)
  \to
  \Functor^{\minimal{\Ck}}_2 (\seqname{A}^2, E^2, C^2)
\end{multline}

The minimal coordinate $\Ck$ model space with neighborhoods in the
$\Ck$-atlas $\seqname{A}^i$ of $E^i$ in the coordinate space $C^i$ is
$\Functor^{\minimal{\Ck}}_2 (\seqname{A}^1, E^i, C^i)$.

The coordinate mapping associated with the
$\seqname{E}^1$-$\seqname{E}^2$ $\Ck$-atlas morphism
$(\funcname{f}_0,\funcname{f}_1)$ of $\seqname{A}^1$ to $\seqname{A}^2$
in the coordinate spaces $\seqname{C}^1$, $\seqname{C}^2$ is
$
  \funcname{f}_1
  \maps
  \Functor^{\minimal{\Ck}}_2 (\seqname{A}^1, E^1, C^1)
  \to
  \Functor^{\minimal{\Ck}}_2 (\seqname{A}^2, E^2, C^2)
$.
If it is a model function then the it is also the
coordinate m-atlas morphism associated with the
$\seqname{E}^1$-$\seqname{E}^2$ $\Ck$-atlas morphism
$(\funcname{f}_0,\funcname{f}_1)$ of $\seqname{A}^1$ to $\seqname{A}^2$
in the coordinate spaces $\seqname{C}^1$, $\seqname{C}^2$.
\end{definition}

\begin{lemma}[Coordinate model spaces associated with $\mathrm{C^k}$-atlases]
\label{lem:ismodC}
Let $\seqname{A}$ be a $\Ck$-atlas of $E$ in the coordinate space $C$.
Then $\Functor^{\minimal{\Ck}}_2 (\seqname{A}, E, C)$ is a $\Ck$ linear
model space iff $\union{\pi_2[\seqname{A}]}$ contains a ball.

\begin{proof}
$\Functor^{\minimal{\Ck}}_2 (\seqname{A}, E, C)$ satisfies the
conditions for a model space.
\begin{enumerate}
\item Since $\pi_2[\seqname{A}]$ is an open cover of
$\union{\pi_2[\seqname{A}]}$, the set of finite intersections is also an
open cover.

\item Finite intersections of finite intersections are finite intersections

\item Restrictions of continuous function are continuous

\item If $f \maps A \to B$ is a morphism of
$\Functor^{\minimal{\Ck}}_2 (\seqname{A}, E, C)$,
$A, A', B, B'$ model meighborhoods of
$\Functor^{\minimal{\Ck}}_2 (\seqname{A}, E, C)$,
$A' \subseteq A$, $B' \subseteq B$ and $f[A'] \subseteq B'$ then since
$f \maps A \to B$ is a morphism it is a restriction of a transition
function between its restrictions to sets in $\pi_2[\seqname{A}]$ and
its restrictions are also, hence morphisms, and thus
$f \restriction_{A'} A' \to B'$ is a morphism.

\item
If $(U,V,\phi) \in \seqname{A}$ then $\Id_V = \phi \compose \phi^{-1}A$
is a transition function and hence a morphism of
$\Functor^{\minimal{\Ck}}_2 (\seqname{A}, E, C)$. If $A, A'$ objects of
$\bigl (\Functor^{\minimal{\Ck}}_2 (\seqname{A}, E, C) \bigr )$ and
$A' \subseteq A$ then the inclusion map
$\funcname{i} \maps A' \hookrightarrow A$ is a restriction of an
identity morphism of $\Functor^{\minimal{\Ck}}_2 (\seqname{A}, E, C)$
and hence a morphism.

\item Restricted sheaf condition: let
\begin{enumerate}
\item $U_\alpha,V_\alpha$, $\alpha \prec \Alpha$, be objects of
$\pi_2 \bigl ( \Functor^{\minimal{\Ck}}_2 (\seqname{A}, E, C) \bigr )$
\item $\funcname{f}_\alpha \maps U_\alpha \to V_\alpha$ be morphisms of
$\pi_2 \bigl ( \Functor^{\minimal{\Ck}}_2 (\seqname{A}, E, C) \bigr )$
\item $U \defeq \union[\alpha \prec \Alpha]{U_\alpha}$
\item $V \defeq \union[\alpha \prec \Alpha]{V_\alpha}$
\item $\funcname{f} \maps U \to V$ be continuous and
$
  \uquant
    {{\alpha \prec \Alpha},{x \in U_\alpha}}
    {\funcname{f}(x) = \funcname{f}_\alpha(x)}
$
\end{enumerate}
Then $\funcname{f}$ is $\Ck$ and hence a morphism of
$\pi_2 \bigl ( \Functor^{\minimal{\Ck}}_2 (\seqname{A}, E, C) \bigr )$.
\end{enumerate}

If $\union{\pi_2[\seqname{A}]}$ is a linear space, then by
\pagecref{def:lin}\!, $\union{\pi_2[\seqname{A}]}$ it contains a ball.
Conversly, $\union{\pi_2[\seqname{A}]}$ is locally connected so
if $\union{\pi_2[\seqname{A}]}$ contains a ball then
the conditions of \pagecref{def:lin} are met.
\end{proof}

Let $\seqname{A}^i$, $i=1,2$, be a semi-maximal $\Ck$-atlas of $E^i$ in
the coordinate space $C^i$, $\funcname{f}_0 \maps E^1 \to E^2$ a
continuous function, $\funcname{f}_1 \maps C^1 \to C^2$ a $\Ck$
function and $(\funcname{f}_0, \funcname{f}_1)$ a $\Ck$-atlas morphism
from  $\seqname{A}^1$ to $\seqname{A}^2$. Then
$
  \funcname{f}_1
  \maps
  \Functor^{\minimal{\Ck}}_2 (\seqname{A}^1, E^1, C^1)
  \to
  \Functor^{\minimal{\Ck}}_2 (\seqname{A}^2, E^2, C^2)
$
is well defined.

\begin{proof}
Let $v^1 \in \Functor^{\minimal{\Ck}}_2 (\seqname{A}^1, E^1, C^1)$,
$(U^i,V^i,\phi^i) \in \seqname{A}^i$, $i=1,2$, be a chart with
$u^1 \in U^1$,  $\phi^1(u^1) = v^1$ and $\funcname{f}_0(u^1) \in U^2$.
Then there are open sets
$U'^1 \subseteq I \defeq U^1 \cap \funcname{f}_0^{-1}[U^2]$,
$V'^1 \subseteq V^1$, $U'^2 \subseteq U^2$, $V'^2 \subseteq V^2$,
$\hat{V'^2} \subseteq C^2$ and a $\Ck$ diffeomorphism
$\hat{\funcname{f}} \maps \hat{V'^2} \toiso V'^2$ such that
{
  \showlabelsinline
  \crefrange{eq:xin}{eq:fphigeqgphi} on
}
\cpagerefrange{eq:xin}{eq:fphigeqgphi} in \pagecref{def:M-ATLmorph}
hold with $x \defeq u^1$. Since $\seqname{A}^2$ is semi-maximal,
$(U'^2, \hat{V'^2}, \hat{\funcname{f}} \compose \phi^2)$ is a chart of
$\seqname{A}^2$ and by \cref{eq:fphixin}
$\funcname{f}_1(v^1) \in \hat{V'^2}$.
\end{proof}
\end{lemma}

\begin{definition}[Model spaces associated with $\mathrm{C^k}$-atlases]
Let $\seqname{A}^i$, $i=1,2$ be $\Ck$-atlases of $E^i$ in the coordinate
spaces $C^i$, $\funcname{f}_0 \maps E^1 \to E^2$ a continuous function
and $\funcname{f}_1 \maps C^1 \to C^2$ a $\Ck$ function. Then

\begin{multline}
\Functor^{\minimal{\Ck}}_1 (\seqname{A}^i, E^i, C^i)
\defeq \\
\minimal{\Mod}
  \left (
    E^i,
    \pi_1[\seqname{A}^i],
    \set
      {\phi'^{-1} \compose \phi}%
      [
        \equant
          {
            {(U,V,\phi) \in \seqname{A}^i},
            {(U',V',\phi') \in \seqname{A}^i}
          }
          {
            U \cap U' \ne \emptyset
          }
      ]
  \right )
\end{multline}
\begin{multline}
\Functor^{\minimal{\Ck}}_1
  \bigl (
    (
      \funcname{f}_0,
      \funcname{f}_1
    ),
    (\seqname{A}^1, E^1, C^1),
    (\seqname{A}^2, E^2, C^2)
  \bigr )
\defeq \\
  \funcname{f}_0
  \maps
  \Functor^{\minimal{\Ck}}_1 (\seqname{A}^1, E^1, C^1)
  \to
  \Functor^{\minimal{\Ck}}_1 (\seqname{A}^2, E^2, C^2)
\end{multline}
The minimal $\Ck$ model space with neighborhoods in the $\Ck$-atlas
$\seqname{A}^i$ of $E^i$ in the coordinate space $C^i$ is
$\Functor^{\minimal{\Ck}}_1 (\seqname{A}^1, E^i, C^i)$.

The mapping associated with the $\seqname{E}^1$-$\seqname{E}^2$
$\Ck$-atlas morphism $(\funcname{f}_0,\funcname{f}_1)$ of
$\seqname{A}^1$ to $\seqname{A}^2$ in the coordinate spaces
$\seqname{C}^1$, $\seqname{C}^2$ is
$
  \funcname{f}_0
  \maps
  \Functor^{\minimal{\Ck}}_1 (\seqname{A}^1, E^1, C^1)
  \to
  \Functor^{\minimal{\Ck}}_1 (\seqname{A}^2, E^2, C^2)
$.
If it is a model function then the it is also the m-atlas morphism
associated with the $\seqname{E}^1$-$\seqname{E}^2$ $\Ck$-atlas morphism
$(\funcname{f}_0,\funcname{f}_1)$ of $\seqname{A}^1$ to $\seqname{A}^2$
in the coordinate spaces $\seqname{C}^1$, $\seqname{C}^2$.
\end{definition}

\begin{lemma}[Model spaces associated with $\mathrm{C^k}$-atlases]
\label{lem:ismodE}
Let $\seqname{A}$ be a $\Ck$-atlas of $E$ in the coordinate
space $C$. Then $\Functor^{\minimal{\Ck}}_1 (\seqname{A}, E, C)$ is a
model space.

\begin{proof}
\Pagecref{lem:minmod}
\end{proof}
\end{lemma}

\begin{theorem}[Functors from $\mathrm{C^k}$ atlases to model spaces]
Let $\seqname{E}$ be a set of topological spaces and $\seqname{C}$ a set
of linear spaces. Then $\Functor^{\minimal{\Ck}}_1$ is a functor from
$\Atl^\Ck(\seqname{E}, \seqname{C})$ to $\Triv{\seqname{E}}$,
$\Functor^{\minimal{\Ck}}_1$ is a functor from
$\full{\Atl^\Ck}(\seqname{E}, \seqname{C})$ to $\Triv{\seqname{E}}$,
$\Functor^{\minimal{\Ck}}_2$ is a functor from
$\maximal[S-]{\Atl^\Ck}(\seqname{E}, \seqname{C})$ to
$\optriv[\Ck-]{\seqname{C}}$ and $\Functor^{\minimal{\Ck}}_2$ is a
functor from $\full{\Atl^\Ck}(\seqname{E}, \seqname{C})$ to
$\Triv[\Ck-]{\seqname{C}}$.

\begin{proof}
Let $\seqname{o}^i \defeq (\seqname{A}^i, E^i, C^i)$, $i \in [1,3]$, be
objects in $\Atl^\Ck(\seqname{E}, \seqname{C})$
and let \\
$
m^i \defeq
  \bigl (
    (\funcname{f}^i_0, \funcname{f}^i_1),
    o^i,
    o^{i+1}
  \bigr )
$
be morphisms in $\Atl^\Ck(\seqname{E}, \seqname{C})$.

$
  \Functor^{\minimal{\Ck}}_1 \maps
  \Atl^\Ck(\seqname{E}, \seqname{C}) \to
  \Triv{\seqname{E}}
$:
\begin{enumerate}
\item $\Functor (f \maps A \to B) \maps \Functor(A) \to \Functor (B)$:
\newline
$
\Functor^{\minimal{\Ck}}_1(\seqname{m}^i) =
\funcname{f}^i_0 \maps \Functor^{\minimal{\Ck}}_1(\seqname{o}^i) \to
\Functor^{\minimal{\Ck}}_1(\seqname{o}^{i+1})
$
\item $\Functor (g \compose f) = \Functor (g) \compose \Functor (f)$: \newline
$
\Functor^{\minimal{\Ck}}_1(\funcseqname{m}^2 \compose[A] \funcseqname{m}^1)
=
\newline
\Functor^{\minimal{\Ck}}_1
  \bigl (
    (
      \funcname{f}^2_0 \compose \funcname{f}^1_0,
      \funcname{f}^2_1 \compose \funcname{f}^1_1
    )
    (\seqname{A}^1, E^1, C^1),
    (\seqname{A}^3, E^3, C^3)
  \bigr )
=
\\
\funcname{f}^2_0 \compose \funcname{f}^1_0 \maps
  \Functor^{\minimal{\Ck}}_1(\seqname{o}^1) \to
  \Functor^{\minimal{\Ck}}_1(\seqname{o}^3)
=
\\
\bigl (
  \funcname{f}^2_0 \maps
  \Functor^{\minimal{\Ck}}_1(\seqname{o}^2) \to
  \Functor^{\minimal{\Ck}}_1(\seqname{o}^3)
  \bigr )
\compose
\bigl (
  \funcname{f}^1_0 \maps
  \Functor^{\minimal{\Ck}}_1(\seqname{o}^1) \to
  \Functor^{\minimal{\Ck}}_1(\seqname{o}^2)
\bigr )
=
\\
\Functor^{\minimal{\Ck}}_1
  \bigl (
    (
      \funcname{f}^2_0,
      \funcname{f}^2_1
    )
    (\seqname{A}^2, E^2, C^2),
    (\seqname{A}^3, E^3, C^3)
  \bigr )
\compose \\
\Functor^{\minimal{\Ck}}_1
  \bigl (
    (
      \funcname{f}^1_0,
      \funcname{f}^1_1
    )
    (\seqname{A}^1, E^1, C^1),
    (\seqname{A}^2, E^2, C^2)
  \bigr )
=
\\
\Functor^{\minimal{\Ck}}_1(\funcseqname{m}^2)
\compose
\Functor^{\minimal{\Ck}}_1(\funcseqname{m}^1)
$
\item $\Functor (\Id_A) = \Id_{\Functor (A)}$:
\begin{enumerate}
\item
$
\Functor^{\minimal{\Ck}}_1(\Id_{o^i}) =
\Functor^{\minimal{\Ck}}_1
  \bigl (
    (\Id_{E^i}, \Id_{C^i}),
    (\seqname{A}^i, E^i, C^i),
    (\seqname{A}^i, E^i, C^i)
  \bigr )
=
\\
\Id_{E^i} \maps \Functor^{\minimal{\Ck}}_1(\seqname{o}^i) \to
\Functor^{\minimal{\Ck}}_1(\seqname{o}^i)
$
\item
$
\Id_{\Functor^{\minimal{\Ck}}_1}(o^i) =
\Id_{E^i} \maps \Functor^{\minimal{\Ck}}_1(\seqname{o}^i) \to
\Functor^{\minimal{\Ck}}_1(\seqname{o}^i)
$
\end{enumerate}
\end{enumerate}

The proof for
$
\Functor^{\minimal{\Ck}}_1 \maps
\full{\Atl^\Ck}(\seqname{E}, \seqname{C}) \to \Triv{\seqname{E}}
$
is identical.

$
  \Functor^{\minimal{\Ck}}_2 \maps
  \Atl^\Ck(\seqname{E}, \seqname{C}) \to
  \optriv[\Ck-]{\seqname{C}}
$:

\begin{enumerate}
\item $\Functor (f \maps A \to B) \maps \Functor(A) \to \Functor (B)$:
\\
$
\Functor^{\minimal{\Ck}}_2(\seqname{m}^i) =
\funcname{f}^i_1 \maps \Functor^{\minimal{\Ck}}_2(\seqname{o}^i) \to
\Functor^{\minimal{\Ck}}_2(\seqname{o}^{i+1})
$
\item $\Functor (g \compose f) = \Functor (g) \compose \Functor (f)$:
$
\Functor^{\minimal{\Ck}}_2(\funcseqname{m}^2 \compose[A] \funcseqname{m}^1)
=
\\
\Functor^{\minimal{\Ck}}_2
  \bigl (
    (
      \funcname{f}^2_0 \compose \funcname{f}^1_0,
      \funcname{f}^2_1 \compose \funcname{f}^1_1
    )
    (\seqname{A}^1, E^1, C^1),
    (\seqname{A}^3, E^3, C^3)
  \bigr )
=
\\
\funcname{f}^2_0 \compose \funcname{f}^1_0 \maps
  \Functor^{\minimal{\Ck}}_2(\seqname{o}^1) \to
  \Functor^{\minimal{\Ck}}_2(\seqname{o}^3)
=
\\
\bigl (
  \funcname{f}^2_1 \maps
  \Functor^{\minimal{\Ck}}_2(\seqname{o}^2) \to
  \Functor^{\minimal{\Ck}}_2(\seqname{o}^3)
  \bigr )
\compose
\bigl (
  \funcname{f}^1_1 \maps
  \Functor^{\minimal{\Ck}}_2(\seqname{o}^1) \to
  \Functor^{\minimal{\Ck}}_2(\seqname{o}^2)
\bigr )
=
\\
\Functor^{\minimal{\Ck}}_2
  \bigl (
    (
      \funcname{f}^2_0,
      \funcname{f}^2_1
    )
    (\seqname{A}^2, E^2, C^2),
    (\seqname{A}^3, E^3, C^3)
  \bigr )
\compose \\
\Functor^{\minimal{\Ck}}_2
  \bigl (
    (
      \funcname{f}^1_0,
      \funcname{f}^1_1
    )
    (\seqname{A}^1, E^1, C^1),
    (\seqname{A}^2, E^2, C^2)
  \bigr )
=
\\
\Functor^{\minimal{\Ck}}_2(\funcseqname{m}^2)
\compose
\Functor^{\minimal{\Ck}}_2(\funcseqname{m}^1)
$
\item $\Functor (\Id_A) = \Id_{\Functor (A)}$:
\begin{enumerate}
\item
$
\Functor^{\minimal{\Ck}}_2(\Id_{\seqname{o}^i}) =
\Functor^{\minimal{\Ck}}_2
  \bigl (
    (\Id_{E^i}, \Id_{C^i}),
    (\seqname{A}^i, E^i, C^i),
    (\seqname{A}^i, E^i, C^i)
  \bigr )
=
\\
\Id_{C^i} \maps \Functor^{\minimal{\Ck}}_2(\seqname{o}^i) \to
\Functor^{\minimal{\Ck}}_2(\seqname{o}^i)
$
\item
$
\Id_{\Functor^{\minimal{\Ck}}_2}(\seqname{o}^i) =
\Id_{C^i} \maps \Functor^{\minimal{\Ck}}_2(\seqname{o}^i) \to
\Functor^{\minimal{\Ck}}_2(\seqname{o}^i)
$
\end{enumerate}
\end{enumerate}

The proof for
$
  \Functor^{\minimal{\Ck}}_2 \maps
  \full{\Atl^\Ck}(\seqname{E}, \seqname{C}) \to
  \Triv[\Ck-]{\seqname{C}}
$
is identical.
\end{proof}
\end{theorem}

\subsection {$\Ck$ manifolds}
\label{sub:ckman}
Conventionally a manifold is different from its atlases, but
$\maximal{\Atl^\Ck}(\seqname{E}, \seqname{C})$ in \pagecref{def:Ck-ATLcat}
encourages treating them on an equal footing. All of the results for
maximal atlases carry directly over to results for manifolds.

\begin{remark}
The variations of the $\LCS^\Ck$ constructors in
\pagecref{def:CktoLCS} are very similar, as are the variations of
the $\Functor^\Ck_{\Man,\LCS}$ functors; likewise the variations of the
$\Functor^\Ck_{\LCS,\Man}$ functors in  \pagecref{def:LCStoCk} are very
similar.

The corresponding results and proofs in \pagecref{the:CktoLCS}
are likewise very similar, as are the corresponding results and proofs
in \pagecref{the:LCStoCk}.
\end{remark}

\begin{definition}[$\mathrm{C^k}$ manifolds]
\label{def:man}
Let $E$ be a topological space, $C$ a linear space and $\seqname{A}$ a
maximal\footnote{
Requiring that the atlas be full would eliminate some pathologies.}
$\Ck$-atlas of $E$ in the coordinate space $C$. Then $(E, C,
\seqname{A})$ is a $\Ck$ manifold.

Let $\seqname{E}$ be a set of topological spaces and $\seqname{C}$ be a
set of linear spaces. Then
\begin{equation}
\Man^\Ck_\Ob(\seqname{E}, \seqname{C}) \defeq
\maximal{\Atl^\Ck}(\seqname{E}, \seqname{C})
\end{equation}
\begin{remark}
The manifold $(E, C, \seqname{A})$ corresponds to the object
$(\seqname{A}, E, C)$.
\end{remark}
\end{definition}

\begin{definition}[$\mathrm{C^k}$ manifold morphisms]
\label{def:manmorph}
Let $S^i \defeq (E^i,C^i,\seqname{A}^i)$, $i=1,2$, be $\Ck$ manifolds
and $(\funcname{f}_0 \maps E^1 \to E^2, \funcname{f}_1 \maps C^1 \to C^2)$
be an $E^1$-$E^2$ $\Ck$-morphism of $\seqname{A}^1$ to $\seqname{A}^2$
in the coordinate spaces $C^1$, $C^2$. Then
$(\funcname{f}_0, \funcname{f}_1)$ is a $\Ck$ morphism of $S^1$ to
$S^2$.

Let $E^i$, $i=1,2$ be topological spaces and $C^i$ be linear spaces. Then
\begin{multline}
\Man^\Ck_\Ar(E^1, C^1, E^2, C^2) \defeq
\Atl^\Ck_\Ar(E^1, C^1, E^2, C^2)
\end{multline}

Let $\seqname{E}$ be a set of topological spaces and $\seqname{C}$ a
set of linear spaces. Then
\begin{equation}
\Man^\Ck_\Ar(\seqname{E}, \seqname{C}) \defeq
\maximal{\Atl^\Ck_\Ar}(\seqname{E}, \seqname{C})
\end{equation}
\begin{equation}
\Man^\Ck(\seqname{E}, \seqname{C}) \defeq
\maximal{\Atl^\Ck}(\seqname{E}, \seqname{C})
\end{equation}
\end{definition}

\begin{theorem}[Categories of $\mathrm{C^k}$ manifolds]
Let $\seqname{E}$ be a set of topological spaces and $\seqname{C}$ a set
of linear spaces. Then $\Man^\Ck(\seqname{E}, \seqname{C})$ is a
category and the identity morphism of $(\seqname{A}, E, C)$ is an
identity morphism.

\begin{proof}
{
  \showlabelsinline
  The result follows directly from \cref{def:man,def:manmorph} above and
}
\pagecref{lem:Ck-ATLiscat}.
\end{proof}
\end{theorem}

\begin{definition}[Functors from $\mathrm{C^k}$ manifolds to Local 2-$\emptyset$ coordinate spaces]
\label{def:CktoLCS}
Let $\seqname{E}$ be a set of topological spaces,
$\seqname{C}$ be a set of $\Ck$ linear spaces,
$\catname{E}$ and $\hat{\catname{E}}$ model categories,
$\catname{C}$ and $\hat{\catname{C}}$ $\Ck$ linear model categories,
$\catseqname{M} \defeq (\catname{E}, \catname{C})$,
$
  \hat{\catseqname{M}} \defeq
  (\hat{\catname{E}}, \hat{\catname{C}})
$,
$
  \Triv{\catseqname{M}} \defeq
  \Bigl ( \Trivcat{\seqname{E}}, \Trivcat[\Ck-]{\seqname{C}} \Bigr )
$.
$
  \optriv{\catseqname{M}} \defeq
  \Bigl ( \Trivcat{\seqname{E}}, \optriv[\Ck-]{\seqname{C}} \Bigr )
$.
$
  \Triv{\catseqname{M}} \SUBCAT[full-] \catseqname{M}
  \SUBCAT[full-] \hat{\catseqname{M}}
$,
$\optriv{\catseqname{M}} \SUBCAT \hat{\catseqname{M}}$.

Let $E^i \in \seqname{E}$, $i=1,2$,
$C^i \in \seqname{C}$,
$\seqname{M}^i \defeq \left ( \Triv{E^i}, \Triv[\Ck-]{C^i} \right )$, \\
$
  \Singcat{\catseqname{M}^i} \defeq
  \left (
    \singcat{\seqname{M}^i_0},
    \singcat[\Ck-]{\seqname{M}^i_1}
  \right )
$,
$\seqname{A}^i$ a maximal $\Ck$-atlas of $E^i$ in $C^i$,
$
  \Functor^{\minimal{\Ck}}_1 (\seqname{A}^i, E^i, C^i) \objin
  \hat{\catname{E}}
$,
$
  \Functor^{\minimal{\Ck}}_2 (\seqname{A}^i, E^i, C^i) \objin
  \hat{\catname{C}}
$,
$(E^i, C^i, \seqname{A}^i)$ a $\Ck$ manifold,
$\seqname{S}^i \defeq (\seqname{A}^i, E^i, C^i)$,
$
  \minimal{\seqname{M}^i} \defeq \\
    \Bigl (
      \Functor^{\minimal{\Ck}}_1 (\seqname{A}^i, E^i, C^i),
      \Functor^{\minimal{\Ck}}_2 (\seqname{A}^i, E^i, C^i)
    \Bigr )
$, \\
$
  \Singcat{\Triv{\catseqname{M}^i}} \defeq
  (
    \singcat{\Triv{\seqname{M}^i_0}},
    \singcat[\Ck-]{\Triv{\seqname{M}^i_1}}
  )
$,
$\funcname{f}_0 \maps E^1 \to E^2$ continuous and
$\funcname{f}_1 \maps C^1 \to C^2$ $\Ck$.
Then
\begin{equation}
\Functor^\Ck_{\Man,\LCS} \seqname{S}^i
\defeq
  \left (
    \Singcat{\catseqname{M}^i},
    \seqname{M}^i,
    \seqname{A}^i,
    \emptyset,
    \emptyset
  \right )
\end{equation}
\begin{equation}
\Functor^{\Ck,\catseqname{M}}_{\Man,\LCS} \seqname{S}^i
\defeq
  \left (
    \catseqname{M},
    \seqname{M}^i,
    \seqname{A}^i,
    \emptyset,
    \emptyset
  \right )
\end{equation}
\begin{equation}
\Functor^{\minimal{\Ck}}_{\Man,\LCS} \seqname{S}^i
\defeq
  \left (
    \Singcat{\Triv{\catseqname{M}^i}},
    \Triv{\seqname{M}^i},
    \seqname{A}^i,
    \emptyset,
    \emptyset
  \right )
\end{equation}
\begin{equation}
\Functor^{\minimal{\Ck},\hat{\catseqname{M}}}_{\Man,\LCS} \seqname{S}^i
\defeq
  \left (
    \hat{\catseqname{M}},
    \minimal{\seqname{M}^i},
    \seqname{A}^i,
    \emptyset,
    \emptyset
  \right )
\end{equation}
\begin{multline}
\Functor^\Ck_{\Man,\LCS}
  \bigl (
    (
      \funcname{f}_0,
      \funcname{f}_1
    ),
    \seqname{S}^1,
    \seqname{S}^2
  \bigr )
\defeq
  \Bigl (
    (
      \funcname{f}_0,
      \funcname{f}_1
    ),
    \Functor^\Ck_{\Man,\LCS} \seqname{S}^1,
    \Functor^\Ck_{\Man,\LCS} \seqname{S}^2
  \Bigr )
\end{multline}
\begin{multline}
\Functor^{\Ck,\catseqname{M}}_{\Man,\LCS}
  \bigl (
    (
      \funcname{f}_0,
      \funcname{f}_1
    ),
    \seqname{S}^1,
    \seqname{S}^2
  \bigr )
\defeq
  \Bigl (
    (
      \funcname{f}_0,
      \funcname{f}_1
    ),
    \Functor^{\Ck,\catseqname{M}}_{\Man,\LCS} \seqname{S}^1,
    \Functor^{\Ck,\catseqname{M}}_{\Man,\LCS} \seqname{S}^2
  \Bigr )
\end{multline}
\begin{multline}
\Functor^{\minimal{\Ck}}_{\Man,\LCS}
  \bigl (
    (
      \funcname{f}_0,
      \funcname{f}_1
    ),
    \seqname{S}^1,
    \seqname{S}^2
  \bigr )
\defeq
  \Biggl (
    \bigl (
      \funcname{f}_0,
      \funcname{f}_1
    \bigr ),
    \Functor^{\minimal{\Ck}}_{\Man,\LCS} \seqname{S}^1,
    \Functor^{\minimal{\Ck}}_{\Man,\LCS} \seqname{S}^2
  \Biggr )
\end{multline}
\begin{multline}
\Functor^{\minimal{\Ck},\catseqname{M}}_{\Man,\LCS}
  \bigl (
    (
      \funcname{f}_0,
      \funcname{f}_1
    ),
    \seqname{S}^1,
    \seqname{S}^2
  \bigr )
\defeq
  \Biggl (
    \bigl (
      \funcname{f}_0,
      \funcname{f}_1
    \bigr ),
    \Functor^{\minimal{\Ck},\catseqname{M}}_{\Man,\LCS} \seqname{S}^1,
    \Functor^{\minimal{\Ck},\catseqname{M}}_{\Man,\LCS} \seqname{S}^2
  \Biggr )
\end{multline}

\begin{multline}
\LCS^\Ck_\Ob(\seqname{E}, \seqname{C}) \defeq
\set
  {
    \Functor^\Ck_{\Man,\LCS} (\seqname{A}, E, C)
  }%
  [
    E \in \seqname{E},
    C \in \seqname{C},
    {\maximal{\isAtl^\Ck_\Ob}(\seqname{A}, E, C)}
  ]*
\end{multline}
\begin{multline}
\LCS^\Ck_\Ar(\seqname{E}, \seqname{C}) \defeq
\set
  {
    \Functor^\Ck_{\Man,\LCS}
      \bigl (
        (
          \funcname{f}_0,
          \funcname{f}_1
        ),
        (
          \seqname{A}^1,
          E^1,
          C^1
        ),
        (
          \seqname{A}^2,
          E^2,
          C^2
        )
      \bigr )
  }%
  [
    E^i \in \seqname{E},
    C^i \in \seqname{C},
    {
      \maximal{\isAtl^\Ck_\Ar}
        (
          \seqname{A}^1, E^1, C^1,
          \seqname{A}^2, E^2, C^2,
          \funcname{f}_0,\funcname{f}_1
        )
    }
  ]*
\end{multline}
\begin{equation}
\LCS^\Ck(\seqname{E}, \seqname{C}) \defeq
\Bigl (
  \LCS^\Ck_\Ob(\seqname{E}, \seqname{C}),
  \LCS^\Ck_\Ar(\seqname{E}, \seqname{C}),
  \compose[A]
\Bigr )
\end{equation}
\begin{multline}
\LCS^{\Ck,\catseqname{M}}_\Ob(\seqname{E}, \seqname{C}) \defeq
\set
  {
    \Functor^{\Ck,\catseqname{M}}_{\Man,\LCS} (\seqname{A}, E, C)
  }%
  [
    E \in \seqname{E},
    C \in \seqname{C},
    {\maximal{\isAtl^\Ck_\Ob}(\seqname{A}, E, C)}
  ]*
\end{multline}
\begin{multline}
\LCS^{\Ck,\catseqname{M}}_\Ar(\seqname{E}, \seqname{C}) \defeq
\set
  {
    \Functor^{\Ck,\catseqname{M}}_{\Man,\LCS}
      \bigl (
        (
          \funcname{f}_0,
          \funcname{f}_1
        ),
        (
          \seqname{A}^1,
          E^1,
          C^1
        ),
        (
          \seqname{A}^2,
          E^2,
          C^2
        )
      \bigr )
  }%
  [
    E^i \in \seqname{E},
    C^i \in \seqname{C},
    {
      \maximal{\isAtl^\Ck_\Ar}
        (
          \seqname{A}^1, E^1, C^1,
          \seqname{A}^2, E^2, C^2,
          \funcname{f}_0,\funcname{f}_1
        )
    }
  ]*
\end{multline}
\begin{equation}
\LCS^{\Ck,\catseqname{M}}(\seqname{E}, \seqname{C}) \defeq
\Bigl (
  \LCS^{\Ck,\catseqname{M}}_\Ob(\seqname{E}, \seqname{C}),
  \LCS^{\Ck,\catseqname{M}}_\Ar(\seqname{E}, \seqname{C}),
  \compose[A]
\Bigr )
\end{equation}
\begin{multline}
\LCS^{\minimal{\Ck}}_\Ob(\seqname{E}, \seqname{C}) \defeq
\set
  {
    \Functor^{\minimal{\Ck}}_{\Man,\LCS} \seqname{S}
  }%
  [
    {\seqname{S} \defeq (\seqname{A}, E, C)},
    E \in \seqname{E},
    C \in \seqname{C},
    {\maximal{\isAtl^\Ck_\Ob}(\seqname{A}, E, C)}
  ]*
\end{multline}
\begin{multline}
\LCS^{\minimal{\Ck}}_\Ar(\seqname{E}, \seqname{C}) \defeq \\
\set
  {
    {
      \Functor^{\minimal{\Ck}}_{\Man,\LCS}
        \bigl (
          (
            \funcname{f}_0,
            \funcname{f}_1
          ),
          (
            \seqname{A}^1,
            E^1,
            C^1
          ),
          (
              \seqname{A}^2,
              E^2,
              C^2
          )
        \bigr )
    }
  }%
  [
    E^i \in \seqname{E},
    C^i \in \seqname{C},
    {
      \maximal{\isAtl^\Ck_\Ar}
        (
          \seqname{A}^1, E^1, C^1,
          \seqname{A}^2, E^2, C^2,
          \funcname{f}_0, \funcname{f}_1
        )
    }
  ]*
\end{multline}
\begin{equation}
\LCS^{\minimal{\Ck}}(\seqname{E}, \seqname{C}) \defeq
\Bigl (
  \LCS^{\minimal{\Ck}}_\Ob(\seqname{E}, \seqname{C}),
  \LCS^{\minimal{\Ck}}_\Ar(\seqname{E}, \seqname{C}),
  \compose[A]
\Bigr )
\end{equation}
\begin{multline}
\LCS^{\minimal{\Ck},\hat{\catseqname{M}}}_\Ob(\seqname{E}, \seqname{C}) \defeq
\set
  {
    \Functor^{\minimal{\Ck}}_{\Man,\LCS} \seqname{S}
  }%
  [
    {\seqname{S} \defeq (\seqname{A}, E, C)},
    E \in \seqname{E},
    C \in \seqname{C},
    {\maximal{\isAtl^\Ck_\Ob}(\seqname{A}, E, C)}
  ]*
\end{multline}
\begin{multline}
\LCS^{\minimal{\Ck},\hat{\catseqname{M}}}_\Ar(\seqname{E}, \seqname{C}) \defeq \\
\set
  {
    {
      \Functor^{\minimal{\Ck},\hat{\catseqname{M}}}_{\Man,\LCS}
        \bigl (
          (
            \funcname{f}_0,
            \funcname{f}_1
          ),
          (
            \seqname{A}^1,
            E^1,
            C^1
          ),
          (
              \seqname{A}^2,
              E^2,
              C^2
          )
        \bigr )
    }
  }%
  [
    E^i \in \seqname{E},
    C^i \in \seqname{C},
    {
      \maximal{\isAtl^\Ck_\Ar}
        (
          \seqname{A}^1, E^1, C^1,
          \seqname{A}^2, E^2, C^2,
          \funcname{f}_0, \funcname{f}_1
        )
    }
  ]*
\end{multline}
\begin{equation}
\LCS^{\minimal{\Ck},\hat{\catseqname{M}}}(\seqname{E}, \seqname{C}) \defeq
\Bigl (
  \LCS^{\minimal{\Ck},\hat{\catseqname{M}}}_\Ob(\seqname{E}, \seqname{C}),
  \LCS^{\minimal{\Ck},\hat{\catseqname{M}}}_\Ar(\seqname{E}, \seqname{C}),
  \compose[A]
\Bigr )
\end{equation}
\end{definition}

\begin{theorem}[Functors from $\mathrm{C^k}$ manifolds to Local 2-$\emptyset$ coordinate spaces]
\label{the:CktoLCS}
Let $\seqname{E}$ be a set of topological spaces,
$\seqname{C}$ be a set of $\Ck$ linear spaces,
$\catname{E}$ and $\hat{\catname{E}}$ model categories,
$\catname{C}$ and $\hat{\catname{C}}$ $\Ck$ linear model categories,
$\catseqname{M} \defeq (\catname{E}, \catname{C})$,
$
  \hat{\catseqname{M}} \defeq
  (\hat{\catname{E}}, \hat{\catname{C}})
$,
$
  \Triv{\catseqname{M}} \defeq
  \Bigl ( \Trivcat{\seqname{E}}, \Trivcat[\Ck-]{\seqname{C}} \Bigr )
$.
$
  \optriv{\catseqname{M}} \defeq
  \Bigl ( \Trivcat{\seqname{E}}, \optriv[\Ck-]{\seqname{C}} \Bigr )
$.
$
  \Triv{\catseqname{M}} \SUBCAT \catseqname{M}
  \SUBCAT \hat{\catseqname{M}}
$,
$\optriv{\catseqname{M}} \SUBCAT \hat{\catseqname{M}}$.

Then
$\LCS^\Ck(\seqname{E}, \seqname{C})$,
$\LCS^{\Ck,\catseqname{M}}(\seqname{E}, \seqname{C})$,
$\LCS^{\minimal{\Ck}}(\seqname{E}, \seqname{C})$ and
$\LCS^{\minimal{\Ck},\hat{\catseqname{M}}}(\seqname{E}, \seqname{C})$
are categories and the identity morphism of
$
   \seqname{L}^\alpha \defeq
    \bigl (
      (\catname{E}, \catname{C}),
      (\seqname{E}^\alpha, \seqname{C}^\alpha),
      \seqname{A}^\alpha,
      \emptyset,
      \emptyset
    \bigr )
$
is that given in \pagecref{def:LCSmorphAr}:
$
  \Id_{\seqname{L}^\alpha}
  \defeq
  \Bigl (
    \bigl ( \Id_{\seqname{E}^\alpha}, \Id_{\seqname{C}^\alpha} \bigr ),
    \seqname{L}^\alpha,
    \seqname{L}^\alpha
  \Bigr )
$.

\begin{proof}
$\LCS^\Ck(\seqname{E}, \seqname{C})$,
$\LCS^{\Ck,\catseqname{M}}(\seqname{E}, \seqname{C})$,
$\LCS^{\minimal{\Ck}}(\seqname{E}, \seqname{C})$ and
$\LCS^{\minimal{\Ck},\hat{\catseqname{M}}}(\seqname{E}, \seqname{C})$
satisfy the definition of a category:
\begin{enumerate}
\item Composition: \\
Let
$
  m^i \defeq
  \Functor^{\Ck}_{\Man,\LCS}
    \bigl (
      (
        \funcname{f}^i_0,
        \funcname{f}^i_1
      ),
      \seqname{S}^i,
      \seqname{S}^{i+1}
    \bigr )
$
be morphisms of $\LCS^{\Ck}(\seqname{E}, \seqname{C})$.
$
  m^2 \compose[A] m^1 =
  \Functor^{\Ck}_{\Man,\LCS}
    \bigl (
      (
        \funcname{f}^2_0 \compose \funcname{f}^1_0,
        \funcname{f}^2_1 \compose \funcname{f}^2_1
      ),
      \seqname{S}^1,
      \seqname{S}^3
    \bigr )
$.
By \pagecref{lem:M-ATLmorph},
$
  (
    \funcname{f}^2_0 \compose \funcname{f}^1_0,
    \funcname{f}^2_1 \compose \funcname{f}^2_1
  )
$
is an m-atlas morphism.
\item Associativity:
Composition is associative by \pagecref{lem:atlcomp}.
\item Unit:
$\Id_{\seqname{L}^\alpha}$ is an identity morphism by \cref{lem:atlcomp}.
\end{enumerate}

The same proof applies to $\LCS^{\Ck,\catseqname{M}}(\seqname{E}, \seqname{C})$,
$\LCS^{\minimal{\Ck}}(\seqname{E}, \seqname{C})$ and
$\LCS^{\minimal{\Ck},\hat{\catseqname{M}}}(\seqname{E}, \seqname{C})$.
\end{proof}

Let
$E^i \in \seqname{E}$, $i \in [1,3]$, $C^i \in \seqname{C}$,
$\seqname{M}^i \defeq \Bigl ( \Triv{E^i}, \Triv[\Ck-]{C^i} \Bigr )$,
$
  \Singcat{\catseqname{M}^i} \defeq
  (
    \singcat{\seqname{M}^i_0},
    \singcat[\Ck-]{\seqname{M}^i_1}
  )
$,
$\seqname{A}^i$ a maximal $\Ck$-atlas of $E^i$ in $C^i$,
$
  \Functor^{\minimal{\Ck}}_1 (\seqname{A}^i, E^i, C^i) \objin
  \hat{\catname{E}}
$,
$
  \Functor^{\minimal{\Ck}}_2 (\seqname{A}^i, E^i, C^i) \objin
  \hat{\catname{C}}
$,
$(E^i, C^i, \seqname{A}^i)$ a $\Ck$ manifold,
$\seqname{S}^i \defeq (\seqname{A}^i,E^i, C^i)$,
$
  \minimal{\seqname{M}^i} \defeq \\
    \Bigl (
      \Functor^{\minimal{\Ck}}_1 (\seqname{A}^i, E^i, C^i),
      \Functor^{\minimal{\Ck}}_2 (\seqname{A}^i, E^i, C^i)
    \Bigr )
$,
$
  \Singcat{\Triv{\catseqname{M}^i}} \defeq
  (
    \singcat{\Triv{\seqname{M}^i_0}},
    \singcat[\Ck-]{\Triv{\seqname{M}^i_1}}
  )
$,
$
\seqname{L}^i \defeq
\Functor^{\Ck}_{\Man,\LCS} \seqname{S}^i =
\Bigl (
  \bigl (
    \Trivcat{E^i}, \Trivcat[\Ck-]{C^i}
  \bigr ),
  \bigl (
    \Triv{E^i}, \Triv[\Ck-]{C^i}
  \bigr ),
  \seqname{A}^i,
  \emptyset,
  \emptyset
\Bigr )
$,
$
\seqname{L}^{i,\catseqname{M}} \defeq
\Functor^{\Ck,\catseqname{M}}_{\Man,\LCS} \seqname{S}^i = \\
\Bigl (
  \catseqname{M},
  \bigl (
    \Triv{E^i}, \Triv[\Ck-]{C^i}
  \bigr ),
  \seqname{A}^i,
  \emptyset,
  \emptyset
\Bigr )
$,
$
\Triv{\seqname{L}^i} \defeq \\
\Bigl (
  \singcat{\Triv{\catseqname{M}^i}},
  \bigl (
    \Triv{E^i}, \optriv[\Ck-]{C^i}
  \bigr ),
  \seqname{A}^i,
  \emptyset,
  \emptyset
\Bigr )
$,
$
\Triv{\seqname{L}^{i,\hat{\catseqname{M}}}} \defeq
\Bigl (
  \hat{\catseqname{M}},
  \bigl (
    \Triv{E^i}, \optriv[\Ck-]{C^i}
  \bigr ),
  \seqname{A}^i,
  \emptyset,
  \emptyset
\Bigr )
$,
$\funcseqname{f}^i \defeq (\funcname{f}^i_0, \funcname{f}^i_1)$,
$\funcname{f}^i_0 \maps E^i \to E^{i+1}$ continuous and
$\funcname{f}^i_1 \maps C^i \to C^{i+1}$ $\Ck$.

$\Functor^\Ck_{\Man,\LCS} \seqname{S}^i$ is a local
$\Singcat{\catseqname{M}^i}$-$\emptyset$ coordinate space and
$\Functor^{\Ck,\catseqname{M}}_{\Man,\LCS} \seqname{S}^i$ is a local
$\catseqname{M}$-$\emptyset$ coordinate space.

\begin{proof}
$\Functor^{\Ck,\catseqname{M}}_{\Man,\LCS} \seqname{S}^i$
satisfies the conditions of \pagecref{def:LCS}:
\begin{enumerate}
\item Model categories:
$\Triv{E^i}$ and $\Triv[\Ck-]{C^i}$ are model spaces
\item
$
\biggl (
  \catseqname{M},
  \Bigl ( \Triv{E^i}, \Triv[\Ck-]{C^i} \Bigr )
  \emptyset,
  \emptyset
\biggr )
$
satisfies the conditions of \pagecref{def:pre}
\item Maximal m-atlas:
$\seqname{A}^i$ is a maximal $\Ck$ atlas by hypothesis and a
maximal atlas of $\Triv{E^1}$ in the coordinate space $\Triv[\Ck-]{C^i}$
by construction.
\item Constraint functions:
There are no adjunct functions so there are no constraint functions.
\end{enumerate}

$\Singcat{\catseqname{M}^i} \SUBCAT \catseqname{M}$, so
$\Functor^{\Ck,\catseqname{M}}_{\Man,\LCS} \seqname{S}^i$ is a local
$\catseqname{M}$-$\emptyset$ coordinate space by \pagecref{lem:LCS}.
\end{proof}

$\Functor^{\minimal{\Ck}}_{\Man,\LCS} \seqname{S}^i$ is a local
$\Singcat{\Triv{\catseqname{M}^i}}$-$\emptyset$ coordinate space and
$
  \Functor^{\minimal{\Ck},\hat{\catseqname{M}}}_{\Man,\LCS}
    \seqname{S}^i
$
is a local $\hat{\catseqname{M}}$-$\emptyset$ coordinate space.

\begin{proof}
$\Functor^{\minimal{\Ck}}_{\Man,\LCS} \seqname{S}^i$
satisfies the conditions of \pagecref{def:LCS}:
\begin{enumerate}
\item Model categories:
$\Triv{E^1}$ and $\optriv[\Ck-]{C^i}$ are model spaces
\item
$
  \Bigl (
    \bigl ( \Triv{E^i}, \optriv[\Ck-]{C^i} \bigr ), \emptyset, \emptyset
  \Bigr )
$
satisfies the conditions of \pagecref{def:pre}
\item Maximal m-atlas:
$\seqname{A}^i$ is a maximal $\Ck$ atlas by hypothesis and a
maximal atlas of $\Triv{E^1}$ in the coordinate space
$\optriv[\Ck-]{C^i}$ by construction.
\item Constraint functions:
There are no adjunct functions so there are no constraint functions.
\end{enumerate}

$\Singcat{\Triv{\catseqname{M}^i}} \SUBCAT \hat{\catseqname{M}}$, so
$\Functor^{\Ck,\hat{\catseqname{M}}}_{\Man,\LCS} \seqname{S}^i$ is a local
$\hat{\catseqname{M}}$-$\emptyset$ coordinate space by \cref{lem:LCS}.
\end{proof}

$(\funcname{f}^1_0, \funcname{f}^1_1)$ is a morphism from
$\seqname{L}^1$ to $\seqname{L}^2$ and a strict morphism from
$\seqname{L}^{1,\catseqname{M}}$ to $\seqname{L}^{2,\catseqname{M}}$.

\begin{proof}
It satisfies the conditions of \pagecref{def:LCSmorph}\!:
\begin{enumerate}
\item Prestructure morphism:
$\funcseqname{f}^1$ $\Sigma$-commutes with $\emptyset,\emptyset$.
$\funcname{f}^1_0$ is continuous, hence a morphism of $\catname{E}$.
$\funcname{f}^1_1$ is $\Ck$, hence a morphism of $\catname{C}$.

\item m-atlas morphism:
$(\funcname{f}^1_0, \funcname{f}^1_1)$ is a $\Triv{E^1}$-$\Triv{E^2}$
m-atlas morphism of $\seqname{A}^1$ to $\seqname{A}^2$ in the coordinate
spaces $\Triv[\Ck-]{C^1}$, $\Triv[\Ck-]{C^i}$.

$(\funcname{f}^1_0, \funcname{f}^1_1)$ is a strict
$\singcat{\Triv{E^1}}$-$\singcat{\Triv{E^2}}$-%
$\Triv{E^1}$-$\Triv{E^2}$-%
$\singcat[\Ck-]{\Triv[\Ck-]{C^1}}$-$\singcat[\Ck-]{\Triv[\Ck-]{C^1}}$
m-atlas morphism of $\seqname{A}^1$ to $\seqname{A}^2$ in the coordinate
spaces $\Triv[\Ck-]{C^1}$, $\Triv[\Ck-]{C^i}$. by
\pagecref{lem:Ck-ATLmorph}.
\end{enumerate}
\end{proof}

$
\Functor^{\minimal{\Ck}}_{\Man,\LCS}
  \bigl (
    (\funcname{f}^1_0, \funcname{f}^1_1), \seqname{S}^1, \seqname{S}^2
  \bigr )
$
is a morphism from $\Triv{\seqname{L}^1}$ to
$\Triv{\seqname{L}^{2,\catseqname{M}}}$ and \\
$
\Functor^{\minimal{\Ck}}_{\Man,\LCS}
  \bigl (
    (\funcname{f}^1_0, \funcname{f}^1_1), \seqname{S}^1, \seqname{S}^2
  \bigr )
$
is a strict morphism from
$\Triv{\seqname{L}^{1,\hat{\catseqname{M}}}}$ to
$\Triv{\seqname{L}^{2,\hat{\catseqname{M}}}}$.

\begin{proof}
It satisfies the conditions of \pagecref{def:LCSmorph}\!:
\begin{enumerate}
\item Prestructure morphism:
$(\funcname{f}^1_0, \funcname{f}^1_1)$ $\Sigma$-commutes with $\emptyset,\emptyset$.
\item m-atlas morphism:
$(\funcname{f}^1_0, \funcname{f}^1_1)$ is an m-atlas morphism from
$(\Triv{E^1}, \optriv[\Ck-]{C^1})$ to
$(\Triv{E^2}, \optriv[\Ck-]{C^2})$.
\end{enumerate}
\end{proof}

$\Functor^{\Ck,\catseqname{M}}_{\Man,\LCS}$ is a functor from
$\Man^\Ck(\seqname{E}, \seqname{C})$ to
$\LCS^{\Ck,\catseqname{M}}(\seqname{E}. \seqname{C})$.

\begin{proof}
$\Functor^{\Ck,\catseqname{M}}_{\Man,\LCS}$ satisfies the definition of a functor:
\begin{enumerate}
\item F(f: A to B) = F(f): F(A) to  F(B):
\begin{enumerate}
\item
$
\Functor^{\Ck,\catseqname{M}}_{\Man,\LCS}
  \bigl (
    (
      \funcseqname{f}^1_0,
      \funcseqname{f}^1_1
    ),
    \seqname{S}^1,
    \seqname{S}^2
  \bigr )
=
  \bigl (
    (
      \funcseqname{f}^1_0,
      \funcseqname{f}^1_1
    ),
    \seqname{L}^i,
    \seqname{L}^i
  \bigr )
$
\item
$
\Functor^{\Ck,\catseqname{M}}_{\Man,\LCS} \seqname{S}^i
=
  \left (
    \Bigl ( \Triv{E^i}, \Triv[\Ck-]{C^i} \Bigr ),
    \seqname{A}^i,
    \emptyset,
    \emptyset
  \right )
=
  \seqname{L}^i
$
\end{enumerate}
\item Composition:
$
  \Functor^{\Ck,\catseqname{M}}_{\Man,\LCS}
    \Bigl (
      \bigl (
        (\funcseqname{f}^2_0, \funcseqname{f}^2_1),
        \seqname{S}^2,
        \seqname{S}^3
      \bigr )
      \compose[A]
      \bigl (
        (\funcseqname{f}^1_0, \funcseqname{f}^1_1),
        \seqname{S}^1,
        \seqname{S}^2
      \bigr )
    \Bigr )
  = \\
  \Functor^{\Ck,\catseqname{M}}_{\Man,\LCS}
    \Bigl (
      \bigl (
        (
          \funcseqname{f}^2_0 \compose \funcseqname{f}^1_0,
          \funcseqname{f}^2_1 \compose \funcseqname{f}^1_1
        ),
        \seqname{S}^1,
        \seqname{S}^3
      \bigr )
    \Bigr )
  = \\
  \bigl (
    (
      \funcseqname{f}^2_0 \compose \funcseqname{f}^1_0,
      \funcseqname{f}^2_1 \compose \funcseqname{f}^1_1
    ),
    \seqname{L}^1,
    \seqname{L}^3
  \bigr )
  = \\
  \Functor^{\Ck,\catseqname{M}}_{\Man,\LCS}
    \bigl (
      (
        \funcseqname{f}^2_0,
        \funcseqname{f}^2_1
      ),
      \seqname{S}^2,
      \seqname{S}^3
    \bigr )
  \compose[A]
  \Functor^{\Ck,\catseqname{M}}_{\Man,\LCS}
    \bigl (
      (
        \funcseqname{f}^1_0,
        \funcseqname{f}^1_1
      ),
      \seqname{S}^1,
      \seqname{S}^2
    \bigr )
$
\item Identity: \newline
$
\Functor^{\Ck,\catseqname{M}}_{\Man,\LCS} (\Id_{\seqname{S}^i}) =
\Functor^{\Ck,\catseqname{M}}_{\Man,\LCS}
  \Bigl (
    (\Id_{E^i}, \Id_{C^i}),
    \seqname{S}^i,
    \seqname{S}^i
  \Bigr ) =\\
  \biggl (
   \Bigl ( \Id_{\Triv{E^i}}, \Id_{\Triv[\Ck-]{C^i}} \Bigr ),
    \seqname{L}^1,
    \seqname{L}^1
  \biggr )
= \\
\Id_{\seqname{L}^1} = \Id_{\Functor^{\Ck,\catseqname{M}}_{\Man,\LCS}(\seqname{S}^i)}
$
\end{enumerate}
\end{proof}

$\Functor^{\minimal{\Ck}}_{\Man,\LCS}$ is a functor from
$\Man^\Ck(\seqname{E}, \seqname{C})$ to
$\LCS^{\minimal{\Ck}}(\seqname{E}. \seqname{C})$.

\begin{proof}
$\Functor^{\minimal{\Ck}}_{\Man,\LCS}$ satisfies the definition of a functor:
\begin{enumerate}
\item F(f: A to B) = F(f): F(A) to  F(B):
\begin{enumerate}
\item
$
\Functor^{\minimal{\Ck}}_{\Man,\LCS}
  \bigl (
    (
      \funcseqname{f}^1_0,
      \funcseqname{f}^1_1
    ),
    \seqname{S}^1,
    \seqname{S}^2
  \bigr )
=
  \bigl (
    (
      \funcseqname{f}^1_0,
      \funcseqname{f}^1_1
    ),
    \seqname{L}^i,
    \seqname{L}^i
  \bigr )
$
\item
$
\Functor^{\minimal{\Ck}}_{\Man,\LCS} \seqname{S}^i
=
  \left (
    \Bigl ( \Triv{E^i}, \optriv[\Ck-]{C^i} \Bigr ),
    \seqname{A}^i,
    \emptyset,
    \emptyset
  \right )
=
  \seqname{L}^i
$
\end{enumerate}
\item Composition:
$
  \Functor^{\minimal{\Ck}}_{\Man,\LCS}
    \Bigl (
      \bigl (
        (\funcseqname{f}^2_0, \funcseqname{f}^2_1),
        \seqname{S}^2,
        \seqname{S}^3
      \bigr )
      \compose[A]
      \bigl (
        (\funcseqname{f}^1_0, \funcseqname{f}^1_1),
        \seqname{S}^1,
        \seqname{S}^2
      \bigr )
    \Bigr )
  = \\
  \Functor^{\minimal{\Ck}}_{\Man,\LCS}
    \Bigl (
      \bigl (
        (
          \funcseqname{f}^2_0 \compose \funcseqname{f}^1_0,
          \funcseqname{f}^2_1 \compose \funcseqname{f}^1_1
        ),
        \seqname{S}^1,
        \seqname{S}^3
      \bigr )
    \Bigr )
  = \\
  \biggl (
    (
      \funcseqname{f}^2_0 \compose \funcseqname{f}^1_0,
      \funcseqname{f}^2_1 \compose \funcseqname{f}^1_1
    ),
    \hat{\seqname{L}}^1,
    \hat{\seqname{L}}^3
  \biggr )
  = \\
  \Functor^{\minimal{\Ck}}_{\Man,\LCS}
    \bigl (
      (
        \funcseqname{f}^2_0,
        \funcseqname{f}^2_1
      ),
      \seqname{S}^2,
      \seqname{S}^3
    \bigr )
  \compose[A]
  \Functor^{\minimal{\Ck}}_{\Man,\LCS}
    \bigl (
      (
        \funcseqname{f}^1_0,
        \funcseqname{f}^1_1
      ),
      \seqname{S}^1,
      \seqname{S}^2
    \bigr )
$
\item Identity: \newline
$
\Functor^{\minimal{\Ck}}_{\Man,\LCS} (\Id_{\seqname{S}^i}) =
\Functor^{\minimal{\Ck}}_{\Man,\LCS}
  \bigl (
    (\Id_{E^i}, \Id_{C^i}),
    \seqname{S}^i,
    \seqname{S}^i
  \bigr )
= \\
  \biggl (
   \bigl ( \Id_{\Triv{E^i}}, \Id_{\Triv[\Ck-]{C^i}} \bigr ),
    \Bigl (
      \bigl ( \Triv{E^i}, \Triv[\Ck-]{C^i} \bigr ),
      \seqname{A}^i,
      \emptyset,
      \emptyset
    \Bigr ),
    \Bigl (
      \bigl ( \Triv{E^i}, \Triv[\Ck-]{C^i} \bigr ),
      \seqname{A}^i,
      \emptyset,
      \emptyset
    \Bigr )
  \biggr )
= \\
\Id_
  {
    \Bigl (
      \bigl ( \Triv{E^i}, \Triv[\Ck-]{C^i} \bigr ),
      \seqname{A}^i,
      \emptyset,
      \emptyset
    \Bigr )
  }
\defeq
\Id_{\Functor^{\minimal{\Ck}}_{\Man,\LCS}(\seqname{S}^i)}
$
\end{enumerate}
\end{proof}
\end{theorem}

\begin{definition}[Functor from Local 2-$\emptyset$ coordinate spaces to $\mathrm{C^k}$ manifolds]
\label{def:LCStoCk}
Let $\seqname{C^i} \defeq (C^i, \catname{C}^i)$,
$i=1,2$ be a $\Ck$ linear model space, $\seqname{E}^i$ a
model space,
$\catname{M}^i_0$ a model category,
$\catname{M}^i_1$ a $\Ck$ linear model category,
$\catseqname{M}^i \defeq (\catname{M}^i_0, \catname{M}^i_1)$,
$\seqname{M}^i \defeq (\seqname{E}^i, \seqname{C}^i)$,
$\seqname{M}^i \seqin \catseqname{M}^i$,
$\seqname{A}^i$ maximal m-atlases of $\seqname{E}^i$ in $\seqname{C}^i$,
$\funcname{f}_0 \maps \seqname{E}^1 \to \seqname{E}^2$,
$\funcname{f}_1 \maps \seqname{C}^1 \to \seqname{C}^2)$,
$\funcname{f} \defeq (\funcname{f}_0, \funcname{f}_1)$ and
$
  \seqname{L}^i \defeq
  (
    \catseqname{M}^i,
    \seqname{M}^i,
    \seqname{A}^i,
    \emptyset,
    \emptyset
  )
$
a local 2-$\emptyset$ coordinate spaces. Then
\begin{equation}
\Functor^\Ck_{\LCS,\Man} \seqname{L}^i \defeq
(E^i, C^i, \seqname{A}^i)
\end{equation}
\begin{equation}
\Functor^\Ck_{\LCS,\Man}
  (
    \funcseqname{f},
    \seqname{L}^1,
    \seqname{L}^2
  )
\defeq
  \bigl (
    \funcseqname{f},
    (E^1, C^1, \seqname{A}^1),
    (E^2, C^2, \seqname{A}^2)
  \bigr )
\end{equation}
\end{definition}

\begin{theorem}[Functor from Local 2-$\emptyset$ coordinate spaces to $\mathrm{C^k}$ manifolds]
\label{the:LCStoCk}
Let $\seqname{C^i} \defeq (C^i, \catname{C}^i)$,
$i=1,2$, be $\Ck$ linear model spaces,
$\seqname{E}^i \defeq (E^i, \catname{E}^i)$ model spaces,
$\seqname{A}^i$ maximal $\Ck$-atlases of $\seqname{E}^i$ in $\seqname{C}^i$,
$\seqname{M}^i \defeq (\seqname{E}^i, \seqname{C}^i)$,
$\catseqname{M}^i \defeq \Singcat{\seqname{M}^i}$,
$\funcname{f}^i_0 \maps \seqname{E}^i \to \seqname{E}^{i+1}$
$\funcname{f}^i_1 \maps \seqname{C}^i \to \seqname{C}^{i+1})$ and
$
  \seqname{L}^i \defeq
  (
    \catseqname{M}^i,
    \seqname{M}^i,
    \seqname{A}^i,
    \emptyset,
    \emptyset
  )
$
a local 2-$\emptyset$ coordinate space. Then
\begin{enumerate}
\item
$
  \Functor^\Ck_{\LCS,\Man} \seqname{L}^i
  \defeq (E^i, \seqname{C}^i, \seqname{A}^i)
$
is a $\Ck$ manifold.
\begin{proof}
$(E^i, \seqname{C}^i, \seqname{A}^i)$ satisfies the conditions of
\pagecref{def:man}:
$(E^i$ is a topological space, $\seqname{C}^i$ is a $\Ck$ linear space and
$\seqname{A}^i$ is a maximal $\Ck$-atlas $\seqname{E}^i$ in $\seqname{C}^i$.
\end{proof}
\item
$\Functor^\Ck_{\LCS,\Man}$ is a functor from
$\LCS^{\Ck,\catseqname{M}}(\seqname{E}. \seqname{C})$ to
$\Man^\Ck(\seqname{E}, \seqname{C})$.
\begin{proof}
$\Functor^\Ck_{\LCS,\Man}$ satisfies the definition of a functor:
Let $\seqname{S}^i=(\seqname{A}^i, E^i, \seqname{C}^i)$.
\begin{enumerate}
\item $\Functor(f \maps A \to B) \maps \Functor(A) \to \Functor(B)$: \newline
\begin{enumerate}
\item
$
\Functor^\Ck_{\LCS,\Man}
  \bigl
    (
      \funcname{f}^i_0,
      \funcname{f}^i_1
    ),
    \seqname{L}^i,
    \seqname{L}^{i+1}
  \bigr )
=
  \bigl
    (
      \funcname{f}^i_0,
      \funcname{f}^i_1
    ),
    \seqname{S}^i,
    \seqname{S}^{i+1}
  \bigr )
$
\item
$
\Functor^\Ck_{\LCS,\Man} (\seqname{L}^i) = \seqname{S}^i
$
\end{enumerate}
\item Composition: \newline
$
\Functor^\Ck_{\LCS,\Man}
  \Bigl (
    \bigl (
      (
        \funcname{f}^2_0,
        \funcname{f}^2_1
      ),
      \seqname{L}^2,
      \seqname{L}^3
    \bigr )
  \compose[A]
    \bigl (
      (
        \funcname{f}^1_0,
        \funcname{f}^1_1
      ),
      \seqname{L}^1,
      \seqname{L}^2
    \bigr )
  \Bigr )
= \\
\Functor^\Ck_{\LCS,\Man}
  \Bigl (
    \bigl (
      (
        \funcname{f}^2_0 \compose \funcname{f}^1_0,
        \funcname{f}^2_1 \compose \funcname{f}^1_1
      ),
      \seqname{L}^1,
      \seqname{L}^3
    \bigr )
  \Bigr )
= \\
  \Bigl (
    \bigl
      (
        \funcname{f}^2_0 \compose \funcname{f}^1_0,
        \funcname{f}^2_1 \compose \funcname{f}^1_1
      ),
      \seqname{S}^1,
      \seqname{S}^3
    \bigr )
  \Bigr )
= \\
  \Bigl (
    \bigl (
      (
        \funcname{f}^2_0,
        \funcname{f}^2_1
      ),
      \seqname{S}^2,
      \seqname{S}^3
    \bigr )
  \Bigr )
  \compose[A]x
  \Bigl (
    \bigl (
      (
        \funcname{f}^1_0,
        \funcname{f}^1_1
      ),
      \seqname{S}^1,
      \seqname{S}^2
    \bigr )
  \Bigr )
= \\
\Functor^\Ck_{\LCS,\Man}
  \Bigl (
    \bigl (
      (
        \funcname{f}^2_0,
        \funcname{f}^2_1
      ),
      \seqname{L}^2,
      \seqname{L}^3
    \bigr )
  \Bigr )
\compose[A]
\Functor^\Ck_{\LCS,\Man}
  \Bigl (
    \bigl (
      (
        \funcname{f}^1_0,
        \funcname{f}^1_1
      ),
      \seqname{L}^1,
      \seqname{L}^2
    \bigr )
  \Bigr )
$
\item Identity: \newline
$
\Functor^\Ck_{\LCS,\Man} (\Id_{L^i})
=
\Functor^\Ck_{\LCS,\Man}
  \bigl
    (\Id_{E^i}, \Id_{\seqname{C}^i)}, \seqname{L}^i, \seqname{L}^i
  \bigr )
= \\
  \bigl
    (\Id_{E^i}, \Id_{\seqname{C}^i}), \seqname{S}^i, \seqname{S}^i
  \bigr )
=
\Id_{\Functor^\Ck_{\LCS,\Man}} (L^i)
$
\end{enumerate}
\end{proof}
\item
$\Functor^\Ck_{\LCS,\Man} \compose \Functor^{\Ck,\catseqname{M}}_{\Man,\LCS} = \Id$
\begin{proof}
$
\Functor^\Ck_{\LCS,\Man} \compose \Functor^{\Ck,\catseqname{M}}_{\Man,\LCS}
  (\seqname{S}^i) =
\Functor^\Ck_{\LCS,\Man} (\seqname{L}^i) =
\seqname{S}^i
$
\end{proof}
\item
$\Functor^{\Ck,\catseqname{M}}_{\Man,\LCS} \compose \Functor^\Ck_{\LCS,\Man} = \Id$
\begin{proof}
$
\Functor^{\Ck,\catseqname{M}}_{\Man,\LCS} \compose \Functor^\Ck_{\LCS,\Man}
  (\seqname{L}^i) =
\Functor^{\Ck,\catseqname{M}}_{\Man,\LCS} (\seqname{S}^i) =
\seqname{L}^i
$
\end{proof}
\end{enumerate}
\end{theorem}

\section{Equivalence of fiber bundles}
\label{sec:bun}
For fiber bundles\footnote{
  The literature has several definitions of fiber bundle. This paper
  uses one chosen for clarity of exposition. It differs from
  \cite[p.~8]{TopFib} in that, e.g., it uses the machinery of maximal
  atlases rather than equivalence classes of coordinate bundles, the
  nomenclature differs in several minor regards.
},
the adjunct spaces are the base space $\seqname{X} = (X,\catname{X})$,
the fiber $\seqname{Y} = (Y, \catname{Y})$
and the group $G$;
the category of the coordinate space is the category of Cartesian
products $\{U \times Y \mid U \objin \catname{X}\}$ of model
neighborhoods in the base space with the entire fiber, with morphisms
$\funcname{t} \maps U \times Y \to U \times Y$ that preserve the
fibers, i.e., $\pi_1 \compose \funcname{t} = \pi_1$, and are
generated by the group action on the fiber
(Equation~\eqref{eq:transitiongroup}).

The sole adjunct functions are the projection $\pi \maps E \onto X$, the
group operation and the group action on the fiber.

This section defines bundle atlases, fiber bundles, local coordinate
spaces equivalent to fiber bundles, categories of them and functors, and
gives basic results

\begin{definition}[Trivial group category of groups]
Let $\seqname{G}$ be a set of topological groups. The trivial group
category of $\seqname{G}$, abbreviated $\Trivcat[\mathrm{group}-]{\seqname{G}}$, is
the category of all continuous homomorphisms between groups of
$\seqname{G}$. By abuse of language it will be shortened to
$\Triv{\seqname{G}}$ when the meaning is clear from context.
\end{definition}

\begin{definition}[Group actions]
Let $Y$ be a topological space, $G$ a topological group, $\rho$ an
effective group action of $G$ on $Y$, $y \in Y$ and $g \in G$. Then
$y \star g \defeq \rho(y,g)$.

Let $X$ be a topological space and $x \in X$. Then
$(x, y) \star g \defeq (x, y \star g)$.

A $\star$ with a subscript or superscript refers to the group action
$\rho$ with the same subscript or superscript.

\begin{remark}
This notation is only used when it is clear from context what the group
action is.
\end{remark}
\end{definition}

\begin{definition}[Protobundles]
\label{def:ProtoBun}
Let $\seqname{B} \defeq (E, X, Y, G, \pi, \rho)$, where $E$, $X$ and $Y$
are topological spaces, $G$ a topological group, $\pi \maps E \onto X$ a
continuous surjection and $\rho$ an efective group action of $G$ on $Y$.
Then $\seqname{B}$ is a protobundle.
\begin{remark}
While this definition does not itself require $E$ to have a local product
structure nor require $\pi$
to have the Covering Homotopy Property, only those protobundles having
an atlas are of interest, and for them \pagecref{def:BunAtl} imposes
additional constraints.
\end{remark}
\end{definition}

\subsection{$G$-$\rho$ model spaces}
\label{sub:grhomod}
\begin{definition}[$G$-$\rho$-model spaces]
\label{def:Grhomod}
Let $\seqname{XY} \defeq (X \times Y, \catname{XY})$ be a model space,
$G$ a topological group and $\rho$ an effective group action of $G$ on
$Y$ such that the objects of $\catname{XY}$ are products of open sets
with $Y$ and the morphisms are fiber-preserving automorphisms generated
by the group action, i.e.,
\begin{equation}
\uquant
  {
    U \objin \catname{XY}
  }
  {
    \equant
      {
        V \in \op{X}
      }
      {
        U = V \times Y
      }
  }
\end{equation}
\begin{equation}
\label{eq:fVV}
\uquant
  {
   \funcname{f} \arin \catname{XY} \maps V \times Y \toiso V' \times Y
  }
  {
    V = V'
  }
\end{equation}
\begin{equation}
\label{eq:fGrho}
\uquant
  {
   \funcname{f} \arin \catname{XY} \maps V \times Y \toiso V \times Y
  }
  {
    \equant
    {
      \funcname{g} \in G^V
    }
    {
    \uquant
      {
        {(x,y) \in V \times Y}
      }
      {
        \funcname{f}(x,y) =
          \bigl ( x, y \star \funcname{g}(x) \bigr )
      }
    }
  }
\end{equation}

Then $\seqname{XY}$ is a $G$-$\rho$ model space of $X \times Y$,
abbreviated \\
$
\mathrm{isG{\rho}}
(
  \seqname{XY}.
  Y,
  G,
  \rho
)
$,
$G_{\seqname{XY},\funcname{f}} \defeq \range(\funcname{g})$
and
$
  G_\seqname{XY} \defeq
  \union%
    [{\funcname{f} \arin \catname{XY}}]
    {G_{\seqname{XY},\funcname{f}}}
$.
\end{definition}

\begin{lemma}[$G$-$\rho$-model spaces]
Let $\seqname{XY}$ be a $G$-$\rho$ model space of $X \times Y$ and
$\funcname{f} \arin \catname{XY} \maps V \times Y \toiso V \times Y$.

The function $\funcname{g}$ in \cref{eq:fGrho} is unique.

\begin{proof}
The group action is effective.
\end{proof}

$G_\seqname{XY}$ is unique.

\begin{proof}
The function $\funcname{g}$ is unique.
\end{proof}
\end{lemma}

\begin{definition}[Morphisms of $G$-$\rho$-model spaces]
\label{def:Grhomorph}
Let
$X^i.Y^i$, $i=1,2$, be topological spaces,
$G^i$ a topological group,
$\rho^i$ an effective group action on $Y^i$ and
$\seqname{XY}^i \defeq (X^i \times Y^i, \catname{XY}^i)$ a
$G^i$-$\rho^i$ model space of $X^i \times Y^i$.
Then a model function
$\funcname{f}_C \maps \seqname{XY}^1 \to \seqname{XY}^2$ is a
$G^1$-$G^2$-$\rho^1$-$\rho^2$ morphism of $X^1 \times Y^1$ to $X^2 \times Y^2$,
abbreviated
$
\mathrm{isG{\rho}morph}
(
  \seqname{XY}^1,
  G^1,
  \rho^1,
  \seqname{XY}^2,
  G^2,
  \rho^2,
  \funcname{f}_C
)
$
iff it preserves the group action, i.e., there is a continuous
homomorphism $\funcname{f}_G \maps G^1 \to G^2$ such that
\fullcref{fig:Grho} is commutative, i.e., \cref{eq:Grho} holds.
\end{definition}

\begin{figure}
\[ \bfig
\node xyg1(0,1000)[\bigl ( (x^1,y^1), g^1 \bigr )]
\node xys1(0,0)[(x^1, y^1 \star g^1)]
\node xyg2(1000,1000)[\bigl ( (x^2,y^2), g^2 \bigr )]
\node xys2(1000,0)[(x^2, y^2 \star g^2)]
\arrow |r|[xyg1`xyg2;\funcname{f}_C \times \funcname{f}_G]
\arrow |a|[xyg1`xys1;\star^1]
\arrow |r|[xys1`xys2;\funcname{f}_C]
\arrow |a|[xyg2`xys2;\star^2]
\efig \]
\caption{Preserving group action}
\label{fig:Grho}
\end{figure}
\begin{equation}
\showlabelsinline
\label{eq:Grho}
\equant
  {
    \funcname{f}_G \maps G_1 \to G_2
  }
  {
    \uquant
    {
      {
        (x,y) \in X^1 \times Y^1
      },
      {
        g \in G^1
      }
    }
    {
      \funcname{f}_C \bigl ( (x,y) \star^1 g \bigr )
      =
      \funcname{f}_C \bigl ( (x,y) \bigr ) \star^2 \funcname{f}_G (g)
    }
  }
\end{equation}

\begin{lemma}[Morphisms of $G$-$\rho$-model spaces]
\label{lem:Grhomorph}
Let $X^i.Y^i$, $i=1,2,3$, be topological spaces, $G^i$ a topological
group, $\rho^i$ an effective group action on $Y^i$,
$\seqname{XY}^i \defeq (X^i \times Y^i, \catname{XY}^i)$ a
$G^i$-$\rho^i$ model space of $X^i \times Y^i$ and
$\funcname{f}^i_C \maps X^i \times Y^i \to X^{i+1} \times Y^{i+1}$ a
$G^i$-$G^{i+1}$-$\rho^i$-$\rho^{i+1}$ morphism of $X^i \times Y^i$ to
$X^{i+1} \times Y^{i+1}$.

The function $\funcname{f}_G$ in \cref{eq:Grho} is unique.

\begin{proof}
The group action is effective.
\end{proof}

$\funcname{f}^i_C$ preserves fibers, i.e.,
$
  \pi_1 \bigl ( \funcname{f}^i_C(x,y) \bigr ) =
  \pi_1 \bigl ( \funcname{f}^i_C(x,y') \bigr )
$
for $x$ in $X^i$ and $y,y'$ in $ Y^i$,
\begin{proof}
Since $\rho^i$ is effective, there is a $g \in G^i$ such that
$y' = y \star^i g$ and thus by \cref{def:Grhomorph}
\begin{equation}
\begin{split}
  \pi_1 \bigl ( \funcname{f}^i_C(x,y') \bigr )
  & =
\\
  \pi_1 \bigl ( \funcname{f}^i_C(x,y \star^i g) \bigr )
  & =
\\
  \pi_1 \bigl ( \funcname{f}^i_C(x,y) \star^{i+1} \funcname{f}^i_G(g) \bigr )
  & =
\\
  \pi_1 \bigl ( \funcname{f}^i_C(x,y) \bigr )
\end{split}
\end{equation}
\end{proof}

There exists a unique function $\funcname{f}^i_X \maps X^i \to X^{i+1}$
such that
$\funcname{f}^i_X \compose \pi_1 = \pi_1 \compose \funcname{f}^i_C$.

\begin{proof}
Define
$\funcname{f}^i_X(x) = \pi_1 \bigl ( \funcname{f}^i_C(x,y) \bigr )$,
where $y$ is an arbitrary point of $Y^i$. It does not depend on the
choice of $y$ because $\funcname{f}^i_C$ preserves fibers.
\end{proof}

\begin{remark}
$\funcname{f}^i_C$ may be a twisted product:  there need not exist
$\funcname{f}^i_Y \maps Y^i \to Y^{i+1}$ such that
$\funcname{f}^i_C = \funcname{f}^i_X \times \funcname{f}^i_Y$.
\end{remark}

$\funcname{f}^2_C \compose \funcname{f}^1_C$ is a
$G^1$-$G^3$-$\rho^1$-$\rho^3$ morphism of $X^1 \times Y^1$ to $X^3 \times Y^3$.

\begin{proof}

\begin{figure}
\[ \bfig
\node xyg1(0,1500)[\bigl ( (x^1,y^1), g^1 \bigr )]
\node xys1(0,0)[(x^1, y^1 \star g^1)]
\node xyg2(1000,1000)[\bigl ( (x^2,y^2), g^2 \bigr )]
\node xys2(1000,0)[(x^2, y^2 \star g^2)]
\node xyg3(2000,1500)[\bigl ( (x^3,y^3), g^3 \bigr )]
\node xys3(2000,0)[(x^3, y^3 \star g^3)]
\arrow |r|[xyg1`xyg3;\funcname{f}^2_C \compose \funcname{f}^1_C \times%
  \funcname{f}^2_G \compose \funcname{f}^1_G]
\arrow |r|[xyg1`xyg2;\funcname{f}^1_C \times \funcname{f}^1_G]
\arrow |a|[xyg1`xys1;\star^1]
\arrow |r|[xys1`xys2;\funcname{f}^1_C]
\arrow |a|[xyg2`xys2;\star^2]
\arrow |r|[xyg2`xyg3;\funcname{f}^2_C \times \funcname{f}^2_G]
\arrow |r|[xys2`xys3;\funcname{f}^1_C]
\arrow |a|[xyg3`xys3;\star^3]
\efig \]
\caption{Preserving group actions}
\label{fig:Grho2}
\end{figure}

Let $\funcname{f}^i_G \maps G^i \to G^{i+1}$ be a continuous
homomorphism such that
\begin{equation}
\showlabelsinline
\label{eq:Grho2}
    \uquant
    {
      {
        (x^i,y^i) \in X^i \times Y^i
      },
      {
        g^i \in G^i
      }
    }
    {
      \funcname{f}^i_C \bigl ( (x^i,y^i) \star g^i \bigr )
      =
      \funcname{f}^i_C \bigl ( (x^i,y^i) \bigr ) \star \funcname{f}^i_G (g^i)
    }
\end{equation}
Let $(x^1,y^1) \in X^1 \times Y^1$ and $g^1 \in G^1$. Then \fullcref{fig:Grho2} is
commutative and
\begin{equation}
\begin{split}
& \funcname{f}^2_C \compose \funcname{f}^1_C
  \bigl ( (x^1,y^1) \star g^1 \bigr ) = \\*
& \funcname{f}^2_C
  \Bigl (
    \funcname{f}^1_C \bigl ( (x^1,y^1) \bigr ) \star
    \funcname{f}^1_G (g^1)
  \Bigr ) = \\*
& \funcname{f}^2_C
  \Bigl (
    \funcname{f}^1_C \bigl ( (x^1,y^1) \bigr )
  \Bigr )
  \star\funcname{f}^2_G \compose \funcname{f}^1_G (g^1) = \\*
& \funcname{f}^2_C \compose \funcname{f}^1_C \bigl ( (x^1,y^1) \bigr )
  \star\funcname{f}^2_G \compose \funcname{f}^1_G (g^1)
\end{split}
\end{equation}
\end{proof}
\end{lemma}

\begin{definition}[Categories of $G$-$\rho$-model spaces]
\label{def:Grhomodcat}
Let $X^\alpha.Y^\alpha$, $\alpha \prec \Alpha$, be topological spaces,
$G^\alpha$ a topological group, $\rho^\alpha$ an effective group action
on $Y^\alpha$,
$
  \seqname{XY}^\alpha \defeq \\
  (X^\alpha \times Y^\alpha, \catname{XY}^\alpha)
$
a $G^\alpha$-$\rho^\alpha$ model space of $X^\alpha \times Y^\alpha$,
$
  \catseqname{XYG\rho}_\Ob \defeq \\
  \set
    {(X^\alpha, Y^\alpha, G^\alpha, \rho^\alpha)}%
    [\alpha \prec \Alpha]
$,
$
  \catseqname{XYG\rho}_\Ar \defeq
  \set
    {\funcname{f}_C \maps X^\alpha \times Y^\alpha \to X^\beta \times Y^\beta}%
    [
      {
        \mathrm{isG{\rho}morph}
        (
          \seqname{XY}^\alpha,
          G^\alpha,
          \rho^\alpha,
          \seqname{XY}^\beta,
          G^\beta,
          \rho^\beta,
          \funcname{f}_C
        )
      },
      \alpha \prec \Alpha,
      \beta \prec \Alpha
    ]*
$
and
$
  \catseqname{XYG\rho} \defeq \\
  (
    \catseqname{XYG\rho}_\Ob,
    \catseqname{XYG\rho}_\Ar
  )
$.
Then any subcategory of $\catseqname{XYG\rho}$ is a $G$-$\rho$ model category.
\end{definition}

\begin{definition}[Trivial $G$-$\rho$-model spaces]
\label{def:Grhomodtriv}
Let $X$ and $Y$ be topological spaces, $G$ a topological group, $\rho$
an effective group action of $G$ on $Y$ and $\catname{XY}$ the category
of all products of open subsets of $X$ with $Y$ and all homeomorphisms
induced by the group action, i.e,

\begin{equation}
\Ob(\catname{XY})
\defeq
\set
  {V \times Y}%
  [V \in \op{X}]
\end{equation}
\begin{multline}
\Ar(\catname{XY})
\defeq
\set
  {
   \funcname{f} \maps V \times Y \toiso V \times Y
  }%
  [
    {V \in \op{X}},
    {\pi_1 \compose \funcname{f} = \pi_1},
    {
      \equant
        {
         \funcname{g} \in G^V,
        }
        {
          \uquant
            {
              {x \in V},
              {y \in Y}
            }
            {
               f(x, y) = \bigl ( x, y \star \funcname{g}(x) \bigr )
            }
        }
    }
  ]*
\end{multline}
Then the trivial $G$-$\rho$ model space of $X,Y$, abbreviated
$\Triv[G-\rho-]{X,Y}$, is $(X \times Y, \catname{XY})$ and $\Triv[G-\rho-]{X,Y}$ is
a trivial $G$-$\rho$ model space of $X,Y$.

\begin{remark}
Let $G'$ be a topological group and $\rho'$ an effective group action of
$G'$ on $Y$ such that $\Triv[G-\rho-]{X,Y} = \Triv[G'-\rho'-]{X,Y}$.
Although $G'$ must be isomorphic to $G$. it need not have the same
topology.
\end{remark}

The identity morphism of $\Triv[G-\rho-]{X,Y}$ is
$\Id_{\Triv[G-\rho-]{X,Y}} \defeq \Id_{X \times Y}$.

Let
$
  \seqname{B}^\alpha \defeq
  (E^\alpha, X^\alpha, Y^\alpha, G^\alpha, \pi^\alpha, \rho^\alpha)
$,
$\alpha \prec \Alpha$, be a protobundle and
$\seqname{B} \defeq \set{{\seqname{B}^\alpha}}[\alpha \prec \Alpha]$
be a set of protobundles.

The trivial coordinate model category of $\seqname{B}$,
$\Trivcat[\Bun-]{\seqname{B}}$, is the category with objects all trivial
$G_\alpha$-$\rho_\alpha$ model spaces of $X_\alpha, Y_\alpha$,
$\alpha \prec \Alpha$ and morphisms all continuous functions compatible
with the group action:

\begin{equation}
\Trivcat[\Bun-]{\seqname{B}}_\Ob \defeq
\set
  {
    \Triv[G-\rho-]{X,Y}
  }%
  [
    {
      (
        E,
        X,
        Y,
        G,
        \pi,
        \rho
      )
      \in \seqname{B}
    }
  ]
\end{equation}
\begin{multline}
\Trivcat[\Bun-]{\seqname{B}}_\Ar \defeq
\set
  {
    \funcname{f}_C \maps
    X^1 \times Y^1 \to
    X^2 \times Y^2
  }%
  [
    {
      \equant
        {
          {
            (
              E^i,
              X^i,
              Y^i,
              G^i,
              \pi^i,
              \rho^i
            )
            \in \seqname{B}
          },
          \funcname{f}_G \maps G^1 \to G^2
        }
        {
          \mathrm{isG{\rho}morph}
          (
            \Triv[G-\rho-]{X^1,Y^1},
            G^1,
            \rho^1,
            \Triv[G-\rho-]{X^2,Y^2},
            G^2,
            \rho^2,
            \funcname{f}^C
          )
        }
    }
  ]*
\end{multline}
\begin{equation}
\Trivcat[\Bun-]{\seqname{B}} \defeq
\bigl (
  \Trivcat[\Bun-]{\seqname{B}}_\Ob,
  \Trivcat[\Bun-]{\seqname{B}}_\Ar
\bigr )
\end{equation}

\begin{remark}
The morphisms $\funcname{f}_C$ may be twisted products:  there need not
exist $\funcname{f}_Y \maps Y^1 \to Y^2$ such that
$\funcname{f}_C = \funcname{f}_X \times \funcname{f}_Y$.
\end{remark}

The trivial product coordinate model category of $\seqname{B}$,
$\Trivcat[\BunProd-]{\seqname{B}}$, is the category with objects all
trivial $G^\alpha$-$\rho^\alpha$ model spaces of $X^\alpha, Y^\alpha$,
$\alpha \prec \Alpha$ and morphisms all products of continuous functions
compatible with the group action:

\begin{equation}
\Trivcat[\BunProd-]{\seqname{B}}_\Ob \defeq \Trivcat[\Bun-]{\seqname{B}}_\Ob
\end{equation}
\begin{multline}
\Trivcat[\BunProd-]{\seqname{B}}_\Ar \defeq
\set
  {
    \funcname{f}_C \defeq  \funcname{f}_X \times \funcname{f}_Y \maps
    X^1 \times Y^1 \to
    X^2 \times Y^2
  }%
  [
    {
      \equant
        {
          {
            (
              E^i,
              X^i,
              Y^i,
              G^i,
              \pi^i,
              \rho^i
            )
            \in \seqname{B}
          },
          \funcname{f}_G \maps G^1 \to G^2
        }
        {
          \mathrm{isG{\rho}morph}
          (
            \Triv[G-\rho-]{X^1,Y^1},
            G^1,
            \rho^1,
            \Triv[G-\rho-]{X^2,Y^2},
            G^2,
            \rho^2,
            \funcname{f}_C
          )
        }
    }
  ]*
\end{multline}
\begin{equation}
\Trivcat[\BunProd-]{\seqname{B}} \defeq
\Bigl (
  \Trivcat[\BunProd-]{\seqname{B}}_\Ob,
  \Trivcat[\BunProd-]{\seqname{B}}_\Ar
\Bigr )
\end{equation}

Any subcategory of $\Trivcat[\Bun-]{\seqname{B}}$ is a trivial coordinate
model category and any subcategory of $\Trivcat[\BunProd-]{\seqname{B}}$ is
a trivial product coordinate model category.
\end{definition}

\begin{lemma}[The trivial coordinate model category of $\seqname{B}$ is a category]
Let
$
  \seqname{B}^\alpha \defeq
  (E^\alpha, X^\alpha, Y^\alpha, G^\alpha, \pi^\alpha, \rho^\alpha)
$,
$\alpha \prec \Alpha$, be a protobundle and
$\seqname{B} \defeq \set{{\seqname{B}^\alpha}}[\alpha \prec \Alpha]$
be a set of protobundles.

Then $\Trivcat[\Bun-]{\seqname{B}}$ is a category and the identity
morphism for object $\seqname{B}^\alpha$ is
$\Id_{X^\alpha \times Y^\alpha}$.

\begin{proof}
$\Trivcat[\Bun-]{\seqname{B}}$ satisfies the definition of a category:
\begin{enumerate}
\item Composition: \newline
The composition of $G$-$\rho$ morphisms is a $G$-$\rho$ morphism by
\pagecref{lem:Grhomorph}.
\item Associativity: \newline
Morphisms are simply functions and composition of morphisms is simply
composition of functions.
\item Unit: \newline
The identity morphisms are simply identity functions and composition of
morphisms is simply composition of functions.
\end{enumerate}
\end{proof}

$\Trivcat[\BunProd-]{\seqname{B}}$ is a category and the identity
morphism for object $\seqname{B}^\alpha$ is
$\Id_{X^\alpha \times Y^\alpha}$.

\begin{proof}
$\Trivcat[\BunProd-]{\seqname{B}}$ satisfies the definition of a category:
\begin{enumerate}
\item Composition: \newline
The composition of $G$-$\rho$ morphisms is a $G$-$\rho$ morphism by
\pagecref{lem:Grhomorph}. Let
$\Triv[G^i-\rho^i-]{X^i,Y^i} \objin\Trivcat[\BunProd-]{\seqname{B}},$
$i=1,2,3$,
$
  \funcname{f}^i_C =
  \funcname{f}^i_X \times \funcname{f}^i_Y \maps
  X^i \times Y^i \to X^{i+1} \times Y^{i+1}
  \arin  \Trivcat[\BunProd-]{\seqname{B}}
$.
Then
$
  \funcname{f}^{i+1}_C \compose \funcname{f}^i_C =
  (\funcname{f}^{i+1}_X \compose \funcname{f}^i_X)
  \times
  (\funcname{f}^{i+1}_Y \compose \funcname{f}^i_Y)
$.
\item Associativity: \newline
Morphisms are simply functions and composition of morphisms is simply
composition of functions.
\item Unit: \newline
The identity morphisms are simply identity functions and composition of
morphisms is simply composition of functions.
\end{enumerate}
\end{proof}
\end{lemma}

\begin{lemma}[The trivial $G$-$\rho$ model space of $X,Y$ is a $G$-$\rho$ model space of $X,Y$]
Let $X$ and $Y$ be topological spaces, $G$ a topological group and
$\rho$ an effective group action on $Y$. Then $\Triv[G-\rho-]{X,Y}$ is
a $G$-$\rho$ model space of $X \times Y$.

\begin{proof}
$\Triv[G-\rho-]{X,Y}$ satisfies the conditions of \pagecref{def:Grhomod}\!:
\begin{enumerate}
\item Product with fiber:
\[
\uquant
  {
    U \objin \pi_2(\Triv[G-\rho-]{X,Y})
  }
  {
    \equant
      {
        V \in \op{X}
      }
      {
        U = V \times Y
      }
  }
\]
by \cref{def:Grhomodtriv}.
\item Generated by $\rho$:
Let $\funcname{f} \arin \pi_2(\Triv[G-\rho-]{X,Y}) \maps V \times Y \toiso V \times Y$.
Then
\[
  \equant
    {
      \funcname{g} \in G^V
    }
    {
    \uquant
      {
        {(x,y) \in V \times Y}
      }
      {
        \funcname{f}(x,y) =
          \bigl ( x, y \star \funcname{g}(x) \bigr )
      }
    }
\]
by \cref{def:Grhomodtriv}.
\end{enumerate}
\end{proof}
\end{lemma}

\subsection{$G$-$\rho$-nearly commutative diagrams}
\label{sub:grhoncd}
Let $X$ and $Y$ be topological spaces, $G$ a topological group, $\rho$
an effective group action on $Y$,
$\seqname{C}=(C,\catname{C}) \defeq \Triv[G-\rho-]{X,Y}$, and $D$ a
tree with two branches, whose nodes are topological spaces $U_i$ and
$V^j$ and whose links are continuous functions
$\funcname{f}_i \maps U_i \to U_{i+1}$ and
$\funcname{f}'_j \maps U_j \to U_{j+1}$ between the spaces:

\begin{align*}
D = \{
  &
    \funcname{f}_0 \maps U_0 = V_0 \to U_1,
    \dotsc,
    \funcname{f}_{m - 1} \maps U_{m - 1} \to U_m,
\\
  &
    \funcname{f}'_0 \maps  U_0 = V_0 \to V_1,
    \dotsc,
    \funcname{f}'_{m-1} \maps V_{m-1} \to V_n
\}
\end{align*}
with $U_0 = V_0$, $U_m \subseteq C$ open and
$V_n \subseteq C$ open, as shown in \pagecref{fig:NCDb}.

\begin{definition}[$G$-$\rho$-nearly commutative diagrams]
$D$ is nearly commutative in $X,Y,\pi,\rho$ iff $D$ is nearly
commutative in category $\catname{C}$.
\end{definition}

\begin{definition}[$G$-$\rho$-nearly commutative diagrams at a point]
\label{def:grhoNCD}
Let $\catname{C}$ and $D$ be as above and $x$ be an element of the
initial node. $D$ is nearly commutative in $X,Y,\pi,\rho$ at $x$ iff
$D$ is nearly commutative in $\catname{C}$ at $x$.
\end{definition}

\begin{definition}[$G$-$\rho$-locally nearly commutative diagrams]
Let $\catname{C}$ and $D$ be as above. $D$ is locally nearly commutative
in $X,Y,\pi,\rho$ iff $D$ is locally nearly commutative in
$\catname{C}$.
\end{definition}

\subsection{Bundle charts}
\label{sub:buncharts}
\begin{definition}[$Y$-$\pi$-bundle charts]
\label{def:Ypichart}
Let $E,X,Y$ be topological spaces and \\
 $\pi \maps E \to X$.
A $Y$-$\pi$-bundle chart $(U, V \times Y, \phi)$ of $E$ in the
coordinate space $X \times Y$ consists of

\begin{enumerate}
\item An open set $U \subseteq E$, known as a coordinate patch
\item An open set $V \times Y \subseteq X \times Y$
\item a homeomorphism $\phi \maps U \toiso V \times Y$, known as a
coordinate function, that preserves fibers. i.e.,
$\pi_1 \compose \phi = \pi$.
\end{enumerate}
\end{definition}

\begin{lemma}[Properties of projection]
\label{def:piepi}
Let $(U, V \times Y, \phi)$ be a $Y$-$\pi$-bundle chart of $E$ in the
coordinate space $X \times Y$ and $v \in V$. Then

$\pi \restriction_U$ is a surjection.
\begin{proof}
Let $v \in V$, $y \in Y$ and
$u \defeq \phi^{-1} \bigl ( (v,y) \bigr ) \in U$. Then $\pi(u) = v$.
\end{proof}
$\pi^{-1}[\set{v}]$ is homeomorphic to $Y$.
\begin{proof}
$\phi$ and $\phi^{-1}$ are homeomorphisms, so their restrictions are
homeomorphisms and thus $\phi^{-1}[\set{v}]$ is homeomorphic to
$\set{v} \times Y$, which is homeomorphic to $Y$.
\end{proof}
\end{lemma}

\begin{definition}[$Y$-$\pi$ subcharts]
\label{def:Ypisub}
Let $(U, V \times Y, \phi)$ be a $Y$-$\pi$-bundle chart of $E$ in the
coordinate space $X \times Y$ and $U' \subseteq U$ open.
Then \\
$(U', V' \times Y, \phi') \defeq (U', \phi[U'], \phi \restriction_{U'})$
is a subchart of $(U, V \times Y, \phi)$.
\end{definition}

\begin{lemma}[$Y$-$\pi$ subcharts]
Let $(U, V \times Y, \phi)$ be a $Y$-$\pi$-bundle chart of $E$ in the
coordinate space $X \times Y$ and $(U', V' \times Y, \phi')$ a subchart
of $(U, V \times Y, \phi)$. Then $(U', V' \times Y, \phi')$ is a
$Y$-$\pi$-bundle chart of $E$ in the coordinate space $X \times Y$.
\begin{proof}
$(U', V' \times Y, \phi')$ satisfies the conditions of \cref{def:Ypichart}:
\begin{enumerate}
\item $U' \subseteq E$ is open.
\item $\phi[U']$ os open since $\phi$ is a homeomorphism.
\item $\phi \restriction_{U'} \maps U' \to V' \times Y$ is the
restriction of a homeomorphism and thus a homeomorphism.
$\phi \restriction_{U'}$ preserves fibers because $\phi$ does.
\end{enumerate}
\end{proof}
\end{lemma}

\begin{definition}[$G$-$\rho$-compatibility]
Let $E,X,Y$ be topological spaces, $G$ a topological group,
$\rho \maps Y \times G \to Y$ an effective right action of
$G$ on $Y$, $\pi \maps E \onto X$ surjective
and $(U^\mu, V^\mu \times Y, \phi^\mu)$,
$(U^\nu, V^\nu \times Y, \phi^\nu)$ $Y$-$\pi$-bundle charts.
$(U^\mu, V^\mu \times Y, \phi^\mu)$ and
$(U^\nu, V^\nu \times Y, \phi^\nu)$ are $G$-$\rho$-compatible if either
\begin{enumerate}
\item $U^\mu$ and $U^\nu$ are disjoint
\item
The transition function
$\funcname{t}^\mu_\nu=\phi^\mu \compose {\phi^\nu}^{-1} \restriction_{\phi^\nu[U^\mu \cap U^\nu]}$
is generated by the group action, i.e., there is a continuous function
$\funcname{g}^\mu_\nu \maps \pi_1[\phi^\nu[U^\mu \cap U^\nu]] \to G$
such that
\begin{equation}
\showlabelsinline
\label{eq:transitiongroup}
\uquant
  {{(x,y) \in \phi^\nu[U^\mu \cap U^\nu]}}
  {
    \funcname{t}^\mu_\nu(x,y) =
    \bigl ( x, y \star \funcname{g}^\mu_\nu(x) \bigr )
  }
\end{equation}
\end{enumerate}
\end{definition}

\begin{lemma}[Symmetry of $G$-$\rho$ compatibility]
Let $(U^\mu, V^\mu \times Y, \phi^\mu)$ and \\
$(U^\nu, V^\nu \times Y, \phi^\nu)$ be $Y$-$\pi$-bundle charts of $E$ in the
coordinate space $X \times Y$. Then
$(U^\mu, V^\mu \times Y, \phi^\mu)$ is $G$-$\rho$-compatible with
$(U^\nu, V^\nu \times Y, \phi^\nu)$ iff
$(U^\nu, V^\nu \times Y, \phi^\nu)$ is $G$-$\rho$-compatible with
$(U^\mu, V^\mu \times Y, \phi^\mu)$.
\begin{proof}
It suffices to prove the result in one direction. If
$
  (U^\mu, V^\mu \times Y, \phi^\mu)
  \cap
  (U^\nu, V^\nu \times Y, \phi^\nu)
$
then
$
  (U^\nu, V^\nu \times Y, \phi^\nu)
  \cap
  (U^\mu, V^\mu \times Y, \phi^\mu)
$.
Otherwise, let
$\funcname{g}^\mu_\nu \maps \pi_1[\phi^\nu[U^\mu \cap U^\nu]] \to G$ be
a continuous function such that
\[
\uquant
  {{(x,y) \in \phi^\nu[U^\mu \cap U^\nu]}}
  {
    \funcname{t}^\mu_\nu(x,y) =
    \bigl ( x, y \star \funcname{g}^\mu_\nu(x) \bigr )
  }
\]
and the inverse transition function is also generated by the group action:
\[
\uquant
  {{(x,y) \in \phi^\nu[U^\mu \cap U^\nu]}}
  {
    \funcname{t}^\nu_\mu(x,y) =
    \bigl ( x, y \star {\funcname{g}^\mu_\nu(x)}^{-1} \bigr )
  }
\]
\end{proof}
\end{lemma}

\begin{lemma}[$G$-$\rho$-compatibility of subcharts]
Let $(U^i, V^i \times Y, \phi^i)$ be a $Y$-$\pi$-bundle chart of $E$ in
the coordinate space $X \times Y$, $(U'^i, V'^i \times Y, \phi'^i)$ a
subchart and $(U^1, V^1, \phi^1)$ be $G$-$\rho$-compatible with
$(U^2, V^2, \phi^2)$.  Then $(U'^1, V'^1 \times Y, \phi'^1)$ is
$G$-$\rho$-compatible with $(U'^2, V'^2 \times Y, \phi'^2)$.

\begin{proof}
If $U^1 \cap U^2 = \emptyset$ then $U'^1 \cap U'^2 = \emptyset$. If
$U'^1 \cap U'^2 = \emptyset$ then $(U'^1, V'^1 \times Y, \phi'^1)$ is
$G$-$\rho$-compatible with $(U'^2, V'^2 \times Y, \phi'^2)$. Otherwise,
the transition function
$
  t^1_2 \defeq
  \phi^1 \compose \phi^{2-1} \restriction_{\phi^2[U^1 \cap U^2]}
$
is generated by the group action and hence
$
  t^1_2 \restriction_{\phi^2[U'^1 \cap U'^2]}
  \maps
  \phi^2[U'^1 \cap U'^2]
  \toiso
  \phi^1[U'^1 \cap U'^2]
$
is generated by the group action.
\end{proof}
\end{lemma}

\begin{corollary}[$G$-$\rho$-compatibility with subcharts]
Let $(U, V \times Y, \phi)$, be a $Y$-$\pi$-bundle chart of $E$ in the
coordinate space $X \times Y$ and $(U', V' \times Y, \phi')$ a subchart.
Then $(U', V' \times Y, \phi')$ is $G$-$\rho$-compatible with
$(U, V \times Y, \phi)$.

\begin{proof}
$(U, V \times Y, \phi)$ is $G$-$\rho$-compatible with itself and is a
subchart of itself,
\end{proof}
\end{corollary}

\begin{definition}[Covering by $Y$-$\pi$-bundle charts]
Let $\seqname{A}$ be a set of $Y$-$\pi$-bundle
charts of $E$ in the coordinate space $X \times Y$.
$\seqname{A}$ covers $E$ iff $E = \union{{\pi_1[\seqname{A}]}}$.
\end{definition}

\begin{lemma}
Let $\seqname{A}$ be a set of $Y$-$\pi$-bundle charts of $E$ in the
coordinate space $X \times Y$ that covers $E$ and $x \in X$. Then
$\pi^{-1}[\set{x}]$ is homeomorphic to $Y$.
\begin{proof}
Since $\seqname{A}$ covers $E$, there is a chart $(U, V, \phi)$ in
$\seqname{A}$ containing $x$.
Then $\pi^{-1}[\set{x}]$ is homeomorphic to $Y$ by \pagecref{def:piepi}.
\end{proof}
\end{lemma}

\subsection{Bundle atlases}
\label{sub:bunatlases}
A set of charts can be atlases for different fiber bundles even if it is
for the same total model space, base space and fiber. In order to
aggregate atlases into categories, there must be a way to distinguish
them. Including the spaces\footnote{
The spaces are redundant, but convenient.},
group and group action in the definitions of the categories serves the
purpose.

\begin{definition}[Bundle atlases]
\label{def:BunAtl}
Let $\seqname{B} \defeq (E, X, Y, G, \pi, \rho)$, be a protobundle. Then
$\seqname{A}$ is a bundle atlas of $B$, abbreviated
$\isAtl^\Bun_\Ob(\seqname{A}, \seqname{B})$ and $\seqname{A}$ is a
$\pi$-$G$-$\rho$-bundle atlas of $E$ in the coordinate space
$X \times Y$, abbreviated
$\isAtl^\Bun_\Ob(\seqname{A}, E, X, Y, G, \pi, \rho)$, iff it consists
of a set of mutually $G$-$\rho$-compatible $Y$-$\pi$-bundle charts of
$E$ in the coordinate space $X \times Y$ that covers $E$\footnote{
  There is no need to introduce the concept of a full
  $\pi$-$G$-$\rho$-bundle atlas because a $\pi$-$G$-$\rho$-bundle
  atlas is automatically full.
}.

By abuse of language we write $U \in \seqname{A}$ for
$U \in \pi_1[\seqname{A}]$.

\begin{remark}
The definition of a $\pi$-$G$-$\rho$-bundle atlas is by design similar
to the definition of a coordinate bundle in \cite[p.~7]{TopFib}, but
there are significant differences. This paper will use the term bundle
atlas to avoid confusion.
\end{remark}

Let
$
  \seqname{B}^\alpha \defeq
  (E^\alpha, X^\alpha, Y^\alpha, G^\alpha, \pi^\alpha, \rho^\alpha)
$,
$\alpha \prec \Alpha$, be a protobundle and
$\seqname{B} \defeq \set{{\seqname{B}^\alpha}}[\alpha \prec \Alpha]$
be a set of protobundles.
Then
\begin{multline}
\Atl^\Bun_\Ob(\seqname{B}^\alpha) \defeq
  \set
  {{
    (
      \seqname{A},
      \seqname{B}^\alpha
    )
  }}%
  [
    {{
      \isAtl^\Bun_\Ob
        (
          \seqname{A},
          E^\alpha,
          X^\alpha,
          Y^\alpha,
          G^\alpha,
          \pi^\alpha,
          \rho^\alpha
        )
    }}
  ]
\end{multline}

\begin{equation}
\Atl^\Bun_\Ob \seqname{B} \defeq
  \union[\alpha \prec \Alpha]{\Atl^\Bun_\Ob(\seqname{B}^\alpha)}
\end{equation}
\end{definition}

\begin{lemma}[Bundle atlases]
Let $E,X,Y$ be topological spaces, $G$ a topological group,
$\pi \maps E \onto X$ surjective and $\rho \maps Y \times G \to Y$ an
effective right action of $G$ on $Y$. Then $\seqname{A}$ is a
$\pi$-$G$-$\rho$-bundle atlas of $E$ in the coordinate space
$X \times Y$ iff it is an m-atlas of $\Triv{E}$ in the coordinate space
$\Triv[G-\rho-]{X,Y}$ and every coordinate function preserves fibers.

\begin{proof}
If $\seqname{A}$ is a $\pi$-$G$-$\rho$-bundle atlas of $E$ in the
coordinate space $X \times Y$ then
\begin{enumerate}
\item Every chart in $\seqname{A}$ is a $Y$-$\pi$-bundle chart of $E$ in
the coordinate space $X \times Y$, and hence its coordinate function
preserves fibers.
\item The charts in $\seqname{A}$ are mutually $G$-$\rho$-compatible;
hence the transition functions are generated by the group action and are
morphisms of $\Triv[G-\rho-]{X,Y}$.
\end{enumerate}

If $\seqname{A}$ is an m-atlas of $\Triv{E}$ in the coordinate space
$\Triv[G-\rho-]{X,Y}$ then the transition functions are generated by the
group action and and thus the charts are mutually $G$-$\rho$-compatible.

If every coordinate function preserves fibers then the m-charts of
$\seqname{A}$ are $Y$-$\pi$-bundle charts.
\end{proof}
\end{lemma}

\begin{definition}[Compatibility of charts with bundle atlases]
A $Y$-$\pi$-bundle chart $(U, V \times Y, \phi)$ of $E$ in the
coordinate space $X \times Y$ is $G$-$\rho$-compatible with a
$\pi$-$G$-$\rho$-bundle atlas $\seqname{A}$ iff it is
$G$-$\rho$-compatible with every chart in the atlas.
\end{definition}

\begin{lemma}[Compatibility of subcharts with bundle atlases]
\label{lem:BunCompat}
Let $\seqname{A}$ be a $\pi$-$G$-$\rho$-bundle atlas of $E$ in the
coordinate space $X \times Y$ and
$\seqname{C}^1 = (U, V \times Y, \phi)$ a $Y$-$\pi$-bundle chart in
$\seqname{A}$. Then any subchart of $\seqname{C}^1$ is
$G$-$\rho$-compatible with $\seqname{A}$.

\begin{proof}
Let $\seqname{C}'^1 = (U'^1, V'^1 \times Y, \phi'^1)$ be a subchart of
$\seqname{C}^1$ and $\seqname{C}^2 = (U^2, V^2 \times Y, \phi^2)$
another chart in $\seqname{A}$.

\begin{enumerate}
\item If $U^1 \cap U^2 = \emptyset$, then $U' \cap U^2 = \emptyset$.
\item If $U' \cap U^2 = \emptyset$ then $\seqname{C}'^1$ is
$G$-$\rho$-compatible with $\seqname{C}^2$.
\item Otherwise the transition function
$
  t^1_2 \defeq
  \phi^1 \compose \phi^{2-1} \restriction_{\phi^1[U^1 \cap U^2]}
$
is generated by the group action and thus
$t^1_2 \restriction_{\phi^2[U'^1 \cap U^2]}$ is generated by the group
action.
\end{enumerate}
\end{proof}
\end{lemma}

\begin{lemma}[Extensions of bundle atlases]
\label{def:piGrho-atl:extensions}
Let $\seqname{A}$ be a $\pi$-$G$-$\rho$-atlas of $\seqname{E}$ in the
coordinate space $X \times Y$ and $(U_i,V_i,\phi_i)$, $i=1,2$ be
$\pi$-$G$-$\rho$-charts of $E$ in the coordinate space $X \times Y$
$G$-$\rho$-compatible with $\seqname{A}$ in the coordinate space $X
\times Y$.  Then $(U_1,V_1,\phi_1)$ is $G$-$\rho$-compatible with
$(U_2,V_2,\phi_2)$ in the coordinate space $X \times Y$.

\begin{proof}
If $U_1 \cap U_2 = \emptyset$ then $(U_1,V_1,\phi_1)$ is
$G$-$\rho$-compatible with $(U_2,V_2,\phi_2)$. Otherwise,
$
  \phi_2 \compose \phi^{-1}_1 \restriction_{\phi_1[U_1 \cap U_2]}
  \maps
  \phi_1[U_1 \cap U_2] \toiso
  \phi_2[U_1 \cap U_2]
$
is a homeomorphism.  It remains to show that
$\phi_2 \compose \phi^{-1}_1 \restriction_{\phi_1[U_1 \cap U_2]}$ is
generated by the group action.  Let $(U'_\alpha,V'_\alpha,\phi'_\alpha)$,
$\alpha \prec \Alpha$, be charts in $\seqname{A}$ such that
$U_1 \cap U_2 \subseteq \union[\alpha \prec \Alpha]{U'_\alpha}$ and
$U_1 \cap U_2 \cap U'_\alpha \neq \emptyset$, $\alpha \prec \Alpha$.
Since the charts are $G$-$\rho$-compatible with
$(U'_\alpha,V'_\alpha,\phi'_\alpha)$,
$
\phi_2 \compose \phi'^{-1}_\alpha \restriction_{U_1 \cap U_2 \cap
U'_\alpha}
$
and
$\phi'_\alpha \compose \phi^{-1}_1 \restriction_{U^1 \cap U^2 \cap
U'_\alpha}
$
are generated by the group action and thus
$\phi^2 \compose \phi^{-1}_1 = \phi^2 \compose \phi'^{-1}_\alpha \compose
\phi'_\alpha \compose \phi^{-1}_1$
is generated by the group action.
\end{proof}
\end{lemma}

\begin{definition}[Maximal bundle atlases]
Let $\seqname{A}$ be a $\pi$-$G$-$\rho$-bundle atlas of $E$ in the
coordinate space $X \times Y$. $\seqname{A}$ is a maximal
$\pi$-$G$-$\rho$-bundle atlas, abbreviated
$
  \maximal{\isAtl^\Bun_\Ob}
    (
      \seqname{A},
      E,
      X,
      Y,
      G,
      \pi,
      \rho
    )
$,
iff it cannot be extended by adding an additional $G$-$\rho$-compatible
$Y$-$\pi$-bundle chart.

$\seqname{A}$ is a semi-maximal $\pi$-$G$-$\rho$-bundle atlas
of $\seqname{E}$ in the coordinate space $\seqname{C}$, abbreviated
$
  \maximal[S-]{\isAtl^\Bun_\Ob}
    (
      \seqname{A},
      E,
      X,
      Y,
      G,
      \pi,
      \rho
    )
$,
iff whenever $(U, V \times Y, \phi) \in \seqname{A}$,
$U' \subseteq U, V' \times Y \subseteq V \times Y$ and
$V'' \times Y \subseteq X \times Y$ are open sets,
$\phi[U'] = V' \times Y$,
$\phi' \maps V' \times Y \toiso V'' \times Y$ is a fiber
preserving homeomorphism generated by the group action then
$(U', V'' \times Y, \phi' \compose \phi) \in \seqname{A}$.

Let
$
  \seqname{B}^\alpha \defeq
  (E^\alpha, X^\alpha, Y^\alpha, G^\alpha, \pi^\alpha, \rho^\alpha)
$,
$\alpha \prec \Alpha$, be a protobundle and
$\seqname{B} \defeq \set{{\seqname{B}^\alpha}}[\alpha \prec \Alpha]$
be a set of protobundles.
Then
\begin{multline}
\maximal{\Atl^\Bun_\Ob}(\seqname{B}^\alpha) \defeq
  \set
  {{
    (
      \seqname{A},
      \seqname{B}^\alpha
    )
  }}%
  [
    {{
      \maximal{\isAtl^\Bun_\Ob}
        (
          \seqname{A},
          E^\alpha,
          X^\alpha,
          Y^\alpha,
          G^\alpha,
          \pi^\alpha,
          \rho^\alpha
        )
    }}
  ]
\end{multline}

\begin{equation}
\maximal{\Atl^\Bun_\Ob} \seqname{B} \defeq
  \union[\alpha \prec \Alpha]{\maximal{\Atl^\Bun_\Ob}(\seqname{B}^\alpha)}
\end{equation}

\begin{multline}
\maximal[S-]{\Atl^\Bun_\Ob}(\seqname{B}^\alpha) \defeq
  \set
  {{
    (
      \seqname{A},
      \seqname{B}^\alpha
    )
  }}%
  [
    {{
      \maximal[S-]{\isAtl^\Bun_\Ob}
        (
          \seqname{A},
          E^\alpha,
          X^\alpha,
          Y^\alpha,
          G^\alpha,
          \pi^\alpha,
          \rho^\alpha
        )
    }}
  ]
\end{multline}

\begin{equation}
\maximal[S-]{\Atl^\Bun_\Ob} \seqname{B} \defeq
  \union[\alpha \prec \Alpha]{\maximal[S-]{\Atl^\Bun_\Ob}(\seqname{B}^\alpha)}
\end{equation}
\end{definition}

\begin{lemma}[Maximal $\pi$-$G$-$\rho$-bundle atlases are semi-maximal $\pi$-$G$-$\rho$-bundle atlases]
Let $E,X,Y$ be topological spaces, $G$ a topological group,
$\pi \maps E \onto X$ surjective, $\rho \maps Y \times G \to Y$ an
effective right action of $G$ on $Y$ and $\seqname{A}$ a maximal
$\pi$-$G$-$\rho$-bundle atlas of $E$ in the coordinate space $X \times Y$.
Then $\seqname{A}$ is a semi-maximal $\pi$-$G$-$\rho$-bundle atlas of $E$ in the
coordinate space $X \times Y$.

\begin{proof}
Let $(U, V \times Y, \phi) \in \seqname{A}$,
$U' \subseteq U, V' \times Y \subseteq V \times Y$ and
$V'' \times Y \subseteq X \times Y$ be open sets and
$\phi[U'] = V' \times Y$, $\phi' \maps V' \times Y \toiso V'' \times Y$
be a fiber preserving homeomorphism generated by the group action.
$(U', V', \phi)$ is a subchart of $(U,V,\phi)$ and by
\pagecref{lem:BunCompat} is $G$-$\rho$-compatible with the charts of
$\seqname{A}$.  Since $\phi'$ is a fiber preserving homeomorphism
generated by the group action, $(U', V'' \times Y, \phi' \compose \phi)$
is $G$-$\rho$-compatible with the charts of $\seqname{A}$. Since
$\seqname{A}$ is maximal, $(U', V'', \phi' \compose \phi)$ is a chart of
$\seqname{A}$.
\end{proof}
\end{lemma}

\begin{theorem}[Existence and uniqueness of maximal $\pi$-$G$-$\rho$-bundle atlases]
Let $\seqname{B} \defeq (E, X, Y, G, \pi, \rho)$ be a protobundle and
$\seqname{A}$ $\pi$-$G$-$\rho$-bundle atlas of $E$ in the coordinate
space $X \times Y$. Then there exists a unique maximal
$\pi$-$G$-$\rho$-bundle atlas
$\maximal{\Atlas^\Bun}(\seqname{A},\seqname{B})$ of $E$ in the
coordinate space $X \times Y$ $G$-$\rho$-compatible with $\seqname{A}$.

\begin{proof}
Let $\seqname{P}$ be the set of all $\pi$-$G$-$\rho$-bundle atlases $E$
in the coordinate space $X \times Y$ containing $\seqname{A}$ and
$G$-$\rho$ compatible in the coordinate space $X \times Y$ with
$\seqname{A}$. Let $\maximal{\seqname{P}}$ be a maximal chain of
$\seqname{P}$. Then $A'=\union{\maximal{\seqname{P}}}$ is a maximal
$\pi$-$G$-$\rho$-bundle atlas of $\seqname{E}$ in the coordinate space
$X \times Y$ $G$-$\rho$ compatible with $\seqname{A}$. Uniqueness
follows from \pagecref{def:piGrho-atl:extensions}.
\end{proof}
\end{theorem}

\begin{lemma}[Existence and uniqueness of projection for atlases in $G$-$\rho$-model spaces]
\label{lem:mGrho}
Let $E$, $X$ and $Y$ be topological spaces, $G$ a topological group,
$\rho$ an effective action of $G$ on $Y$,
$\seqname{C} \defeq (X \times Y, \catname{XY})$ a $G$-$\rho$ model space
of $X \times Y$ and $\seqname{A}$ an m-atlas of $E$ in the coordinate
model space $\seqname{C}$.  Then there exists a unique function
$\pi \maps \seqname{E} \to X$ such that for any chart $(U, V, \phi)$ in
$\seqname{A}$, $\pi \restriction_U = \pi_1 \compose \phi$.  If
$\seqname{A}$ is full then $\pi$ is surjective.

\begin{proof}
Let $(U, V, \phi)$ be an arbitrary chart in $\seqname{A}$ and define
$\pi(e \in U) = \pi_1 \compose \phi(e)$. $\pi(e)$ does not depend on the
choice of chart because the morphisms of a $G$-$\rho$ model space
preserve fibers. $\pi$ is continuous because it is continuous on each
coordinate patch.

Let $x \in X$. If $\seqname{A}$ is full then there exists a chart
$(U, V, \phi)$ in $\seqname{A}$ such that $x \in \pi_1[V]$. Let $u$ be an
arbitrary point in $\phi^{-1}[\set{x} \times Y]$. Then
$\pi(u) = \pi_1(\phi(u)) = x$.
\end{proof}
\end{lemma}
\subsection{Bundle atlas morphisms and functors}
\label{sub:bunatlmorph}
This section defines categories of bundle atlases,
$\Atl^\Bun \seqname{B}$. and constructs functors from them to categories
of m-atlases, $\Trivcat[\Bun-]{\seqname{B}}$ and
$\Trivcat[\BunProd-]{\seqname{B}}$.
It only constructs reverse functors for
$\Trivcat[\BunProd-]{\seqname{B}}$.

\begin{definition}[Bundle-atlas morphisms]
\label{def:Bun-ATLmorph}
Let
$
  \seqname{B}^i
  \defeq
  (
    E^i,
    X^i,
    Y^i,
    G^i,
    \pi^i,
    \rho^i
  )
$,
$i=1,2$, be a protobundle and let $\seqname{A}^i$ be a
$\pi^i$-$G^i$-$\rho^i$-atlas of $E^i$ in the coordinate space
$C^i = X^i \times Y^i$. Then
$
  \funcseqname{f}
  \defeq
  (
    \funcname{f}_E \maps E^1 \to E^2,
    \funcname{f}_X \maps X^1 \to X^2,
    \funcname{f}_Y \maps Y^1 \to Y^2,
    \funcname{f}_G \maps G^1 \to G^2
  )
$
is a $\seqname{B^1}$-$\seqname{B^2}$ bundle-atlas morphism from
$\seqname{A}^1$ to $\seqname{A}^2$, abbreviated \\
$
  \isAtl^\Bun_\Ar
    (
      \seqname{A^1}, \seqname{B^1},
      \seqname{A^2}, \seqname{B^2},
      \funcseqname{f}
    )
$,
iff
\begin{enumerate}
\item all four functions are continuous
\item $\funcname{f}_G$ is a homomorphism
\item $\funcseqname{f}$ commutes with $\pi^i$ and $\rho^i$, i.e.,
\begin{enumerate}
\item $\pi^2 \compose \funcname{f}_E = \funcname{f}_X \compose \pi^1$
\item
$
  \uquant
    {{y \in Y^1},{g \in G^1}}
    {{
      \funcname{f}_Y(y) \star^2 \funcname{f}_G(g) =
      \funcname{f}_Y(y \star^1 g)
    }}
$
\end{enumerate}
\item for any
$(U^1, V^1, \phi^1 \maps U^1 \toiso V^1) \in \seqname{A}^1$,
$(U^2, V^2, \phi^2 \maps U^2 \toiso V^2) \in \seqname{A}^2$,
the diagram
$D \defeq (\{I \defeq U^1 \cap \funcname{f}_0^{-1}[U^2], V^1, E^2, U^2, V^2 \}$,
$\{ \funcname{f}_0, \phi^2, \phi^1, \funcname{f}_1 \})$
is locally nearly commutative in $X,Y,\pi,\rho$.
\end{enumerate}

If $\seqname{A}^1$ and $\seqname{A}^2$ are maximal atlases then
$\funcseqname{f}$ is also a maximal $\seqname{B^1}$-$\seqname{B^2}$
bundle-atlas morphism from $\seqname{A}^1$ to $\seqname{A}^2$,
abbreviated
$
  \maximal{\isAtl^\Bun_\Ar}
    (
      \seqname{A^1}, \seqname{B^1},
      \seqname{A^2}, \seqname{B^2},
      \funcseqname{f}
    )
$

The identity morphism of $(\seqname{A}^i, \seqname{B}^i)$ is
\begin{equation}
\Id_{(\seqname{A}^i, \seqname{B}^i)} \defeq
\Bigl (
\bigl (
  \Id_{E^i},
  \Id_{X^i},
  \Id_{Y^i},
  \Id_{ G^i}
\bigr )
\bigl ( \seqname{A}^i, \seqname{B}^i\bigr ),
\bigl ( \seqname{A}^i, \seqname{B}^i\bigr )
\Bigr )
\end{equation}
This nomenclature will be justified later.

Let
$
  \seqname{B}^i \defeq
  (
    E^i,
    X^i,
    Y^i,
    G^i,
    \pi^i,
    \rho^i
  )
$,
$i=1,2$, be a protobundle. Then

\begin{equation}
\Atl^\Bun_\Ar
  (
    \seqname{B}^1,
    \seqname{B}^2
  )
\defeq
\set
  {
    \bigl (
      \funcseqname{f},
      (\seqname{A}^1, \seqname{B}^1),
      (\seqname{A}^2, \seqname{B}^2)
    \bigr )
  }%
  [
    {
      \isAtl^\Bun_\Ar
        (
          \seqname{A}^1, \seqname{B}^1,
          \seqname{A}^2, \seqname{B}^2,
          \funcseqname{f}
        )
    }
  ]
\end{equation}

\begin{equation}
\maximal{\Atl^\Bun_\Ar}
  (
    \seqname{B}^1,
    \seqname{B}^2
  )
\defeq
\set
  {
    \bigl (
      \funcseqname{f},
      (\seqname{A}^1, \seqname{B}^1),
      (\seqname{A}^2, \seqname{B}^2)
    \bigr )
  }%
  [
    {
      \maximal{\isAtl^\Bun_\Ar}
        (
          \seqname{A}^1, \seqname{B}^1,
          \seqname{A}^2, \seqname{B}^2,
          \funcseqname{f}
        )
    }
  ]
\end{equation}

\begin{equation}
\maximal[S-]{\Atl^\Bun_\Ar}
  (
    \seqname{B}^1,
    \seqname{B}^2
  )
\defeq
\set
  {
    \bigl (
      \funcseqname{f},
      (\seqname{A}^1, \seqname{B}^1),
      (\seqname{A}^2, \seqname{B}^2)
    \bigr )
  }%
  [
    {
      \maximal[S-]{\isAtl^\Bun_\Ar}
        (
          \seqname{A}^1, \seqname{B}^1,
          \seqname{A}^2, \seqname{B}^2,
          \funcseqname{f}
        )
    }
  ]
\end{equation}

\begin{equation}
\Atl^\Bun
  (
    \seqname{B}^i
  )
\defeq
  \bigl (
    \Atl^\Bun_\Ob
      (
        \seqname{B}^i
      ),
    \Atl^\Bun_\Ar
      (
        \seqname{B}^i,
        \seqname{B}^i
      ),
    \compose[A]
  \bigr )
\end{equation}
\begin{equation}
\Atl^\Bun
  (
    \seqname{B}^i
  )
\defeq
  \bigl (
    \Atl^\Bun_\Ob
      (
        \seqname{B}^i
      ),
    \Atl^\Bun_\Ar
      (
        \seqname{B}^i,
        \seqname{B}^i
      ),
    \compose[A]
  \bigr )
\end{equation}
\end{definition}

\begin{lemma}[Bundle-atlas morphisms]
\label{lem:Bun-ATLmorph}
Let
$
  \seqname{B}^i
  \defeq
  (
    E^i,
    X^i,
    Y^i,
    G^i,
    \pi^i,
    \rho^i
  )
$, $i=1,2$, be a protobundle and let
$\seqname{A}^i$ be a $\pi^i$-$G^i$-$\rho^i$-atlas of $E^i$ in the coordinate space
$C^i = X^i \times Y^i$. Then
$
  \funcseqname{f}
  \defeq
  (
    \funcname{f}_E \maps E^1 \to E^2,
    \funcname{f}_X \maps X^1 \to X^2,
    \funcname{f}_Y \maps Y^1 \to Y^2,
    \funcname{f}_G \maps G^1 \to G^2
  )
$
is a $\seqname{B^1}$-$\seqname{B^2}$ bundle-atlas morphism from $A^1$ to
$A^2$ iff $\funcname{f}_X \times \funcname{f}_Y$ is a
$G^1$-$G^2$-$\rho^1$-$\rho^2$ morphism from
$\Triv[G^1-\rho^1-]{X^1,Y^1}$ to $\Triv[G^2-\rho^2-]{X^2,Y^2}$ and
$(\funcname{f}_E, \funcname{f}_X \times \funcname{f}_Y)$ is a
$\Triv{E^1}$-$\Triv{E^2}$ m-atlas morphism from $A^1$ to $A^2$ in the
coordinate spaces $\Triv[G^1-\rho^1-]{X^1,Y^1}$,
$\Triv[G^2-\rho^2-]{X^2,Y^2}$.

\begin{proof}
If $\funcseqname{f}$ is a $\seqname{B^1}$-$\seqname{B^2}$ bundle-atlas
morphism then $\funcname{f}_X \times \funcname{f}_Y$ is a model function
and $\funcname{f}_G$ is the function asserted to exist in \cref{eq:Grho}
(preservation of group action) of \pagecref{def:Grhomorph}, so
$\funcname{f}_X \times \funcname{f}_Y$ is a $G^1$-$G^2$-$\pi^1$-$\pi^2$
morphism.

A diagram is locally nearly commutative in $X^2,Y^2,\pi^2,\rho^2$ iff it
is m-locally nearly commutative in $\pi_2 \bigl ( \Triv[G^2-\rho^2-]{X^2,Y^2} \bigr )$,
thus $(\funcname{f}_E, \funcname{f}_X \times \funcname{f}_Y)$ is a
m-atlas morphism in the coordinate space $\Triv[G^2-\rho^2-]{X^2,Y^2}$,
so $(\funcname{f}_E, \funcname{f}_X \times \funcname{f}_Y)$ is an
m-atlas morphism from $A^1$ to $A^2$ in the coordinate space
$\Triv[G^2-\rho^2-]{X^2,Y^2}$.

If $\funcname{f}_X \times \funcname{f}_Y$ is a
$G^1$-$G^2$-$\pi^1$-$\pi^2$ morphism from $\Triv[G^1-\rho^1-]{X^1,Y^1}$
to $\Triv[G^2-\rho^2-]{X^2,Y^2}$ then $\funcseqname{f}$ commutes with
$\rho^i$.

If $(\funcname{f}_E, \funcname{f}_X \times \funcname{f}_Y)$ is a
$\Triv{E}^1$-$\Triv{E}^2$ m-atlas morphism from $A^1$ to $A^2$ in the
coordinate space $\Triv[G^2-\rho^2-]{X^2,Y^2}$ then $\funcseqname{f}$
commutes with $\pi^i$.

A diagram is locally nearly commutative in $X^2,Y^2,\pi^2,\rho^2$ iff it
is m-locally nearly commutative in $\pi_2(\Triv[G^2-\rho^2-]{X^2,Y^2})$,
thus $(\funcname{f}_E, \funcname{f}_X \times \funcname{f}_Y)$ is a
m-atlas morphism in the coordinate space $\Triv[G^2-\rho^2-]{X^2,Y^2}$,
so $\funcseqname{f}$ is a $\seqname{B^1}$-$\seqname{B^2}$ bundle-atlas
morphism from $A^1$ to $A^2$.
\end{proof}
\end{lemma}

\begin{corollary}[Bundle-atlas morphisms]
\label{cor:Bun-ATLmorph}
Let
$
  \seqname{B}^i
  \defeq
  (
    E^i,
    X^i,
    Y^i,
    G^i,
    \pi^i,
    \rho^i
  )
$,
$i=1,2,3$, be a protobundle,
$\seqname{A}^i$ be a $\pi^i$-$G^i$-$\rho^i$-atlas of $E^i$ in the coordinate space
$C^i = X^i \times Y^i$ and
$
  \funcseqname{f}^i
  \defeq
  (
    \funcname{f}^i_E \maps E^i \to E^{i+1},
    \funcname{f}^i_X \maps X^i \to X^{i+1},
    \funcname{f}^i_Y \maps Y^i \to Y^{i+1},
    \funcname{f}^i_G \maps G^i \to G^{i+1}
  )
$
a $\seqname{B^1}$-$\seqname{B^2}$ bundle-atlas morphism from $A^i$ to $A^{i+i}$.

$\funcseqname{f}^2 \compose[()] \funcseqname{f}^1$ is a bundle-atlas
morphism from $A^1$ to $A^3$.

\begin{proof}
$
  \funcname{f}^1_X \times \funcname{f} 1_Y \compose
  \funcname{f}^2_X \times \funcname{f}^2_Y
$
is a $G^1$-$G^3$-$\pi^1$-$\pi^3$ morphism from
$\Triv[G^1-\rho^1-]{X^1,Y^1}$ to $\Triv[G^3-\rho^3-]{X^3,Y^3}$ by
\pagecref{lem:M-ATLmorph} and
$
  (\funcname{f}^1_E, \funcname{f}^1_X \times \funcname{f}^1_Y)
  \compose[()]
  (\funcname{f}^2_E, \funcname{f}^2_X \times \funcname{f}^2_Y)
$
is a $\Triv{E}^1$-$\Triv{E}^3$ m-atlas morphism from $A^1$ to $A^3$ in
the coordinate space $\Triv[G^3-\rho^3-]{X^3,Y^3}$ by
\pagecref{lem:Grhomorph}.
\end{proof}
\end{corollary}

\begin{definition}[Categories of bundle atlases]
\label{def:BonATLcat}

Let
$
  \seqname{B}^\alpha \defeq
  (E^\alpha, X^\alpha, Y^\alpha, G^\alpha, \pi^\alpha, \rho^\alpha)
$,
$\alpha \prec \Alpha$, be a protobundle and
$\seqname{B} \defeq \set{{\seqname{B}^\alpha}}[\alpha \prec \Alpha]$
be a set of protobundles.
Then

\begin{equation}
\Atl^\Bun_\Ar \seqname{B}
\defeq
\union
  [
    \seqname{B^\mu \in \seqname{B}},
    \seqname{B}^\nu \in \seqname{B}
  ]
  {
    \Atl^\Bun_\Ar
      (
        \seqname{B}_\mu,
        \seqname{B}_\nu
      )
  }
\end{equation}
\begin{equation}
\Atl^\Bun \seqname{B}
\defeq
  \bigl (
    \Atl^\Bun_\Ob
      (
        \seqname{B}
      ),
    \Atl^\Bun_\Ar
      (
        \seqname{B}
      ),
    \compose[A]
  \bigr )
\end{equation}

\begin{equation}
\maximal{\Atl^\Bun_\Ar} \seqname{B}
\defeq
\union
  [
    \seqname{B^\mu \in \seqname{B}},
    \seqname{B}^\nu \in \seqname{B}
  ]
  {
    \maximal{\Atl^\Bun_\Ar}
      (
        \seqname{B}_\mu,
        \seqname{B}_\nu
      )
  }
\end{equation}
\begin{equation}
\maximal{\Atl^\Bun} \seqname{B}
\defeq
  \bigl (
    \maximal{\Atl^\Bun_\Ob}
      (
        \seqname{B}
      ),
    \maximal[S-]{\Atl^\Bun_\Ar}
      (
        \seqname{B}
      ),
    \compose[A]
  \bigr )
\end{equation}

\begin{equation}
\maximal[S-]{\Atl^\Bun_\Ar} \seqname{B}
\defeq
\union
  [
    \seqname{B^\mu \in \seqname{B}},
    \seqname{B}^\nu \in \seqname{B}
  ]
  {
    \maximal[S-]{\Atl^\Bun_\Ar}
      (
        \seqname{B}_\mu,
        \seqname{B}_\nu
      )
  }
\end{equation}
\begin{equation}
\maximal[S-]{\Atl^\Bun} \seqname{B}
\defeq
  \bigl (
    \maximal[S-]{\Atl^\Bun_\Ob}
      (
        \seqname{B}
      ),
    \maximal[S-]{\Atl^\Bun_\Ar}
      (
        \seqname{B}
      ),
    \compose[A]
  \bigr )
\end{equation}
\end{definition}

\begin{lemma}[$\Atl^\Bun \seqname{B}$ is a category]
\label{lem:BunATLiscat}
Let
$
  \seqname{B}^\alpha \defeq
  (E^\alpha, X^\alpha, Y^\alpha, G^\alpha, \pi^\alpha, \rho^\alpha)
$,
$\alpha \prec \Alpha$, be a protobundle and
$\seqname{B} \defeq \set{{\seqname{B}^\alpha}}[\alpha \prec \Alpha]$
be a set of protobundles.
Then $\Atl^\Bun \seqname{B}$ is a category

Let
$
  (\seqname{A}^\alpha, \seqname{B}^\alpha)
  \in
  \Atl^\Bun_\Ob \seqname{B}
$.
Then $\Id_{(\seqname{A}^\alpha, \seqname{B}^\alpha)}$ is the
identity morphism for \\
$(\seqname{A}^\alpha, \seqname{B}^\alpha)$.

\begin{proof}
Let $(\seqname{A}^i,\seqname{B}^i)$, $i=1,2,3$
be objects of $\Atl^\Bun \seqname{B}$ and
let \\
$
  \funcname{m}^i \defeq
  \bigl (
    \funcseqname{f}^i,
    (\seqname{A}^i,\seqname{B}^i),
    (\seqname{A}^{i+1},\seqname{B}^{i+1})
  \bigr )
$
be morphisms of $\Atl^\Bun \seqname{B}$. Then
\begin{enumerate}
\item Composition: \newline
$
  \bigl (
    \funcname{m}^2 \compose[()] \funcname{m}^2,
    (\seqname{A}^1,\seqname{E}^1,\seqname{C}^1),
    (\seqname{A}^3,\seqname{E}^3,\seqname{C}^3)
  \bigr )
$
is a morphism of $\Atl^\Bun \seqname{B}$ by \pagecref{cor:Bun-ATLmorph}.
\item Associativity: \newline
Composition is associative by \pagecref{lem:atlcomp}.
\item Identity: \newline
$\Id_{(\seqname{A}^i, \seqname{E}^i, \seqname{C}^i}$ is an identity
morphism by \cref{lem:atlcomp}.
\end{enumerate}
\end{proof}
\end{lemma}

\begin{definition}[Functor from Bundle atlases to m-atlases]
\label{def:FuncBuntoM}
Let \\
$
  \seqname{B}^i
  \defeq
  (
    E^i,
    X^i,
    Y^i,
    G^i,
    \pi^i,
    \rho^i
  )
$,
$i=1,2$, be a protobundle, let $\seqname{A}^i$ be a
$\pi^i$-$G^i$-$\rho^i$-atlas of $E^i$ in the coordinate space
$C^i = X^i \times Y^i$ and let
$
  \funcseqname{f} \\
  \defeq
  (
    \funcname{f}_E \maps E^1 \to E^2,
    \funcname{f}_X \maps X^1 \to X^2,
    \funcname{f}_Y \maps Y^1 \to Y^2,
    \funcname{f}_G \maps G^1 \to G^2
  )
$
be a $\seqname{B^1}$-$\seqname{B^2}$ bundle-atlas morphism from $A^1$ to
$A^2$. Then
\begin{equation}
\Functor^\Bun_{\Bun,\M} (\seqname{A}^i, \seqname{B}^i)
\defeq
(\seqname{A}^i, \Triv{E}^i, \Triv[G^i-\rho^i-]{X^i,Y^i})
\end{equation}
\begin{multline}
\Functor^\Bun_{\Bun,\M}
  \bigl (
    \funcseqname{f},
    (\seqname{A}^1, \seqname{B}^1),
    (\seqname{A}^2, \seqname{B}^2)
  \bigr )
\defeq \\
  \bigl (
    (\funcseqname{f}_E, \funcseqname{f}_X \times \funcseqname{f}_Y),
    (\seqname{A}^1, \Triv{E^1}, \Triv[G^1-\rho^1-]{C^1}),
    (\seqname{A}^2, \Triv{E^2}, \Triv[G^2-\rho^2-]{C^2})
  \bigr )
\end{multline}
\end{definition}

\begin{theorem}[Functor from Bundle atlases to m-atlases]
\label{the:FuncBuntoM}
Let $\seqname{E}$ be a set of topological spaces.
$
  \seqname{B}^\alpha \defeq
  (E^\alpha \in \seqname{E}, X^\alpha, Y^\alpha, G^\alpha, \pi^\alpha, \rho^\alpha)
$,
$\alpha \prec \Alpha$, be a protobundle,
$\seqname{B} \defeq \set{{\seqname{B}^\alpha}}[\alpha \prec \Alpha]$
be a set of protobundles and
$
  \seqname{C}^\alpha \defeq
  \Triv[G^\alpha-\rho^\alpha-]{(X^\alpha, Y^\alpha)}
$.
Then $\Functor^\Bun_{\Bun,\M}$ is a functor from
$\Atl^\Bun(\seqname{B})$ to
$\Atl(\Trivcat{\seqname{E}}, \Trivcat[\Bun-]{\seqname{B}}$ and a functor
from $\Atl^\Bun(\seqname{B})$ to \\
$\Atl(\Trivcat{\seqname{E}}, \Trivcat[\BunProd-]{\seqname{B}}$.

\begin{proof}
Let $\seqname{o}^i \defeq (\seqname{A}^i, B^i)$, $i \in [1,3]$, be objects of
$\Atl^\Bun(\seqname{B})$ and
$
\seqname{m}^i \defeq
  \bigl (
    \funcseqname{f}^i,
    \seqname{o}^i,
    \seqname{o}^{i+1}
  \bigr )
$, $i=1,2$, be morphisms.

$\Functor^\Bun_{\Bun,\M} (\seqname{m}^i)$ is a morphism from
$\Functor^\Bun_{\Bun,\M} \seqname{o}^i$ to $\Functor^\Bun_{\Bun,\M} \seqname{o}^{i+1}$:

\begin{equation}
\begin{split}
  &
\Functor^\Bun_{\Bun,\M} (\seqname{m}^i) =
\\
  &
\Functor^\Bun_{\Bun,\M}
  \bigl (
    \funcseqname{f}^i,
    (\seqname{A}^i, \seqname{B}^i),
    (\seqname{A}^{i+1}, \seqname{B}^{i+1})
  \bigr )
=
\\
  &
  \bigl (
    (
      \funcseqname{f}^i_E,
      \funcseqname{f}^i_X \times \funcseqname{f}^i_Y
    ),
    (\seqname{A}^i, \Triv{E}^i, \Triv[G-\rho-]{C^i}),
    (\seqname{A}^{i+1}, \Triv{E}^{i+1}, \Triv[G-\rho-]{C^{i+1}})
  \bigr )
\end{split}
\end{equation}

$\Functor^\Bun_{\Bun,\M}$ maps identity functions to identity functions:

\begin{equation}
\begin{split}
  &
\Functor^\Bun_{\Bun,\M} \Id_{(\seqname{A}^i, \seqname{B}^i)} =
\\
  &
\Functor^\Bun_{\Bun,\M}
  \bigl (
    (\Id_{E^i}, \Id_{X^i}, \Id_{Y^i}, \Id_{G^i}),
    (\seqname{A}^i, \seqname{B}^i),
    (\seqname{A}^i, \seqname{B}^i)
  \bigr ) =
\\
  &
  \bigl (
    (\Id_{\Triv{E}^i}, \Id_{C^i}),
    (\seqname{A}^i, \Triv{E}^i, \Triv[G-\rho-]{C^i}),
    (\seqname{A}^i, \Triv{E}^i, \Triv[G-\rho-]{C^i})
  \bigr ) =
\\
  &
  \bigl (
    (\Id_{\Triv{E}^i}, \Id_{C^i}),
    \Functor^\Bun_{\Bun,\M} (\seqname{A}^i, \seqname{B}^i),
    \Functor^\Bun_{\Bun,\M} (\seqname{A}^i, \seqname{B}^i)
  \bigr ) =
\\
  &
\Id_{\Functor^\Bun_{\Bun,\M} (\seqname{A}^i, \seqname{B}^i)}
\end{split}
\end{equation}

$
  \Functor^\Bun_{\Bun,\M} (m^2) \compose[A] \Functor^\Bun_{\Bun,\M} (m^1)
  =
  \Functor^\Bun_{\Bun,\M} (m^2 \compose[A] m^1)
$:

\begin{enumerate}
\item
$
  m^2 \compose[A] m^1 =
  \bigl (
    (
      \funcname{f}^2_0 \compose \funcname{f}^1_0,
      \funcname{f}^2_X \compose \funcname{f}^1_X,
      \funcname{f}^2_Y \compose \funcname{f}^1_Y,
      \funcname{f}^2_G \compose \funcname{f}^1_G
    ),
    (\seqname{A}^1, \seqname{B}^1),
    (\seqname{A}^3, \seqname{B}^3)
  \bigr )
  $
\item
$
  \Functor^\Bun_{\Bun,\M}
    \bigl (
       (\seqname{A}^i, \seqname{B}^i)
    \bigr )
  =
  (\seqname{A}^i, \Triv{E}^i, \Triv[G^i-\rho^i-]{X^i,Y^i})
$
\item
$
  \Functor^\Bun_{\Bun,\M} (m^i)
    = \\
    \bigl (
      (\funcname{f}^i_0, \funcname{f}^i_1),
      (\seqname{A}^i, \Triv{E}^i, \Triv[G^i-\rho^i-]{X^i,Y^i}),
      (\seqname{A}^{i+1}, \Triv{E}^{i+1}, \Triv[G^i-\rho^i-]{C^{i+1}})
    \bigr )
$
\item
$
  \Functor^\Bun_{\Bun,\M} (m^2) \compose[A] \Functor^\Bun_{\Bun,\M} (m^1)
  = \\
    \bigl (
      (
        \funcname{f}_0^2 \compose \funcname{f}_0^1,
        \funcname{f}_1^2 \compose \funcname{f}_1^2
      ),
      (\seqname{A}^1, \Triv{E}^1, \Triv[G^1-\rho^1-]{C^1}),
      (\seqname{A}^3, \Triv{E}^3, \Triv[G^3-\rho^3-]{C^3})
    \bigr )
$
\item
$
  \Functor^\Bun_{\Bun,\M} (m^2 \compose m^1)
  = \\
    \bigl (
      (
        \funcname{f}_0^2 \compose \funcname{f}_0^1,
        \funcname{f}_1^2 \compose \funcname{f}_1^2
      ),
      (\seqname{A}^1, \Triv{E}^1, \Triv[G^1-\rho^1-]{C^1}),
      (\seqname{A}^3, \Triv{E}^3, \Triv[G^3-\rho^3-]{C^3})
    \bigr )
$
\end{enumerate}

\end{proof}
\end{theorem}

\begin{lemma}[Base space functions derived from bundle-atlas morphisms]
\label{lem:BunAtlbasemorph}
Let $\seqname{E}^i$, $i=1,2$, be a model space, $X^i$, $Y^i$ topological
spaces, $G^i$ a topological group, $\rho^i$ an effective action of $G^i$
on $Y^i$, $\seqname{C}^i \defeq (X^i \times Y^i, \catname{XY}^i)$ a
$G^i$-$\rho^i$ model space of $X^i \times Y^i$,
$\funcname{f}_C \maps \seqname{C}^1 \to \seqname{C}^2$ a
$G^1$-$G^2$-$\rho^1$-$\rho^2$ morphism of
$X^1 \times Y^1$ to $X^2 \times Y^2$, i.e., a model function that
preserves group action, $\seqname{A}^i$ an m-atlas of $E^i$ in the
coordinate model space $\seqname{C}^i$ and
$
  \funcseqname{f} \defeq
  (
    \funcname{f}_E \maps \seqname{E^1} \to \seqname{E}^2,
    \funcname{f}_C \maps \seqname{C}^1 \to \seqname{C}^2
  )
$
an $\seqname{E}^1$-$\seqname{E}^2$ m-atlas morphism of $\seqname{A}^1$
to $\seqname{A}^2$ in the coordinate spaces $\seqname{C}^1$,
$\seqname{C}^2$.

Then there exists a unique function $\funcname{f}_X \maps X^1 \to X^2$
such that $\funcname{f}_X \compose \pi_1 = \pi_1 \compose \funcname{f}_C$.

\begin{proof}
Let $x$ in $X^1$, $y,y'$ in $Y^1$, Then
$\funcname{f}_C$ preserves fibers, i.e.,
$
\pi_1 \bigl ( \funcname{f}_C(x,y) \bigr ) =
\pi_1 \bigl ( \funcname{f}_C(x,y') \bigr )
$.
by \pagecref{lem:Grhomorph}. Define
$\funcname{f}_X(x) \defeq \pi_1 \bigl ( \funcname{f}_C(x,y) \bigr )$
\end{proof}

\begin{remark}
There need not exist $\funcname{f}_Y \maps Y^1 \to Y^2$ such that
$\funcname{f}_C = \funcname{f}_X \times \funcname{f}_Y$.
\end{remark}
\end{lemma}

\begin{definition}[Functor from m-atlases to Bundle atlases]
\label{def:FuncMtoBun}
Let $\catseqname{E}$ be a trivial model category, $\catname{C}$ be a
trivial product coordinate model category with objects
$
\set
  {
    \seqname{C}^\alpha \defeq
    (X^\alpha \times Y^\alpha, \catname{XY}^\alpha)
  }%
  [\alpha \prec \Alpha]
$,
$G$ a group valued function on $\Ob(\seqname{C})$ and $\rho$ a function
valued function on $\Ob(\seqname{C})$ such that for every
$
  \seqname{C}^\alpha \defeq
  (X^\alpha \times Y^\alpha, \catname{XY}^\alpha) \objin \catname{C}
$,
$\rho(\seqname{C}^\alpha)$ is an effective action of
$G(\seqname{C}^\alpha)$ on $Y^\alpha$ and
$
\Triv%
  [
    G(\seqname{C}^\alpha)-\rho(\seqname{C}^\alpha)-
  ]
  {X^\alpha, Y^\alpha} =
  \seqname{C}^\alpha
$.

Let $\seqname{E}^i \defeq (E^i, \catname{E}^i) \objin \catseqname{E}$,
$i=1,2$, be a trivial model space,
$
  \seqname{C}^i \defeq (X^i \times Y^i, \catname{XY}^i)
  \objin \catname{C}
$
$G^i \defeq G(\seqname{C}^i)$, $\rho^i \defeq \rho(\seqname{C}^i)$,
$\seqname{A}^i$ a full m-atlas of $\seqname{E}^i$ in the coordinate
space $\seqname{C}^i$,
$\pi^i \maps E^i \onto X^i$ the unique function asserted in
\cref{lem:mGrho} and
$
  \seqname{B}^i \defeq \\
    (
      E^i,
      X^i,
      Y^i,
      G^i,
      \pi^i,
      \rho^i
    )
$.
Then define
\begin{equation}
\Functor^\Bun_{\M-G-\rho,\Bun}
(\seqname{A}^i, \seqname{E}^i, \seqname{C}^i)
\defeq
(\seqname{A}^i, \seqname{B}^i)
\end{equation}

Let
$
  \funcseqname{f} \defeq
  (
    \funcname{f}_E \maps \seqname{E^1} \to \seqname{E}^2,
    \funcname{f}_C \defeq \funcname{f}_X \times \funcname{f}_Y \maps
    \seqname{C}^1 \to \seqname{C}^2
  )
$
be an $\seqname{E}^1$-$\seqname{E}^2$ m-atlas morphism of
$\seqname{A}^1$ to $\seqname{A}^2$ in the coordinate spaces
$\seqname{C}^1, \seqname{C}^2$ that preserves the group action,
$\funcname{f}_G \maps G^1 \to G^2$ the unique function asserted to exist
in \cref{eq:Grho} (preservation of group action) of
\pagecref{def:Grhomorph}.

Then
\begin{multline}
\Functor^\Bun_{\M-G-\rho,\Bun}
  \bigl (
    \funcseqname{f},
    (\seqname{A}^1, \seqname{E^1}, \seqname{C}^1),
    (\seqname{A}^2, \seqname{E^2}, \seqname{C}^2)
  \bigr )
\defeq \\
  \bigl (
    (\funcname{f}_E, \funcname{f}_X, \funcname{f}_Y, \funcname{f}_G),
    (\seqname{A}^1, \seqname{B}^1),
    (\seqname{A}^2, \seqname{B}^2)
  \bigr )
\end{multline}
\end{definition}

\begin{theorem}[Functor from m-atlases to bundle atlases]
\label{the:FuncMtoBun}
Let $\catname{C}$ be a trivial product coordinate model category,
$\catname{E}$ a model category,
$G$ a group valued function on $\Ob(\seqname{C})$ and $\rho$ a function
valued function on $\Ob(\seqname{C})$ such that for every
$
  \seqname{C}^\alpha \defeq
  (X^\alpha \times Y^\alpha, \catname{XY}^\alpha) \objin \catname{C}
$,
$\rho(\seqname{C}^\alpha)$ is an effective action of
$G(\seqname{C}^\alpha)$ on $Y^\alpha$ and
$
\Triv%
  [
    G(\seqname{C}^\alpha)-\rho(\seqname{C}^\alpha)-
  ]
  {X^\alpha, Y^\alpha} =
  \seqname{C}^\alpha
$.

Let $\pi$ be the unique function valued function on
$\full{\Atl_\Ob}(\catname{E}, \catname{C})$ such that for every
$
\bigl (
  \seqname{A}^\alpha,
  (E^\alpha,\catname{E}^\alpha),
  (C^\alpha, \catname{C}^\alpha)
  \bigr )
  \in \full{\Atl_\Ob}(\catname{E}, \catname{C})
$ and $(U, V, \phi) \in  \seqname{A}^\alpha$,
$\pi_1 \compose \phi = \pi(\seqname{A}^\alpha) \restriction_U$
and let
$
  \seqname{B} \defeq
  \set
  {
    \seqname{B}^\alpha \defeq
      (
        E^\alpha,
        X^\alpha,
        Y^\alpha,
        G(\seqname{C}^\alpha),
        \pi(\seqname{A}^\alpha),
        \rho(\seqname{C}^\alpha)
      )
  }%
  [
    {
      \bigl (
        \seqname{A}^\alpha,
        (E^\alpha,\catname{E}^\alpha),
        (C^\alpha \defeq X^\alpha \times Y^\alpha, \catname{C}^\alpha)
      \bigr ) \in \full{\Atl_\Ob}(\catname{E}, \catname{C})
    }
  ]*
$.

Then $\Functor^\Bun_{\M-G-\rho,\Bun}$ is a functor from
$\Atl(\Trivcat{\seqname{E}}, \Trivcat[\BunProd-]{\seqname{B}}$
to $\Atl^\Bun(\seqname{B})$.

\begin{proof}
Let
$
  \seqname{o}^i \defeq
  (
    \seqname{A}^i,
    \seqname{E}^i,
    \seqname{C}^i \defeq (X^i \times Y^i, \catname{XY}^i)
  )
  \objin \catname{M})
$,
$i \in [1,3]$,
$G^i \defeq G(\seqname{o}^i)$, \\
$\pi^i \defeq \pi(\seqname{o}^i)$,
$\rho^i \defeq \rho(\seqname{o}^i)$,
$\seqname{A}^i$ an m-atlas of $E^i$ in the coordinate model space $\seqname{C}^i$, \\
$
  \seqname{B}^i \defeq
    (
      E^i,
      X^i,
      Y^i,
      G^i,
      \pi^i,
      \rho^i
    )
$,
$
  \seqname{m}^i \defeq
    \bigl (
      \funcseqname{f}^i,
      \seqname{o}^i,
      \seqname{o}^{i+1}
    \bigr )
    \arin \catname{M}
$, $i=1,2$,
an $E^i$-$E^{i+1}$ m-atlas morphism of $\seqname{A}^i$ to
$\seqname{A}^{i+1}$ in the coordinate spaces
$\seqname{C}^i, \seqname{C}^{i+1}$ that preserves the group action,
$\funcname{f}_G \maps G^i \to G^{i+1}$ the unique function asserted to
exist in \cref{eq:Grho} (preservation of group action) of
\pagecref{def:Grhomorph} and $\pi^i \maps E^i \onto X^i$ the unique
function asserted in \cref{lem:mGrho}.

Let $\funcname{f}^i_G \maps G^i \to G^{i+1}$ be the function asserted by
\pagecref{eq:Grho}; let
$\funcname{f}^i_1 = \funcname{f}^i_X \times \funcname{f}^i_Y$ be the
decomposition given by \pagecref{def:Grhomodtriv}.

$\Functor^\Bun_{\M-G-\rho,\Bun} (\seqname{m}^i)$ is a morphism from
$\Functor^\Bun_{\M-G-\rho,\Bun} \seqname{o}^i$ to \\
$\Functor^\Bun_{\M-G-\rho,\Bun} \seqname{o}^{i+1}$:

\begin{equation}
\begin{split}
  &
\Functor^\Bun_{\M-G-\rho,\Bun} (\seqname{m}^i) =
\\
  &
\Functor^\Bun_{\M-G-\rho,\Bun}
  (
    \funcseqname{f}^i,
    \seqname{o}^i,
    \seqname{o}^{i+1}
  )
=
\\
  &
  \bigl (
    (
      \funcname{f}^i_0,
      \funcname{f}^i_X,
      \funcname{f}^i_Y,
      \funcname{f}^i_G
    ),
    (\seqname{A}^i, \seqname{B}^i),
    (\seqname{A}^{i+1}, \seqname{B}^{i+1})
  \bigr )
=
\\
  &
  \bigl (
    (
      \funcname{f}^i_0,
      \funcname{f}^i_X,
      \funcname{f}^i_Y,
      \funcname{f}^i_G
    ),
    \Functor^\Bun_{\M-G-\rho,\Bun} \seqname{o}^i,
    \Functor^\Bun_{\M-G-\rho,\Bun} \seqname{o}^{i+1}
  \bigr )
\end{split}
\end{equation}

$\Functor^\Bun_{\M-G-\rho,\Bun}$ maps identity functions to identity
functions:

\begin{equation}
\begin{split}
  &
\Functor^\Bun_{\M-G-\rho,\Bun} \Id_{\seqname{o}^i} =
\\
  &
\Functor^\Bun_{\M-G-\rho,\Bun}
  \bigl (
    (\Id_{E^i}, \Id_{C^i}),
    \seqname{o}^i,
    \seqname{o}^i
  \bigr ) =
\\
  &
  \bigl (
    (\Id_{E^i}, \Id_{X^i}, \Id_{Y^i}, \Id_{G^i}),
    (\seqname{A}^i, \seqname{B}^i),
    (\seqname{A}^i, \seqname{B}^i)
  \bigr ) =
\\
  &
  \bigl (
    (\Id_{E^i}, \Id_{X^i}, \Id_{Y^i}, \Id_{G^i}),
    \Functor^\Bun_{\M-G-\rho,\Bun} \seqname{o}^i,
    \Functor^\Bun_{\M-G-\rho,\Bun} \seqname{o}^i
  \bigr ) =
\\
  &
\Id_{\Functor^\Bun_{\M-G-\rho,\Bun} \seqname{o}^i}
\end{split}
\end{equation}

$
  \Functor^\Bun_{\M-G-\rho,\Bun} (m^2) \compose[A] \Functor^\Bun_{\M-G-\rho,\Bun} (m^1)
  =
  \Functor^\Bun_{\M-G-\rho,\Bun} (m^2 \compose[A] m^1)
$:

\begin{enumerate}
\item
$
  m^2 \compose[A] m^1 =
  \bigl (
    (
      \funcname{f}_0^2 \compose \funcname{f}_0^1,
      \funcname{f}_1^2 \compose \funcname{f}_1^1
    ),
    \seqname{o}^1,
    \seqname{o}^3
  \bigr )
  $
\item
$
  \Functor^\Bun_{\M-G-\rho,\Bun}
    \bigl (
       \seqname{o}^i
    \bigr )
  =
  (\seqname{A}^i, \seqname{B}^i)
$
\item
$
  \Functor^\Bun_{\M-G-\rho,\Bun} (m^i)
    = \\
    \bigl (
      (
        \funcname{f}^i_0,
        \funcname{f}^i_X,
        \funcname{f}^i_Y,
        \funcname{f}^i_G
      ),
      (\seqname{A}^i, \seqname{B}^i),
      (\seqname{A}^{i+1}, \seqname{B}^{i+1})
    \bigr )
$
\item
$
  \Functor^\Bun_{\M-G-\rho,\Bun} (m^2) \compose[A] \Functor^\Bun_{\M-G-\rho,\Bun} (m^1)
  = \\
    \bigl (
      (
        \funcname{f}^2_0 \compose \funcname{f}^1_0,
        \funcname{f}^2_X \compose \funcname{f}^1_X,
        \funcname{f}^2_Y \compose \funcname{f}^1_Y,
        \funcname{f}^2_G \compose \funcname{f}^1_G
      ),
      (\seqname{A}^1, \seqname{B}^1),
      (\seqname{A}^3, \seqname{B}^3)
    \bigr )
$
\item
$
  \Functor^\Bun_{\M-G-\rho,\Bun} (m^2 \compose m^1)
  = \\
    \bigl (
      (
        \funcname{f}^2_0 \compose \funcname{f}^1_0,
        \funcname{f}^2_X \compose \funcname{f}^1_X,
        \funcname{f}^2_Y \compose \funcname{f}^1_Y,
        \funcname{f}^2_G \compose \funcname{f}^1_G
      ),
      (\seqname{A}^1, \seqname{B}^1),
      (\seqname{A}^3, \seqname{B}^3)
    \bigr )
$
\end{enumerate}
\end{proof}

$\Functor^\Bun_{\Bun,\M} \compose \Functor^\Bun_{\M-G-\rho,\Bun} = \Id$.
\begin{proof}
Expanding the definitions, we have
\begin{enumerate}
\item
$
  \Functor^\Bun_{\M-G-\rho,\Bun}
  (\seqname{A}^i, \seqname{E}^i, \seqname{C}^i) =
  (\seqname{A}^i, \seqname{B}^i)
$
\item
$
  \Functor^\Bun_{\Bun,\M} (\seqname{A}^i, \seqname{B}^i) =
  (\seqname{A}^i, \Triv{E}^i, \Triv[G^i-\rho^i-]{X^i,Y^i})
  (\seqname{A}^i, \seqname{E}^i, \seqname{C}^i)
$,
since by \cref{def:FuncMtoBun}, $\seqname{E}^i$ and $\seqname{C}^i$ are trivial.
\item
$
  \Functor^\Bun_{\M-G-\rho,\Bun}
    \bigl (
      (
        \funcseqname{f}^1_E,
        \funcseqname{f}^1_X \times \funcseqname{f}^1_Y
      ),
      (\seqname{A}^1, \seqname{E^1}, \seqname{C}^1),
      (\seqname{A}^2, \seqname{E^2}, \seqname{C}^2)
    \bigr ) = \\
    \bigl (
      (
        \funcname{f}^1_E,
        \funcname{f}^1_X,
        \funcname{f}^1_Y,
        \funcname{f}^1_G
      ),
      (\seqname{A}^1, \seqname{B}^1),
      (\seqname{A}^2, \seqname{B}^2)
    \bigr )
$
\item
$
  \Functor^\Bun_{\Bun,\M}
    \bigl (
      (
        \funcname{f}^1_E,
        \funcname{f}^1_X,
        \funcname{f}^1_Y,
        \funcname{f}^1_G
      ),
      (\seqname{A}^1, \seqname{B}^1),
      (\seqname{A}^2, \seqname{B}^2)
    \bigr ) = \\
    \bigl (
      (
        \funcseqname{f}^1_E,
        \funcseqname{f}^1_X \times \funcseqname{f}^1_Y
      ),
      (\seqname{A}^1, \Triv{E^1}, \Triv[G^1-\rho^1-]{C^1}),
      (\seqname{A}^2, \Triv{E^2}, \Triv[G^2-\rho^2-]{C^2})
    \bigr ) = \\
    \bigl (
      (
        \funcseqname{f}^1_E,
        \funcseqname{f}^1_X \times \funcseqname{f}^1_Y
      ),
      (\seqname{A}^1, \seqname{E}^1, \seqname{C}^1),
      (\seqname{A}^2, \seqname{E}^2, \seqname{C}^2)
    \bigr )
$,
since by \cref{def:FuncMtoBun}, $\seqname{E}^i$ and $\seqname{C}^i$ are trivial.
\end{enumerate}
\end{proof}
\end{theorem}

\subsection{Associated model spaces and functors}
\label{sub:assmodB}
\begin{definition}[Coordinate model spaces associated with bundle atlases]
\label{def:ismodBF}
Let \\
$
  \seqname{B}^i
  \defeq
  (
    E^i,
    X^i,
    Y^i,
    G^i,
    \pi^i,
    \rho^i
  )
$,
$\seqname{A}^i$ be a $\pi^i$-$G^i$-$\rho^i$-bundle
atlas of $E^i$ in the coordinate space $C^i = X^i \times Y^i$ and
$
  \funcseqname{f} \\
  \defeq
  (
    \funcname{f}_E \maps E^1 \to E^2,
    \funcname{f}_X \maps X^1 \to X^2,
    \funcname{f}_Y \maps Y^1 \to Y^2,
    \funcname{f}_G \maps G^1 \to G^2
  )
$
such that \\
$
  \isAtl^\Bun_\Ob
  (
    \seqname{A^i},
    E^i,
    X^i,
    Y^i,
    G^i,
    \pi^i,
    \rho^i
  )
$
and
$
  \isAtl^\Bun_\Ar
    (
      \seqname{A^1}, \seqname{B^1},
      \seqname{A^2}, \seqname{B^2},
      \funcseqname{f}
    )
$.
Then

\begin{multline}
\Functor^\Bun_2 (\seqname{A}^i, \seqname{B}^i)
\defeq \\
\minimal{\Mod}
  \Biggl (
    C^i,
    \pi_2[\seqname{A}^i],
    \set
      {\phi' \compose \phi^{-1}}%
      [
        \equant
          {
            {(U,V,\phi) \in \seqname{A}^i},
            {(U',V',\phi') \in \seqname{A}^i}
          }
          {
            U \cap U' \ne \emptyset
          }
      ]
  \Biggr )
\end{multline}
\begin{equation}
\Functor^\Bun_2
  \bigl (
    \funcseqname{f},
    (\seqname{A}^1, \seqname{B}^1),
    (\seqname{A}^2, \seqname{B}^2)
  \bigr )
\defeq
  \funcname{f}_X \times \funcname{f}_Y
  \maps
  \Functor^\Bun_2 (\seqname{A}^1, \seqname{B}^1)
  \to
  \Functor^\Bun_2 (\seqname{A}^2, \seqname{B}^2)
\end{equation}

The minimal $G$-$\rho$ coordinate model space with neighborhoods in
$\seqname{A}^i$ is $\Functor^\Bun_2 (\seqname{A}^i, \seqname{B}^i)$.

The coordinate mapping associated with the
$\seqname{B^1}$-$\seqname{B^2}$ bundle-atlas morphism
$\funcseqname{f}$ from $\seqname{A^1}$ to $\seqname{A^2}$ is
$
  \funcname{f}_X \times \funcname{f}_Y
  \maps
  \Functor^\Bun_2 (\seqname{A}^1, \seqname{B}^1)
  \to
  \Functor^\Bun_2 (\seqname{A}^2, \seqname{B}^2)
$.
If it is a model function then it is also the coordinate
$G^1$-$G^2$-$\rho^1$-$\rho^2$ morphism associated with the
$\seqname{B^1}$-$\seqname{B^2}$ bundle-atlas morphism $\funcseqname{f}$
from $\seqname{A^1}$ to $\seqname{A^2}$.
\end{definition}

\begin{lemma}[Coordinate model spaces associated with bundle atlases]
\label{lem:ismodBF}
Let
$
  \seqname{B}
  \defeq
  (
    E,
    X,
    Y,
    G,
    \pi,
    \rho
  )
$
and let $\seqname{A}$ be a $\pi$-$G$-$\rho$-bundle atlas of $E$ in the
coordinate space $C = X \times Y$. Then
$\Functor^\Bun_2 (\seqname{A}, \seqname{B})$ is a model space.
\end{lemma}

\begin{proof}
$\Functor^\Bun_2 (\seqname{A}, \seqname{B})$ satisfies the
conditions for a model space.
for a model space. Let
$
  \catname{C} \defeq
  \pi_2 \bigl ( \Functor^\Bun_2 (\seqname{A}, \seqname{B}) \bigr )
$.
\begin{enumerate}
\item Since $\pi_2[\seqname{A}]$ is an open cover of
$\union{\pi_2[\seqname{A}]}$, the set of finite intersections is also an
open cover.
\item Finite intersections of finite intersections are finite intersections
\item Restrictions of continuous function are continuous
\item If $f \maps A \to B$ is a morphism of
$\Functor^\Bun_2 (\seqname{A}, \seqname{B})$ $A, A', B, B'$ objects of
$\catname{C}$ $A' \subseteq A$, $B' \subseteq B$ and $f[A'] \subseteq
B'$ then since $f \maps A \to B$ is a morphism it is a restriction of a
transition function between its restrictions to sets in
$\pi_2[\seqname{A}]$ and its restrictions are also, hence morphisms, and
thus $f \restriction_{A'} A' \to B'$ is a morphism.
\item
If $(U,V,\phi) \in \seqname{A}$ then $\Id_V = \phi \compose \phi^{-1}A$
is a transition function and hence a morphism of
$\Functor^\Bun_2 (\seqname{A}, \seqname{B})$. If $A, A'$ objects of
$\catname{C}$ and
$A' \subseteq A$ then the inclusion map
$\funcname{i} \maps A' \hookrightarrow A$ is a restriction of an
identity morphism of $\Functor^\Bun_2 (\seqname{A}, \seqname{B})$ and
hence a morphism.

\item Restricted sheaf condition: let
\begin{enumerate}
\item $U_\alpha,V_\alpha$, $\alpha \prec \Alpha$, be objects of
$\catname{C}$ \item $\funcname{f}_\alpha \maps U_\alpha \to V_\alpha$ be
morphisms of $\catname{C}$
\item $U \defeq \union[\alpha \prec \Alpha]{U_\alpha}$
\item $V \defeq \union[\alpha \prec \Alpha]{V_\alpha}$
\item $\funcname{f} \maps U \to V$ be continuous and
$
  \uquant
    {{\alpha \prec \Alpha},{x \in U_\alpha}}
    {\funcname{f}(x) = \funcname{f}_\alpha(x)}
$
\end{enumerate}
Then $\funcname{f}$ is generated by the group action and hence a
morphism of $\catname{C}$
\end{enumerate}
\end{proof}

\begin{definition}[Model spaces associated with bundle atlases]
\label{def:ismodT}
Let \\
$
  \seqname{B}^i
  \defeq
  (
    E^i,
    X^i,
    Y^i,
    G^i,
    \pi^i,
    \rho^i
  )
$,
$i=1,2$, and $\seqname{A}^i$ be a $\pi^i$-$G^i$-$\rho^i$-bundle atlas of
$E^i$ in the coordinate space $C^i = X^i \times Y^i$.   Then

\begin{multline}
\Functor^\Bun_1 (\seqname{B}^i, \seqname{A}^i)
\defeq \\
\minimal{\Mod}
  \bigl (
    C^i,
    \pi_1[\seqname{A}^i],
    \set
      {\phi'^{-1} \compose \phi}%
      [
        \equant
          {
            {(U,V,\phi) \in \seqname{A}^i},
            {(U',V',\phi') \in \seqname{A}^i}
          }
          {
            U \cap U' \ne \emptyset
          }
      ]
  \bigr )
\end{multline}
\begin{equation}
\Functor^\Bun_1
  \bigl (
    \funcseqname{f},
    (\seqname{A}^1, \seqname{B}^1),
    (\seqname{A}^2, \seqname{B}^2)
  \bigr )
\defeq
  \funcname{f}_E
  \maps
  \Functor^\Bun_1 (\seqname{A}^1, \seqname{B}^1)
  \to
  \Functor^\Bun_1 (\seqname{A}^2, \seqname{B}^2)
\end{equation}

The minimal model spaces with neighborhoods in the
atlas $\seqname{A}^i$ is \\
$\Functor^\Bun_1 (E^i, X^i, Y^i, \pi^i, G^i, \rho^i, \seqname{A}^i)$.

The mapping associated with the $\seqname{B^1}$-$\seqname{B^2}$
bundle-atlas morphism $\funcseqname{f}$ from $\seqname{A^1}$ to
$\seqname{A^2}$ is
$
  \funcname{f}_E
  \maps
  \Functor^\Bun_1 (\seqname{A}^1, \seqname{B}^1)
  \to
  \Functor^\Bun_1 (\seqname{A}^2, \seqname{B}^2)
$.
If it is a model function then it is also the m-atlas morphism
associated with the $\seqname{B^1}$-$\seqname{B^2}$ bundle-atlas
morphism $\funcseqname{f}$ from $\seqname{A^1}$ to $\seqname{A^2}$.
\end{definition}

\begin{lemma}[Model spaces associated with bundle atlases]
\label{lem:ismodT}
Let \\
$
  \seqname{B}
  \defeq
  (
    E,
    X,
    Y,
    G,
    \pi,
    \rho
  )
$
and let $\seqname{A}$ be a $\pi$-$G$-$\rho$-bundle atlas of $E$ in the
coordinate space $C = X \times Y$. Then
$\Functor^\Bun_1 (\seqname{A}, \seqname{B})$ is a model space.

\begin{proof}
\Pagecref{lem:minmod}
\end{proof}
\end{lemma}

\begin{theorem}[Functors from bundle atlases to model spaces]
Let
$
  \seqname{B}^\alpha \defeq \\
  (E^\alpha, X^\alpha, Y^\alpha, G^\alpha, \pi^\alpha, \rho^\alpha)
$,
$\alpha \prec \Alpha$, be a protobundle, and
$\seqname{B} \defeq \set{{\seqname{B}^\alpha}}[\alpha \prec \Alpha]$
be a set of protobundles. Then
$\Functor^\Bun_1$ is a functor from
$\Atl^\Bun \seqname{B}$ to $\Triv{\seqname{E}}$ and
$\Functor^\Bun_2$ is a functor from
$\Atl^\Bun \seqname{B}$ to $\Trivcat[\Bun-]{\seqname{B}}$.

\begin{proof}
Let $\seqname{o}^i \defeq (\seqname{A}^i, E^i, C^i)$, $i \in [1,3]$, be
objects in
$
\Bun
  (
    \seqname{E},
    \seqname{X},
    \seqname{Y},
    \funcseqname{\pi},
    \funcseqname{G},
    \funcseqname{\rho}
  )
$
and let
$
m^i \defeq
  \bigl (
    (\funcname{f}^i_0, \funcname{f}^i_1),
    o^i,
    o^{i+1}
  \bigr )
$
be morphisms in $\Atl^\Bun(\seqname{E}, \seqname{C})$.

$
  \Functor^\Bun_1 \maps
  \Atl^\Bun \seqname{B} \to
  \Triv{\seqname{E}}
$:
\begin{enumerate}
\item $\Functor (f \maps A \to B) \maps \Functor(A) \to \Functor (B)$:
\newline
$
\Functor^\Bun_1(\seqname{m}^i) =
\funcname{f}^i_0 \maps \Functor^\Bun_1(\seqname{o}^i) \to
\Functor^\Bun_1(\seqname{o}^{i+1})
$
\item $\Functor (g \compose f) = \Functor (g) \compose \Functor (f)$: \newline
$
\Functor^\Bun_1(\funcseqname{m}^2 \compose[A] \funcseqname{m}^1)
=
\\
\Functor^\Bun_1
  \bigl (
    (
      \funcname{f}^2_0 \compose \funcname{f}^1_0,
      \funcname{f}^2_1 \compose \funcname{f}^1_1
    )
    (\seqname{A}^1, E^1, C^1),
    (\seqname{A}^3, E^3, C^3)
  \bigr )
=
\\
\funcname{f}^2_0 \compose \funcname{f}^1_0 \maps
  \Functor^\Bun_1(\seqname{o}^1) \to
  \Functor^\Bun_1(\seqname{o}^3)
=
\\
\bigl (
  \funcname{f}^2_0 \maps
  \Functor^\Bun_1(\seqname{o}^2) \to
  \Functor^\Bun_1(\seqname{o}^3)
  \bigr )
\compose
\bigl (
  \funcname{f}^1_0 \maps
  \Functor^\Bun_1(\seqname{o}^1) \to
  \Functor^\Bun_1(\seqname{o}^2)
\bigr )
=
\\
\Functor^\Bun_1
  \bigl (
    (
      \funcname{f}^2_0,
      \funcname{f}^2_1
    )
    (\seqname{A}^2, E^2, C^2),
    (\seqname{A}^3, E^3, C^3)
  \bigr )
\compose \\
\Functor^\Bun_1
  \bigl (
    (
      \funcname{f}^1_0,
      \funcname{f}^1_1
    )
    (\seqname{A}^1, E^1, C^1),
    (\seqname{A}^2, E^2, C^2)
  \bigr )
=
\\
\Functor^\Bun_1(\funcseqname{m}^2)
\compose
\Functor^\Bun_1(\funcseqname{m}^1)
$
\item $\Functor (\Id_A) = \Id_{\Functor (A)}$:
\begin{enumerate}
\item
$
\Functor^\Bun_1(\Id_{o^i}) =
\Functor^\Bun_1
  \bigl (
    (\Id_{E^i}, \Id_{C^i}),
    (\seqname{A}^i, E^i, C^i),
    (\seqname{A}^i, E^i, C^i)
  \bigr )
=
\\
\Id_{E^i} \maps \Functor^\Bun_1(\seqname{o}^i) \to
\Functor^\Bun_1(\seqname{o}^i)
$
\item
$
\Id_{\Functor^\Bun_1}(o^i) =
\Id_{E^i} \maps \Functor^\Bun_1(\seqname{o}^i) \to
\Functor^\Bun_1(\seqname{o}^i)
$
\end{enumerate}
\end{enumerate}

$
  \Functor^\Bun_2 \maps
  \Atl^\Bun \seqname{B} \to
  \Trivcat[\Bun-]{\seqname{B}}
$:

\begin{enumerate}
\item $\Functor (f \maps A \to B) \maps \Functor(A) \to \Functor (B)$:
\newline
$
\Functor^\Bun_2(\seqname{m}^i) =
\funcname{f}^i_1 \maps \Functor^\Bun_2(\seqname{o}^i) \to
\Functor^\Bun_2(\seqname{o}^{i+1})
$
\item $\Functor (g \compose f) = \Functor (g) \compose \Functor (f)$:
\newline
$
\Functor^\Bun_2(\funcseqname{m}^2 \compose[A] \funcseqname{m}^1)
=
\\
\Functor^\Bun_2
  \bigl (
    (
      \funcname{f}^2_0 \compose \funcname{f}^1_0,
      \funcname{f}^2_1 \compose \funcname{f}^1_1
    )
    (\seqname{A}^1, E^1, C^1),
    (\seqname{A}^3, E^3, C^3)
  \bigr )
=
\\
\funcname{f}^2_0 \compose \funcname{f}^1_0 \maps
  \Functor^\Bun_2(\seqname{o}^1) \to
  \Functor^\Bun_2(\seqname{o}^3)
=
\\
\bigl (
  \funcname{f}^2_1 \maps
  \Functor^\Bun_2(\seqname{o}^2) \to
  \Functor^\Bun_2(\seqname{o}^3)
  \bigr )
\compose
\bigl (
  \funcname{f}^1_1 \maps
  \Functor^\Bun_2(\seqname{o}^1) \to
  \Functor^\Bun_2(\seqname{o}^2)
\bigr )
=
\\
\Functor^\Bun_2
  \bigl (
    (
      \funcname{f}^2_0,
      \funcname{f}^2_1
    )
    (\seqname{A}^2, E^2, C^2),
    (\seqname{A}^3, E^3, C^3)
  \bigr )
\compose \\
\Functor^\Bun_2
  \bigl (
    (
      \funcname{f}^1_0,
      \funcname{f}^1_1
    )
    (\seqname{A}^1, E^1, C^1),
    (\seqname{A}^2, E^2, C^2)
  \bigr )
=
\\
\Functor^\Bun_2(\funcseqname{m}^2)
\compose
\Functor^\Bun_2(\funcseqname{m}^1)
$
\item $\Functor (\Id_A) = \Id_{\Functor (A)}$:
\begin{enumerate}
\item
$
\Functor^\Bun_2(\Id_{\seqname{o}^i}) =
\Functor^\Bun_2
  \bigl (
    (\Id_{E^i}, \Id_{C^i}),
    (\seqname{A}^i, E^i, C^i),
    (\seqname{A}^i, E^i, C^i)
  \bigr )
=
\\
\Id_{C^i} \maps \Functor^\Bun_2(\seqname{o}^i) \to
\Functor^\Bun_2(\seqname{o}^i)
$
\item
$
\Id_{\Functor^\Bun_2}(\seqname{o}^i) =
\Id_{C^i} \maps \Functor^\Bun_2(\seqname{o}^i) \to
\Functor^\Bun_2(\seqname{o}^i)
$
\end{enumerate}
\end{enumerate}
\end{proof}
\end{theorem}

\subsection{Fiber bundles}
\label{sub:fib}
Conventionally a fiber bundle is different from its atlases, but
\pagecref{def:BonATLcat} encourages treating them on an equal footing. All
of the results for maximal bundle atlases carry directly over to results
for fiber bundles.

\begin{definition}[fiber bundles]
\label{def:Bun}
Let $E$, $X$ and $Y$ be topological spaces, $\pi \maps E \onto X$ surjective,
$G$ a topological group,
$\rho \maps Y \times G \to Y$ an effective right action of
$G$ on $Y$ and $\seqname{A}$
a maximal $\pi$-$G$-$\rho$-bundle atlas of $E$ in the coordinate space $X \times Y$.
Then $(E, X, Y, \pi, G, \rho, \seqname{A})$ is a fiber bundle.

Let
$
  \seqname{B}
  \defeq
  \set
    {{
      \seqname{B}^\alpha \defeq
      (
        E^\alpha,
        X^\alpha,
        Y^\alpha,
        G^\alpha,
        \pi^\alpha,
        \rho^\alpha
      )
    }}%
    [
      \alpha \prec \Alpha
    ]
$,
where $E^\alpha, X^\alpha, Y^\alpha$ are topological spaces,
$G^\alpha$ a topological group,
$\pi^\alpha \maps E^\alpha \onto X^\alpha$ surjective and
$\rho^\alpha \maps Y^\alpha \times G^\alpha \to Y^\alpha$
an effective right action of $G^\alpha$ on $Y^\alpha$. Then
\begin{multline}
\Bun_\Ob \seqname{B} \defeq \maximal{\Atl^\Bun_\Ob} \seqname{B}
\end{multline}
\end{definition}

\begin{definition}[Bundle maps]
\label{def:Bunmorph}
Let \\
$
  \seqname{B}
  \defeq
  \set
    {{
      \seqname{B}^\alpha \defeq
      (
        E^\alpha,
        X^\alpha,
        Y^\alpha,
        G^\alpha,
        \pi^\alpha,
        \rho^\alpha
      )
    }}%
    [
      \alpha \prec \Alpha
    ]
$,
where $E^\alpha, X^\alpha, Y^\alpha$ are topological spaces,
$G^\alpha$ a topological group,
$\pi^\alpha \maps E^\alpha \onto X^\alpha$ surjective and
$\rho^\alpha \maps Y^\alpha \times G^\alpha \to Y^\alpha$
an effective right action of $G^\alpha$ on $Y^\alpha$. Then
\begin{multline}
\Bun_\Ar \seqname{B} \defeq \maximal{\Atl^\Bun_\Ar} \seqname{B}
\end{multline}
\begin{multline}
\Bun \seqname{B}
\defeq
  \bigl (
    \Bun_\Ob \seqname{B},
    \Bun_\Ar \seqname{B},
    \compose[A]
  \bigr )
\end{multline}

Let $(\seqname{A}^i,\seqname{B}^i) \in \Bun_\Ob \seqname{B}$, $i=1,2$.
Then
$
  \funcseqname{f}
  \defeq \\
  (
    \funcname{f}_E \maps E^1 \to E^2,
    \funcname{f}_X \maps X^1 \to X^2,
    \funcname{f}_Y \maps Y^1 \to Y^2,
    \funcname{f}_G \maps G^1 \to G^2
  )
$
is a bundle map from $(\seqname{A}^1,\seqname{B}^1)$ to
$(\seqname{A}^1,\seqname{B}^1)$ iff it is a bundle-atlas morphism from
$A^1$ to $A^2$. The identity morphism for
$(\seqname{A}^i,\seqname{B}^i)$ is
\begin{equation}
\Id_{(\seqname{A}^i,\seqname{B}^i)} \defeq
\bigl (
  (\Id_{E^i}, \Id_{X^i}, \Id_{Y^i}, \Id_{G^i}),
  (\seqname{A}^i,\seqname{B}^i),
  (\seqname{A}^i,\seqname{B}^i)
\bigr )
\end{equation}
\end{definition}

\begin{theorem}[Categories of fiber bundles]
\label{the:Bun}
Let \\
$
  \seqname{B}
  \defeq
  \set
    {{
      \seqname{B}^\alpha \defeq
      (
        E^\alpha,
        X^\alpha,
        Y^\alpha,
        G^\alpha,
        \pi^\alpha,
        \rho^\alpha
      )
    }}%
    [
      \alpha \prec \Alpha
    ]
$,
where $E^\alpha, X^\alpha, Y^\alpha$ are topological spaces,
$G^\alpha$ a topological group,
$\pi^\alpha \maps E^\alpha \onto X^\alpha$ surjective and
$\rho^\alpha \maps Y^\alpha \times G^\alpha \to Y^\alpha$
an effective right action of $G^\alpha$ on $Y^\alpha$. Then
$\Bun \seqname{B}$ is a category and $Id_{B^\alpha}$ is the identity
morphism for $B^\alpha$.

\begin{proof}
{
  \showlabelsinline
  The result follows directly from \fullcref{def:Bun},
}
\fullcref{def:Bunmorph} and \pagecref{lem:BunATLiscat}.
\end{proof}
\end{theorem}

\begin {definition}[Functor from fiber bundles to Local Coordinate Spaces]
\label{def:BunLCS}
Let $\seqname{E}$, $\seqname{X}$ and $\seqname{Y}$ be sets of
topological spaces, $\seqname{G}$ be a set of topological groups,
$
  \seqname{B}^\alpha \defeq
  (
    E^\alpha \in \seqname{E},
    X^\alpha \in \seqname{X},
    Y^\alpha \in \seqname{Y},
    G^\alpha \in \seqname{Y},
    \pi^\alpha,
    \rho^\alpha
  )
$,
$\alpha \prec \Alpha$, be a protobundle,
$\seqname{B} \defeq \set{{\seqname{B}^\alpha}}[\alpha \prec \Alpha]$
be a set of protobundles,
$C^\alpha \defeq X^\alpha \times Y^\alpha$,
$\catname{E}$, $\catname{X}$, $\catname{Y}$ and $\catname{G}$ be model
categories,
$\catseqname{XYG\rho}$ be a $G$-$\rho$-model category,
$
  \catseqname{M} \defeq
   \Bigl  (
      \catname{E},
      \catseqname{XYG\rho},
      \catname{X},
      \catname{Y},
      \catname{G}
   \Bigr  )
$,
$
  \Triv{\catseqname{M}} \defeq \\
    \Bigl  (
      \Trivcat{\seqname{E}},
      \Trivcat[\Bun-]{\seqname{B}},
      \Trivcat{\seqname{X}},
      \Trivcat{\seqname{Y}},
      \Trivcat{\seqname{G}}
    \Bigr  )
$,
$\Triv{\catseqname{M}} \SUBCAT[full-] \catseqname{M}$,
$
  \seqname{B}^i \defeq
  (
    E^i,
    X^i,
    Y^i,
    G^i,
    \pi^i,
    \rho^i
  )
  \in \seqname{B}
$
$i=1,2$,
with group operation $\star^i \maps G^i \times G^i \onto G^i$,
$\seqname{A}^i$ a maximal $\pi^i$-$G^i$-$\rho^i$-bundle atlas of $E^i$
in the coordinate space $C^i \defeq X^i \times Y^i$,
$\pi^i_X \defeq \pi_1 \maps X^i \times Y^i \onto X^i$,
$\pi^i_Y \defeq \pi_2 \maps X^i \times Y^i \onto Y^i$,
$
  \seqname{F}^i \defeq
  (
    \pi^i,
    \pi^i_X,
    \pi^i_Y,
    \star^i,
    \rho^i
  )
$,
$
  \Singcat{\Triv{\catseqname{M}^i}} \defeq \\
    \Bigl  (
      \singcat{\Triv{E^i}},
      \singcat{\Triv[G^i-\rho^i-]{X^i,Y^i}},
      \singcat{X^i},
      \singcat{Y^i},
      \singcat{G^i},
    \Bigr  )
$,
$
  \seqname{M}^i \defeq
    \Bigl (
      \Triv{E^i},
      \Triv[G^i-\rho^i-]{X^i,Y^i},
      X^i,
      Y^i,
      G^i
    \Bigr )
$,
$
  \Sigma \defeq \\
  \bigl (
      (0,2),
      (1,2),
      (1,3),
      (4,4,4),
      (3,4,3)
  \bigr )
$,
$
  \seqname{L}^i \defeq
  \left (
    \Singcat{\Triv{\catseqname{M}^i}},
    \seqname{M}^i,
    \seqname{A}^i,
    \seqname{F}^i,
    \Sigma
  \right )
$,
$
  \seqname{L}^{i,\catseqname{M}} \defeq \\
  (
    \catseqname{M},
    \seqname{M}^i,
    \seqname{A}^i,
    \seqname{F}^i,
    \Sigma
  )
$
and
$
   \funcseqname{f} \defeq \\
     (
       \funcname{f}_E \maps E^1 \to E^2,
       \funcname{f}_X \maps X^1 \to X^2,
       \funcname{f}_Y \maps Y^1 \to Y^2,
       \funcname{f}_G \maps G^1 \to G^2
     )
$
a $\seqname{B^1}$-$\seqname{B^2}$ bundle-atlas morphism from
$\seqname{A}^1$ to $\seqname{A}^2$. Then
\begin{multline}
\Functor^\Bun_{\Fib,\LCS}
  (\seqname{A}^i, \seqname{B}^i) \defeq
\seqname{L}^i
\end{multline}
\begin{multline}
\Functor^\Bun_{\Fib,\LCS}
  \bigl (
    \funcseqname{f},
    (
      \seqname{A}^1,
      \seqname{B}^1
    ),
    (
      \seqname{A}^2,
      \seqname{B}^2
    )
  \bigr )
\defeq
\Bigl (
  \bigl (
    \funcname{f}_E,
    \funcname{f}_X \times \funcname{f}_Y,
    \funcname{f}_X,
    \funcname{f}_Y,
    \funcname{f}_G
  \bigr ),
  \seqname{L}^1,
  \seqname{L}^2
\Bigr )
\end{multline}

\begin{multline}
\Functor^{\Bun,\catseqname{M}}_{\Fib,\LCS}
  (\seqname{A}^i, \seqname{B}^i) \defeq
\seqname{L}^{i,\catseqname{M}}
\end{multline}
\begin{multline}
\Functor^{\Bun,\catseqname{M}}_{\Fib,\LCS}
  \bigl (
    \funcseqname{f},
    (
      \seqname{A}^1,
      \seqname{B}^1
    ),
    (
      \seqname{A}^2,
      \seqname{B}^2
    )
  \bigr )
\defeq \\
\Bigl (
  \bigl (
    \funcname{f}_E,
    \funcname{f}_X \times \funcname{f}_Y,
    \funcname{f}_X,
    \funcname{f}_Y,
    \funcname{f}_G
  \bigr ),
  \seqname{L}^{1,\catseqname{M}},
  \seqname{L}^{2,\catseqname{M}}
\Bigr )
\end{multline}

\begin{multline}
\LCS^\Bun_\Ob \seqname{B}
\defeq
\set
  {
    \Functor^\Bun_{\Fib,\LCS}
      \bigl (
        \seqname{A},
        (
          E^\alpha,
          X^\alpha,
          Y^\alpha,
          G^\alpha,
          \pi^\alpha,
          \rho^\alpha
        )
      \bigr )
  }%
  [
    {
      (
        E^\alpha,
        X^\alpha,
        Y^\alpha,
        G^\alpha,
        \pi^\alpha,
        \rho^\alpha
      )
      \in
      \seqname{B}
    },
    {
      \maximal{\isAtl^\Bun_\Ob}
        (
          \seqname{A},
          E^\alpha,
          X^\alpha,
          Y^\alpha,
          G^\alpha,
          \pi^\alpha,
          \rho^\alpha
        )
    }
  ]*
\end{multline}
\begin{multline}
\LCS^\Bun_\Ar \seqname{B}
\defeq
\set
  {
    \Functor^\Bun_{\Fib,\LCS}
      \bigl (
        \funcseqname{f},
        (\seqname{A}^1, \seqname{B}^1),
        (\seqname{A}^2, \seqname{B}^2)
      \bigr )
  }%
  [
    {\seqname{B}^i \in \seqname{B}},
    {
      \maximal{\isAtl^\Bun_\Ar}
        (
          \seqname{A}^1,
          \seqname{B}^1,
          \seqname{A}^2,
          \seqname{B}^2,
          \funcseqname{f}
        )
    }
  ]*
\end{multline}
\begin{multline}
\LCS^\Bun \seqname{B}
\defeq
  \bigl (
    \LCS^\Bun_\Ob \seqname{B},
    \LCS^\Bun_\Ar \seqname{B},
    \compose[A]
  \bigr )
\end{multline}

\begin{multline}
\LCS^{\Bun,\catseqname{M}}_\Ob \seqname{B}
\defeq
\set
  {
    \Functor^{\Bun,\catseqname{M}}_{\Fib,\LCS}
      \bigl (
        \seqname{A},
        (
          E^\alpha,
          X^\alpha,
          Y^\alpha,
          G^\alpha,
          \pi^\alpha,
          \rho^\alpha
        )
      \bigr )
  }%
  [
    {
      (
        E^\alpha,
        X^\alpha,
        Y^\alpha,
        G^\alpha,
        \pi^\alpha,
        \rho^\alpha
      )
      \in
      \seqname{B}
    },
    {
      \maximal{\isAtl^\Bun_\Ob}
        (
          \seqname{A},
          E^\alpha,
          X^\alpha,
          Y^\alpha,
          G^\alpha,
          \pi^\alpha,
          \rho^\alpha
        )
    }
  ]*
\end{multline}
\begin{multline}
\LCS^{\Bun,\catseqname{M}}_\Ar \seqname{B}
\defeq
\set
  {
    \Functor^{\Bun,\catseqname{M}}_{\Fib,\LCS}
      \bigl (
        \funcseqname{f},
        (\seqname{A}^1, \seqname{B}^1),
        (\seqname{A}^2, \seqname{B}^2)
      \bigr )
  }%
  [
    {\seqname{B}^i \in \seqname{B}},
    {
      \maximal{\isAtl^\Bun_\Ar}
        (
          \seqname{A}^1,
          \seqname{B}^1,
          \seqname{A}^2,
          \seqname{B}^2,
          \funcseqname{f}
        )
    }
  ]*
\end{multline}
\begin{multline}
\LCS^{\Bun,\catseqname{M}} \seqname{B}
\defeq
  \bigl (
    \LCS^{\Bun,\catseqname{M}}_\Ob \seqname{B},
    \LCS^{\Bun,\catseqname{M}}_\Ar \seqname{B},
    \compose[A]
  \bigr )
\end{multline}
\end{definition}

\begin{theorem}[Functor from fiber bundles to Local Coordinate Spaces]
\label{the:BunLCS}
Let $\seqname{E}$, $\seqname{X}$ and $\seqname{Y}$ be sets of
topological spaces, $\seqname{G}$ be a set of topological groups,
$
  \seqname{B}^\alpha \defeq
  (
    E^\alpha \in \seqname{E},
    X^\alpha \in \seqname{X},
    Y^\alpha \in \seqname{Y},
    G^\alpha \in \seqname{Y},
    \pi^\alpha,
    \rho^\alpha
  )
$,
$\alpha \prec \Alpha$, be a protobundle,
$\seqname{B} \defeq \set{{\seqname{B}^\alpha}}[\alpha \prec \Alpha]$
be a set of protobundles,
$C^\alpha \defeq X^\alpha \times Y^\alpha$,
$\catname{E}$, $\catname{X}$, $\catname{Y}$ and $\catname{G}$ be model
categories,
$\catseqname{XYG\rho}$ be a $G$-$\rho$-model category,
$
  \catseqname{M} \defeq
   \Bigl  (
      \catname{E},
      \catseqname{XYG\rho},
      \catname{X},
      \catname{Y},
      \catname{G}
   \Bigr  )
$,
$
  \Triv{\catseqname{M}} \defeq \\
    \Bigl  (
      \Trivcat{\seqname{E}},
      \Trivcat[\Bun-]{\seqname{B}},
      \Trivcat{\seqname{X}},
      \Trivcat{\seqname{Y}},
      \Trivcat{\seqname{G}}
    \Bigr  )
$,
$\Triv{\catseqname{M}} \SUBCAT[full-] \catseqname{M}$,
$
  \seqname{B}^i \defeq
  (
    E^i,
    X^i,
    Y^i,
    G^i,
    \pi^i,
    \rho^i
  )
  \in \seqname{B}
$,
$i=1,2$,
with group operation $\star^i \maps G^i \times G^i \onto G^i$,
$\seqname{A}^i$ a maximal $\pi^i$-$G^i$-$\rho^i$-bundle atlas of $E^i$
in the coordinate space $C^i \defeq X^i \times Y^i$,
$\seqname{B}^i \defeq (E^i, X^i, Y^i, G^i, \pi^i, \rho^i)$,
$\pi^i_X \defeq \pi_1 \maps X^i \times Y^i \onto X^i$,
$\pi^i_Y \defeq \pi_2 \maps X^i \times Y^i \onto Y^i$,
$
  \seqname{F}^i \defeq
  (
    \pi^i,
    \pi^i_X,
    \pi^i_Y,
    \star^i,
    \rho^i
  )
$,
$
  \Singcat{\Triv{\catseqname{M}^i}} \defeq
    \Bigl  (
      \singcat{\Triv{E^i}},
      \singcat{\Triv[G^i-\rho^i-]{X^i,Y^i}},
      \singcat{X^i},
      \singcat{Y^i},
      \singcat{G^i},
    \Bigr  )
$,
$
  \seqname{M}^i \defeq
    \Bigl (
      \Triv{E^i},
      \Triv[G^i-\rho^i-]{X^i,Y^i},
      X^i,
      Y^i,
      G^i
    \Bigr )
$,
$
  \Sigma \defeq \\
  \bigl (
      (0,2),
      (1,2),
      (1,3),
      (4,4,4),
      (3,4,3)
  \bigr )
$,
$
  \seqname{L}^i \defeq
  \left (
    \Singcat{\Triv{\catseqname{M}^i}},
    \seqname{M}^i,
    \seqname{A}^i,
    \seqname{F}^i,
    \Sigma
  \right )
$
$
  \seqname{L}^{i,\catseqname{M}} \defeq \\
  (
    \catseqname{M},
    \seqname{M}^i,
    \seqname{A}^i,
    \seqname{F}^i,
    \Sigma
  )
$
and
$
   \funcseqname{f}^i \defeq \\
     (
       \funcname{f}^i_E \maps E^i \to E^{i+1},
       \funcname{f}^i_X \maps X^i \to X^{i+1},
       \funcname{f}^i_Y \maps Y^i \to Y^{i+1},
       \funcname{f}^i_G \maps G^i \to G^{i+1}
     )
$,
$i=1,2$, a bundle map from $(\seqname{A}^1,\seqname{B}^1)$ to
$(\seqname{A}^2,\seqname{B}^2)$. Then:

$\seqname{L}^i$ is a local $\Singcat{\Triv{\catseqname{M}^i}}$-$\Sigma$
coordinate space and $\seqname{L}^{i,\catseqname{M}}$ is a local
$\catseqname{M}$-$\Sigma$ coordinate space.

\begin{proof}
They satisfy the criteria in \pagecref{def:LCS}:
\begin{enumerate}
\item $\Singcat{\Triv{\catseqname{M}^i}}$ and $\catseqname{M}$ are
sequences of categories by construction.

\item
$
  \seqname{M}^i \seqin
  \Singcat{\Triv{\catseqname{M}^i}} \SUBCAT[full-]
  \Triv{\catseqname{M}}
$
by construction and $\Triv{\catseqname{M}} \SUBCAT[full-]\catseqname{M}$
by hypothesis.

\item $\singcat{\Triv{E^i}}$ is a model category by
\pagecref{def:trivmod}, $\singcat{\Triv[G^i-\rho^i-]{X^i,Y^i}}$ is a
model category by \pagecref{def:Grhomodtriv}, $\catname{E}$ is a model
category by hypothesis and $\catseqname{XYG\rho}$ is a model category by
\pagecref{def:Grhomodcat}.

\item
$
  \left (
    \Singcat{\Triv{\catseqname{M}^i}},
    \seqname{M}^i,
    \Sigma,
    \seqname{F}^i
  \right )
$
and $(\catseqname{M}, \seqname{M}^i, \Sigma, \seqname{F}^i)$
are prestructures:
$\seqname{F}^i$ has $\Singcat{\Triv{\catseqname{M}^i}}$-signatures
$\Sigma$ and $\seqname{F}^i$ has $\catseqname{M}$-signatures $\Sigma$.

\item $\seqname{A}^i$ is a maximal m-atlas of $\Triv{E^i}$ in
$\Triv[G^i-\rho^i-]{X^i,Y^i}$ by \pagecref{lem:Bun-ATLmorph}.
\item There are no constraint functions.
\end{enumerate}
\end{proof}

$\LCS^\Bun \seqname{B}$ and $\LCS^{\Bun,\catseqname{M}} \seqname{B}$ are
categories and the identity morphism for $\seqname{L}^i$ is
$
  \Id_{\seqname{L}^i} \defeq
  \Bigl (
    \bigl (
      \Id_{E^i}
      \Id_{X^i}
      \Id_{X^i \times Y^i}
      \Id_{Y^i}
      \Id_{G^i}
    \bigr ),
    \seqname{L}^i,
    \seqname{L}^i
  \Bigr )
$.

\begin{proof}
$\LCS^\Bun \seqname{B}$ and $\LCS^{\Bun,\catseqname{M}}\seqname{B}$
satisfy the definition of a category:
\begin{enumerate}
\item Composition: \newline
$\funcseqname{f}^2 \compose \funcseqname{f}^1$ is a bundle map from
$(\seqname{A}^1,\seqname{B}^1)$ to $(\seqname{A}^1,\seqname{B}^1)$ by
\pagecref{cor:Bun-ATLmorph}.
Then \\
$
  \Functor^\Bun_{\Fib,\LCS}
    \bigl (
      \funcseqname{f}^2,
      (\seqname{A}^2, \seqname{B}^2),
      (\seqname{A}^3, \seqname{B}^3)
    \bigr )
  \compose
  \Functor^\Bun_{\Fib,\LCS}
    \bigl (
      \funcseqname{f}^1,
      (\seqname{A}^1, \seqname{B}^1),
      (\seqname{A}^2, \seqname{B}^2)
    \bigr )
  = \\
  \Functor^\Bun_{\Fib,\LCS}
    \bigl (
      \funcseqname{f}^2 \compose \funcseqname{f}^1,
      (\seqname{A}^1, \seqname{B}^1),
      (\seqname{A}^3, \seqname{B}^3)
    \bigr )
$
and \\
$
  \Functor^{\Bun,\catseqname{M}}_{\Fib,\LCS}
    \bigl (
      \funcseqname{f}^2,
      (\seqname{A}^2, \seqname{B}^2),
      (\seqname{A}^3, \seqname{B}^3)
    \bigr )
  \compose
  \Functor^{\Bun,\catseqname{M}}_{\Fib,\LCS}
    \bigl (
      \funcseqname{f}^1,
      (\seqname{A}^1, \seqname{B}^1),
      (\seqname{A}^2, \seqname{B}^2)
    \bigr )
  = \\
  \Functor^{\Bun,\catseqname{M}}_{\Fib,\LCS}
    \bigl (
      \funcseqname{f}^2 \compose \funcseqname{f}^1,
      (\seqname{A}^1, \seqname{B}^1),
      (\seqname{A}^3, \seqname{B}^3)
    \bigr )
$.

\item Associativity: \newline
Composition is associative by \pagecref{lem:atlcomp}.
\item Unit: \newline
$\Id_{\seqname{L}^i}$ is an identity morphism by \pagecref{lem:atlcomp}.
\end{enumerate}
\end{proof}

$
\Functor^\Bun_{\Fib,\LCS}
  \bigl (
    \funcseqname{f}^i,
    (\seqname{A}^i, \seqname{B}^i),
    (\seqname{A}^{i+1}, \seqname{B}^{i+1})
  \bigr )
$
is a morphism from $\LCS^\Bun \seqname{B}^i$ to \\
$\LCS^\Bun \seqname{B}^{i+1}$ and
$
\Functor^{\Bun,\catseqname{M}}_{\Fib,\LCS}
  \bigl (
    \funcseqname{f}^i,
    (\seqname{A}^i, \seqname{B}^i),
    (\seqname{A}^{i+1}, \seqname{B}^{i+1})
  \bigr )
$
is a strict morphism from
$\LCS^{\Bun,\catseqname{M}} \seqname{B}^i$ to
$\LCS^{\Bun,\catseqname{M}} \seqname{B}^{i+1}$.

\begin{proof}
It satisfies the conditions of \pagecref{def:LCSmorph}\!:
\begin{enumerate}
\item Prestructure morphism:
\newline
$\funcseqname{f}^1$ $\Sigma$-commutes with
$\funcseqname{F}^1$, $\funcseqname{F}^2$.
\item m-atlas morphism:
\newline
$(\funcname{f}_0, \funcname{f}_1)$ is an m-atlas morphism from
$(\Triv{E^1}, \Triv[G^1-\rho^1-]{C^1})$ to \\
$(\Triv{E^2}, \Triv[G^2-\rho^2-]{C^2})$ by \pagecref{lem:Bun-ATLmorph}.
\end{enumerate}
\end{proof}

$\Functor^\Bun_{\Fib,\LCS}$ is a functor from $\Bun \seqname{B}$ to
$\LCS^\Bun \seqname{B}$ and
$\Functor^{\Bun,\catseqname{M}}_{\Fib,\LCS}$ is a functor from
$\Bun \seqname{B}$ to $\LCS^{\Bun,\catseqname{M}} \seqname{B}$.

\begin{proof}
$\Functor^\Bun_{\Fib,\LCS}$ and
$\Functor^{\Bun,\catseqname{M}}_{\Fib,\LCS}$ satisfy the definition of a
functor:
\begin{enumerate}
\item F(f: A to B) = F(f): F(A) to  F(B):
\begin{enumerate}
\item
$
\Functor^\Bun_{\Fib,\LCS}
  \bigl (
    \funcseqname{f}^i
    (\seqname{A}^i, \seqname{B}^i),
    (\seqname{A}^{i+1}, \seqname{B}^{i+1})
  \bigr )
= \\
  \bigl (
    (
      \funcname{f}^i_E,
      \funcname{f}^i_X \times \funcname{f}^i_Y,
      \funcname{f}^i_X,
      \funcname{f}^i_Y,
      \funcname{f}^i_G
    ),
    \seqname{L}^i,
    \seqname{L}^{i+1}
  \bigr ) \\
\Functor^{\Bun,\catseqname{M}}_{\Fib,\LCS}
  \bigl (
    \funcseqname{f}^i
    (\seqname{A}^i, \seqname{B}^i),
    (\seqname{A}^{i+1}, \seqname{B}^{i+1})
  \bigr )
= \\
  \bigl (
    (
      \funcname{f}^i_E,
      \funcname{f}^i_X \times \funcname{f}^i_Y,
      \funcname{f}^i_X,
      \funcname{f}^i_Y,
      \funcname{f}^i_G
    ),
    \seqname{L}^{i,\catseqname{M}},
    \seqname{L}^{i+1,\catseqname{M}}
  \bigr )
$
\item
$
 \Functor^\Bun_{\Fib,\LCS} (\seqname{A}^i, \seqname{B}^i)i =
  \seqname{L}^i \\
 \Functor^{\Bun,\catseqname{M}}_{\Fib,\LCS}(\seqname{A}^i, \seqname{B}^i) =
  \seqname{L}^{i,\catseqname{M}}
$
\end{enumerate}

\item Composition: F(f g) = F(f) F(g) \\
This follows from the proof above that $\LCS^\Bun \seqname{B}$ and
$\LCS^{\Bun,\catseqname{M}} \seqname{B}$ are categories.

\item Identity: \newline
$
\Functor^\Bun_{\Fib,\LCS} (\Id_{(\seqname{A}^i, \seqname{B}^i)}) = \\
\Functor^\Bun_{\Fib,\LCS}
  \Bigl (
    \bigl ( \Id_{E^i}, \Id_{X^i}, \Id_{Y^i}, \Id_{G^i} \bigr ),
    (\seqname{A}^i, \seqname{B}^i),
    (\seqname{A}^i, \seqname{B}^i)
  \Bigr ) = \\
  \Bigl (
    \bigl (
      \Id_{E^i},
      \Id_{X^i} \times \Id_{Y^i},
      \Triv{X^i},
      \Triv{Y^i},
      \Triv{G^i}
    \bigr ),
    \seqname{L}^i,
    \seqname{L}^i
  \Bigr ) = \\
\Id_{\seqname{L}^i} = \\
\Id_{\Functor^\Bun_{\Fib,\LCS} (\seqname{A}^i, \seqname{B}^i)} \\
\Functor^{\Bun,\catseqname{M}}_{\Fib,\LCS}
  (\Id_{(\seqname{A}^i, \seqname{B}^i)}) = \\
\Functor^{\Bun,\catseqname{M}}_{\Fib,\LCS}
  \Bigl (
    \bigl ( \Id_{E^i}, \Id_{X^i}, \Id_{Y^i}, \Id_{G^i} \bigr ),
    (\seqname{A}^i, \seqname{B}^i),
    (\seqname{A}^i, \seqname{B}^i)
  \Bigr ) = \\
  \Bigl (
    \bigl (
      \Id_{E^i},
      \Id_{X^i} \times \Id_{Y^i},
      \Triv{X^i},
      \Triv{Y^i},
      \Triv{G^i}
    \bigr ),
    \seqname{L}^{i,\catseqname{M}},
    \seqname{L}^{i,\catseqname{M}}
  \Bigr ) = \\
\Id_{\seqname{L}^{i,\catseqname{M}}} = \\
\Id_{\Functor^{\Bun,\catseqname{M}}_{\Fib,\LCS}
  (\seqname{A}^i, \seqname{B}^i)}
$
\end{enumerate}
\end{proof}
\end{theorem}

\begin{definition}[Functor from local coordinate spaces to fiber bundles]
\label{def:LCStoBun}
Let $\catname{E}$ be a model category,
$\catseqname{XYG\rho}$ a $G$-$\rho$-model category,
$\catname{X}$, $\catname{Y}$ and $\catname{G}$ topological categories,
$\catseqname{M}$ a sequence of categories,
$\seqname{M}^i$, $i=1,2$, a sequence,
$\seqname{E}^i \objin \catname{E}$,
$\seqname{C}^i = (C^i,\catseqname{C}^i) \objin \catseqname{XYG\rho}$,
$X^i \objin \catname{X}$,
$Y^i \objin \catname{Y}$,
$G^i \objin \catname{G}$,
$\pi^i \maps E^i \onto X^i$,
$\pi^i_X \maps C^i \onto X^i$,
$\pi^i_Y \maps C^i \onto Y^i$,
$\star \maps G^i \times G^i \onto G^i$,
$\rho^i \maps Y^i \times G^i \onto Y^i$,
$\Singcat{\Triv{\catseqname{M}^i}}$ a sequence of categories,
$
  \seqname{L}^i \defeq
  \Biggl (
    \Singcat{\Triv{\catseqname{M}^i}},
    \seqname{M}^i,
    \seqname{A}^i,
    \seqname{F}^i,
    \Sigma
  \Biggr )
$
and
$
  \seqname{L}^{i,\catseqname{M}} \defeq
  (
    \catseqname{M}^i,
    \seqname{M}^i,
    \seqname{A}^i,
    \seqname{F}^i,
    \Sigma
  )
$
local coordinate spaces,
$\seqname{B}^i \defeq (E^i, X^i, Y^i, G^i, \pi^i, \rho^i)$
and
$\funcseqname{f} \defeq (\funcname{f}_\beta, \beta \prec \length(\seqname{M}^1)$
a morphism from $\seqname{L}^1$ to $\seqname{L}^2$,
satisfying
\begin{enumerate}
\item
$
  \head(\catseqname{M},5) =
  \Bigl  (
    \catname{E},
    \catseqname{XYG\rho},
    \catname{X},
    \catname{Y},
    \catname{G}
  \Bigr  )
$
\item
$
  \head \Biggl ( \Singcat{\Triv{\catseqname{M}^i}},5 \Biggr )  =
  \Biggl  (
    \singcat{\Triv{E^i}},
    \singcat{\Triv[G^i-\rho^i-]{X^i,Y^i}},
    \singcat{X^i},
    \singcat{Y^i},
    \singcat{G^i}
  \Biggr  )
$
\item
$
  \head(\seqname{M}^i,5) =
  (
    \seqname{E}^i,
    \seqname{C}^i,
    X^i,
    Y^i,
    G^i
  )
$
\item $\seqname{M}^i \seqin \catseqname{M}$
\item
$
  \head(\seqname{F}^i,5) =
  (
    \pi^i,
    \pi^i_X,
    \pi^i_Y,
    \star^i,
    \rho^i
  )
$
\item $C^i = X^i \times Y^i$
\item $G^i$ is a topological group with group operation $\star^i$.
Subsequent references to $G^i$ should be read as referring to
the group rather than just the underlying topological space.
\item $\pi^i$ is surjective.
\item $\pi^i_X = \pi_1 \maps X^i \times Y^i \onto X^i$.
\item $\pi^i_Y = \pi_2 \maps X^i \times Y^i \onto Y^i$.
\item $\rho^i$ is an effective right action of $G^i$ on $Y^i$.
\item
$
  \head(\Sigma,5) =
  \bigl ( (0,2), (1,2), (1,3), (4,4,4), (3,4,3) \bigr )
$
\end{enumerate}
\begin{remark}
Due to the commutation requirement, specifying
$\catname{M}_4 = \Trivcat[\mathrm{group-}]{\seqname{G}}$ is not
necessary in order to ensure that $\funcname{f}_G$ is a homomorphism.
\end{remark}
Then
\begin{equation}
\Functor^\Bun_{\LCS,\Fib} L^i \defeq
  (E^i, X^i, Y^i, \pi^i, G^i, \rho^i, \seqname{A}^i)
\end{equation}
\begin{equation}
\Functor^\Bun_{\LCS,\Fib}
  \bigl (
    \funcseqname{f},
    (\seqname{A}^1, \seqname{B}^1),
    (\seqname{A}^2, \seqname{B}^2)
  \bigr )
\defeq
  \bigl (
    (\funcname{f}_0, \funcname{f}_2, \funcname{f}_3, \funcname{f}_4),
    \Functor^\Bun_{\LCS,\Fib} L^1,
    \Functor^\Bun_{\LCS,\Fib} L^2
  \bigr )
\end{equation}
\end{definition}

\begin{theorem}[Functor from local coordinate spaces to fiber bundles]
\label{the:LCStoBun}
Let $\catname{E}$ be a model category,
$\catseqname{XYG\rho}$ a $G$-$\rho$-model category,
$\catname{X}$, $\catname{Y}$ and $\catname{G}$ topological categories,
$\catseqname{M}$ a sequence of categories,
$
  \seqname{E}^\alpha \defeq (E^\alpha,\catname{E}^\alpha)
  \objin
  \catname{E}
$,
$\alpha \prec \Alpha$,
$
  \seqname{C}^\alpha \defeq (C^\alpha,\catname{C}^\alpha)
  \objin
  \catname{C}
$,
$X^\alpha \objin \catname{X}$,
$Y^\alpha \objin \catname{Y}$,
$G^\alpha \objin \catname{G}$ a topological group with group operation
$\star^\alpha$, \\
$\pi^\alpha \maps E^\alpha \onto X^\alpha$ surjective and
$\pi^\alpha_X \maps C^\alpha \onto X^\alpha$,
$\pi^\alpha_Y \maps C^\alpha \onto Y^\alpha$,
$\star^\alpha \maps G^\alpha \times G^\alpha \onto G^\alpha$,
$\rho^\alpha \maps Y^\alpha \times G^\alpha \onto Y^\alpha$ is
an effective right action of $G^\alpha$ on $Y^\alpha$,
$\Singcat{\Triv{\catseqname{M}^\alpha}}$ a sequence of categories,
$\seqname{M}^\alpha$ a sequence,
$
  \seqname{L}^\alpha \defeq
  \Biggl (
    \Singcat{\Triv{\catseqname{M}^\alpha}},
    \seqname{M}^\alpha,
    \seqname{A}^\alpha,
    \seqname{F}^\alpha,
    \Sigma
  \Biggr )
$
and
$
  \seqname{L}^{\alpha,\catseqname{M}} \defeq \\
  (
    \catseqname{M}^\alpha,
    \seqname{M}^\alpha,
    \seqname{A}^\alpha,
    \seqname{F}^\alpha,
    \Sigma
  )
$
local coordinate spaces,
$
  \seqname{B}^\alpha \defeq
  (
    E^\alpha,
    X^\alpha,
    Y^\alpha,
    G^\alpha,
    \pi^\alpha,
    \rho^\alpha
  )
$
and
$\seqname{B} \defeq \set{{ \seqname{B}^\alpha}}[\alpha \prec \Alpha]$,
satisfying
\begin{enumerate}
\item
$
  \head(\catseqname{M},5) =
  \Bigl  (
    \catname{E},
    \catseqname{XYG\rho},
    \catname{X},
    \catname{Y},
    \catname{G}
  \Bigr  )
$
\item
$
  \head \Biggl ( \Singcat{\Triv{\catseqname{M}^\alpha}},5 \Biggr ) =
  \Biggl  (
    \singcat{\Triv{E^\alpha}},
    \singcat{\Triv[G^\alpha-\rho^\alpha-]{X^\alpha,Y^\alpha}},
    \singcat{X^\alpha},
    \singcat{Y^\alpha},
    \singcat{G^\alpha}
  \Biggr  )
$
\item
$
  \head(\seqname{M}^\alpha,5) =
  (
    \seqname{E}^\alpha,
    \seqname{C}^\alpha,
    X^\alpha,
    Y^\alpha,
    G^\alpha
  )
$
\item $\seqname{M}^\alpha \seqin \catseqname{M}$
\item
$
  \head(\seqname{F}^\alpha,5) =
  (
    \pi^\alpha,
    \pi^\alpha_X,
    \pi^\alpha_Y,
    \star^\alpha,
    \rho^\alpha
  )
$
\item $C^\alpha = X^\alpha \times Y^\alpha$
\item $G^\alpha$ is a topological group with group operation $\star^\alpha$.
Subsequent references to $G^\alpha$ should be read as referring to
the group rather than just the underlying topological space.
\item $\pi^\alpha$ is surjective.
\item $\pi^\alpha_X = \pi_1 \maps X^\alpha \times Y^\alpha \onto X^\alpha$.
\item $\pi^\alpha_Y = \pi_2 \maps X^\alpha \times Y^\alpha \onto Y^\alpha$.
\item $\rho^\alpha$ is an effective right action of $G^\alpha$ on $Y^\alpha$.
\item
$
  \head(\Sigma,5) =
  \bigl ( (0,2), (1,2), (1,3), (4,4,4), (3,4,3) \bigr )
$
\end{enumerate}

Let $\alpha^i \prec \Alpha$, $i=1,2$,
$\seqname{M}^i \defeq \seqname{M}^{\alpha^i}$,
$\seqname{E}^i = (E^i,\catseqname{E}^i) \defeq \seqname{E}^{\alpha^i}$,
$\seqname{C}^i = (C^i,\catseqname{C}^i) \defeq \seqname{C}^{\alpha^i}$,
$X^i \defeq X^{\alpha^i}$,
$Y^i \defeq Y^{\alpha^i}$,
$G^i \defeq G^{\alpha^i}$,
$\pi^i \defeq \pi^{\alpha^i}$,
$\pi^i_X \defeq \pi^{\alpha^i}_X$,
$\pi^i_Y \defeq \pi^{\alpha^i}_Y$,
$\star^i \defeq \star^{\alpha^i}$, \\
$\rho^i \defeq \rho^{\alpha^i}$,
$
  \Singcat{\Triv{\catseqname{M}^i}} \defeq
  \Singcat{\Triv{\catseqname{M}^{\alpha^i}}}
$,
$\seqname{L}^i \defeq \seqname{L}^{\alpha^i}$,
$
  \seqname{L}^{i,\catseqname{M}} \defeq
  \seqname{L}^{\alpha^i,\catseqname{M}}
$,
$\seqname{B}^i \defeq \seqname{B}^{\alpha^i}$,
$
  \funcseqname{f}^i \defeq \\
  (\funcname{f}^i_\beta, \beta \prec \length(\seqname{M}^i)
$
a morphism from $\seqname{L}^1$ to $\seqname{L}^2$,
$\funcname{f}_E \defeq \funcname{f}_0$,
$\funcname{f}_C \defeq \funcname{f}_1$,
$\funcname{f}_X \defeq \funcname{f}_2$,
$\funcname{f}_Y \defeq \funcname{f}_3$ and
$\funcname{f}_G \defeq \funcname{f}_4$.
Then:

$\Functor^\Bun_{\LCS,\Fib} \seqname{L}^i$ and
$\Functor^\Bun_{\LCS,\Fib} \seqname{L}^{\alpha^i,\catseqname{M}}$ are
fiber bundles.

\begin{proof}
$\Functor^\Bun_{\LCS,\Fib} \seqname{L}^i$
$\Functor^\Bun_{\LCS,\Fib} \seqname{L}^{i,\catseqname{M}}$
satisfy the conditions of
\pagecref{def:Bun} by hypothesis.
\end{proof}

$\funcname{f}_C = \funcname{f}_X \times \funcname{f}_Y$
\begin{proof}
Since $\funcname{f}$ is a morphism from $L^1$ to $L^2$, it $\Sigma$
commutes with $\seqname{F}^1$, $\seqname{F}^2$.
$
  \pi_1 \compose \funcname{f}_C =
  \pi^2_X \compose \funcname{f}_C =
  \funcname{f}_X \compose \pi^2_X =
  \funcname{f}_X \compose \pi_1
$ and
$
  \pi_2 \compose \funcname{f}_C =
  \pi^2_Y \compose \funcname{f}_C =
  \funcname{f}_Y \compose \pi^2_Y =
  \funcname{f}_Y \compose \pi_2
$.
\end{proof}

$
\Functor^\Bun_{\LCS,\Fib}
  (
    \funcseqname{f},
    \seqname{L}^1,
    \seqname{L}^2
  )
$
and
$
\Functor^\Bun_{\LCS,\Fib}
  (
    \funcseqname{f},
    \seqname{L}^{1,\catseqname{M}},
    \seqname{L}^{2,\catseqname{M}}
  )
$
are bundle maps.

\begin{proof}
$
\Functor^\Bun_{\LCS,\Fib}
  (
    \funcseqname{f},
    \seqname{L}^1,
    \seqname{L}^2
  )
$
and
$
\Functor^\Bun_{\LCS,\Fib}
  (
    \funcseqname{f},
    \seqname{L}^{1,\catseqname{M}},
    \seqname{L}^{2,\catseqname{M}}
  )
$
satisfy the conditions of \pagecref{def:Bun-ATLmorph}
\begin{enumerate}
\item $\funcname{f}_E$, $\funcname{f}_X$, $\funcname{f}_Y$ and
$\funcname{f}_G$ are continuous
\item $\funcname{f}_G$ is a homomorphism due to commutation
\item $\funcseqname{f}$ commutes with $\pi^i$ and $\rho^i$
\item for any
$(U^1, V^1, \phi^1 \maps U^1 \toiso V^1) \in \seqname{A}^1$,
$(U^2, V^2, \phi^2 \maps U^2 \toiso V^2) \in \seqname{A}^2$,
the diagram
$\{ \funcname{f}_0, \phi^2, \phi^1, \funcname{f}_1 \})$
is locally nearly commutative in $X,Y,\pi,\rho$ because
$\catseqname{XYG\rho}$ is a $G$-$\rho$-model category by hypothesis,
$\singcat{\Triv[G^\alpha-\rho^\alpha-]{X^\alpha,Y^\alpha}}$ is a
$G$-$\rho$-model category by construction and
$(\funcname{f}_0, \funcname{f}_1)$ is an m-atlas morphism
from $\seqname{A}^1$ to $\seqname{A}^2$,
\end{enumerate}
$\Functor^\Bun_{\LCS,\Fib} \funcseqname{f}$ satisfies the conditions of
\pagecref{def:Bunmorph}.
\end{proof}

$\Functor^\Bun_{\LCS,\Fib}$ is a functor from $\LCS^\Bun \seqname{B}$ to
$\Bun \seqname{B}$ and is a functor from \\
$\LCS^{\Bun,\catseqname{M}} \seqname{B}$ to $\Bun \seqname{B}$.

\begin{proof}
$\Functor^\Bun_{\LCS,\Fib}$ satisfies the definition of a functor:
\begin{enumerate}
\item F(f: A to B) = F(f): F(A) to  F(B):
\begin{enumerate}
\item
$
  \Functor^\Bun_{\LCS,\Fib}
    (
      \funcseqname{f}^1,
      \seqname{L}^1,
      \seqname{L}^2
    \bigr )
  = \\
    \bigl (
      (
        \funcname{f}^1_E,
        \funcname{f}^1_X,
        \funcname{f}^1_Y,
        \funcname{f}^1_G
      ),
      \seqname{B}^1,
      \seqname{B}^2
    \bigr )
$
\item
$
  \Functor^\Bun_{\LCS,\Fib}
    (\catseqname{M}^i, \seqname{M}^i, \seqname{A}^i, \seqname{F}^i, \Sigma)
  =
    \seqname{B}^i
$
\end{enumerate}
\item Composition: \newline
$
  \Functor^\Bun_{\LCS,\Fib}
    \Bigl (
      \bigl (
        \funcseqname{f}^2,
        \seqname{L}^2,
        \seqname{L}^3
      \bigr )
      \compose[A]
      \bigl (
        \funcseqname{f}^1,
        \seqname{L}^1,
        \seqname{L}^2
      \bigr )
    \Bigr )
  = \\
  \Functor^\Bun_{\LCS,\Fib}
    \Bigl (
       \bigl (
         \funcseqname{f}^2 \compose[()] \funcseqname{f}^1,
         \seqname{L}^1,
         \seqname{L}^3
       \bigr )
    \Bigr )
  = \\
  \bigl (
    (
      \funcname{f}^2_E \compose \funcname{f}^2_2,
      \funcname{f}^2_X \compose \funcname{f}^1_X,
      \funcname{f}^2_Y \compose \funcname{f}^1_Y,
      \funcname{f}^2_G \compose\funcname{f}^1_G
    ),
    \seqname{B}^1,
    \seqname{B}^3
  \bigr )
  = \\
  \Functor^\Bun_{\LCS,\Fib}
    \bigl (
      \funcseqname{f}^2,
      \seqname{L}^2,
      \seqname{L}^3
    \bigr )
  \compose[A]
  \Functor^\Bun_{\LCS,\Fib}
    \bigl (
      \funcseqname{f}^1,
      \seqname{L}^1,
      \seqname{L}^2
    \bigr )
$
\item Identity: \newline
$
\Functor^\Bun_{\LCS,\Fib} (\Id_{\seqname{L}^i}) =
\Functor^\Bun_{\LCS,\Fib}
  \Bigl (
    \Id_{\seqname{M}^i},
    \seqname{L}^i,
    \seqname{L}^i
  \Bigr ) = \\
  \biggl (
    \bigl (
      \Id_{E^i},
      \Id_{X^i},
      \Id_{Y^i},
      \Id_{G^i}
    \bigr ),
    (\seqname{A}^i,\seqname{B}^i),
    (\seqname{A}^i,\seqname{B}^i)
  \biggr )
=
\Id_{(\seqname{A}^i,\seqname{B}^i)} =
\Id_{\Functor^\Bun_{\LCS,\Fib} \seqname{L}^i}
$
\end{enumerate}

The proof does not depend on the category, so it applies to
$\LCS^{\Bun,\catseqname{M}} \seqname{B}$ as well.
\end{proof}

$
  \Functor^\Bun_{\LCS,\Fib}
  \compose
  \Functor^\Bun_{\Fib,\LCS}
  = \Id
$
and
$
  \Functor^\Bun_{\LCS,\Fib}
  \compose
  \Functor^{\Bun,\catseqname{M}}_{\Fib,\LCS}
  = \Id
$.
\begin{proof}
Expanding the definitions, we have
\begin{enumerate}
\item
$
  \Functor^{\Bun,\catseqname{M}}_{\Fib,\LCS}
    (\seqname{A}^i, \seqname{B}^i) =
  \seqname{L}^i
$
\item
$
  \Functor^\Bun_{\LCS,\Fib} \seqname{L}^i =
  (\seqname{A}^i, \seqname{B}^i)
$
\item
$
  \Functor^{\Bun,\catseqname{M}}_{\Fib,\LCS}
    \bigl (
      (
        \funcname{f}_E,
        \funcname{f}_X,
        \funcname{f}_Y,
        \funcname{f}_G
      ),
      (\seqname{A}^1, \seqname{B}^1),
      (\seqname{A}^2, \seqname{B}^2)
    \bigr ) = \\
    \bigl (
      (
        \funcname{f}_E,
        \funcname{f}_X \times \funcname{f}_Y,
        \funcname{f}_X,
        \funcname{f}_Y,
        \funcname{f}_G
      ),
      \seqname{L}^1,
      \seqname{L}^2
    \bigr )
$
\item
$
  \Functor^\Bun_{\LCS,\Fib}
    \bigl (
      (
        \funcname{f}_E,
        \funcname{f}_X \times \funcname{f}_Y,
        \funcname{f}_X,
        \funcname{f}_Y,
        \funcname{f}_G
      ),
      \seqname{L}^1,
      \seqname{L}^2
    \bigr ) = \\
    \bigl (
      (
        \funcname{f}_E,
        \funcname{f}_X,
        \funcname{f}_Y,
        \funcname{f}_G
      ),
      (\seqname{A}^1, \seqname{B}^1),
      (\seqname{A}^2, \seqname{B}^2)
    \bigr )
$
\end{enumerate}

The proof does not depend on the category, so it applies to
$\LCS^{\Bun,\catseqname{M}} \seqname{B}$ as well.
\end{proof}

$
  \Functor^\Bun_{\Fib,\LCS}
  \compose
  \Functor^\Bun_{\LCS,\Fib}
$
is the identity functor on $\LCS^\Bun \seqname{B}$ and
$
  \Functor^{\Bun,\catseqname{M}}_{\Fib,\LCS}
  \compose
  \Functor^\Bun_{\LCS,\Fib}
$
is the identity functor on $\LCS^{\Bun,\catseqname{M}}$.

\begin{proof}
Expanding the definitions, we have
\begin{enumerate}
\item
$
  \Functor^\Bun_{\LCS,\Fib} \seqname{L}^i =
  (\seqname{A}^i, \seqname{B}^i)
$
\item
$
  \Functor^\Bun_{\Fib,\LCS} (\seqname{A}^i, \seqname{B}^i) =
  \seqname{L}^i
$
\item
$
  \Functor^\Bun_{\Fib,\LCS}
    \bigl (
      (
        \funcname{f}^1_E,
        \funcname{f}^1_X,
        \funcname{f}^1_Y,
        \funcname{f}^1_G
      ),
      (\seqname{A}^1, \seqname{B}^1),
      (\seqname{A}^2, \seqname{B}^2)
    \bigr ) = \\
    \bigl (
      (
        \funcname{f}^1_E,
        \funcname{f}^1_X \times \funcname{f}^1_Y,
        \funcname{f}^1_X,
        \funcname{f}^1_Y,
        \funcname{f}^1_G
      ),
      \seqname{L}^1,
      \seqname{L}^2
    \bigr )
$
\item
$
  \Functor^\Bun_{\LCS,\Fib}
    \bigl (
      (
        \funcname{f}^1_E,
        \funcname{f}^1_X \times \funcname{f}^1_Y,
        \funcname{f}^1_X,
        \funcname{f}^1_Y,
        \funcname{f}^1_G
      ),
      \seqname{L}^1,
      \seqname{L}^2
    \bigr ) = \\
    \bigl (
      (
        \funcname{f}^1_E,
        \funcname{f}^1_X,
        \funcname{f}^1_Y,
        \funcname{f}^1_G
      ),
      (\seqname{A}^1, \seqname{B}^1),
      (\seqname{A}^2, \seqname{B}^2)
    \bigr )
$
\end{enumerate}

The same proof applies for $\Functor^{\Bun,\catseqname{M}}_{\Fib,\LCS}$
with $\seqname{L}^{i,\catseqname{M}}$ in place of $\seqname{L}^i$.
\end{proof}
\end{theorem}

\section{Future directions}
\label{sec:fut}
If this paradigm proves useful, it can be extended to include a set of
admissible functions on the model neighborhoods of the charts,
possibly using the language of sheaves. That might be desirable for
coordinate spaces more general than Fr\'echet spaces.

Further work is needed to determine whether it is productive to
allow a local coordinate space to have more than one atlas, e.g., for
more than one bundle structure on the same base space.

The extension of paracompactness to model spaces is intended to be
useful for partitions of unity on fiber bundles. Further work is
needed to determine whether that is actually the case.

The definitions given here include some fairly strong conditions,
e.g., AOC. Further work is needed to determine whether they should be
relaxed for applications beyond manifolds and fiber bundles.

Further work is needed to determine whether the concepts of
category-based atlases\footnote{As opposed to pseudogroup based}
of model spaces and prestructures have general utility.

If the concept of nearly commutative diagrams proves useful, further
work is needed to determine whether a more general definition has
utility.

Further work is needed to devise a definition of constraints that
expresses global properties, e.g., compactness, and is both clear and
rigorous.

Further work is needed to determine conditions for mappings associated
with atlas morphisms to be model functions.

This paper uses the language of category theory as an organizing
principle, but defines various notions concretely with sets. It may be
desirable to abstract away some of the details, in the spirit of, e.g.,
topoi.

\end{document}